\documentclass[11pt,a4paper,oneside]{amsbook}
\usepackage[T1]{fontenc}
\usepackage[utf8]{inputenc}
\usepackage{mathpazo}                 % Palatino text and maths, as in Vistoli
\usepackage[scaled=0.92]{helvet}
\usepackage{microtype}
\usepackage[a4paper,textwidth=6.0in,textheight=8.8in,
            top=1.25in,headsep=0.3in,footskip=0.4in]{geometry}
\usepackage{amsmath,amssymb,amsthm,mathtools}
\usepackage{mathrsfs}
\numberwithin{equation}{section}
\numberwithin{figure}{chapter}
\numberwithin{table}{chapter}

\usepackage{graphicx}
\usepackage{xcolor}
\usepackage{tikz}
\usepackage{tikz-cd}
\usetikzlibrary{arrows.meta,decorations.markings,calc,positioning,
                shapes.geometric,backgrounds,3d}
\usepackage{float}
\usepackage[font=small,labelfont=bf,labelsep=period,skip=6pt]{caption}
\usepackage{subcaption}
\usepackage{booktabs}
\usepackage{array}
\usepackage[most]{tcolorbox}

\usepackage{enumitem}
\setlist[enumerate]{leftmargin=*,itemsep=2pt,topsep=3pt,label=(\roman*)}
\setlist[itemize]{leftmargin=*,itemsep=2pt,topsep=3pt,label=--}

\usepackage[square,numbers]{natbib}

\usepackage[colorlinks,linkcolor=black!70!blue,citecolor=black!60!blue,
            urlcolor=black!60!blue]{hyperref}
\usepackage[capitalize,nameinlink,noabbrev]{cleveref}

\newtheoremstyle{vplain}{7pt plus 2pt}{7pt plus 2pt}{\itshape}{0pt}{\bfseries}{.}{ }{}
\newtheoremstyle{vdefinition}{7pt plus 2pt}{7pt plus 2pt}{\normalfont}{0pt}{\bfseries}{.}{ }{}
\newtheoremstyle{vremark}{7pt plus 2pt}{7pt plus 2pt}{\normalfont}{0pt}{\itshape}{.}{ }{}
\theoremstyle{vplain}
\newtheorem{theorem}{Theorem}[section]
\newtheorem{proposition}[theorem]{Proposition}
\newtheorem{lemma}[theorem]{Lemma}
\newtheorem{corollary}[theorem]{Corollary}

\theoremstyle{vdefinition}
\newtheorem{definition}[theorem]{Definition}
\newtheorem{example}[theorem]{Example}
\newtheorem{construction}[theorem]{Construction}

\theoremstyle{vremark}
\newtheorem{remark}[theorem]{Remark}

\crefname{theorem}{Theorem}{Theorems}
\crefname{proposition}{Proposition}{Propositions}
\crefname{lemma}{Lemma}{Lemmas}
\crefname{corollary}{Corollary}{Corollaries}
\crefname{definition}{Definition}{Definitions}
\crefname{example}{Example}{Examples}
\crefname{remark}{Remark}{Remarks}
\crefname{construction}{Construction}{Constructions}
\crefname{equation}{}{}
\crefname{figure}{Figure}{Figures}
\crefname{section}{Section}{Sections}
\crefname{subsection}{Section}{Sections}
\crefname{chapter}{Chapter}{Chapters}
\crefname{appendix}{Appendix}{Appendices}

\definecolor{intuitionbg}{RGB}{243,246,250}
\definecolor{intuitionrule}{RGB}{110,130,165}
\newtcolorbox{intuition}[1][]{
  enhanced, breakable,
  colback=intuitionbg, colframe=intuitionbg,
  boxrule=0pt, arc=1.5pt, outer arc=1.5pt,
  left=10pt, right=9pt, top=6pt, bottom=6pt,
  boxsep=0pt,
  borderline west={2pt}{0pt}{intuitionrule},
  before skip=8pt plus 2pt, after skip=8pt plus 2pt,
  fontupper=\small,
  before upper={\textsc{\small Intuition\ifx\\#1\\\else\ (#1)\fi.}\enspace}
}
\newtcolorbox{informal}[1]{
  enhanced, breakable,
  colback=intuitionbg, colframe=intuitionbg,
  boxrule=0pt, arc=1.5pt, outer arc=1.5pt,
  left=10pt, right=9pt, top=6pt, bottom=6pt, boxsep=0pt,
  borderline west={2pt}{0pt}{intuitionrule},
  before skip=8pt plus 2pt, after skip=8pt plus 2pt,
  fontupper=\small,
  before upper={\textsc{\small #1.}\enspace}
}

\newcommand{\N}{\mathbb{N}}
\newcommand{\Z}{\mathbb{Z}}
\newcommand{\Q}{\mathbb{Q}}
\newcommand{\R}{\mathbb{R}}
\newcommand{\F}{\mathcal{F}}                  % functions on the vertex set
\newcommand{\PG}{\mathcal{P}}                 % polar grid
\newcommand{\cA}{\mathcal{A}}
\newcommand{\cC}{\mathcal{C}}
\newcommand{\cD}{\mathcal{D}}
\newcommand{\cO}{\mathcal{O}}
\newcommand{\cU}{\mathcal{U}}
\newcommand{\Xs}{\mathbf{X}}                  % a coordinate system
\newcommand{\Nset}{\N\text{-}\mathbf{Set}}
\newcommand{\CoordSet}{\mathbf{Coord}\N\text{-}\mathbf{Set}}
\newcommand{\DGA}{\mathbf{DGA}}
\newcommand{\Set}{\mathbf{Set}}

\newcommand{\Top}{\mathbf{Top}}
\newcommand{\PSh}{\mathbf{PSh}}

\newcommand{\Bx}{\square}                     % the box category
\newcommand{\Bxinj}{\square_{\mathrm{inj}}}

\newcommand{\dR}{\mathrm{dR}}
\newcommand{\HdR}[1]{H^{#1}_{\dR}}
\newcommand{\HN}[1]{H^{#1}_{N}}
\DeclareMathOperator{\Coord}{Coord}
\DeclareMathOperator{\Hom}{Hom}
\DeclareMathOperator{\End}{End}
\DeclareMathOperator{\Fun}{Fun}
\DeclareMathOperator{\Ext}{Ext}
\DeclareMathOperator{\Tor}{Tor}
\DeclareMathOperator{\im}{im}

\DeclareMathOperator{\Sym}{Sym}
\DeclareMathOperator{\Cay}{Cay}
\DeclareMathOperator{\Fix}{Fix}
\DeclareMathOperator{\stab}{stab}
\DeclareMathOperator{\Tot}{Tot}
\DeclareMathOperator{\Der}{Der}
\DeclareMathOperator{\Coind}{Coind}
\DeclareMathOperator{\Spec}{Spec}

\DeclareMathOperator{\Nerve}{N}
\newcommand{\id}{\mathrm{id}}
\newcommand{\op}{\mathrm{op}}
\newcommand{\cube}[2]{[#1;#2]}                % the cube [x;I]
\newcommand{\Qc}{Q}                           % (semi-)cubical set of an N-set
\newcommand{\Cc}{C_{\Bx}}                     % cubical (co)chains
\newcommand{\Dc}{D_{\Bx}}
\newcommand{\Nc}{N_{\Bx}}
\newcommand{\Cs}{C^{\bullet}_{\Bx}}
\newcommand{\Ds}{D^{\bullet}_{\Bx}}
\newcommand{\Ns}{N^{\bullet}_{\Bx}}
\newcommand{\Fp}{F}                           % face poset F(X)
\newcommand{\Lam}{\textstyle\bigwedge}

\newcommand{\defeq}{\coloneqq}
\newcommand{\isoto}{\xrightarrow{\ \sim\ }}

\newcommand{\vone}{\mathbf{1}}

\definecolor{cA}{RGB}{201,58,58}
\definecolor{cB}{RGB}{38,86,186}
\definecolor{cC}{RGB}{16,124,96}
\definecolor{cD}{RGB}{200,120,10}
\definecolor{cE}{RGB}{124,52,166}
\definecolor{polN}{RGB}{201,58,58}
\definecolor{polS}{RGB}{38,86,186}
\definecolor{polE}{RGB}{124,52,166}
\definecolor{e1col}{RGB}{198,62,62}
\definecolor{e2col}{RGB}{40,88,188}
\definecolor{fillcol}{RGB}{128,96,180}
\definecolor{gedge}{RGB}{60,80,120}
\tikzset{
  ar/.style={line width=1.05pt,line cap=round,shorten >=6pt,shorten <=6pt,
    postaction={decorate,decoration={markings,
      mark=at position 0.58 with {\arrow[line width=1.05pt]{Stealth[length=5pt,width=5pt]}}}}},
  lp/.style={line width=0.85pt,-{Stealth[length=4pt,width=4pt]}},
  dlp/.style={line width=1.05pt,gray!55,densely dashed,-{Stealth[length=5pt,width=5pt]}},
  vtx/.style={circle,fill=black!80,inner sep=0pt,minimum size=4.6pt},
  gvtx/.style={circle,fill=gedge!80,inner sep=0pt,minimum size=3pt},
  gline/.style={gedge!70,line width=0.55pt},
  cedge/.style={line cap=round,shorten >=2.3pt,shorten <=2.3pt,
    postaction={decorate,decoration={markings,
      mark=at position 0.60 with {\arrow[line width=1.05pt]{Stealth[length=4.5pt,width=4.5pt]}}}}},
  callout/.style={black,line width=1.4pt,line cap=round,
    shorten >=2.3pt,shorten <=2.3pt,
    postaction={decorate,decoration={markings,
      mark=at position 0.60 with {\arrow[line width=1.2pt]{Stealth[length=4.5pt]}}}}},
  poleloop/.style={line width=0.75pt,-{Stealth[length=3pt]}},
  every picture/.style={line join=round},
}

\makeatletter
\renewcommand{\@secnumfont}{\bfseries}
\makeatother
\definecolor{cF}{RGB}{190,60,140}
\definecolor{edgec}{RGB}{60,60,60}
\definecolor{fillc}{RGB}{225,232,244}
\definecolor{idc}{RGB}{201,58,58}
\tikzset{
  idar/.style={idc,line width=1.8pt,line cap=round,shorten >=6pt,shorten <=6pt,
    postaction={decorate,decoration={markings,
      mark=at position 0.80 with {\arrow[line width=1.8pt]{Stealth[length=6pt,width=6pt]}}}}},
  grar/.style={edgec,line width=0.85pt,line cap=round,shorten >=6pt,shorten <=6pt,
    postaction={decorate,decoration={markings,
      mark=at position 0.58 with {\arrow[line width=0.85pt]{Stealth[length=5pt,width=5pt]}}}}},
}
\newcommand{\figscale}{1}

\title{Discrete geometry and de Rham cohomology of graphs}
\author{Utkarsh Bajaj}
\address{University of Oxford}
\date{Trinity Term 2026}

\begin{document}

\frontmatter
\begin{titlepage}
  \centering
  \vspace*{1.5cm}
  {\LARGE\bfseries Discrete geometry and\\[4pt] de Rham cohomology of graphs\par}
  \vspace{1.8cm}
  {\large Utkarsh Bajaj\par}
  \vspace{1.6cm}
  \includegraphics[width=4cm]{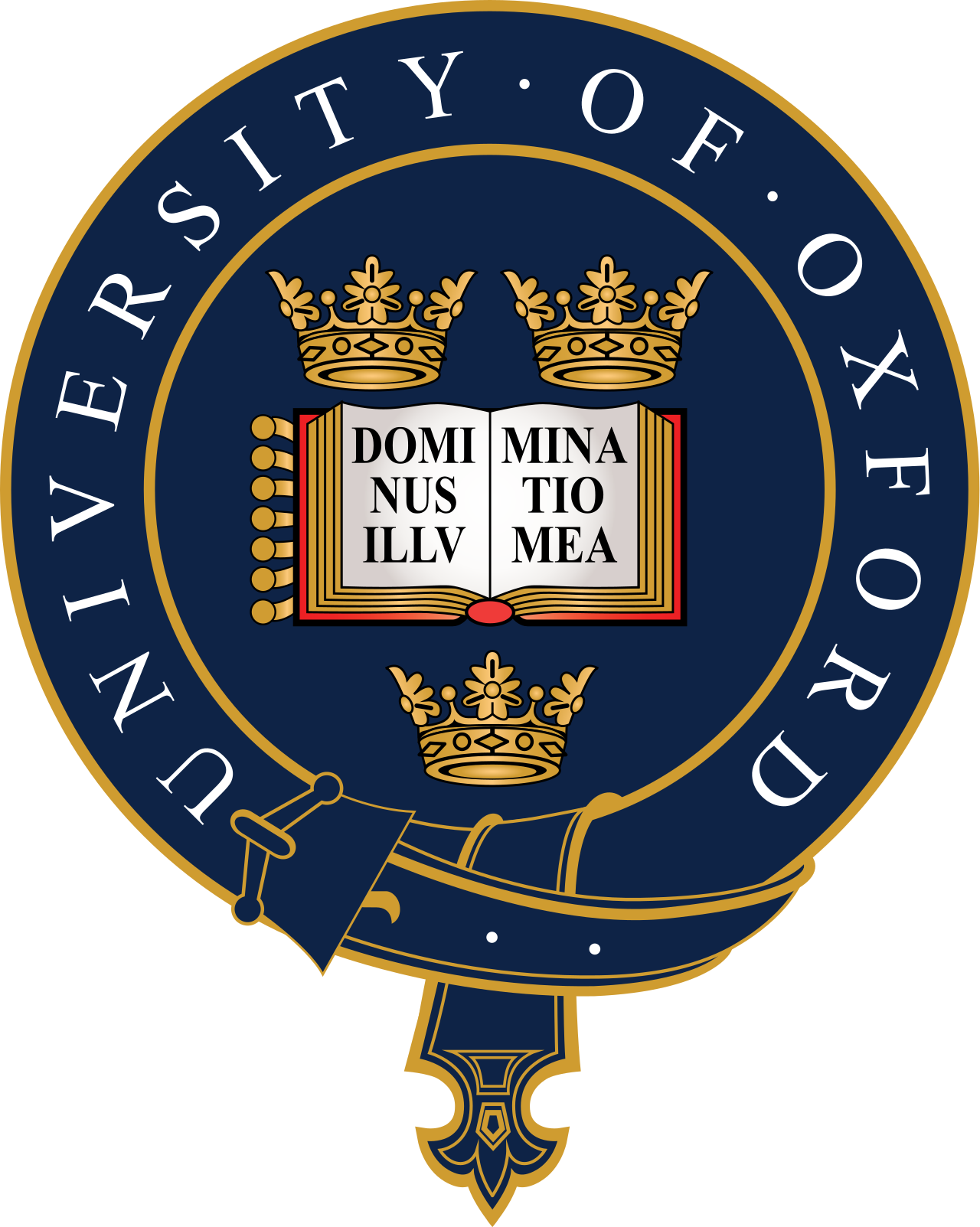}\par
  \vspace{1.6cm}
  {\normalsize A dissertation submitted for the degree of\\[2pt]
   \textit{MSc in Mathematics and Foundations of Computer Science}\par}
  \vspace{1.4cm}
  {\normalsize University of Oxford\\[2pt] Trinity Term 2026\par}
  \vfill
\end{titlepage}

\setcounter{tocdepth}{1}
\tableofcontents

\mainmatter
\setcounter{chapter}{-1}
\chapter{Introduction}
\label{ch:intro}

%Claude note: rewritten a second time, at your request, in a shorter and more
%restrained form. The earlier version's claims about a "de Rham cohomology
%spectrum" and about stacks are not proved in the body and appear only as
%directions for further work in the Conclusion.

This dissertation develops a de Rham theory for graphs equipped with
coordinate systems. A graph is treated as a space on which calculus is set
up as it is on a manifold, by choosing local coordinates, differentiating in
those coordinates, and gluing.

A \emph{coordinate system} on a graph $G=(V,E)$ is a finite family
$X_1,\dots,X_n\colon V\to V$ of pairwise commuting maps, each of which sends
a vertex to a neighbour or fixes it. Each $X_i$ defines a difference
operator $\Delta_i f=f\circ X_i-f$ on functions, and the $\Delta_i$ commute
because the $X_i$ do. Commuting difference operators and the exterior
algebra on symbols $dX_1,\dots,dX_n$ give a cochain complex
\[
  0\to\Omega^0(G,\Xs)\to\Omega^1(G,\Xs)\to\dots\to\Omega^n(G,\Xs)\to 0,
\]
the discrete de Rham complex, whose cohomology $\HdR{\bullet}(G,\Xs)$ is the
de Rham cohomology of the coordinatized graph $(G,\Xs)$. The cycle graph with
its rotation has the cohomology of the circle, and the discrete torus
$C_{N_1}\times\dots\times C_{N_d}$ with its $d$ shifts has the cohomology of
the $d$-torus (\cref{sec:first-computations}). 
 
A coordinate system of rank $n$ on $G = (V,E)$ is equivalently the action of a monoid $\mathbb{N}^{n}$ on $V$. 
We use the language of $\mathbb{N}-$sets throughout the thesis. The de Rham complex of an $\mathbb{N}-$set, as a graded algebra,
carries a non-commutative product that makes it into a differential graded algebra$-$we show 
that the non-commutativity is deeply linked to underlying discreteness of the space. 

Our first major contribution is a collection of results about discrete geometry which have analogs in the smooth world.
We prove the discrete Poincar\'e lemma, which says that the discrete space $\mathbb{Z}^{n}$ $-$ viewed as the infinite "grid" graph - has de Rham cohomology
0 in higher degrees i.e. is contractible. This is proven using homotopy invariance for $\N-$sets. We also show that for \textbf{any} finite graph 
admitting a global coordinate system, the associated de Rham cohomology 
is that of a torus. This result is shown to be a discrete analogue of Atiyah-Bott's localization theorem in equivariant cohomology.

However, we will show that this is not really the "correct" cohomology to consider. 
A coordinate system on a bare set $V$ \textit{induces} a directed graph, where the edges are pairs
of vertices connected by the action of a coordinate. And for most finite graphs, a single coordinate
system on the graph induces a directed graph that is not isomorphic to the underlying graph. We fix this by
introducing the notion of a cubical set, and the de Rham complex of a cubical set is identical to cubical cochain complex. 
Then, the "correct" de Rham complex will be the normalized cubical cochain complex, where the degenerate cubes
will account for the mismatch between the directed graph induced by the coordinate system and the underlying graph.

We show how the language of cubical sets fit naturally to describe the geometry of discrete spaces. Differential forms
are shown to be functions on the cubes of a cubical set, and the cubical cochain complex carries a cup product, 
diagonal to the Serre diagonal map  on cubes, which is shown to be isomorphic to the product structure on the de Rham complex. 
We also use cubical sets to define an atlas of $\mathbb{N}-$sets by gluing their cubical sets along cubical subsets. 
Cubical sets have a natural topology  on them$-$the Alexandrov topology$-$and discrete differential forms are described as sheaves on this space. 
The second major contribution of this thesis: the de Rham cohomology of an atlas of $\mathbb{N}-$sets is shown to be isomorphic to the sheaf cohomology
of the constant sheaf $\underline\Z$, similar to continuous case. 

One feature separates this construction from the existing cohomology theories of graphs.
Those theories are intrinsic: a graph determines its cohomology. 
Here the cohomology depends on the atlas of coordinate patches on the graph, and that dependence is deliberate. 
A coordinate chart on a graph determines a cell structure on the part of the graph it covers, 
with that part as the 1-skeleton; what a chart supplies is the cells above dimension one, 
and different charts supply different ones. A finite graph admits only finitely many coordinate systems, 
so the set of cohomologies attached to a finite graph is itself a finite invariant of the graph. We define this to 
be \textit{de Rham cohomology spectrum} of a finite graph. This is our third major contribution.

Overall, the thesis 
Proofs that consist mainly of sign bookkeeping are collected in
\cref{app:proofs}. Informal explanations are set apart
in tinted boxes, and explain the intuition behind the mathematics. We assume familiarity with differential
geometry and algebraic topology at the level of Bott and Tu~\citep{bott-tu1982} and
Hatcher~\citep{hatcher2002}, with homological algebra at the level of
Weibel~\citep{weibel1994}, and with basic category theory.

\chapter{Motivation}
\label{ch:motivation}

Differential geometry is the study of spaces on which one can
differentiate. The spaces in question are manifolds, and the reason one can
differentiate on them is that they are glued from copies of $\R^n$, on
which differentiation is already understood. Every geometric concept, from
vector fields and differential forms to the exterior derivative and de Rham
cohomology, is built on that local model and on the compatibility of the
charts where they overlap \citep{bott-tu1982}.

Discrete spaces have no such local model, and calculus on them is not part
of the standard curriculum. Consider the identity
\begin{equation}\label{eq:ftc}
  \sum_{i=1}^{N}\bigl[f(i)-f(i-1)\bigr]=f(N)-f(0),
\end{equation}
valid for any function $f\colon\N\to\Z$, or indeed for functions into any
abelian group. Read the left-hand side as the integral of a derivative over
an interval and the right-hand side as an evaluation on the boundary of that
interval, and \eqref{eq:ftc} is the fundamental theorem of calculus. Knill
takes this identity as the starting point of his quantum calculus programme
\citep{knill-web}, which uses graphs as the geometric context for a discrete
theory and establishes graph-theoretic counterparts of a range of classical
theorems, among them a Green--Stokes theorem, a Gauss--Bonnet theorem and a
Poincar\'e--Hopf theorem \citep{knill2011,knill2012a,knill2012b}.

\section{Calculus on \texorpdfstring{$\N$}{N}}
\label{sec:calculus-on-N}

%Claude note: the two notes to yourself at the start and end of this section
%("Me: say some stuff...", "We definitely should do something here") and the
%"Editor's note" left by an earlier AI pass have been removed and replaced by
%the two paragraphs framing the section. Please still check the references
%GKP and Kelley-Peterson yourself before submission, as the editor's note
%asked.

Calculus on $\N$ is the calculus of finite differences. This section
records the part of it that is used later, and fixes the point of view that
a single discrete direction is completely understood, so that the work of
the dissertation lies in combining several commuting directions.

Let $\F$ denote the space of all functions $f\colon\N\to\Z$. The forward
difference operator $\Delta\colon\F\to\F$ and the summation operator
$\Sigma\colon\F\to\F$ are
\[
  (\Delta f)(n)=f(n+1)-f(n),\qquad (\Sigma f)(n)=\sum_{k=0}^{n-1}f(k),
\]
with the convention that the constant of summation is fixed by
$(\Sigma f)(0)=0$. A direct verification gives the discrete fundamental
theorem
\[
  \Delta(\Sigma f)=f,\qquad \Sigma(\Delta f)(n)=f(n)-f(0).
\]
So $\Delta$ and $\Sigma$ play the roles of $d/dx$ and $\int$, and the
second identity is \eqref{eq:ftc}.

The analogy runs further. The classical accounts are Boole~\citep{boole1860} and
Jordan~\citep{jordan1965}, and Graham, Knuth and Patashnik~\citep[Ch.~2]{gkp1994} and
Kelley and Peterson~\citep{kelley-peterson2001} are the modern references. The monomials $x^n$
of the continuous theory, which satisfy $\frac{d}{dx}x^n=nx^{n-1}$, are
replaced by the falling factorials
\[
  x^{(n)}=x(x-1)\cdots(x-n+1),\qquad x^{(0)}=1,
\]
which satisfy $\Delta x^{(n)}=n\,x^{(n-1)}$. The function $2^x$ plays the
role of $e^x$, since $\Delta 2^x=2^x$. The Newton series
\[
  f(n)=\sum_{k=0}^{m}\binom{n}{k}(\Delta^k f)(0)+R_m(n)
\]
replaces the Taylor series, with a remainder $R_m$ that vanishes identically
when $f$ is a polynomial of degree at most $m$. The $z$-transform
stands in for the Laplace transform, and under it $\Delta$ becomes
multiplication by $z-1$ \citep{jordan1965,kelley-peterson2001}. A
$q$-deformed version, in which $\Delta$ is replaced by the Jackson
derivative, is the subject of Kac and Cheung~\citep{kac-cheung2002}, and $q$-de Rham
complexes appear in arithmetic geometry \citep{scholze2017}. We do not
pursue the $q$-deformation here.
%Claude note: the draft said "we shall meet it again in Chapter 8"; there is
%no Chapter 8, so I have replaced the forward reference by this sentence.

The dissertation lifts this one-variable calculus to the orbits of several
commuting direction fields on a graph. The summation operator reappears
whenever a closed form is integrated along coordinate lines, as in the
discrete Poincar\'e lemma of \cref{sec:cartesian-grid}, and the fundamental
theorem \eqref{eq:ftc} is the statement that $\Delta$ has a right inverse,
which is why the first cohomology of a line vanishes (\cref{sec:cartesian-grid}).

\section{Graphs glued from graphs}
\label{sec:glued}

The aim is partly the same as Knill's, to develop calculus on discrete
spaces with graphs as the underlying geometric objects. In addition, a graph
is to be regarded as glued from smaller graphs in the way a manifold is
glued from charts, and this second requirement determines the shape of the
theory.

It also places the construction in a familiar framework. Much of modern
geometry is concerned with notions of space general enough to include
manifolds, algebraic varieties and schemes, and Grothendieck's formulation
is that a space is a sheaf on a site, determined by local pieces subject to
a descent condition \citep{vistoli2005}. Here the local pieces are
$\N$-sets, sets with $n$ commuting self-maps (\cref{sec:N-sets}), each
carrying a semi-cubical set whose $p$-cubes are a vertex together with $p$
distinct directions. The overlaps along which we glue are sub-cubical sets,
and the glued object is a semi-cubical set whose cubes form a poset with
the Alexandrov topology (\cref{sec:what-to-glue}). Differential forms are
functions on the points of this space, and the Mayer--Vietoris sequence
and the \v{C}ech--de Rham complex of \cref{sec:sheaf-cohomology} are
consequences of this sheaf-theoretic description.
%Claude note: the paragraph "Edit later. Once I have settled what the objects
%of the site are..." has been replaced by this description of what the body
%of the dissertation actually does. The claim about differentiable stacks
%is not made in the body and has been dropped here; see the Conclusion.

\section{Discrete derivatives in algebraic geometry}
\label{sec:alg-geom}

A need for a discrete theory of calculus also arises in algebraic geometry,
where differentiation enters through the Zariski cotangent space
\citep{vakil}. Let $X$ be a scheme, let $p\in X$, and let $\mathfrak m$ be
the maximal ideal of the local ring $\cO_{X,p}$. The Zariski cotangent
space at $p$ is the vector space $\mathfrak m/\mathfrak m^2$ over the
residue field $\kappa(p)=\cO_{X,p}/\mathfrak m$. The definition is modelled
on the smooth case. For a manifold $M$ with sheaf $\cO$ of smooth functions
the same quotient recovers the cotangent space $T_p^*M$, and passing to the
quotient by $\mathfrak m^2$ discards second-order information and retains
the first-order part.

The algebraic content of the definition is the Leibniz rule.

\begin{definition}[Derivation]
Let $\varphi\colon R\to A$ be a ring map and let $M$ be an $A$-module. An
$R$-derivation of $A$ into $M$ is an additive map $D\colon A\to M$ which
annihilates $\varphi(R)$ and satisfies the Leibniz rule
$D(ab)=a\,D(b)+b\,D(a)$.
\end{definition}

A derivation $D\colon A\to M$ satisfies $D(\mathfrak a^2)\subseteq\mathfrak
a M$ for every ideal $\mathfrak a$ of $A$, by the Leibniz rule. If
$\mathfrak a$ annihilates $M$, as $\mathfrak m$ does when $M$ is the residue
field, then $D$ vanishes on $\mathfrak m^2$ and its restriction to
$\mathfrak m$ factors through $\mathfrak m/\mathfrak m^2$. The Zariski
cotangent space is thus the object representing derivations into the
residue field. The classical derivative is defined analytically, as a
limit, but the Leibniz rule that governs it is purely algebraic, and this is
what the definition retains.

In a discrete setting the analogue of the derivative is a difference
operator, and a difference operator does not satisfy the Leibniz rule. It
satisfies a twisted form of it, the twist being by the shift that defines
the operator \citep{dmh1994a,dmh1994b}. We prove this as
\cref{prop:twisted-leibniz}. Nothing about schemes is inherently
continuous, $\Spec\Z$ and schemes over finite fields being examples, and
yet the definition of the cotangent space is imported from the continuous
case. In \cref{sec:geometric-pullback} we show that for functions on a
finite set the quotient $\mathfrak m/\mathfrak m^2$ is zero, because there
is no second-order information to discard, and the cotangent space of an
$\N$-set has to be defined differently.

\chapter{Coordinate systems and \texorpdfstring{$\N$}{N}-sets}
\label{ch:coordinates}

Cohomology connects the continuous and the discrete. For a topological
space that admits a triangulation, simplicial and singular cohomology
agree, and for a smooth manifold de Rham's theorem identifies the
cohomology of the complex of smooth differential forms with singular
cohomology \citep{bott-tu1982,hatcher2002}. The same invariant is thus
computed once from an analytic object and once from a combinatorial one.

If a graph $G$ is regarded as a discrete model of a manifold $M$, a
cohomology attached to $G$ may be expected to agree with that of $M$. The
cycle graph on $N$ vertices should have the cohomology of the circle and
the discrete torus $C_M\times C_N$ that of $T^2$. This is verified for the
construction below in \cref{sec:first-computations}.

This chapter defines the de Rham complex of a graph with a coordinate
system. Existing cohomology theories of graphs are defined in terms of
combinatorial structure, cliques or paths. Here a graph is treated as a
$\textit{space}$, and the theory is required to have the properties of a cohomology
theory of spaces, homotopy invariance and a Mayer--Vietoris sequence
\citep{brown1962}. These are proved in \cref{ch:homotopy} and
\cref{ch:gluing}.

We first recall the simplicial cohomology of a graph.

\section{Graphs as simplicial complexes}
\label{sec:simplicial}

Throughout this section $G=(V,E)$ is a simple graph, so that two vertices
are joined by at most one edge, and $\F=\{f\colon V\to\Z\}$ is the ring of
integer-valued functions on the vertices. When $V$ is finite the functions
with values in a ring $R$, denoted $\F_R$, are $\F\otimes_\Z R$.

A straightforward way to view $G$ as a space is to regard it as a
one-dimensional simplicial complex, with the vertices as $0$-cells and the
edges as $1$-cells. Choosing an orientation of each edge, that is, regarding
$G$ as a directed graph and writing $e^-$ and $e^+$ for the tail and head
of $e$, we obtain the cochain complex dual to the simplicial chain complex,
\[
  0\to\F_V\xrightarrow{\ d\ }\F_E\to 0,\qquad (df)(e)=f(e^+)-f(e^-),
\]
where $\F_V$ and $\F_E$ are the functions on the vertices and on the edges.
This is the classical setting of Eckmann~\citep{eckmann1945}.

\begin{proposition}\label{prop:simplicial-H0}
Let $G$ be a finite graph with $c$ connected components. Then
$H^0(G)\cong\Z^c$ and $H^1(G)\cong\Z^{|E|-|V|+c}$.
\end{proposition}

\begin{proof}
A function $f$ lies in $\ker d$ precisely when $f(u)=f(v)$ for every edge
$uv$, that is, when $f$ is constant on each connected component. Such
functions are freely determined by their $c$ values, one on each component,
so $H^0(G)\cong\Z^c$.

For $H^1$, let $k$ be any field. Tensoring with $k$ gives the cochain
complex $0\to\F_V\otimes k\to\F_E\otimes k\to 0$ of $k$-valued functions,
and the same argument gives $\ker d\cong k^c$, hence
$\im d\cong k^{|V|-c}$ and $H^1(G;k)\cong k^{|E|-|V|+c}$. Since the
complex consists of free abelian groups, the universal coefficient theorem
gives
\[
  H^0(G;k)\cong \bigl(H^0(G)\otimes k\bigr)\oplus\Tor_1^\Z\bigl(H^1(G),k\bigr).
\]
The left-hand side is $k^c\cong H^0(G)\otimes k$, so
$\Tor_1^\Z(H^1(G),k)=0$ for every field $k$. If $H^1(G)$ had a summand
$\Z/p$ then $\Tor_1^\Z(\Z/p,\Z/p)=\Z/p\neq 0$ would contradict this with
$k=\Z/p$. Hence $H^1(G)$ is a finitely generated torsion-free abelian group,
and its rank is $\dim_\Q H^1(G;\Q)=|E|-|V|+c$.
\end{proof}

\begin{figure}[htb]
\centering
\resizebox{0.72\textwidth}{!}{\input{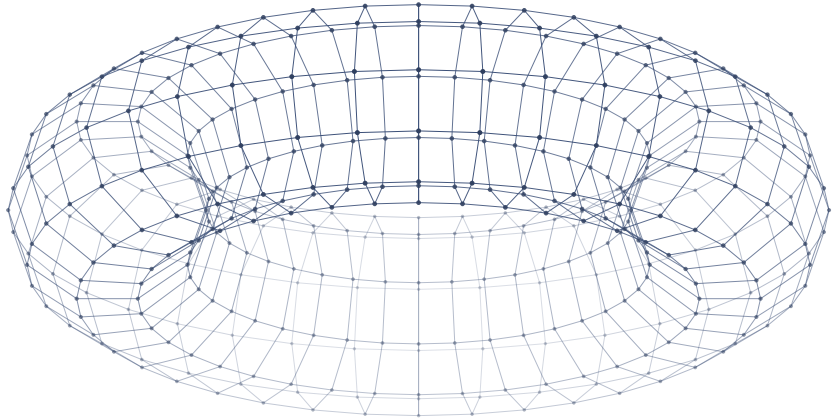}}
\caption{The torus graph $C_{40}\times C_{10}$. As a one-dimensional
simplicial complex it has no cohomology above degree one.}
\label{fig:torus40x10}
\end{figure}

A graph regarded in this way is a one-dimensional complex and has no
cohomology above degree one. The torus graph of \cref{fig:torus40x10} is a
discrete model of a surface, but as a simplicial complex it has $400$
vertices and $800$ edges, so $H^1$ has rank $401$ and $H^2=0$.

The rest of the chapter constructs a complex in which a graph can carry
cohomology in higher degrees.

\section{A de Rham complex from coordinates}
\label{sec:deRham}

For a smooth manifold $M$ the de Rham complex is
\[
  0\to\Omega^0(M)\to\Omega^1(M)\to\dots\to\Omega^p(M)\to\cdots,
\]
where $\Omega^p(M)$ is the space of smooth sections of $\Lam^pT^*M$. It is
a complex because partial derivatives commute. In coordinates
$x_1,\dots,x_n$ on a chart,
\[
  df=\sum_i\frac{\partial f}{\partial x_i}\,dx_i
  \quad\Longrightarrow\quad
  d^2f=\sum_{i,j}\frac{\partial^2 f}{\partial x_i\,\partial x_j}\,
  dx_i\wedge dx_j=0,
\]
since $\partial^2/\partial x_i\partial x_j$ is symmetric in $i,j$ while
$dx_i\wedge dx_j$ is antisymmetric.

The algebraic input is therefore a family of coordinates together with
pairwise commuting derivative operators. The rest of the construction is
the exterior algebra.

\subsection{Coordinate fields}
\label{sec:coordinate-fields}

For a vertex $v\in V$ we write $N(v)=\{w\in V:\ vw\in E\}$ for its
neighbourhood and $N[v]=N(v)\cup\{v\}$ for its closed neighbourhood. A
coordinate field is a choice of direction at each vertex.

\begin{definition}[Coordinate field]\label{def:coordinate-field}
A \emph{coordinate field} on $G$ is a map $X\colon V\to V$ with
$X(v)\in N[v]$ for every $v\in V$. It is \emph{strict} if $X(v)\in N(v)$
for every non-isolated vertex $v$. The set of coordinate fields on $G$ is
denoted $\Coord(G)$.
\end{definition}

\begin{example}\label{ex:fields}
On the cycle graph $C_N$ with vertex set $\Z/N$ the rotation $X(i)=i+1$ is
a strict coordinate field. On the path graph $1-2-3$ the map $X(1)=2$,
$X(2)=3$, $X(3)=3$ is a coordinate field which is not strict, since it
fixes the endpoint $3$. On the polar grid of \cref{fig:polar-grid}, a
discrete disc with concentric rings and radial spokes, the map $R$ that
moves each vertex one step towards the origin and fixes the origin is a
coordinate field, as is the rotation $\Theta$ of each ring. All three
examples are worked out in \cref{sec:first-computations}. The fixed points
in the second and third examples affect the cohomology, see
\cref{sec:polar-grid}.
\end{example}

Each coordinate field $X$ determines a $\Z$-algebra endomorphism
\[
  \psi_X\colon\F\to\F,\qquad \psi_X(f)=f\circ X,
\]
the pullback of functions along $X$. We call $\psi_X$ the \emph{shift} along
$X$.

\begin{proposition}\label{prop:shift-injective}
The map $\Hom_{\Set}(V,V)\to\End_{\Z\text{-alg}}(\F)$, $X\mapsto\psi_X$, is
injective.
\end{proposition}

\begin{proof}
Suppose $X_1\neq X_2$ and choose $v$ with $X_1(v)\neq X_2(v)$. Let $f$ be
the indicator function of the vertex $X_1(v)$. Then
$\psi_{X_1}(f)(v)=f(X_1(v))=1$ while $\psi_{X_2}(f)(v)=f(X_2(v))=0$.
\end{proof}

\subsection{Difference operators}
\label{sec:difference-operators}

A derivative measures the change of a function in a direction. In the
discrete setting the analogue is the difference between the value of a
function at a vertex and its value at the vertex the direction points to.
This idea goes back to the calculus of finite differences
\citep{boole1860,jordan1965}, and it organizes the noncommutative
differential calculus of Dimakis and M\"uller-Hoissen
\citep{dmh1994a,dmh1994b,dmh1999}, where first-order differential calculi on
a discrete set are shown to correspond to directed graphs on that set. A
first-order calculus on a discrete set requires a choice of directions
between its points, and such a choice is a directed graph.

\begin{definition}[Difference operator]\label{def:difference-operator}
Let $X\in\Coord(G)$. The \emph{difference operator} in the direction $X$ is
the homomorphism of abelian groups
\[
  \Delta_X\colon\F\to\F,\qquad \Delta_X=\psi_X-\id,
\]
so that $(\Delta_Xf)(v)=f(X(v))-f(v)$.
\end{definition}

The operator $\Delta_X$ is not a derivation. It satisfies a Leibniz rule
twisted by the shift.

\begin{proposition}[Twisted Leibniz rule]\label{prop:twisted-leibniz}
For all $f,g\in\F$ and all $X\in\Coord(G)$,
\[
  \Delta_X(fg)=\psi_X(f)\,\Delta_X(g)+\Delta_X(f)\,g.
\]
\end{proposition}

\begin{proof}
At a vertex $v$ the left-hand side is $f(X(v))g(X(v))-f(v)g(v)$. Inserting
$\pm f(X(v))g(v)$ gives
$f(X(v))\bigl(g(X(v))-g(v)\bigr)+\bigl(f(X(v))-f(v)\bigr)g(v)$.
\end{proof}

The twisted Leibniz rule is the ordinary Leibniz rule for a twisted
bimodule. Write $A=\F$ and let $Ae$ be the free $A$-module of rank one on a
symbol $e$, with the two actions
\[
  a_1\cdot(a_2e)=\psi_X(a_1)a_2\,e,\qquad (a_1e)\cdot a_2=a_1a_2\,e.
\]
These commute with one another, so $Ae$ is an $A$-$A$-bimodule in which the
left action is twisted by $\psi_X$. Setting $\delta_X(f)=\Delta_X(f)\,e$,
\cref{prop:twisted-leibniz} says that $\delta_X\in\Der(A,Ae)$. In
\cref{sec:algebra-structure} the differential forms on a graph are made
into a differential graded algebra by a twist of the same kind.

There is a standard ring that packages this information. If
$\alpha\colon R\to R$ is a ring homomorphism then $\alpha-\id$ is an
$\alpha$-derivation, and one can form the Ore extension, or skew polynomial
ring, $R[t;\alpha,\alpha-\id]$. It is the free left $R$-module on
$1,t,t^2,\dots$ subject to $tr=\alpha(r)t+\alpha(r)-r$, and it is the ring
of difference operators in one variable. Ore extensions were introduced in
Ore~\citep{ore1933}, and McConnell and Robson~\citep[Ch.~1]{mcconnell-robson2001} is the standard
modern treatment. This machinery is not needed below.

For $d^2=0$ we need the difference operators of different coordinate fields
to commute.

\begin{proposition}\label{prop:commutator}
For $X,Y\in\Coord(G)$ and $f\in\F$,
\[
  [\Delta_Y,\Delta_X](f)=(\Delta_Y\Delta_X-\Delta_X\Delta_Y)(f)
  =\psi_{X\circ Y}(f)-\psi_{Y\circ X}(f),
\]
that is, $[\Delta_Y,\Delta_X](f)(v)=f(X(Y(v)))-f(Y(X(v)))$.
\end{proposition}

\begin{proof}
Expand $(\psi_Y-\id)(\psi_X-\id)-(\psi_X-\id)(\psi_Y-\id)
=\psi_Y\psi_X-\psi_X\psi_Y$ and note that
$\psi_Y\psi_X(f)=f\circ X\circ Y=\psi_{X\circ Y}(f)$.
\end{proof}

\begin{corollary}\label{cor:commute}
$[\Delta_Y,\Delta_X]=0$ as operators on $\F$ if and only if
$X\circ Y=Y\circ X$ as maps $V\to V$.
\end{corollary}

\begin{proof}
By \cref{prop:commutator} the commutator is $\psi_{X\circ Y}-\psi_{Y\circ X}$,
and by \cref{prop:shift-injective} this vanishes if and only if
$X\circ Y=Y\circ X$.
\end{proof}

The difference operators commute exactly when the coordinate fields commute
as maps of the vertex set, and this is the condition in the following
definition.

\subsection{Coordinate systems and forms}
\label{sec:coordinate-systems}

\begin{definition}[Coordinate system]\label{def:coordinate-system}
A \emph{coordinate system} on $G$ is a finite family
$\Xs=\{X_1,\dots,X_n\}\subseteq\Coord(G)$ of pairwise commuting coordinate
fields, so that $X_i\circ X_j=X_j\circ X_i$ for all $i,j$. The integer $n$
is the \emph{rank} of $\Xs$. We abbreviate $S_i=\psi_{X_i}$ and
$\Delta_i=\Delta_{X_i}=S_i-\id$.
\end{definition}

Equivalently, a coordinate system of rank $n$ is an action of the free
commutative monoid $\N^n$ on $V$ in which the generators act by coordinate
fields. When the $X_i$ are invertible it is an action of $\Z^n$. This
reformulation is the subject of \cref{sec:N-sets}, and from
\cref{ch:homotopy} onwards it is the one we use.

Recall the smooth construction of forms \citep[\S1]{bott-tu1982}. On $\R^n$
one lets $\Lam^\bullet$ be the exterior algebra over $\R$ on symbols
$dx_1,\dots,dx_n$, with basis the products $dx_{i_1}\cdots dx_{i_p}$ for
$i_1<\dots<i_p$. The smooth $p$-forms are
$\Omega^p(\R^n)=C^\infty(\R^n)\otimes_\R\Lam^p(\R^n)^*$, and $d$ is
defined on functions by $df=\sum_i(\partial f/\partial x_i)\,dx_i$ and
extended coefficientwise. We repeat the construction over $\Z$.

Fix a coordinate system $\Xs$ of rank $n$ and let
$\Lam^\bullet\Z^n=\bigoplus_{p=0}^n\Lam^p\Z^n$ be the exterior algebra on a
free abelian group of rank $n$, whose generators we name $dX_1,\dots,dX_n$
to match the smooth notation. For $I=\{i_1<\dots<i_p\}\subseteq[n]$ we write
$dX_I=dX_{i_1}\wedge\dots\wedge dX_{i_p}$.

\begin{definition}[Coordinate $p$-forms]\label{def:forms}
The \emph{coordinate $p$-forms} associated with $\Xs$ form the $\F$-module
\[
  \Omega^p(G,\Xs)=\F\otimes_\Z\Lam^p\Z^n\cong\Lam^p_{\F}(\F^n),
\]
free of rank $\binom{n}{p}$ over $\F$ with basis $\{dX_I:\ |I|=p\}$. We
write $\Omega^\bullet(G,\Xs)=\bigoplus_{p=0}^n\Omega^p(G,\Xs)$, and
$\Omega^p(G,\Xs;R)=\F_R\otimes_\Z\Lam^p\Z^n$ when we want coefficients in a
commutative ring $R$.
\end{definition}

\begin{definition}[Exterior derivative]\label{def:d}
The \emph{exterior derivative} is the $\Z$-linear map
$d\colon\Omega^p(G,\Xs)\to\Omega^{p+1}(G,\Xs)$ given by
$d=\sum_{j=1}^n\Delta_j\otimes(dX_j\wedge-)$, that is,
\[
  d\Bigl(\sum_{|I|=p}f_I\,dX_I\Bigr)=\sum_{|I|=p}\sum_{j=1}^n
  \Delta_j(f_I)\,dX_j\wedge dX_I.
\]
\end{definition}

\begin{proposition}\label{prop:d2}
$d^2=0$.
\end{proposition}

\begin{proof}
For $\omega=f\,dX_I$,
\[
  d^2\omega=\sum_{k=1}^n\sum_{j=1}^n\Delta_k\Delta_j(f)\,
  dX_k\wedge dX_j\wedge dX_I.
\]
By \cref{cor:commute} the coefficient $\Delta_k\Delta_j(f)$ is symmetric in
$j$ and $k$, while $dX_k\wedge dX_j$ is antisymmetric. The off-diagonal terms
cancel in pairs and the diagonal terms vanish because $dX_j\wedge dX_j=0$.
\end{proof}

\begin{definition}[Discrete de Rham complex]\label{def:deRham}
For a graph $G$ with a coordinate system $\Xs$ of rank $n$, the
\emph{discrete de Rham complex} is the cochain complex
\[
  0\to\Omega^0(G,\Xs)\xrightarrow{d}\Omega^1(G,\Xs)\xrightarrow{d}\dots
  \xrightarrow{d}\Omega^n(G,\Xs)\to 0,
\]
and its cohomology is written $\HdR{p}(G,\Xs)$, or $\HdR{p}(G,\Xs;R)$ when
the coefficients are taken in $R$.
\end{definition}

The complex has length the rank of the coordinate system, which is not an
invariant of $G$. If $V$ is finite then $\Omega^p(G,\Xs)$ is a free
$\Z$-module of rank $|V|\binom{n}{p}$, with the following consequence.

\begin{theorem}[Euler rigidity]\label{thm:euler}
Let $G$ be a finite graph with a coordinate system $\Xs$ of rank $n\geq 1$,
and let $k$ be a field. Then the Euler characteristic of the de Rham complex
vanishes,
\[
  \chi_{\dR}(G,\Xs)=\sum_{p=0}^n(-1)^p\dim_k\HdR{p}(G,\Xs;k)=0.
\]
\end{theorem}

\begin{proof}
The Euler characteristic of a finite complex of finite-dimensional vector
spaces can be computed from the complex itself,
\[
  \sum_{p=0}^n(-1)^p\dim_k\HdR{p}(G,\Xs;k)
  =\sum_{p=0}^n(-1)^p\dim_k\Omega^p(G,\Xs;k)
  =|V|\sum_{p=0}^n(-1)^p\binom{n}{p}=0.\qedhere
\]
\end{proof}

In particular no finite graph with a coordinate system has the cohomology
$(\Z,0,\dots,0)$ of a point. \Cref{sec:toral-rigidity} sharpens this: the
cohomology of a finite graph with a coordinate system of rank $n$ is that
of a disjoint union of $n$-tori. Other cohomologies require several
coordinate systems on the same graph, which is the subject of
\cref{ch:gluing}.

\subsection{Cayley graphs}
\label{sec:cayley}

We next determine which graphs admit coordinate systems without fixed
points or collapsing. Call a coordinate system $\Xs$ \emph{bijective} if
every $X_i$ is a bijection with no fixed points. For a bijective coordinate
system the action of $\N^n$ extends to an action of $\Z^n$.

\begin{definition}[Homogeneous coordinate system]\label{def:homogeneous}
A coordinate system $\Xs=\{X_1,\dots,X_n\}$ on $G=(V,E)$ is
\emph{homogeneous} if it is bijective and strict, and if every edge of $G$
is of the form $\{v,X_i(v)\}$ for some $v\in V$ and some $i$.
\end{definition}
%Claude note: the draft defined "homogeneous" as bijective and strict only,
%and stated the theorem below without the last clause. The theorem is false
%without it: a single fixed-point-free involution on K_4 (a perfect matching)
%is bijective and strict, K_4 is connected, but K_4 is not a Cayley graph on
%one generator. The proof in the draft uses the clause ("the edge map is
%bijective because the coordinate system covers the graph"), so I have made
%it part of the definition. Please check that this is what you intend.

\begin{theorem}\label{thm:cayley}
A finite connected graph $G=(V,E)$ admits a homogeneous coordinate system
of rank $n$ if and only if it is isomorphic to the Cayley graph
$\Cay(A,\{a_1,\dots,a_n\})$ of a finite abelian group $A$ with respect to a
generating set of $n$ elements, none of which is the identity.
\end{theorem}

\begin{proof}
Suppose $G=\Cay(A,\{a_1,\dots,a_n\})$. The maps $X_i\colon A\to A$,
$X_i(a)=a+a_i$, are bijective, have no fixed points since $a_i\neq 0$, and
commute because $A$ is abelian. Every edge $\{a,a+a_i\}$ is of the required
form, and the graph is connected because the $a_i$ generate.

Conversely let $\Xs$ be a homogeneous coordinate system on the connected
graph $G$. Let $A$ be the subgroup of the symmetric group $\Sym(V)$
generated by $X_1,\dots,X_n$. It is abelian since the $X_i$ commute. We
claim that $A$ acts freely on $V$. Suppose $g\in A$ fixes some $v\in V$. Any
$w\in V$ can be written $w=h(v)$ for some $h\in A$, because $G$ is connected
and every edge is of the form $\{u,X_i(u)\}$ with $X_i$ invertible. Then
$g(w)=gh(v)=hg(v)=h(v)=w$ by commutativity, so $g=\id$.

Fix $v_0\in V$ and define $\phi\colon\Cay(A,\{X_1,\dots,X_n\})\to G$ on
vertices by $\phi(g)=g(v_0)$ and on edges by
$\{g,gX_i\}\mapsto\{g(v_0),X_i(g(v_0))\}$. The vertex map is injective
because the action is free and surjective because $G$ is connected. The
edge map is well defined since $\{g(v_0),X_i(g(v_0))\}$ is an edge of $G$,
and it is bijective onto $E$ because every edge of $G$ has the form
$\{u,X_i(u)\}$ with $u=g(v_0)$ for a unique $g$. Hence $\phi$ is an
isomorphism of graphs, and it carries the coordinate field $X_i$ of the
Cayley graph to the coordinate field $X_i$ of $G$.
\end{proof}

Thus the graphs carrying a global homogeneous coordinate system are the
Cayley graphs of finite abelian groups, and by \cref{cor:toral} their de
Rham cohomology is that of an $n$-torus, independently of the group $A$.
Graphs that are not of this form are treated by covering them with several
charts in \cref{ch:gluing}.

\subsection{The generalized calculus philosophy}
\label{sec:philosophy}

%Claude note: rewritten and shortened as your comment asked.

\begin{informal}{Philosophy}
The definitions are guided by the view that a coordinate chart on a space
is a region on which a commutative group or monoid acts, with points
labelled by the action and coordinate lines given by the orbits of the
generators. On a graph the acting monoid is $\N^n$, generated by the
commuting coordinate fields, and a chart is a subgraph on which $\N^n$ acts. The philosophy
is that to describe a space, we can identify local regions of symmetry, which can be formally described
as the orbits of a certain action, giving rise to charts. In \cref{ch:gluing} a graph is covered by such charts, as a manifold is
covered by open subsets of $\R^n$. 
\end{informal}

\subsection{The algebra structure}
\label{sec:algebra-structure}

The difference operators $\Delta_i$ do not satisfy the Leibniz rule, and
$\Omega^\bullet(G,\Xs)$ with the $\F$-bilinear product is not a differential
graded algebra. We recall the definitions and then describe the product that
makes it one.

\begin{definition}[Graded algebra]\label{def:graded-algebra}
A \emph{graded $R$-algebra} is a family $\{K^p\}_{p\geq 0}$ of $R$-modules
together with bilinear products $K^p\times K^q\to K^{p+q}$ and an element
$1\in K^0$ making $\bigoplus_pK^p$ an associative unital $R$-algebra. It is
\emph{graded commutative} if $a\cdot b=(-1)^{pq}b\cdot a$ for $a\in K^p$
and $b\in K^q$. A morphism of graded $R$-algebras is a morphism of
$R$-algebras which preserves the grading.
\end{definition}

\begin{definition}[Differential graded algebra]\label{def:dga}
A \emph{differential graded algebra}, or DGA, over $R$ is a graded
$R$-algebra with an additive map $d\colon K^p\to K^{p+1}$ satisfying
$d^2=0$ and $d(a\cdot b)=d(a)\cdot b+(-1)^pa\cdot d(b)$ for $a\in K^p$. A
morphism of DGAs is a morphism of graded algebras which is also a cochain
map.
\end{definition}

The obvious product to try is the $\F_R$-bilinear one,
\[
  \Omega^p(G,\Xs;R)\otimes_{\F_R}\Omega^q(G,\Xs;R)\to\Omega^{p+q}(G,\Xs;R),
  \qquad f\,dX_I\otimes g\,dX_J\mapsto fg\,dX_I\wedge dX_J,
\]
which makes $\Omega^\bullet(G,\Xs;R)$ an associative graded-commutative
$\F_R$-algebra with unit $1\in\Omega^0=\F_R$. Concretely, if
$\omega_1=\sum_{|I|=p}f_I\,dX_I$ and $\omega_2=\sum_{|J|=q}g_J\,dX_J$ then
the coefficient of $dX_K$ in $\omega_1\cdot\omega_2$ at the vertex $v$ is
$\sum\varepsilon(I,J)f_I(v)g_J(v)$, the sum over the pairs $(I,J)$ with
$I\sqcup J=K$, where $\varepsilon(I,J)$ is the sign of the permutation
sorting $I$ followed by $J$. This product does not satisfy the Leibniz rule.

\begin{proposition}\label{prop:leibniz-error}
For $\omega_1=f\,dX_I$ of degree $p$ and $\omega_2=g\,dX_J$,
\[
  d(\omega_1\cdot\omega_2)=d(\omega_1)\cdot\omega_2
  +(-1)^p\omega_1\cdot d(\omega_2)
  +\sum_{i=1}^n\Delta_i(f)\Delta_i(g)\,dX_i\wedge dX_I\wedge dX_J.
\]
\end{proposition}

\begin{proof}
By \cref{prop:twisted-leibniz},
$\Delta_i(fg)=\Delta_i(f)g+f\Delta_i(g)+\Delta_i(f)\Delta_i(g)$. Multiply
by $dX_i\wedge dX_I\wedge dX_J$ and sum over $i$. The first term gives
$d(\omega_1)\cdot\omega_2$, and the second gives
$(-1)^p\omega_1\cdot d(\omega_2)$ after moving $dX_i$ past $dX_I$.
\end{proof}

The error term is a product of two first differences and has no continuous
counterpart. As in \cref{sec:difference-operators}, it is removed by a
twist. For
$J=\{j_1<\dots<j_q\}\subseteq[n]$ write $S_J=S_{j_1}\circ\dots\circ S_{j_q}$
and $S_\emptyset=\id$. The order does not matter because the coordinate
fields commute.

\begin{definition}[Twisted product]\label{def:twisted-product}
The \emph{twisted product} on $\Omega^\bullet(G,\Xs;R)$ is the $\Z$-bilinear
map
\[
  \mu\colon\Omega^p(G,\Xs;R)\otimes_\Z\Omega^q(G,\Xs;R)\to\Omega^{p+q}(G,\Xs;R),
  \qquad f\,dX_I\otimes g\,dX_J\mapsto S_J(f)\,g\,dX_I\wedge dX_J.
\]
\end{definition}

The tensor product is over $\Z$ and not over $\F_R$ because the map is not
$\F_R$-bilinear. Pulling a function out of the left-hand factor shifts it.
The product is associative, since $S_K(S_J(f)g)h=S_{J\cup K}(f)S_K(g)h$
whenever $J\cap K=\emptyset$, and when $J\cap K\neq\emptyset$ both triple
products vanish because $dX_J\wedge dX_K=0$. With unit $1\in\Omega^0$ this
makes $\Omega^\bullet(G,\Xs;R)$ a graded $\Z$-algebra.

\begin{theorem}\label{thm:dga}
$(\Omega^\bullet(G,\Xs;R),\mu,d)$ is a differential graded $\Z$-algebra.
\end{theorem}

\begin{proof}
We check the Leibniz rule on generators $\omega_1=f\,dX_I$ of degree $p$
and $\omega_2=g\,dX_J$. By \cref{prop:twisted-leibniz},
\[
  d(\omega_1\cdot\omega_2)=\sum_{i=1}^n\Delta_i\bigl(S_J(f)g\bigr)\,
  dX_i\wedge dX_I\wedge dX_J,\qquad
  \Delta_i(S_J(f)g)=\Delta_i(S_J(f))\,g+S_iS_J(f)\,\Delta_i(g).
\]
On the other side,
$d(\omega_1)\cdot\omega_2=\sum_iS_J(\Delta_i(f))\,g\,dX_i\wedge dX_I\wedge dX_J$
and
$(-1)^p\omega_1\cdot d(\omega_2)=\sum_iS_{J\cup\{i\}}(f)\,\Delta_i(g)\,
dX_i\wedge dX_I\wedge dX_J$. The two expressions agree term by term. Since
$S_J$ commutes with $\Delta_i$ we have $\Delta_i(S_J(f))=S_J(\Delta_i(f))$,
and $S_iS_J=S_{J\cup\{i\}}$ whenever $i\notin J$, while if $i\in J$ then
$dX_i\wedge dX_J=0$ and both terms vanish.
\end{proof}

This algebra is not graded commutative, since
$f\,dX_I\cdot g\,dX_J=S_J(f)g\,dX_I\wedge dX_J$ while
$g\,dX_J\cdot f\,dX_I=S_I(g)f\,dX_J\wedge dX_I$. In \cref{sec:serre} the
twisted product is identified with the cup product on cubical cochains dual
to the Serre diagonal.

\section{Other cohomology theories for graphs}
\label{sec:other-theories}

\subsection{Clique cohomology}
\label{sec:clique}

The standard device for producing higher-dimensional cells from a graph,
and the one Knill uses throughout \citep{knill2011,knill2012a}, is the
clique complex.

\begin{definition}[Whitney complex]\label{def:whitney}
Let $G=(V,E)$ be a simple graph. The \emph{clique complex}, or Whitney
complex, of $G$ is the abstract simplicial complex on the vertex set $V$
whose $k$-simplices are the $(k+1)$-cliques of $G$, that is, the subsets of
$V$ of size $k+1$ which are pairwise adjacent.
\end{definition}

A $(k+1)$-clique is the $1$-skeleton of a $k$-simplex, so the clique
complex fills in every simplex whose $1$-skeleton is present in $G$. Knill's
curvature, a count of simplices in the unit sphere of a vertex, is defined
in terms of it.

The clique complex ties the dimension of the space to the presence of
complete subgraphs, and a graph modelling a higher-dimensional space need
not contain them. The discrete torus of \cref{fig:torus40x10} is
triangle-free, so its clique complex is its own $1$-skeleton and its clique
cohomology vanishes above degree one. A theory whose cells are simplices
does not detect the squares of a grid, and this is one reason for the
passage to cubes in \cref{ch:cubical}.

\begin{figure}[htb]
\centering
\resizebox{0.62\textwidth}{!}{\input{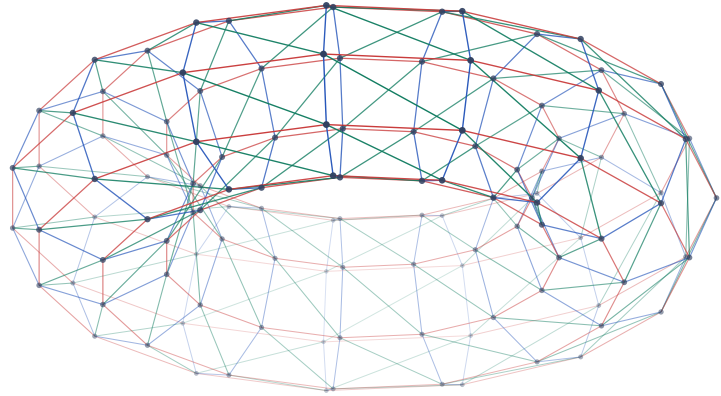}}
\caption{A triangulation of the torus realized as the Cayley graph
$\Cay(\Z/15\times\Z/8,\{(1,0),(0,1),(1,1)\})$. The three generators are
drawn in three colours.}
\label{fig:torus-tri}
\end{figure}

Conversely, the de Rham complex of a single coordinate system does not
detect simplices. Consider the triangulated torus of \cref{fig:torus-tri},
the Cayley graph $\Cay(\Z/n\times\Z/m,\{(1,0),(0,1),(1,1)\})$. Its cliques
are its vertices, edges and triangles, so its clique cohomology is the
simplicial cohomology of the triangulated torus, $(\Z,\Z^2,\Z)$. The
coordinate system given by the three generators is homogeneous of rank $3$,
and by \cref{cor:toral} its de Rham cohomology is that of the $3$-torus,
$(\Z,\Z^3,\Z^3,\Z)$.

The simplicial cohomology is recovered by an atlas of small charts in
place of a single global one. Each triangle carries a chart of rank $2$ in
which one edge of the corresponding square is collapsed
(\cref{sec:two-simplex}), and the atlas of all such charts has the
normalized cohomology of the $2$-torus (\cref{sec:triangulated-torus}).

\subsection{Path homology}
\label{sec:path-homology}

Path homology, introduced by Grigor'yan, Lin, Muranov and Yau
\citep{glmy2012}, is a homology theory for directed graphs built from
paths. An undirected graph is treated as the directed graph in which every
edge is replaced by two opposite arrows. We recall the definitions briefly.

%Claude note: the draft contained a screenshot of the definitions from
%another paper; I have written them out following GLMY.
Let $G=(V,E)$ be a digraph and $R$ a commutative ring. For $p\geq 0$ let
$\Lambda_p$ be the free $R$-module on the symbols $e_{i_0\cdots i_p}$ with
$i_q\in V$, the \emph{elementary $p$-paths}, with boundary
\[
  \partial e_{i_0\cdots i_p}=\sum_{q=0}^p(-1)^q e_{i_0\cdots\widehat{i_q}\cdots i_p},
\]
where the hat denotes omission. A path is \emph{regular} if
$i_q\neq i_{q+1}$ for all $q$, and one passes to the quotient $R_p$ of
$\Lambda_p$ by the non-regular paths, on which $\partial$ is still defined.
A regular path is \emph{allowed} if $i_q\to i_{q+1}$ is an edge for every
$q$, and $A_p\subseteq R_p$ is spanned by the allowed paths. The boundary of
an allowed path need not be allowed, so one sets
\[
  \Omega_p(G)=\{u\in A_p:\ \partial u\in A_{p-1}\},
\]
the $\partial$-invariant paths. Then $\partial\Omega_p\subseteq\Omega_{p-1}$
and the \emph{path homology} of $G$ is the homology of
$(\Omega_\bullet(G),\partial)$. Path cohomology is the cohomology of the
dual complex.

\begin{intuition}
An allowed $p$-path is a walk of length $p$ along the arrows, and its
boundary is the alternating sum of the walks obtained by omitting one
vertex. Omitting an interior vertex $i_q$ produces the step
$i_{q-1}\to i_{q+1}$, which need not be an edge, so the boundary of an
allowed path is in general not allowed. The $\partial$-invariant paths are
the linear combinations of allowed paths in which the non-allowed terms of
the boundary cancel. For a square $a\to b\to d$, $a\to c\to d$ the
combination $e_{abd}-e_{acd}$ is $\partial$-invariant, the two terms
$\pm e_{ad}$ cancelling. The cells detected by path homology are thus
determined by coincidences among walks.
\end{intuition}

Two features distinguish path homology from the theory developed here.
Path homology is an invariant of the digraph, while the cohomology defined
below depends on the coordinate system placed on the graph. The algebraic
setting also differs. The path cochains form a differential graded algebra
which is a quotient of the path algebra $RG$, while the de Rham complex is
a Koszul complex over a polynomial ring (\cref{sec:koszul}). The two
theories appear to agree in low degrees for graphs with few coincidences
among walks, but a comparison is not needed for what follows and is not
attempted. Path homology is in general \textbf{not} the homology of a space
naturally attached to the graph \citep{fu-ivanov2024}, whereas the
normalized cohomology of \cref{ch:cubical} is the cohomology of a cubical
set by construction, which is the singular cohomology of its geometric realization.

\subsection{Knill's exterior derivative}
\label{sec:knill}

%Claude note: written from Knill's papers as your comment requested; please
%check that this is a fair summary of the parts of his work you have in mind.
Knill's calculus on a finite simple graph is the simplicial calculus of its
Whitney complex \citep{knill2013,knill2014}. Functions on the $k$-cliques
are the $k$-forms, the exterior derivative is the simplicial coboundary
after a choice of orientations, and the Dirac operator $D=d+d^*$ and the
Hodge Laplacian $L=D^2$ give a Hodge theory whose harmonic forms compute
the cohomology of the clique complex. More recently Knill has studied
cohomology of open sets of finite topologies \citep{knill2023}. The theory
is combinatorial and intrinsic to the graph, and its cells are those of the
clique complex, with the limitation described in \cref{sec:clique}. It
recovers a large number of theorems of differential geometry, among them
Gauss--Bonnet, Poincar\'e--Hopf and the McKean--Singer formula.

\subsection{Discrete exterior calculus}
\label{sec:dec}

Discrete exterior calculus \citep{hirani2003,dhlm2005} is a numerical
theory. Its objects are simplicial complexes embedded in Euclidean space
together with their circumcentric duals, and its cochains are integrals of
smooth forms over primal and dual cells. The exterior derivative is again
the coboundary, and the theory provides a discrete Hodge star and a
discrete wedge product that converge to the smooth ones under refinement.
It deals with simplicial complexes rather than graphs and requires an
embedding to define the metric operators. The constructions in this
dissertation involve no metric and no embedding.

\section{The category of \texorpdfstring{$\N$}{N}-sets}
\label{sec:N-sets}

\subsection{What is an \texorpdfstring{$\N$}{N}-set?}
\label{sec:what-is-N-set}

A coordinate system $\Xs=(X_1,\dots,X_n)$ on a graph $G=(V,E)$ is an
action of the monoid $\N^n$ on the vertex set $V$, and for most of what
follows only the action is needed. Here $\N=\{0,1,2,\dots\}$, so that
$\N$ is a monoid.

\begin{definition}[$\N$-set]\label{def:N-set}
An \emph{$\N$-set} is a pair $(X,\N^n)$ consisting of a set $X$ with an
action $\N^n\times X\to X$ of the monoid $\N^n$, so that $0\cdot x=x$ and
$(a+b)\cdot x=a\cdot(b\cdot x)$. We write $e_1,\dots,e_n$ for the standard
generators of $\N^n$.
\end{definition}

Giving an action of $\N^n$ on $X$ is the same as giving $n$ pairwise
commuting maps $X\to X$, the actions $x\mapsto e_i\cdot x$ of the
generators, and this is the reformulation we use throughout. We write
$\F_X=\Fun(X,\Z)$ for the functions on $X$, $S_i(f)(x)=f(e_i\cdot x)$ for
the shifts, $S_I=S_{i_1}\circ\dots\circ S_{i_p}$ for $I\subseteq[n]$, and
$\Delta_i=S_i-\id$ for the difference operators, exactly as before. Every
$\N$-set has an \emph{induced directed graph}, with vertex set $X$ and an
arrow $x\to e_i\cdot x$ for every $x$ and every $i$, and forgetting the
directions gives the \emph{induced graph}. A coordinate system on $G$ is an
$\N$-set whose induced graph is a subgraph of $G$ together with loops at the
fixed points.

\begin{definition}[Morphisms of $\N$-sets]\label{def:N-morphism}
A \emph{morphism} $F\colon(X,\N^n)\to(Y,\N^m)$ of $\N$-sets is a pair
$(f,\phi)$ of a set map $f\colon X\to Y$ and a monoid homomorphism
$\phi\colon\N^n\to\N^m$ such that $f(a\cdot x)=\phi(a)\cdot f(x)$.
Composition is componentwise and the identity is $(\id_X,\id_{\N^n})$. This
defines the category $\Nset$.
\end{definition}

\begin{proposition}\label{prop:products}
$\Nset$ has finite products. The product of $(X,\N^n)$ and $(Y,\N^m)$ is
$(X\times Y,\N^{n+m})$ with the action $(a,b)\cdot(x,y)=(a\cdot x,b\cdot y)$,
the projections being $(\pi_X,p)$ and $(\pi_Y,q)$ for the two coordinate
projections $p,q$ of $\N^{n+m}$. The terminal object is $(*,\N^0)$.
\end{proposition}

\begin{proof}
Given $(g,\gamma)\colon(Z,\N^l)\to(X,\N^n)$ and
$(h,\eta)\colon(Z,\N^l)\to(Y,\N^m)$, the pair $((g,h),(\gamma,\eta))$ is the
unique morphism to $(X\times Y,\N^{n+m})$ compatible with the projections,
and it is a morphism of $\N$-sets because
$(g,h)(a\cdot z)=(\gamma(a)\cdot g(z),\eta(a)\cdot h(z))$.
\end{proof}

\subsection{The de Rham complex of an \texorpdfstring{$\N$}{N}-set}
\label{sec:deRham-N-set}

The de Rham complex of an $\N$-set $(X,\N^n)$ is defined exactly as for a
coordinatized graph,
\[
  \Omega^p(X,\N^n)=\F_X\otimes_\Z\Lam^p\Z^n,\qquad
  d=\sum_{j=1}^n\Delta_j\otimes(dX_j\wedge-),
\]
and it is a complex because the monoid is commutative. Its cohomology is
again denoted $\HdR{p}(X,\N^n)$, and the twisted product of
\cref{def:twisted-product} makes $\Omega^\bullet(X,\N^n)$ a DGA.

We show that $\Omega^\bullet$ takes finite products of finite $\N$-sets to
tensor products of DGAs. Recall the tensor product of DGAs
\citep[Definition 22.3.4]{stacks}.
%Claude note: please check this section number against the current Stacks
%project (it moves as the project is edited); citing the tag is safer.

\begin{definition}[Tensor product of DGAs]\label{def:tensor-dga}
Let $(A,d)$ and $(B,d)$ be differential graded algebras over a commutative
ring $R$. Their tensor product is the graded algebra $A\otimes_RB$ with
$(A\otimes B)^k=\bigoplus_{i+j=k}A^i\otimes B^j$, multiplication
$(a\otimes b)(a'\otimes b')=(-1)^{\deg a'\deg b}\,aa'\otimes bb'$, and
differential $d(a\otimes b)=d(a)\otimes b+(-1)^{\deg a}a\otimes d(b)$.
\end{definition}

The underlying cochain complex of $A\otimes_RB$ is the total complex of the
tensor product of the complexes $A$ and $B$.

\begin{proposition}\label{prop:product-dga}
Let $(X,\N^n)$ and $(Y,\N^m)$ be finite $\N$-sets. The map
\[
  \Omega^\bullet(X,\N^n)\otimes_\Z\Omega^\bullet(Y,\N^m)\to
  \Omega^\bullet(X\times Y,\N^{n+m}),\qquad
  f\,dX_I\otimes g\,dY_J\mapsto (f\otimes g)\,dZ_{I\cup J[n]},
\]
where $(f\otimes g)(x,y)=f(x)g(y)$ and $J[n]=\{j+n:\ j\in J\}$, is an
isomorphism of differential graded algebras. Consequently
$\Omega^\bullet(\prod_{i=1}^k(X_i,\N^{n_i}))\cong
\bigotimes_{i=1}^k\Omega^\bullet(X_i,\N^{n_i})$ for finite $\N$-sets
$(X_i,\N^{n_i})$.
\end{proposition}

\begin{proof}
The map $\kappa\colon\F_X\otimes_\Z\F_Y\to\F_{X\times Y}$,
$\kappa(f\otimes g)(x,y)=f(x)g(y)$, is an isomorphism because $X$ and $Y$
are finite. Both sides are free of rank $|X|\cdot|Y|$ and $\kappa$ takes the
basis $\delta_x\otimes\delta_y$ of indicator functions to the basis
$\delta_{(x,y)}$. Also
$\bigoplus_{a+b=p}\Lam^a\Z^n\otimes\Lam^b\Z^m\cong\Lam^p\Z^{n+m}$ via
$dX_I\otimes dY_J\mapsto dZ_{I\cup J[n]}$, which is a bijection on bases.
Tensoring the two isomorphisms gives the stated $\Z$-module isomorphism in
each degree. That it respects the products and the differentials is a
direct expansion of both sides, carried out in \cref{app:product}. The
statement for $k$ factors follows by induction, using the associativity of
the tensor product of DGAs.
\end{proof}

The assignment $(X,\N^n)\mapsto\Omega^\bullet(X,\N^n)$ is not a functor
on all of $\Nset$. The next section identifies the morphisms along which
forms can be pulled back.

\subsection{Functoriality of differential forms}
\label{sec:functoriality}

We ask whether $\Omega^\bullet$ extends to a contravariant functor
$\Nset^\op\to\DGA$. Let $F=(f,\phi)\colon(X,\N^n)\to(Y,\N^m)$ be a
morphism of $\N$-sets. On zero-forms there is an evident pullback,
$F^*\colon\Omega^0(Y,\N^m)\to\Omega^0(X,\N^n)$, $g\mapsto g\circ f$. For
$a=(a_1,\dots,a_n)\in\N^n$ write $S^X_a=(S^X_1)^{a_1}\circ\dots\circ(S^X_n)^{a_n}$
for the pullback along $x\mapsto a\cdot x$.

\begin{proposition}\label{prop:pullback-shift}
$S^X_i\circ F^*=F^*\circ S^Y_{\phi(e_i)}$ and hence
$\Delta^X_i\circ F^*=F^*\circ(S^Y_{\phi(e_i)}-\id)$ as maps
$\Omega^0(Y,\N^m)\to\Omega^0(X,\N^n)$.
\end{proposition}

\begin{proof}
$(S^X_iF^*g)(x)=g(f(e_i\cdot x))=g(\phi(e_i)\cdot f(x))=(F^*S^Y_{\phi(e_i)}g)(x)$,
because $\phi$ is a monoid map and $F$ is a morphism of $\N$-sets.
\end{proof}

It remains to define $F^*$ on one-forms. We first determine the algebraic
constraints and give the geometric interpretation in
\cref{sec:geometric-pullback}. Suppose $F^*(dY_j)=\sum_ih_{ji}\,dX_i$ with
$h_{ji}\in\F_X$. Since
$g\,dY_j=dY_j\cdot g$ and $F^*$ is to be multiplicative, we must have
$F^*(g\,dY_j)=F^*(g)\sum_ih_{ji}\,dX_i$. For $F^*$ to be a cochain map
the square
\[
\begin{tikzcd}[column sep=large]
  \Omega^0(Y,\N^m)\arrow[r,"d_Y"]\arrow[d,"F^*"'] & \Omega^1(Y,\N^m)\arrow[d,"F^*"]\\
  \Omega^0(X,\N^n)\arrow[r,"d_X"] & \Omega^1(X,\N^n)
\end{tikzcd}
\]
must commute, that is, $F^*(d_Yg)=d_X(F^*g)$ for every $g\in\F_Y$. The
left-hand side is
\[
  F^*\Bigl(\sum_{j=1}^m\Delta^Y_j(g)\,dY_j\Bigr)
  =\sum_{i=1}^n\Bigl(\sum_{j=1}^mh_{ji}\,F^*\Delta^Y_j(g)\Bigr)dX_i,
\]
and the right-hand side is $\sum_i\Delta^X_i(F^*g)\,dX_i$. Comparing
coefficients of the basis $dX_1,\dots,dX_n$, the condition is
\begin{equation}\label{eq:chain-condition}
  \Delta^X_i(F^*g)=\sum_{j=1}^mh_{ji}\,F^*\bigl(\Delta^Y_j(g)\bigr)
  \qquad\text{for all }i\in[n]\text{ and }g\in\F_Y.
\end{equation}
By \cref{prop:pullback-shift} the left-hand side is
$F^*\bigl((S^Y_{\phi(e_i)}-\id)g\bigr)$. If $\phi(e_i)=0$ then
$S^Y_{\phi(e_i)}=\id$ and we can take $h_{ji}=0$ for all $j$. If
$\phi(e_i)=e_{j_0}$ is a generator then $S^Y_{\phi(e_i)}-\id=\Delta^Y_{j_0}$
and we can take $h_{j_0i}=1$ and $h_{ji}=0$ otherwise. In general, if
$\phi(e_i)=(a_{i1},\dots,a_{im})$, then
\[
  S^Y_{\phi(e_i)}-\id=\prod_{j=1}^m(\id+\Delta^Y_j)^{a_{ij}}-\id
  =\sum_{j=1}^ma_{ij}\Delta^Y_j+(\text{terms of total degree}\geq 2\text{ in the }\Delta^Y_j),
\]
and the higher-order terms cannot be matched by any choice of the $h_{ji}$,
which only multiply first differences. The simplest instance already fails.

\begin{proposition}\label{prop:no-pullback}
Let $(\Z,\N)$ be $\Z$ with $\N$ acting by translation, and let
$F=(f,\phi)\colon(\Z,\N)\to(\Z,\N)$ be given by $f(x)=2x$ and $\phi(1)=2$.
Then $F$ is a morphism of $\N$-sets, and there is no $h\in\F_\Z$ for which
$F^*(g\,dY)=F^*(g)\,h\,dX$ defines a cochain map.
\end{proposition}

\begin{proof}
%Claude note: the draft omitted this proof as "lengthy"; the following short
%argument was supplied by me. Please check it.
$F$ is a morphism because $f(x+1)=2x+2=\phi(1)+f(x)$. Fix $x\in\Z$ and let
$g$ be the indicator function of $2x+2$. Condition
\eqref{eq:chain-condition} at $x$ reads
$g(2x+2)-g(2x)=h(x)\bigl(g(2x+1)-g(2x)\bigr)$, whose left-hand side is $1$
and whose right-hand side is $0$ for every value of $h(x)$.
\end{proof}

For morphisms sending each generator either to $0$, in which case the
coordinate is said to be collapsed, or to another generator, $F^*$ can be
defined and is a morphism of DGAs.

\begin{definition}[Coordinate morphisms]\label{def:coord-morphism}
A \emph{coordinate morphism} $(X,\N^n)\to(Y,\N^m)$ is a pair $(f,\sigma)$
of a map $f\colon X\to Y$ and a map of index sets
$\sigma\colon[n]\to[m]\cup\{*\}$ such that $(f,\phi_\sigma)$ is a morphism of
$\N$-sets, where $\phi_\sigma(e_i)=e_{\sigma(i)}$ if $\sigma(i)\in[m]$ and
$\phi_\sigma(e_i)=0$ if $\sigma(i)=*$. Coordinate morphisms compose, with
the convention $\sigma(*)=*$, and contain the identities. We write
$\CoordSet$ for the resulting subcategory of $\Nset$ with the same objects.
\end{definition}

We could have written $0$ instead of $*$. The symbol $*$ is chosen to match
the notation for degenerate coordinates in \cref{ch:cubical}, where we show
that coordinate morphisms are closely related to morphisms of cubical sets.

Let $F=(f,\sigma)$ be a coordinate morphism. Since $\Omega^\bullet(Y,\N^m)$
is generated as an algebra by $\F_Y$ and the one-forms $dY_i$, it suffices
to define $F^*$ on these. We set $F^*(g)=g\circ f$ and
\begin{equation}\label{eq:pullback-one-forms}
  F^*(dY_i)=\sum_{j:\ \sigma(j)=i}dX_j,
\end{equation}
an empty sum being $0$. The motivation for this formula is geometric and is
explained in the next section. Multiplicativity then forces
\[
  F^*(dY_I)=F^*(dY_{i_1})\wedge\dots\wedge F^*(dY_{i_p})
  =\sum_{K}\varepsilon_\sigma(K)\,dX_K,
\]
the sum over the $p$-element subsets $K\subseteq[n]$ on which $\sigma$
restricts to a bijection $K\to I$, and $\varepsilon_\sigma(K)$ the sign of
the permutation that $\sigma|_K$ induces relative to the increasing orders on
$K$ and $I$. Finally $F^*(g\,dY_I)=F^*(g)\,F^*(dY_I)$.

\begin{theorem}\label{thm:functor}
With these definitions, $\Omega^\bullet\colon\CoordSet^\op\to\DGA$ is a
functor.
\end{theorem}

The proof consists of sign bookkeeping and is given in \cref{app:functor}.
The one point to record here is that a direction $j$ with $\sigma(j)=*$ is
a direction along which every pulled-back function is constant, since
$\Delta^X_j\circ F^*=F^*\circ(S^Y_0-\id)=0$. Such a direction contributes
nothing to $d_X(F^*\omega)$, which is why collapsing a coordinate is
compatible with the cochain condition.

\subsection{Geometric interpretation of the pullback}
\label{sec:geometric-pullback}

Here, we will use some geometry to justify our choice of $F^{*}$. We can view the 
$1-$differential forms $\Omega^{1}(X,\mathbb{N}^{n})$ as sections on the "discrete cotangent bundle". 
At a point $x$, the cotangent space is spanned by $dX_{1},...,dX_{n}$, indicating that there are $n$ directions.
If we have a morphism of $N-$sets, the pullback of differential forms can then be defined as the pullback of functions on their
respective cotangent spaces. How do we define the cotangent space for a $\mathbb{N}-$set? 

Fix $x\in X$ and let $\mathfrak m_x=\{f\in\F_X:\ f(x)=0\}$ be the ideal of
functions vanishing at $x$. In algebraic and differential geometry the
cotangent space at $x$ is $\mathfrak m_x/\mathfrak m_x^2$, the quotient by
$\mathfrak m_x^2$ discarding second-order information. In the discrete
setting this quotient is zero.

\begin{proposition}\label{prop:m2}
Let $X$ be finite and $\F_X=\Fun(X,\Z)$. Then $\mathfrak m_x^2=\mathfrak m_x$
for every $x\in X$.
\end{proposition}

\begin{proof}
For $y\in X$ let $\delta_y$ be the indicator function of $y$, so that
$\delta_y\delta_z=0$ for $y\neq z$ and $\delta_y^2=\delta_y$. As a
$\Z$-module $\mathfrak m_x$ is free on $\{\delta_y:\ y\neq x\}$, since
$f=\sum_yf(y)\delta_y$. Each generator satisfies
$\delta_y=\delta_y\cdot\delta_y\in\mathfrak m_x^2$.
\end{proof}

The functions on a discrete set are spanned by idempotents, so there is no
second-order information to discard. We use the following substitute.

\begin{definition}[Open neighbourhood]\label{def:open-nbhd}
Let $(X,\N^n)$ be an $\N$-set. The \emph{open neighbourhood} of $x\in X$ is
the set $U_x=\{x,e_1\cdot x,\dots,e_n\cdot x\}$.
\end{definition}

We define the cotangent space $T_{x}^{*} X$ at a point $x$ in the $\mathbb{N}-$set $(X,\N^{n})$ as the set of functions
$f : U_{x} \rightarrow \Z$ such that $f(x) = 0$; this can be thought of as the maximal ideal $\mathfrak m_{x}$ 
of the "stalk" at $x$ (This is not a concrete definition, but it can give the correct intuition. We investigate these ideas
formally in chapter 5 when discussing sheaves). Then, the pullback of one forms is the pullback of sections.

A one-form $g=\sum_ig_i\,dY_i$ on $Y$, evaluated at $y=f(x)$, is the function $\tilde g : U_{y} \rightarrow \Z$
where $\tilde g(y) = 0$, $\tilde g(e^{Y}_{i}\cdot y) = g_{i}(y)$. The pullback of this function on 
$U_{x} = \{x, e_{1}^{X} \cdot x,...,e_{n}^{X} \cdot x\} \subset f^{-1}(U)$, is $F^{*} g(x) = g(f(x)) = 0$ and 
$(F^{*}g)(e^{X}_{i} \cdot x) = g(f(e^{X}_{i} \cdot x)) = g(e^{Y}_{\sigma(i)} \cdot f(x)) = g_{\sigma(i)}(f(x))$. 
This is in the cotangent space $T_{x}^{*}X$ because it vanishes at $x$. Then, we can define the pullback as  
\[
  F^*\Bigl(\sum_{i=1}^mg_i\,dY_i\Bigr)=\sum_{j:\ \sigma(j)\neq *}
  F^*(g_{\sigma(j)})\,dX_j,\qquad\text{so}\qquad
  F^*(dY_i)=\sum_{j:\ \sigma(j)=i}dX_j,
\]
which is formula \eqref{eq:pullback-one-forms}.

\section{First computations}
\label{sec:first-computations}

We begin with the zeroth cohomology, which is the same for every coordinate
system.

\begin{proposition}\label{prop:H0}
Let $\Xs=\{X_1,\dots,X_n\}$ be a coordinate system on $G$. Then
\[
  \HdR{0}(G,\Xs)=\bigcap_{i=1}^n\ker\Delta_i
\]
is the group of functions $V\to\Z$ that are constant on the orbits of the
monoid $\N^n$ acting through $\Xs$, that is, on the connected components
of the induced graph. If there are $r$ such components then
$\HdR{0}(G,\Xs)\cong\Z^r$.
\end{proposition}

\begin{proof}
Since the $dX_i$ form an $\F$-basis of $\Omega^1$, the form
$df=\sum_i\Delta_i(f)\,dX_i$ vanishes if and only if $\Delta_if=0$ for
every $i$, and $\Delta_if=0$ says that $f(X_i(v))=f(v)$ for all $v$. So
$f$ is constant on the classes of the equivalence relation generated by
$v\sim X_i(v)$, which are the connected components of the induced graph.
\end{proof}

The group $\bigcap_i\ker\Delta_i$ consists of the functions that are
constant on the components of the induced graph, which need not be the
components of $G$. A coordinate system that uses only some of the edges of
$G$ has more components than $G$. In \cref{ch:gluing} a graph is covered by
charts and the union of the induced graphs of the charts is required to be
the graph itself (\cref{def:coordinatized-graph}).

\subsection{The cycle graph}
\label{sec:cycle}

\begin{figure}[htb]
\centering
% The cycle graph C_10 with the shift coordinate X(i)=i+1 (arrows).
\begin{tikzpicture}[x=1cm,y=1cm]
\foreach \i in {0,...,9}{
  \pgfmathsetmacro{\ang}{90-36*\i}\pgfmathsetmacro{\angb}{90-36*(\i+1)}
  \draw[cA!85,line width=0.7pt,cedge] (\ang:1.7) -- (\angb:1.7);}
\foreach \i in {0,...,9}{
  \pgfmathsetmacro{\ang}{90-36*\i}
  \node[gvtx] at (\ang:1.7) {};
  \node[font=\scriptsize] at (\ang:2.05) {$\i$};}
\end{tikzpicture}
\caption{The cycle graph $C_{10}$ with the rotation coordinate
$X(i)=i+1$.}
\label{fig:cycle10}
\end{figure}
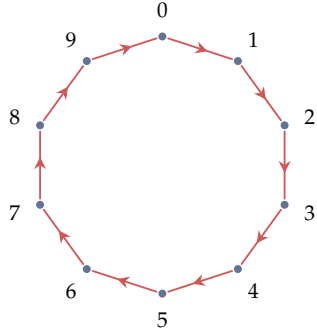

\begin{example}[The cycle graph]\label{ex:cycle}
Let $G=C_N$ be the cycle on $V=\Z/N$ and let $\Xs=\{X\}$ with $X(i)=i+1$
the rotation, as in \cref{fig:cycle10}. Then $\Omega^0=\F$,
$\Omega^1=\F\,dX$, and the complex is
$0\to\F\xrightarrow{\Delta}\F\to 0$ with $\Delta f(i)=f(i+1)-f(i)$. There
is a single orbit, so $\HdR{0}\cong\Z$ by \cref{prop:H0}. The image of
$\Delta$ consists of the functions $g$ with $\sum_ig(i)=0$, since
$\sum_i\Delta f(i)$ telescopes to zero and conversely any such $g$ is
$\Delta f$ for $f(i)=\sum_{k<i}g(k)$. Hence
$\HdR{1}=\F/\im\Delta\cong\Z$ via $g\mapsto\sum_ig(i)$, and
$\HdR{\bullet}(C_N,\Xs)\cong H^\bullet(S^1)$. Here the de Rham complex is
the simplicial cochain complex of $C_N$ with edges oriented by $X$.
\end{example}

\subsection{The torus graph}
\label{sec:torus}

\begin{definition}[Torus graph]\label{def:torus-graph}
Let $N_1,\dots,N_d\geq 3$. The \emph{torus graph} $T=C_{N_1}\times\dots\times
C_{N_d}$ has vertex set $\Z/N_1\times\dots\times\Z/N_d$, with an edge
between $(a_1,\dots,a_d)$ and $(a_1,\dots,a_i+1,\dots,a_d)$ for every $i$.
Its \emph{product shift coordinate system} $\Xs=\{X_1,\dots,X_d\}$ has
$X_i$ adding $1$ in the $i$-th coordinate.
\end{definition}

The torus graph with this coordinate system is the product $\N$-set
$(\Z/N_1,\N)\times\dots\times(\Z/N_d,\N)$ of $d$ cycles with their
rotations, so by \cref{prop:product-dga}
\[
  \Omega^\bullet(T,\Xs)\cong\bigotimes_{i=1}^d\Omega^\bullet(C_{N_i},\{X_i\}).
\]
Recall the K\"unneth formula \citep[Thm.~3.6.3]{weibel1994}. If $P$ and $Q$
are cochain complexes of $R$-modules such that $P^n$ and $d(P^n)$ are flat
for each $n$, there is a short exact sequence
\[
  0\to\bigoplus_{i+j=n}H^i(P)\otimes_RH^j(Q)\to H^n(\Tot(P\otimes_RQ))
  \to\bigoplus_{i+j=n+1}\Tor_1^R\bigl(H^i(P),H^j(Q)\bigr)\to 0.
\]
The terms of a de Rham complex are free $\Z$-modules and so are their
images, submodules of free $\Z$-modules being free, so the hypotheses hold
for any two $\N$-sets.

\begin{theorem}[Cohomology of the discrete torus]\label{thm:torus}
Let $T=C_{N_1}\times\dots\times C_{N_d}$ with the product shift coordinate
system $\Xs$. Then $\HdR{p}(T,\Xs)\cong\Z^{\binom{d}{p}}$ for
$0\leq p\leq d$. The discrete $d$-torus has the cohomology of the smooth
$d$-torus.
\end{theorem}

\begin{proof}
By \cref{ex:cycle} each factor has cohomology $\Z$ in degrees $0$ and $1$
and zero elsewhere. This is free, so all $\Tor_1$ terms vanish and
induction on $d$ gives
\[
  \HdR{p}(T,\Xs)\cong\bigoplus_{i_1+\dots+i_d=p}
  H^{i_1}(\Omega^\bullet(C_{N_1}))\otimes\dots\otimes H^{i_d}(\Omega^\bullet(C_{N_d})).
\]
Each $i_k$ is $0$ or $1$, each summand is $\Z$, and the number of summands
with $\sum i_k=p$ is $\binom{d}{p}$.
\end{proof}

\begin{figure}[htb]
\centering
\resizebox{0.6\textwidth}{!}{\input{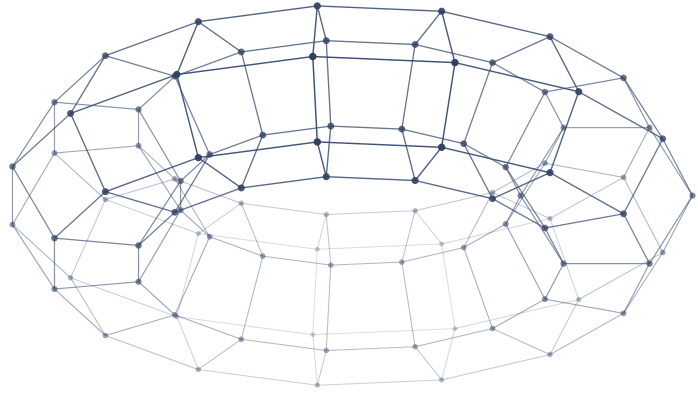}}
\caption{The torus graph $C_{15}\times C_6$.}
\label{fig:torus15x6}
\end{figure}

Each $\Omega^p(G,\Xs)$ is a free $\Z$-module, so the universal coefficient
theorem applies and $\HdR{p}(G,\Xs;R)$ is determined by $\HdR{p}(G,\Xs)$
and $\Tor_1^\Z(\HdR{p+1}(G,\Xs),R)$. For the torus the integral cohomology
is free and the cohomology with coefficients in $R$ is $R^{\binom{d}{p}}$.
We work over $\Z$ from here on.

\subsection{The Cartesian grid}
\label{sec:cartesian-grid}

Our first infinite example is the grid graph on $\Z^n$, with an edge
between $a$ and $a+e_i$ for each standard basis vector $e_i$. It is the
$\N$-set $(\Z^n,\N^n)$ with $e_i\cdot a=a+e_i$, the product of $n$ copies
of the line $(\Z,\N)$. \Cref{prop:product-dga} requires finite sets and
does not apply. The line itself is straightforward: the complex of $(\Z,\N)$
is $\F_\Z\xrightarrow{\Delta}\F_\Z$, and $\Delta$ is surjective, a right
inverse being given by summation as in \cref{sec:calculus-on-N}, so
$\HdR{\bullet}(\Z,\N)=(\Z,0)$.

It is expected for the proof of $(\Z^{n},\N^{n})$ to not be straightforward. We will 
show in section 3 that the de Rham cohomology a $\N-$set is the singular cohomology of 
its associated topological space, and for the $\N-$set $(\Z^{n},\N^{n})$, the associated 
space is $\mathbb{R}^{n}$. And we know that the cohomology of $\R^{n}$ is 0 in higher degrees because it is contractible; this result is the Poincar\'e lemma,
and essentially uses the homotopy invariance theorem for singular cohomology, which is a non-trivial result. 
We want a proof that lies purely in the discrete world, for which we use homotopy invariance.

The case $p=1$ can be done directly, with summation
in place of integration. We give a proof sketch.
%Claude note: the draft stated the homotopy proof in Chapter 3 for (N^n,N^n)
%but announced it here for (Z^n,N^n). The homotopy h(x)=x used there takes
%values in N^n only on the quadrant, so I have separated the two statements.
%The Z^n case is covered by the Shapiro-lemma theorem, which does not need
%finiteness. Please check that this is the intended reading.

\begin{theorem}[Discrete Poincar\'e lemma, degree one]\label{thm:poincare-1}
$\HdR{1}(\Z^n,\N^n)=0$.
\end{theorem}

\begin{proof}
A one-form $\omega=\sum_j\omega_j\,dX_j$ is closed if and only if
$\Delta_i\omega_j=\Delta_j\omega_i$ for all $i,j$, since
$d\omega=\sum_{i<j}(\Delta_i\omega_j-\Delta_j\omega_i)\,dX_i\wedge dX_j$. Given a
closed $\omega$ we integrate along the coordinate axes. Use the summation
convention $\sum_{k=0}^{x-1}a_k=-\sum_{k=x}^{-1}a_k$ for $x<0$, under which
$\sum_{k=0}^{x}a_k-\sum_{k=0}^{x-1}a_k=a_x$ for every $x\in\Z$, and define
\[
  f(x_1,\dots,x_n)=\sum_{j=1}^n\ \sum_{k=0}^{x_j-1}
  \omega_j(x_1,\dots,x_{j-1},k,0,\dots,0).
\]
Then $\Delta_nf=\omega_n$ directly. For $i<n$, the terms with $j<i$ do not
depend on $x_i$, the term $j=i$ contributes $\omega_i(x_1,\dots,x_i,0,\dots,0)$ to the sum,
and each term with $j>i$ contributes
\[
  \sum_{k=0}^{x_j-1}\Delta_i\omega_j(x_1,\dots,x_{j-1},k,0,\dots,0)
  =\sum_{k=0}^{x_j-1}\Delta_j\omega_i(x_1,\dots,x_{j-1},k,0,\dots,0),
\]
using closedness. This sum telescopes to
$\omega_i(x_1,\dots,x_j,0,\dots,0)-\omega_i(x_1,\dots,x_{j-1},0,\dots,0)$.
Summing over $j=i+1,\dots,n$ we get $\omega_i(x_1,\dots,x_n)$. Hence $df=\omega$.
\end{proof}

\subsection{The polar grid}
\label{sec:polar-grid}

\begin{figure}[htb]
\centering
\resizebox{\textwidth}{!}{% (a) the polar grid P(10) truncated at radius 6; (b) the graph induced by the
% coordinate system {R,Theta}: R points inwards (red), Theta rotates (blue),
% and both fix the origin, which is why two loops appear there.
\begin{tikzpicture}[x=1cm,y=1cm]
\def\m{10}\def\K{6}\def\step{0.52}
% ---------------- (a) -------------------------------------------------------
\begin{scope}
  \foreach \a in {1,...,\K}{
    \pgfmathsetmacro{\rr}{\a*\step}
    \foreach \i in {0,...,9}{
      \pgfmathsetmacro{\ta}{360*\i/\m}\pgfmathsetmacro{\tb}{360*(\i+1)/\m}
      \draw[gline] (\ta:\rr) -- (\tb:\rr);}}
  \foreach \i in {0,...,9}{
    \pgfmathsetmacro{\ta}{360*\i/\m}
    \draw[gline] (0,0) -- (\ta:\K*\step);}
  \foreach \a in {1,...,\K}{\pgfmathsetmacro{\rr}{\a*\step}
    \foreach \i in {0,...,9}{\pgfmathsetmacro{\ta}{360*\i/\m}\node[gvtx] at (\ta:\rr) {};}}
  \node[gvtx] at (0,0) {};
  \node[font=\small] at (0,-3.9) {(a) the polar grid $\PG(10)$, drawn to radius $6$};
\end{scope}
% ---------------- (b) -------------------------------------------------------
\begin{scope}[xshift=8.6cm]
  \foreach \a in {1,...,\K}{
    \pgfmathsetmacro{\rr}{\a*\step}
    \foreach \i in {0,...,9}{
      \pgfmathsetmacro{\ta}{360*\i/\m}\pgfmathsetmacro{\tb}{360*(\i+1)/\m}
      \draw[cB!85,line width=0.55pt,cedge] (\ta:\rr) -- (\tb:\rr);}}
  \foreach \a in {2,...,\K}{
    \pgfmathsetmacro{\ro}{\a*\step}\pgfmathsetmacro{\ri}{(\a-1)*\step}
    \foreach \i in {0,...,9}{
      \pgfmathsetmacro{\ta}{360*\i/\m}
      \draw[cA!85,line width=0.55pt,cedge] (\ta:\ro) -- (\ta:\ri);}}
  \foreach \i in {0,...,9}{\pgfmathsetmacro{\ta}{360*\i/\m}
    \draw[cA!85,line width=0.55pt,cedge] (\ta:\step) -- (0,0);}
  \foreach \a in {1,...,\K}{\pgfmathsetmacro{\rr}{\a*\step}
    \foreach \i in {0,...,9}{\pgfmathsetmacro{\ta}{360*\i/\m}\node[gvtx] at (\ta:\rr) {};}}
  \node[gvtx] at (0,0) {};
  % the two loops at the origin
  \draw[poleloop,cA] (-0.22,0.16) arc (150:-160:0.25);
  \draw[poleloop,cB] (0.22,0.16) arc (30:340:0.25);
  \node[font=\small] at (0,-3.9) {(b) the induced graph of $(\PG(10),\{R,\Theta\})$};
  \node[font=\scriptsize,align=left,anchor=west] at (3.5,2.6)
    {\textcolor{cA}{\rule{7pt}{1.6pt}}~$R$ (inwards)\\[2pt]
     \textcolor{cB}{\rule{7pt}{1.6pt}}~$\Theta$ (rotation)\\[2pt]
     two loops at $\mathbf 0$};
\end{scope}
\end{tikzpicture}}
\caption{The polar grid $\PG(10)$ and the directed graph induced by the
coordinate system $\{R,\Theta\}$. Both coordinates fix the origin, which
is why two loops appear there.}
\label{fig:polar-grid}
\end{figure}
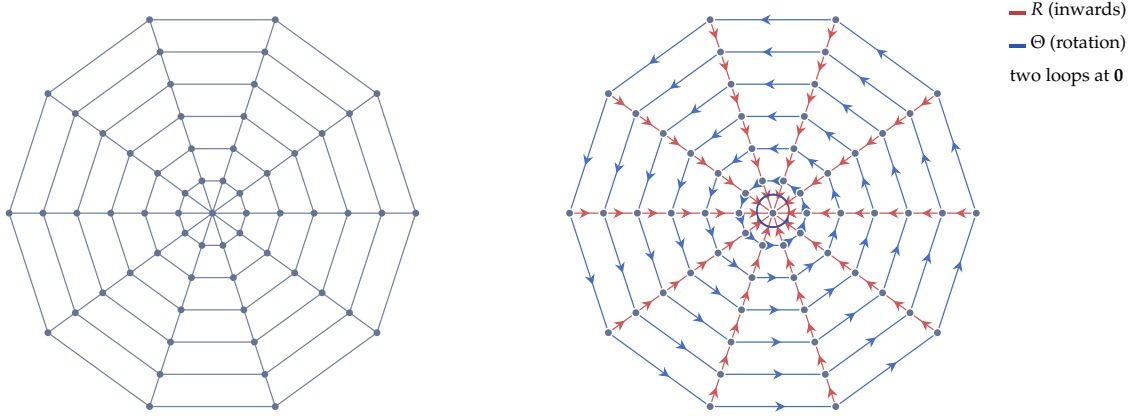

What about a graph modelled on polar coordinates? Let $m\geq 3$. The
\emph{polar grid} $\PG(m)$ has vertex set
$V=(\N\times\Z/m)/\!\sim$, where all the points $(0,i)$ are identified to a
single origin $\mathbf 0$, and edges $\{(r,i),(r+1,i)\}$ and
$\{(r,i),(r,i+1)\}$ for all $r\geq 0$ and $i\in\Z/m$. \Cref{fig:polar-grid}(a)
shows $\PG(10)$ out to radius $6$. The rings are the circles on which the
radial coordinate is constant and the spokes are the lines on which the
angular coordinate is constant.

We introduce the radial and angular coordinate fields. Let
$R(r,i)=(r-1,i)$ for $r>0$ and $R(\mathbf 0)=\mathbf 0$, and let
$\Theta(r,i)=(r,i+1)$ for $r>0$ and $\Theta(\mathbf 0)=\mathbf 0$. Both are
coordinate fields and they commute. Away from the origin
$R\Theta(r,i)=(r-1,i+1)=\Theta R(r,i)$ when $r\geq 2$, both composites send
$(1,i)$ to $\mathbf 0$, and both fix $\mathbf 0$. So $\{R,\Theta\}$ is a
coordinate system of rank $2$.

\begin{proposition}\label{prop:polar-torus}
The de Rham cohomology of $(\PG(m),\{R,\Theta\})$ is $(\Z,\Z^2,\Z)$, the
cohomology of the $2$-torus.
\end{proposition}

\begin{proof}
The complex is $0\to\F\xrightarrow{d^0}\F^2\xrightarrow{d^1}\F\to 0$ with
$d^0f=\Delta_R(f)\,dR+\Delta_\Theta(f)\,d\Theta$ and
$d^1(g\,dR+h\,d\Theta)=(\Delta_R(h)-\Delta_\Theta(g))\,dR\wedge d\Theta$.

\emph{Degree $0$.} Every vertex is carried to the origin by a power of $R$,
so there is a single orbit and $\HdR{0}\cong\Z$ by \cref{prop:H0}.

\emph{Degree $2$.} Both $\Delta_R$ and $\Delta_\Theta$ vanish at the origin,
since $R$ and $\Theta$ fix it, so every element of $\im d^1$ vanishes at
$\mathbf 0$. Conversely let $f\in\F$ with $f(\mathbf 0)=0$, take $g=0$ and
define $h(\mathbf 0)=0$ and $h(r,i)=-\sum_{s=1}^rf(s,i)$ for $r\geq 1$. For
$r\geq 2$ we get $\Delta_R(h)(r,i)=h(r-1,i)-h(r,i)=f(r,i)$, and for $r=1$
we get $\Delta_R(h)(1,i)=h(\mathbf 0)-h(1,i)=f(1,i)$. So
$d^1(h\,d\Theta)=f\,dR\wedge d\Theta$, hence
$\im d^1=\{f:\ f(\mathbf 0)=0\}$ and evaluation at the origin gives
$\HdR{2}\cong\Z$.

\emph{Degree $1$.} A form $g\,dR+h\,d\Theta$ is closed if and only if
$\Delta_R(h)=\Delta_\Theta(g)$, that is,
\[
  h(r-1,i)-h(r,i)=g(r,i+1)-g(r,i)\qquad\text{for all }r\geq 1,\ i\in\Z/m.
\]
We claim that
\[
  \Phi\colon\ker d^1\to\Z^2,\qquad g\,dR+h\,d\Theta\mapsto
  \bigl(g(\mathbf 0),h(\mathbf 0)\bigr)
\]
induces an isomorphism $\HdR{1}\cong\Z^2$. It kills exact forms, since
$\Delta_R(f)(\mathbf 0)=f(R(\mathbf 0))-f(\mathbf 0)=0$ and likewise for
$\Delta_\Theta$. It is surjective, because the forms $\delta_{\mathbf 0}\,dR$
and $1\cdot d\Theta$ are closed and map to $(1,0)$ and $(0,1)$. For
injectivity, suppose $g\,dR+h\,d\Theta$ is closed with
$g(\mathbf 0)=h(\mathbf 0)=0$, and define $f(\mathbf 0)=0$,
$f(r,i)=-\sum_{s=1}^rg(s,i)$ for $r\geq 1$. Then $\Delta_R(f)=g$ by the
computation in degree $2$, using $g(\mathbf 0)=0$ at the origin. It remains
to check $\Delta_\Theta(f)=h$, that is, $f(r,i+1)-f(r,i)=h(r,i)$ for $r\geq 1$,
the case $r=0$ being $0=h(\mathbf 0)$. For $r=1$,
\begin{align*}
  f(1,i+1)-f(1,i)&=g(1,i)-g(1,i+1)=-\Delta_\Theta(g)(1,i)\\
  &=-\Delta_R(h)(1,i)=h(1,i)-h(\mathbf 0)=h(1,i).
\end{align*}
For $r\geq 2$, assuming the identity at $r-1$,
\begin{align*}
  f(r,i+1)-f(r,i)&=\bigl(f(r-1,i+1)-f(r-1,i)\bigr)-\bigl(g(r,i+1)-g(r,i)\bigr)\\
  &=h(r-1,i)-\Delta_\Theta(g)(r,i),
\end{align*}
and $\Delta_\Theta(g)(r,i)=\Delta_R(h)(r,i)=h(r-1,i)-h(r,i)$ by closedness,
so the right-hand side is $h(r,i)$. Hence $g\,dR+h\,d\Theta=df$ is exact.
\end{proof}

The polar grid is a discrete model of the plane, and its de Rham
cohomology is that of a torus rather than of a contractible space. In the
proof, the obstruction to solving the recurrences is at the origin, where
the difference of every function vanishes because both coordinates fix it.

\begin{intuition}[the induced graph]
The de Rham complex of an $\N$-set $(X,\N^n)$ depends only on the induced
directed graph, with an arrow $x\to e_i\cdot x$ for every $x$ and $i$, and
not on the graph on which the coordinate system was defined.
\Cref{fig:polar-grid}(b) shows the induced graph of the polar grid. It
differs from $\PG(m)$ only at the origin, where
$R(\mathbf 0)=\Theta(\mathbf 0)=\mathbf 0$ produces two loops, and the two
classes in $\HdR{1}$ are represented by the forms $\delta_{\mathbf 0}\,dR$
and $d\Theta$, which are supported on these loops. The cohomology computed
above is that of the induced graph with its two loops.
\end{intuition}

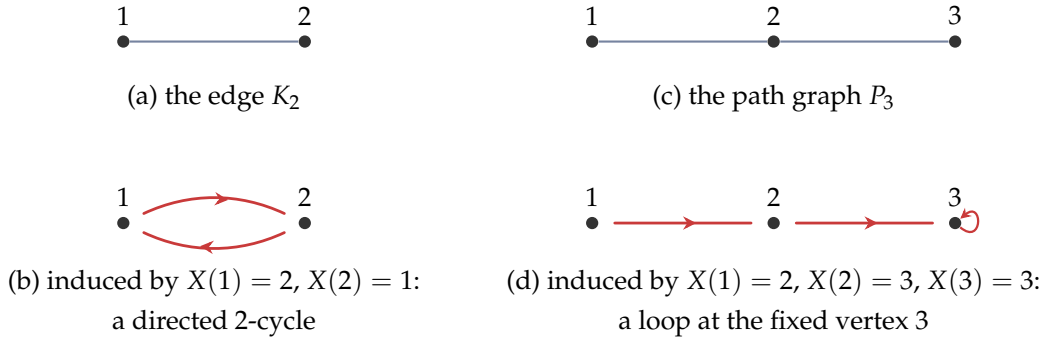
\begin{figure}[htb]
\centering
% The two smallest examples of the "induced graph" phenomenon.
% Left: the edge graph and the graph induced by X(1)=2, X(2)=1.
% Right: the path graph and the graph induced by X(1)=2, X(2)=3, X(3)=3.
\begin{tikzpicture}[x=1cm,y=1cm]
% ---- the edge ----
\begin{scope}
  \node[vtx,label={[font=\small]above:$1$}] (a) at (0,0) {};
  \node[vtx,label={[font=\small]above:$2$}] (b) at (2.4,0) {};
  \draw[gline,line width=0.8pt] (a) -- (b);
  \node[font=\small] at (1.2,-0.75) {(a) the edge $K_2$};
\end{scope}
\begin{scope}[yshift=-2.4cm]
  \node[vtx,label={[font=\small]above:$1$}] (a) at (0,0) {};
  \node[vtx,label={[font=\small]above:$2$}] (b) at (2.4,0) {};
  \draw[ar,cA] (a) to[bend left=25] (b);
  \draw[ar,cA] (b) to[bend left=25] (a);
  \node[font=\small,align=center] at (1.2,-1.05) {(b) induced by $X(1)=2$, $X(2)=1$:\\ a directed $2$-cycle};
\end{scope}
% ---- the path ----
\begin{scope}[xshift=6.2cm]
  \node[vtx,label={[font=\small]above:$1$}] (a) at (0,0) {};
  \node[vtx,label={[font=\small]above:$2$}] (b) at (2.4,0) {};
  \node[vtx,label={[font=\small]above:$3$}] (c) at (4.8,0) {};
  \draw[gline,line width=0.8pt] (a) -- (b) -- (c);
  \node[font=\small] at (2.4,-0.75) {(c) the path graph $P_3$};
\end{scope}
\begin{scope}[xshift=6.2cm,yshift=-2.4cm]
  \node[vtx,label={[font=\small]above:$1$}] (a) at (0,0) {};
  \node[vtx,label={[font=\small]above:$2$}] (b) at (2.4,0) {};
  \node[vtx,label={[font=\small]above:$3$}] (c) at (4.8,0) {};
  \draw[ar,cA] (a) -- (b);
  \draw[ar,cA] (b) -- (c);
  \draw[lp,cA] (c) to[out=-40,in=40,looseness=9] (c);
  \node[font=\small,align=center] at (2.4,-1.05) {(d) induced by $X(1)=2$, $X(2)=3$, $X(3)=3$:\\ a loop at the fixed vertex $3$};
\end{scope}
\end{tikzpicture}
\caption{The two smallest examples of the induced-graph phenomenon. A
coordinate field on the edge produces a $2$-cycle, and a coordinate field
on the path produces a loop at the vertex it fixes.}
\label{fig:small-graphs}
\end{figure}

The same phenomenon occurs for the smallest graphs. The graph with two
vertices $1,2$ and a single edge, with the coordinate field $X(1)=2$,
$X(2)=1$, has induced directed graph a $2$-cycle, \cref{fig:small-graphs}(b),
and its de Rham cohomology is $(\Z,\Z)$, that of the circle rather than of
the edge. In every finite $\N$-set each $X_i$ is defined on all of $V$, so
the coordinates eventually cycle or come to rest, and every fixed point and
every cycle is a feature of the induced graph that the original graph need
not have. This is addressed in \cref{sec:normalized}.

\subsection{The path graph and the square graph}
\label{sec:path-square}

Consider the path graph on the vertices $1,2,3$ with the coordinate field
$X(1)=2$, $X(2)=3$, $X(3)=3$, as in \cref{fig:small-graphs}(c,d). The
induced directed graph has a loop at $3$, and the de Rham cohomology is
$(\Z,\Z)$. Indeed $\Omega^0\to\Omega^1$ is $\Z^3\to\Z^3$,
$(a,b,c)\mapsto(b-a,c-b,0)$, whose cokernel is $\Z$. The class is "carried" by
the loop.

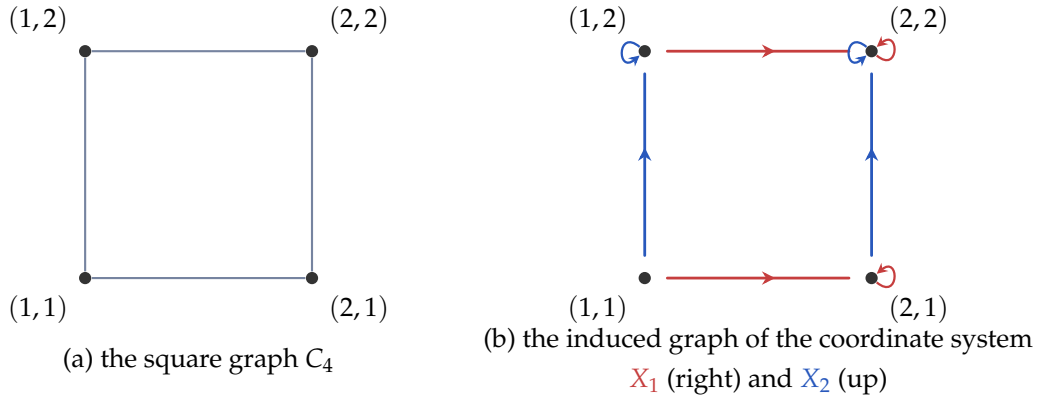
\begin{figure}[htb]
\centering
% The square graph and the rank-two coordinate system on it: e_1 moves right,
% e_2 moves up, and each loops where the move would leave the square.
\begingroup
\begin{tikzpicture}[x=1cm,y=1cm]
\begin{scope}
  \node[vtx] (a) at (0,0) {};  \node[vtx] (b) at (3,0) {};
  \node[vtx] (c) at (0,3) {};  \node[vtx] (d) at (3,3) {};
  \node[font=\small,below left =1pt of a] {$(1,1)$};
  \node[font=\small,below right=1pt of b] {$(2,1)$};
  \node[font=\small,above left =1pt of c] {$(1,2)$};
  \node[font=\small,above right=1pt of d] {$(2,2)$};
  \draw[gline,line width=0.8pt] (a) -- (b) -- (d) -- (c) -- (a);
  \node[font=\small] at (1.5,-1.1) {(a) the square graph $C_4$};
\end{scope}
\begin{scope}[xshift=7.4cm]
  \node[vtx] (a) at (0,0) {};  \node[vtx] (b) at (3,0) {};
  \node[vtx] (c) at (0,3) {};  \node[vtx] (d) at (3,3) {};
  \node[font=\small,below left =1pt of a] {$(1,1)$};
  \node[font=\small,below right=1pt of b] {$(2,1)$};
  \node[font=\small,above left =1pt of c] {$(1,2)$};
  \node[font=\small,above right=1pt of d] {$(2,2)$};
  \draw[ar,e1col] (a) -- (b);
  \draw[ar,e1col] (c) -- (d);
  \draw[ar,e2col] (a) -- (c);
  \draw[ar,e2col] (b) -- (d);
  \draw[lp,e1col] (b) to[out=-40,in=40,looseness=9] (b);
  \draw[lp,e1col] (d) to[out=-40,in=40,looseness=9] (d);
  \draw[lp,e2col] (c) to[out=140,in=220,looseness=9] (c);
  \draw[lp,e2col] (d) to[out=140,in=220,looseness=9] (d);
  \node[font=\small,align=center] at (1.5,-1.1) {(b) the induced graph of the coordinate system\\ $\textcolor{e1col}{X_1}$ (right) and $\textcolor{e2col}{X_2}$ (up)};
\end{scope}
\end{tikzpicture}
\endgroup
\caption{The square graph and its rank-two coordinate system. The
coordinates fix the vertices they cannot leave, and the induced graph has
the cohomology of a torus.}
\label{fig:square-coords}
\end{figure}

The same occurs for a two-dimensional coordinate system on the square
graph $C_4$ with vertices $(1,1),(2,1),(1,2),(2,2)$, \cref{fig:square-coords}.
Regard the square as a subgraph of the grid $\Z^2$ and let $X_1$ move one
step to the right and $X_2$ one step up, with the convention that a
coordinate which would leave the square fixes the vertex. So
$X_1(1,j)=(2,j)$, $X_1(2,j)=(2,j)$, $X_2(i,1)=(i,2)$, $X_2(i,2)=(i,2)$.
The two fields commute, and the induced graph has four loops at the
vertices where a coordinate is forced to stop. The de Rham cohomology of
this $\N$-set is $(\Z,\Z^2,\Z)$, the cohomology of the $2$-torus, the
classes again being supported on the loops.

\begin{intuition}[normalization]
Each $\N$-set $(X,\N^n)$ determines a semi-cubical set, with the points of
$X$ as $0$-cubes, the arrows $x\to e_i\cdot x$ as $1$-cubes, and a base
point together with $p$ distinct directions as a $p$-cube
(\cref{sec:cubes}). The de Rham complex is the cubical cochain complex of
this semi-cubical set. A $1$-cube $x\to e_i\cdot x$ with $e_i\cdot x=x$ has
both faces equal, which is the defining property of a degenerate cube, and
in \cref{sec:normalized} such cubes are \textit{declared} degenerate and the
normalized cochain complex is formed. The normalized cohomology of the
polar grid is that of a disc, of the square graph that of a filled square,
and of the path graph that of a point. The charts of \cref{ch:gluing} are
$\N$-sets with this normalization.
\end{intuition}

\chapter{Homotopy invariance and toral rigidity}
\label{ch:homotopy}

This chapter studies more properties of the de Rham complex. We show that the de Rham complex is
a Koszul complex, and its cohomology is shown to be an $\Ext$ group over a polynomial
ring. For bijective coordinate systems this $\Ext$ group is the group
cohomology of $\Z^n$ with coefficients in the functions on $X$, which is
computed by Shapiro's lemma. In general it is the cohomology of the monoid
$\N^n$, which we identify with the cohomology of the nerve of an action
category. Natural transformations between functors induce homotopies
between nerves, and this gives a homotopy invariance theorem.

\section{Properties of the de Rham complex}
\label{sec:properties}

\subsection{The de Rham complex is a Koszul complex}
\label{sec:koszul}

Recall the Koszul complex \citep[\S4.5]{weibel1994}. If $R$ is a
commutative ring and $x=(x_1,\dots,x_n)$ is a sequence of elements of $R$,
the Koszul complex $K_\bullet(x)$ is the tensor product
$K_\bullet(x_1)\otimes_R\dots\otimes_RK_\bullet(x_n)$ of the two-term
complexes $R\xrightarrow{x_i}R$. Concretely $K_p(x)=\Lam^pR^n$, with basis
$e_I=e_{i_1}\wedge\dots\wedge e_{i_p}$ for $I=\{i_1<\dots<i_p\}$ and
differential
\[
  d(e_{i_1}\wedge\dots\wedge e_{i_p})=\sum_{k=1}^p(-1)^{k+1}x_{i_k}\,
  e_{i_1}\wedge\dots\wedge\widehat{e_{i_k}}\wedge\dots\wedge e_{i_p}.
\]
Its dual $\Hom_R(K_\bullet(x),R)$ is the cochain complex with
$K^p=\Lam^pR^n$ and $d(e_I)=(\sum_ix_ie_i)\wedge e_I$. We need a version
with operators in place of ring elements.

\begin{definition}[Koszul complexes of commuting operators]\label{def:koszul-ops}
Let $M$ be an $R$-module and let $\varphi_1,\dots,\varphi_n\colon M\to M$
be pairwise commuting $R$-linear maps. The \emph{Koszul chain complex}
$K_\bullet(\varphi;M)$ has $K_p=M\otimes_R\Lam^pR^n$ and
\[
  d_p(m\otimes e_{i_1}\wedge\dots\wedge e_{i_p})=\sum_{k=1}^p(-1)^{k+1}
  \varphi_{i_k}(m)\otimes e_{i_1}\wedge\dots\wedge\widehat{e_{i_k}}\wedge\dots\wedge e_{i_p}.
\]
The \emph{Koszul cochain complex} $K^\bullet(\varphi;M)$ has
$K^p=M\otimes_R\Lam^pR^n$ and
\[
  d^p(\alpha\otimes e_I)=\sum_{i=1}^n\varphi_i(\alpha)\otimes e_i\wedge e_I.
\]
\end{definition}

Let $S=R[T_1,\dots,T_n]$ and make $M$ an $S$-module by $T_i\cdot m=\varphi_i(m)$,
which is well defined because the $\varphi_i$ commute.

\begin{proposition}\label{prop:koszul-hom}
There are isomorphisms of complexes of $S$-modules
\[
  K_\bullet(\varphi;M)\cong M\otimes_SK_\bullet(T_1,\dots,T_n;S),\qquad
  K^\bullet(\varphi;M)\cong\Hom_S\bigl(K_\bullet(T_1,\dots,T_n;S),M\bigr).
\]
\end{proposition}

\begin{proof}
The complex $K_\bullet(T;S)$ has $K_p=\Lam^pS^n$ with differential
$d(e_{i_1}\wedge\dots\wedge e_{i_p})=\sum_k(-1)^{k+1}T_{i_k}\,e_{\dots\widehat{i_k}\dots}$.
Tensoring with $M$ over $S$ and using the identification
$M\otimes_S\Lam^pS^n\cong M\otimes_R\Lam^pR^n$ turns $m\otimes T_{i_k}e_J$
into $\varphi_{i_k}(m)\otimes e_J$, which is the differential of
$K_\bullet(\varphi;M)$.

For the second isomorphism, $\Hom_S(\Lam^pS^n,M)\cong M\otimes_R\Lam^pR^n$
via $\vartheta\mapsto\sum_{|I|=p}\vartheta(e_I)\otimes e_I$. Under this
identification the dual differential $\vartheta\mapsto\vartheta\circ d_{p+1}$
sends $\alpha\otimes e_I$, the functional with $\vartheta(e_I)=\alpha$ and
$\vartheta(e_J)=0$ for $J\neq I$, to
$\sum_J\vartheta(d\,e_J)\otimes e_J$. Only those $J=\{j_1<\dots<j_{p+1}\}$
containing $I$ contribute, and if $J=I\cup\{i\}$ with $i=j_k$ then
$\vartheta(d\,e_J)=(-1)^{k+1}T_i\cdot\alpha=(-1)^{k+1}\varphi_i(\alpha)$,
while $e_i\wedge e_I=(-1)^{k-1}e_J$. Hence the image of $\alpha\otimes e_I$
is $\sum_i\varphi_i(\alpha)\otimes e_i\wedge e_I$.
\end{proof}

\begin{proposition}\label{prop:koszul-ext}
With notation as above,
$H_p(K_\bullet(\varphi;M))\cong\Tor^S_p(R,M)$ and
$H^p(K^\bullet(\varphi;M))\cong\Ext^p_S(R,M)$, where $R=S/(T_1,\dots,T_n)$.
\end{proposition}

\begin{proof}
The sequence $(T_1,\dots,T_n)$ is regular in $S$, so the augmented Koszul
complex $K_\bullet(T;S)\to S/(T_1,\dots,T_n)=R$ is a free resolution of
$R$ over $S$ \citep[Cor.~4.5.5]{weibel1994}. Applying $M\otimes_S-$ and
$\Hom_S(-,M)$ and using \cref{prop:koszul-hom} gives the result.
\end{proof}

The de Rham complex of an $\N$-set $(X,\N^n)$ is the Koszul cochain complex
$K^\bullet(\Delta;\F_X)$ of the operators $\Delta_1,\dots,\Delta_n$ over
$R=\Z$, with $e_i$ written $dX_i$. Hence:

\begin{corollary}\label{cor:ext}
For every $\N$-set $(X,\N^n)$, and in particular for every coordinatized
graph, $\HdR{p}(X,\N^n)\cong\Ext^p_{\Z[T_1,\dots,T_n]}(\Z,\F_X)$, where
$T_i$ acts on $\F_X$ by $\Delta_i$ and on $\Z$ by zero.
\end{corollary}

\subsection{Bijective coordinates}
\label{sec:bijective}

Call an $\N$-set $(X,\N^n)$ \emph{bijective} if each generator $e_i$ acts
bijectively on $X$. The action then extends to an action of the group
$\Z^n$, by adjoining the inverses, and $\F_X$ becomes a module over the
group ring $\Z[\Z^n]=\Z[T_1^{\pm1},\dots,T_n^{\pm1}]$ with $T_i$ acting as
the shift $S_i$. The de Rham differential is built from
$\Delta_i=S_i-\id=T_i-1$, so the de Rham complex is the Koszul cochain
complex of the sequence $(T_1-1,\dots,T_n-1)$ acting on $\F_X$. Since
$(T_1-1,\dots,T_n-1)$ is a regular sequence in
$\Z[T_1^{\pm1},\dots,T_n^{\pm1}]$ with quotient $\Z$, the same proof of
\cref{prop:koszul-ext} gives
\begin{equation}\label{eq:group-cohomology}
  \HdR{p}(X,\N^n)\cong\Ext^p_{\Z[\Z^n]}(\Z,\F_X)=H^p(\Z^n;\F_X),
\end{equation}
the group cohomology of $\Z^n$ with coefficients in $\F_X$.

Group cohomology of this kind is computed by Shapiro's lemma
\citep[Lemma~6.3.2]{weibel1994}. Write $X/\Z^n$ for the set of orbits. Since
$\Z^n$ is abelian, points in the same orbit have the same stabilizer, so
each orbit $O$ has a well-defined stabilizer $\stab(O)\leq\Z^n$, a free
abelian group of some rank $r_O$ with $0\leq r_O\leq n$.

\begin{theorem}\label{thm:bijective}
Let $(X,\N^n)$ be a bijective $\N$-set. Then
\[
  \HdR{p}(X,\N^n)\cong\prod_{O\in X/\Z^n}\Z^{\binom{r_O}{p}},
\]
where $r_O$ is the rank of the stabilizer of $O$ and $\binom{r}{p}=0$ for
$p>r$.
\end{theorem}

\begin{proof}
Suppose first that the action is transitive, with stabilizer $H$ of rank
$r$. For any $x_0\in X$ the map $\Z^n/H\to X$, $a+H\mapsto a\cdot x_0$, is a
bijection under which the action of $\Z^n$ becomes translation. So
$\F_X\cong\Fun(\Z^n/H,\Z)$ as $\Z[\Z^n]$-modules. A function on $\Z^n/H$ is
a function on $\Z^n$ constant on the cosets of $H$, which is the description
of a coinduced module: the map
\[
  \psi\colon\Fun(\Z^n/H,\Z)\to\Hom_{\Z H}(\Z[\Z^n],\Z)=\Coind_H^{\Z^n}\Z,
  \qquad \psi(f)(e_a)=f(a+H),
\]
where $\Z$ carries the trivial $H$-action, is an isomorphism of
$\Z[\Z^n]$-modules. It is well defined since
$\psi(f)(e_{h+a})=f(a+H)=h\cdot\psi(f)(e_a)$ for $h\in H$, it is bijective
because a $\Z H$-linear map out of $\Z[\Z^n]$ is determined by its values
on coset representatives, and it is $\Z[\Z^n]$-linear because
$\psi(T_i\cdot f)(e_a)=f(a+e_i+H)=(T_i\cdot\psi(f))(e_a)$. Shapiro's lemma
and \eqref{eq:group-cohomology} now give
\[
  \HdR{p}(X,\N^n)=H^p(\Z^n;\F_X)\cong H^p\bigl(\Z^n;\Coind_H^{\Z^n}\Z\bigr)
  \cong H^p(H;\Z)\cong H^p(\Z^r;\Z)\cong\Z^{\binom{r}{p}},
\]
the last step being the cohomology of the $r$-torus, or equivalently the
Koszul complex of $(T_1-1,\dots,T_r-1)$ on the trivial module $\Z$.

In general $X=\coprod_OO$ is the disjoint union of its orbits, so
$\F_X\cong\prod_O\F_O$ as $\Z[\Z^n]$-modules and
$\Omega^p(X,\N^n)\cong\prod_O\Omega^p(O,\N^n)$ compatibly with the
differentials. Cohomology commutes with products.
\end{proof}

\begin{example}
The grid $(\Z^n,\N^n)$ of \cref{sec:cartesian-grid} is bijective and
transitive with trivial stabilizer, so $\HdR{\bullet}(\Z^n,\N^n)=(\Z,0,\dots,0)$.
This is the discrete Poincar\'e lemma for the grid. The torus graph of
\cref{thm:torus} is bijective and transitive with stabilizer
$N_1\Z\times\dots\times N_d\Z$ of rank $d$, which gives
$\HdR{p}=\Z^{\binom{d}{p}}$ again.
\end{example}

\subsection{Monoid cohomology}
\label{sec:monoid-cohomology}

For bijective coordinates the de Rham cohomology is group cohomology. In
general it is the Eilenberg--MacLane cohomology of the monoid $\N^n$
\citep{cegarra-leech2025}.

\begin{definition}[Monoid cohomology]\label{def:monoid-cohomology}
Let $M$ be a monoid and $A$ a left $\Z M$-module. The \emph{cohomology of
$M$ with coefficients in $A$} is $H^p(M;A)=\Ext^p_{\Z M}(\Z,A)$, where $\Z$
carries the trivial action.
\end{definition}

For $M=\N^n$ the monoid ring $\Z[\N^n]$ is the polynomial ring
$\Z[T_1,\dots,T_n]$ with $T_i=e_i$, and the trivial module is
$\Z=\Z[\N^n]/(T_1-1,\dots,T_n-1)$. The module $\F_X$ of an $\N$-set is a
$\Z[\N^n]$-module with $T_i$ acting by the shift $S_i$, so that
$T_i-1$ acts by $\Delta_i$. Substituting $T_i-1$ for $T_i$ in
\cref{cor:ext} gives the following.

\begin{corollary}\label{cor:monoid}
For every $\N$-set, $\HdR{p}(X,\N^n)\cong H^p(\N^n;\F_X)$.
\end{corollary}

Monoid cohomology, like group cohomology, can be computed from the bar
resolution. The bar resolution of $\Z$ over $\Z M$ is
\[
  \dots\to B_2M\to B_1M\to B_0M\xrightarrow{\ \epsilon\ }\Z\to 0,
\]
where $B_pM$ is the free $\Z M$-module on symbols $[u_1|\dots|u_p]$ with
$u_i\in M$, and the differential is $\sum_{i=0}^p(-1)^id_i$ with
\begin{align*}
  d_0[u_1|\dots|u_p]&=u_1[u_2|\dots|u_p],\\
  d_i[u_1|\dots|u_p]&=[u_1|\dots|u_iu_{i+1}|\dots|u_p]\qquad(0<i<p),\\
  d_p[u_1|\dots|u_p]&=[u_1|\dots|u_{p-1}].
\end{align*}
It is the unnormalized chain complex of a simplicial $\Z M$-module.
Applying $\Hom_{\Z M}(-,A)$ and identifying
$\Hom_{\Z M}(B_pM,A)\cong\Fun(M^p,A)$ gives the cochain complex with
\begin{equation}\label{eq:bar-cochains}
\begin{aligned}
  (\delta f)(m_1,\dots,m_p)={}&m_1\cdot f(m_2,\dots,m_p)
  +\sum_{i=1}^{p-1}(-1)^if(m_1,\dots,m_im_{i+1},\dots,m_p)\\
  &+(-1)^pf(m_1,\dots,m_{p-1}),
\end{aligned}
\end{equation}
whose cohomology is $H^p(M;A)$. The de Rham complex and the bar complex
compute the same groups. The de Rham complex is finite, is built from the
difference operators, and carries the twisted product of
\cref{sec:algebra-structure}. The bar complex is infinite, but it is the
cochain complex of a simplicial set, and this is used in the next section.

\section{Action categories and homotopy invariance}
\label{sec:action-categories}

\subsection{Action categories}
\label{sec:action-cat}

A group acting on a set gives rise to an action groupoid, and many
properties of the action are encoded in the groupoid. Action groupoids are
also called quotient stacks, and they are used extensively in the theory of
differentiable stacks and orbifolds through the notion of a Lie groupoid
\citep{behrend-xu2011}. Monoid actions give categories rather than
groupoids.

\begin{definition}[Action category]\label{def:action-category}
Let a monoid $M$ act on a set $X$. The \emph{action category}
$\cA(X,M)=M\ltimes X$ has object set $X$ and morphisms
$(v,x)\colon x\to v\cdot x$ for $v\in M$ and $x\in X$, so that
$\Hom(x_1,x_2)=\{v\in M:\ v\cdot x_1=x_2\}$. Composition is
$(v',v\cdot x)\circ(v,x)=(v'v,x)$ and the identity of $x$ is $(1,x)$.
\end{definition}

For an $\N$-set we write $\cA(X)=\cA(X,\N^n)$. A morphism
$F=(f,\phi)\colon(X,\N^n)\to(Y,\N^m)$ of $\N$-sets induces a functor
$\cA(F)\colon\cA(X)\to\cA(Y)$, sending the object $x$ to $f(x)$ and the
morphism $(v,x)$ to $(\phi(v),f(x))$. This is a functor if $F$ is \textbf{any} morphism
of $\N$-sets, so every morphism of $\N$-sets induces a functor of action
categories, not only the coordinate morphisms of \cref{def:coord-morphism}.

\subsection{Cohomology of the action category}
\label{sec:action-cohomology}

The cohomology of a category can be defined as the simplicial cohomology of its nerve, which is the
simplicial set $\Nerve\cA(X,M)$ whose $p$-simplices are the strings of $p$
composable arrows. For action categories, we will show that this recovers the monoid
cohomology of \cref{sec:monoid-cohomology}.

\begin{theorem}\label{thm:nerve}
Let a monoid $M$ act on a set $X$ and let $A=\Fun(X,\Z)$ with
$(m\cdot f)(x)=f(m\cdot x)$. Then
$H^p(M;A)\cong H^p\bigl(\Nerve\cA(X,M);\Z\bigr)$ for all $p$.
\end{theorem}

\begin{proof}
A $p$-simplex of $\Nerve\cA(X,M)$ is a string
$x_0\to x_1=m_1\cdot x_0\to\dots\to x_p=m_p\cdot x_{p-1}$, and it is
determined by the pair $[x_0;m_1,\dots,m_p]\in X\times M^p$. The face maps
are
\begin{align*}
  d_0[x_0;m_1,\dots,m_p]&=[m_1\cdot x_0;m_2,\dots,m_p],\\
  d_i[x_0;m_1,\dots,m_p]&=[x_0;m_1,\dots,m_im_{i+1},\dots,m_p]\qquad(0<i<p),\\
  d_p[x_0;m_1,\dots,m_p]&=[x_0;m_1,\dots,m_{p-1}].
\end{align*}
The simplicial cochain complex is $\Fun(X\times M^p,\Z)$ with coboundary
$\delta\varphi=\sum_i(-1)^i\varphi\circ d_i$. We can identify
$\Fun(X\times M^p,\Z)$ with $\Fun(M^p,A)$ by
$f(m_1,\dots,m_p)(x)=\varphi[x;m_1,\dots,m_p]$. Then
\begin{align*}
  (\delta f)(m_1,\dots,m_p)(x)={}&f(m_2,\dots,m_p)(m_1\cdot x)
  +\sum_{i=1}^{p-1}(-1)^if(m_1,\dots,m_im_{i+1},\dots,m_p)(x)\\
  &+(-1)^pf(m_1,\dots,m_{p-1})(x),
\end{align*}
and the first term is $\bigl(m_1\cdot f(m_2,\dots,m_p)\bigr)(x)$. This is
exactly the bar cochain complex \eqref{eq:bar-cochains}.
\end{proof}

\begin{corollary}\label{cor:nerve-deRham}
For every $\N$-set, $\HdR{p}(X,\N^n)\cong H^p(\Nerve\cA(X);\Z)$. In
particular a morphism of $\N$-sets $F\colon(X,\N^n)\to(Y,\N^m)$ induces a
map
\[
  F^\natural\colon\HdR{p}(Y,\N^m)\to\HdR{p}(X,\N^n),
\]
namely the map induced on the cohomology of nerves by $\cA(F)$, and
$(G\circ F)^\natural=F^\natural\circ G^\natural$.
\end{corollary}

\subsection{Homotopy invariance}
\label{sec:homotopy-invariance}

A natural transformation between two functors $\cC\to\cD$ induces a
homotopy between the induced maps of nerves $\Nerve\cC\to\Nerve\cD$, and
homotopic maps induce the same map on cohomology
\citep[\S2]{segal1968}, \citep[\S1, Prop.~2]{quillen1973}. 

Let $F=(f,\phi)$ and $G=(g,\psi)$ be morphisms $(X,\N^n)\to(Y,\N^m)$ of
$\N$-sets. A natural transformation $\eta\colon\cA(F)\Rightarrow\cA(G)$
consists of morphisms $\eta_x\colon f(x)\to g(x)$ in $\cA(Y)$, one for
each $x\in X$. A morphism $f(x)\to g(x)$ is an element $h(x)\in\N^m$ with
$g(x)=h(x)\cdot f(x)$, so $\eta$ is a map $h\colon X\to\N^m$. Naturality
with respect to the morphism $(v,x)\colon x\to v\cdot x$ is the
commutativity of
\[
\begin{tikzcd}[column sep=large,row sep=large]
  f(x)\arrow[r,"h(x)"]\arrow[d,"\phi(v)"'] & g(x)\arrow[d,"\psi(v)"]\\
  f(v\cdot x)\arrow[r,"h(v\cdot x)"'] & g(v\cdot x)
\end{tikzcd}
\]
which says that $h(v\cdot x)+\phi(v)=h(x)+\psi(v)$ in $\N^m$.

\begin{definition}[Homotopy]\label{def:homotopy}
A \emph{homotopy} from $F=(f,\phi)$ to $G=(g,\psi)$ is a map
$h\colon X\to\N^m$ such that $g(x)=h(x)\cdot f(x)$ and
$h(v\cdot x)+\phi(v)=h(x)+\psi(v)$ for all $x\in X$ and $v\in\N^n$. When
$\phi=\psi$ the second condition says that $h$ is constant on orbits.
\end{definition}

\begin{proposition}\label{prop:homotopy}
A homotopy $h$ from $F$ to $G$ defines a natural transformation
$\eta\colon\cA(F)\Rightarrow\cA(G)$ with components $\eta_x=(h(x),f(x))$.
Consequently
$F^\natural=G^\natural$ on de Rham cohomology.
\end{proposition}

\begin{proof}
The composite $\eta_{v\cdot x}\circ\cA(F)(v,x)$ is the morphism
$(h(v\cdot x)+\phi(v),f(x))$ and $\cA(G)(v,x)\circ\eta_x$ is
$(\psi(v)+h(x),f(x))$, so the square commutes exactly when
$h(v\cdot x)+\phi(v)=h(x)+\psi(v)$. The statement about cohomology follows
from \cref{cor:nerve-deRham} and the fact quoted above.
\end{proof}

This notion of homotopy is not symmetric, since $h$ takes values in the
monoid $\N^m$ and cannot be inverted. A homotopy in either direction
suffices for the conclusion that $F$ and $G$ induce the same map on
cohomology. We say that $F\colon X\to Y$ is a
\emph{homotopy equivalence} if there is $G\colon Y\to X$ with homotopies
between $G\circ F$ and $\id_X$, in either direction, and between
$F\circ G$ and $\id_Y$. Then $F^\natural$ is an isomorphism.

\subsection{The discrete Poincar\'e lemma}
\label{sec:poincare}

\begin{theorem}[Discrete Poincar\'e lemma]\label{thm:poincare}
$\HdR{p}(\N^n,\N^n)=0$ for $p\geq 1$, where $\N^n$ acts on itself by
translation.
\end{theorem}

\begin{proof}
Let $P\colon(\N^n,\N^n)\to(*,\N^0)$ be the unique morphism to the terminal
$\N$-set and let $Q\colon(*,\N^0)\to(\N^n,\N^n)$ send $*$ to the origin.
Then $P\circ Q=\id$, and $Q\circ P=(f,\phi)$ with $f(x)=0$ for all $x$ and
$\phi=0$. We claim that $h(x)=x$ is a homotopy from $Q\circ P$ to
$(\id,\id)$. Indeed $h(x)\cdot f(x)=x+0=x$, and
$h(v+x)+\phi(v)=x+v=h(x)+\id(v)$. So
$\HdR{p}(\N^n,\N^n)\cong\HdR{p}(*,\N^0)$, which is $\Z$ in degree $0$ and
zero otherwise.
\end{proof}

\begin{intuition}
The homotopy $h(x)=x$ records that every point $x$ of the quadrant is
reached from the origin by the coordinate action, in $x_1$ steps of $e_1$,
then $x_2$ steps of $e_2$, and so on. It is the discrete counterpart of the
straight-line contraction of $\R^n_{\geq 0}$ to the origin, and the proof
of \cref{thm:poincare-1} integrates along the same path.
\end{intuition}

\section{Toral rigidity}
\label{sec:toral-rigidity}

Let $(X,\N^n)$ be a finite bijective $\N$-set. Its orbits $O_1,\dots,O_c$
are finite, and $O_i\cong\Z^n/\stab(O_i)$ forces $\stab(O_i)$ to have rank
$n$. \Cref{thm:bijective} gives at once:

\begin{lemma}\label{lem:finite-bijective}
Let $(X,\N^n)$ be a finite bijective $\N$-set with $c$ orbits. Then
$\HdR{p}(X,\N^n)\cong\Z^{c\binom{n}{p}}$ for all $p$.
\end{lemma}

To extend this to an arbitrary finite $\N$-set we construct its
\emph{core}, a subset on which the coordinates are bijective, and show
that the $\N$-set contracts onto it. We call the result a localization
theorem because of its formal resemblance to the localization theorem of
Atiyah and Bott~\citep{atiyah-bott1984} in equivariant cohomology. There, for a torus $T$
acting on a compact manifold $M$, the restriction map
$H^\bullet_T(M)\to H^\bullet_T(M^T)$ becomes an isomorphism after inverting
suitable elements of $H^\bullet_T(\mathrm{pt})$, so that the equivariant
cohomology is determined by the fixed locus. The role of the fixed locus is
played here by the eventual image of the coordinate action.

\begin{definition}[Core]\label{def:core}
Let $(X,\N^n)$ be a finite $\N$-set. Its \emph{core} is
$X_\infty=\bigcap_{a\in\N^n}a\cdot X$.
\end{definition}

Write $\vone=(1,\dots,1)\in\N^n$. The sets $N\vone\cdot X$ decrease with $N$
and $X$ is finite, so they stabilize, and since every $a\in\N^n$ satisfies
$a\leq N\vone$ for $N$ large we have $X_\infty=N\vone\cdot X$ for all
sufficiently large $N$. Fix such an $N$. The core is stable under the action,
and the map
\[
  \rho\colon X\to X_\infty,\qquad \rho(x)=N\vone\cdot x,
\]
is a morphism of $\N^n$-sets, since
$\rho(v\cdot x)=(N\vone+v)\cdot x=v\cdot\rho(x)$. On the core the coordinates
are bijective: $e_i\cdot X_\infty=(N\vone+e_i)\cdot X\supseteq(N+1)\vone\cdot X
=X_\infty$, so each $e_i$ acts surjectively on the finite set $X_\infty$
and hence bijectively.

\begin{theorem}[Core localization]\label{thm:core}
Let $(X,\N^n)$ be a finite $\N$-set and let $\iota\colon X_\infty\to X$ be
the inclusion of its core. Then $\iota$ is a homotopy equivalence, and the
restriction of forms induces an isomorphism
\[
  \iota^*\colon\HdR{p}(X,\N^n)\isoto\HdR{p}(X_\infty,\N^n)
  \qquad\text{for every }p.
\]
\end{theorem}

\begin{proof}
Both $\iota$ and $\rho$ are $\N^n$-equivariant maps, so
\cref{rem:two-pullbacks} applies and $\iota^*=\iota^\natural$. The constant
map $h(x)=N\vone$ is a homotopy from $\id_X$ to $\iota\circ\rho$, because
$(\iota\rho)(x)=N\vone\cdot x=h(x)\cdot\id_X(x)$ and $h$ is constant, hence
constant on orbits. The same $h$ restricted to $X_\infty$ is a homotopy
from $\id_{X_\infty}$ to $\rho\circ\iota$. By \cref{prop:homotopy},
$\rho^\natural\iota^\natural=(\iota\rho)^\natural=\id$ and
$\iota^\natural\rho^\natural=(\rho\iota)^\natural=\id$.
\end{proof}

\begin{corollary}[Toral rigidity]\label{cor:toral}
Let $G$ be a finite graph with a coordinate system $\Xs$ of rank $n$, and
let $c$ be the number of connected components of the induced graph. Then
\[
  \HdR{p}(G,\Xs)\cong\Z^{c\binom{n}{p}}\qquad\text{for all }p.
\]
The de Rham cohomology of a finite coordinatized graph is that of $c$
disjoint copies of the $n$-torus, independently of the graph.
\end{corollary}
%Claude note: the draft's Section 3.3 ends with the core localization theorem
%and does not state this corollary, although the introduction announces it.
%I have added it, together with the short argument that the number of orbits
%of the core equals the number of components of the induced graph.

\begin{proof}
By \cref{thm:core} and \cref{lem:finite-bijective},
$\HdR{p}(G,\Xs)\cong\Z^{c'\binom{n}{p}}$ where $c'$ is the number of
$\Z^n$-orbits of the core. It remains to see that $c'=c$. Every $x\in X$
is joined in the induced graph to $\rho(x)=N\vone\cdot x\in X_\infty$, and if
$x$ and $y$ are joined by an edge $y=e_i\cdot x$ then
$\rho(y)=e_i\cdot\rho(x)$ lies in the orbit of $\rho(x)$. So
$x\mapsto\text{orbit of }\rho(x)$ is a well-defined surjection from
components of the induced graph to orbits of the core, and it is injective
because points of the same orbit are joined by edges.
\end{proof}

The Euler characteristic of $\Z^{c\binom{n}{p}}$ is
$c\sum_p(-1)^p\binom{n}{p}=0$, in agreement with \cref{thm:euler}. The
mechanism is the one observed in \cref{sec:polar-grid}. A coordinate field
on a finite set eventually cycles or comes to rest, and the de Rham complex
computes the cohomology of the induced graph, whose cycles and fixed points
account for the toral answer. This result motivates the normalization
introduced in \cref{ch:cubical}, which removes the contribution of the fixed
points, and the atlases of \cref{ch:gluing}, which allow us to deal with more complicated graphs (such as the sphere graph).

\chapter{Cubical sets}
\label{ch:cubical}

The starting point for this chapter, and one forms the central point of the theory from here 
onwards is this: a differential form is a function on cubes.
Under this interpretation the de Rham complex is the
cochain complex of a semi-cubical set. Cube complexes
occur in geometric group theory \citep{sageev1995,wise2012}, and cubical
sets were introduced by Kan~\citep{kan1955} as a combinatorial model for
homotopy theory. Here they serve as models for discrete differential geometry!

\section{Semi-cubical sets}
\label{sec:semi-cubical}

\subsection{\texorpdfstring{$p$}{p}-cubes}
\label{sec:cubes}

Recall that for an $\N$-set $(X,\N^n)$ the $p$-forms are
$\Omega^p(X,\N^n)=\F_X\otimes_\Z\Lam^p\Z^n$, with $\Lam^p\Z^n$ free on the
symbols $dX_I$, $I\subseteq[n]$, $|I|=p$. So a $p$-form is a family
$(\omega_I)_{|I|=p}$ of functions on $X$, which is the same thing as a
single function on the set of pairs $(x,I)$.

\begin{definition}[$p$-cube]\label{def:p-cube}
Let $(X,\N^n)$ be an $\N$-set. A \emph{$p$-cube} is a formal pair
$\cube{x}{I}$ with $x\in X$ and $I\subseteq[n]$, $|I|=p$. We call $x$ the
\emph{base vertex} and $I$ the \emph{index set} of the cube.
\end{definition}

\begin{figure}[htb]
\centering
% What a p-cube [x;I] looks like in X.  (a) a genuine square; (b) a 1-cube
% that is a loop, because e_i x = x; (c) a 2-cube drawn as a triangle,
% because e_i x = e_j x.
\begin{tikzpicture}[x=1cm,y=1cm]
\begin{scope}
  \node[vtx,label={[font=\small]below left:$x$}] (a) at (0,0) {};
  \node[vtx,label={[font=\small]below right:$e_i x$}] (b) at (2.6,0) {};
  \node[vtx,label={[font=\small]above left:$e_j x$}] (c) at (0,2.6) {};
  \node[vtx,label={[font=\small]above right:$e_ie_j x$}] (d) at (2.6,2.6) {};
  \fill[fillcol!14] (0,0) rectangle (2.6,2.6);
  \draw[ar,cA] (a) -- (b); \draw[ar,cA] (c) -- (d);
  \draw[ar,cB] (a) -- (c); \draw[ar,cB] (b) -- (d);
  \node[font=\small,fillcol!70!black] at (1.3,1.3) {$\cube{x}{\{i,j\}}$};
  \node[font=\small] at (1.3,-1.15) {(a) a $2$-cube};
\end{scope}
\begin{scope}[xshift=5.6cm]
  \node[vtx,label={[font=\small]below:$x=e_i x$}] (a) at (1.3,0.6) {};
  \draw[lp,cA] (a) to[out=50,in=130,looseness=12] (a);
  \node[font=\small,fillcol!70!black] at (1.3,2.6) {$\cube{x}{\{i\}}$};
  \node[font=\small] at (1.3,-1.15) {(b) a loop};
\end{scope}
\begin{scope}[xshift=9.4cm]
  \node[vtx,label={[font=\small]below left:$x$}] (a) at (0,0) {};
  \node[vtx,label={[font=\small]below right:$e_i x=e_j x$}] (b) at (2.6,0) {};
  \node[vtx,label={[font=\small]above:$e_ie_j x$}] (d) at (1.3,2.6) {};
  \fill[fillcol!14] (0,0) -- (2.6,0) -- (1.3,2.6) -- cycle;
  \draw[ar,cA] (a) to[bend right=18] (b); \draw[ar,cB] (a) to[bend left=18] (b);
  \draw[ar,cB] (b) to[bend right=18] (d); \draw[ar,cA] (b) to[bend left=18] (d);
  \node[font=\small,fillcol!70!black] at (1.3,0.95) {$\cube{x}{\{i,j\}}$};
  \node[font=\small] at (1.3,-1.15) {(c) a $2$-cube drawn as a triangle};
\end{scope}
\end{tikzpicture}
\caption{A $2$-cube $\cube{x}{\{i,j\}}$ in an $\N$-set, and two degenerate
pictures. When $e_i\cdot x=x$ the $1$-cube is a loop, and when
$e_i\cdot x=e_j\cdot x$ the $2$-cube has a repeated vertex and is drawn as a
triangle with two doubled edges.}
\label{fig:cubes}
\end{figure}
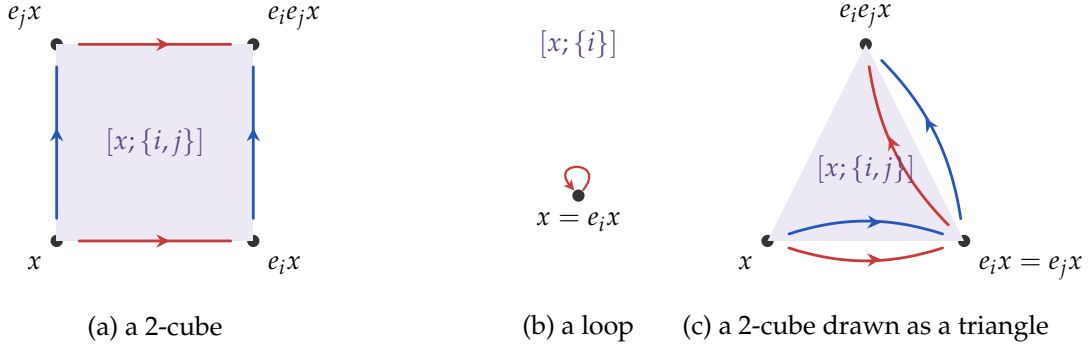

A $p$-cube is a base vertex together with $p$ distinct directions. Because
the coordinate fields commute, the $2^p$ translates $S_J\cdot x$ for
$J\subseteq I$ form the vertices of a cube, as in \cref{fig:cubes}(a). The
pair is said to be "formal" because the actions of the $e_i$ need not be injective, so
the image of a $p$-cube in $X$ need not have $2^p$ distinct vertices. If
$e_i\cdot x=x$ the $1$-cube $\cube{x}{\{i\}}$ is a loop, and if
$e_i\cdot x=e_j\cdot x$ the image of the $2$-cube $\cube{x}{\{i,j\}}$ is a
triangle, as in \cref{fig:cubes}(b,c). The triangle is treated in
\cref{sec:two-simplex}.

We think of a $p$-form $\omega=\sum_{|I|=p}\omega_I\,dX_I$ as the function
on $p$-cubes
\[
  \omega\cube{x}{I}\defeq\omega_I(x).
\]
In this notation $d\omega$ is a function on $(p+1)$-cubes.

\begin{proposition}\label{prop:d-on-cubes}
As a function on $(p+1)$-cubes, with $I=\{i_0<i_1<\dots<i_p\}$,
\[
  (d\omega)\cube{x}{I}=\sum_{j=0}^p(-1)^j\Bigl(
  \omega\cube{e_{i_j}\cdot x}{I\setminus\{i_j\}}-\omega\cube{x}{I\setminus\{i_j\}}\Bigr).
\]
\end{proposition}

\begin{proof}
One writes $d\omega=\sum_{|K|=p+1}(d\omega)_K\,dX_K$ and collects the
terms $\Delta_i(\omega_I)\,dX_i\wedge dX_I$ with $\{i\}\cup I=K$. The
position of $i$ in $K$ produces the sign. The details are in
\cref{app:d-on-cubes}.
\end{proof}

\begin{intuition}[the formula for $d$]
The right-hand side is an alternating sum over the $2(p+1)$ faces of the
cube $\cube{x}{I}$. Removing the direction $i_j$ leaves two opposite faces,
the \emph{front} face $\cube{x}{I\setminus\{i_j\}}$ based at $x$ and the
\emph{back} face $\cube{e_{i_j}\cdot x}{I\setminus\{i_j\}}$ based at
$e_{i_j}\cdot x$, and the formula takes back minus front with the sign
$(-1)^j$. For
$p=0$ the formula reads $(df)\cube{x}{\{i\}}=f(e_i\cdot x)-f(x)=\Delta_if(x)$.
\end{intuition}

It is worth spelling out the similarity to the continuous case of differential forms on 
a manifold, which is what motivated us to make this definition in the first case. If $M$ is a smooth manifold,
a $p$-form is a smooth section of $\Lam^pT^*M$. At every point $x$ the fibre
$(\Lam^pT^*M)_x$ has basis $dx_{i_1}\wedge\dots\wedge dx_{i_p}$, so a
section is a collection of "coefficient" functions $f_I$, and the formal pair
$[x;dx_I]$ may be regarded as an infinitesimal $p$-cube at $x$.
Then it is easy to see that the differential forms are smooth functions on such infinitesimal cubes. We can make this comparison more explicit:

\begin{definition}[Discrete tangent space]\label{def:tangent}
Let $(X,\N^n)$ be an $\N$-set. The \emph{tangent space} at $x\in X$ is
$T_xX=\Z^{\{1,\dots,n\}}\cong\Z^n$, with basis the coordinate fields
$e_1,\dots,e_n$, and the \emph{cotangent space} is $T^*_xX=\Hom_\Z(\Z^n,\Z)$
with dual basis $dX_1,\dots,dX_n$. The tangent bundle is
$TX=\bigsqcup_xT_xX$, and similarly for the cotangent bundle and their
exterior powers $\Lam^pT^*_xX=\Lam^p\Z^n$.
\end{definition}

In this language $\Omega^p(X,\N^n)$ is the set of sections of
$\Lam^pT^*X$, that is, the functions $\omega$ on $X$ with
$\omega(x)\in\Lam^pT_x^*X$ for every $x$, which is precisely
$\bigoplus_{|I|=p}\F_X$ as before. The dictionary between the two
settings is summarized in \cref{tab:dictionary}.

\begin{table}[htb]
\centering
\small
\begin{tabular}{@{}ll@{}}
\toprule
smooth manifold $M$ & $\N$-set $(X,\N^n)$\\
\midrule
point $x\in M$ & vertex $x\in X$\\
tangent space $T_xM$ & $T_xX=\Z^n$\\
basis $\partial/\partial x_1,\dots,\partial/\partial x_n$ & coordinate fields $e_1,\dots,e_n$\\
$\Lam^pT^*_xM$ with basis $dx_I$ & $\Lam^pT^*_xX$ with basis $dX_I$\\
infinitesimal $p$-cube $[x;dx_I]$ & $p$-cube $\cube{x}{I}$\\
$\Omega^p(M)=C^\infty(M,\Lam^pT^*M)$ & $\Omega^p(X,\N^n)=\Fun(X,\Lam^p\Z^n)$\\
$\partial f/\partial x_i$ & $\Delta_if$\\
\bottomrule
\end{tabular}
\caption{The dictionary between smooth and discrete differential forms.}
\label{tab:dictionary}
\end{table}

\subsection{The box category}
\label{sec:box}

\begin{definition}[The box category]\label{def:box}
The \emph{box category} $\Bx$ has as objects the sets $[1]^n=\{0,1\}^n$ for
$n\geq 0$, and its morphisms are generated by the \emph{coface maps}
\[
  \delta^{i}_{\epsilon}\colon[1]^{n-1}\to[1]^n,\quad
  (x_1,\dots,x_{n-1})\mapsto(x_1,\dots,x_{i-1},\epsilon,x_i,\dots,x_{n-1}),
  \quad 1\leq i\leq n,\ \epsilon\in\{0,1\},
\]
which insert $\epsilon$ in position $i$, and the \emph{codegeneracy maps}
\[
  \sigma^i\colon[1]^n\to[1]^{n-1},\qquad
  (x_1,\dots,x_n)\mapsto(x_1,\dots,x_{i-1},x_{i+1},\dots,x_n),
  \qquad 1\leq i\leq n,
\]
which delete the $i$-th coordinate \citep{grandis-mauri2003,husainov2013}.
\end{definition}

Cubical sets can also be endowed with further structure maps called
\emph{connections}, $\gamma^n_{i,0},\gamma^n_{i,1}\colon[1]^n\to[1]^{n-1}$,
which replace the pair $(x_i,x_{i+1})$ by $\max(x_i,x_{i+1})$ or
$\min(x_i,x_{i+1})$ \citep{brown-higgins1981}. We mention them only because
they reappear in \cref{sec:functor-cubical}.

\begin{proposition}[Cocubical identities]\label{prop:cocubical}
The generating maps satisfy
\begin{align*}
  \delta^j_\beta\delta^i_\alpha&=\delta^i_\alpha\delta^{j-1}_\beta && (i<j),\\
  \sigma^j\delta^i_\alpha&=\delta^i_\alpha\sigma^{j-1} && (i<j),\qquad
  \sigma^i\delta^i_\alpha=\id,\qquad
  \sigma^j\delta^i_\alpha=\delta^{i-1}_\alpha\sigma^j\quad(i>j),\\
  \sigma^j\sigma^i&=\sigma^i\sigma^{j+1} && (i\leq j).
\end{align*}
\end{proposition}

\begin{proof}
Each identity is checked by following a tuple through both sides. For the
first, both sides insert $\alpha$ in position $i$ and $\beta$ in position
$j$. For the last, both sides delete the original coordinates $i$ and $j+1$.
The middle ones are similar.
\end{proof}

\begin{definition}[Cubical and semi-cubical objects]\label{def:cubical-set}
A \emph{cubical set} is a functor $\Bx^\op\to\Set$, and more generally a
\emph{cubical object} in a category $\cC$ is a functor $\Bx^\op\to\cC$. Let
$\Bxinj\subseteq\Bx$ be the subcategory with the same objects whose
morphisms are generated by the coface maps alone. A \emph{semi-cubical}, or
pre-cubical, object in $\cC$ is a functor $\Bxinj^\op\to\cC$.
\end{definition}

If $Q$ is a cubical set we write $Q_n=Q([1]^n)$ for the set of $n$-cubes,
$d_i^\epsilon=Q(\delta^i_\epsilon)\colon Q_n\to Q_{n-1}$ for the face maps and
$s_i=Q(\sigma^i)\colon Q_{n-1}\to Q_n$ for the degeneracies. The cubical
identities are the cocubical ones read backwards,
\begin{equation}\label{eq:cubical-identities}
\begin{aligned}
  d_i^\alpha d_j^\beta&=d_{j-1}^\beta d_i^\alpha && (i<j),\\
  d_i^\alpha s_j&=s_{j-1}d_i^\alpha && (i<j),\qquad
  d_i^\alpha s_i=\id,\qquad
  d_i^\alpha s_j=s_jd_{i-1}^\alpha\quad(i>j),\\
  s_is_j&=s_{j+1}s_i && (i\leq j).
\end{aligned}
\end{equation}
The presentation of $\Bx$ by generators and relations used here is the
standard one and can be found in Brown, Higgins and Sivera~\citep{bhs2011} or Cisinski~\citep{cisinski2006};
the variants of $\Bx$ obtained by adding symmetries, connections or
diagonals are due to Brown and Higgins and to Grandis and Mauri. The
conventions we follow, in particular for the differential below, are those
of Doherty et al.~\citep{dkls2024}.

If $A$ is a cubical object in an abelian category, its \emph{unnormalized
chain complex} $Q_\bullet(A)$ has $Q_n(A)=A_n$ and
\begin{equation}\label{eq:cubical-differential}
  \partial_n=\sum_{i=1}^n(-1)^{i-1}(d_i^1-d_i^0)\colon A_n\to A_{n-1}.
\end{equation}

\begin{proposition}\label{prop:cubical-chain}
$\partial_{n-1}\partial_n=0$, so $Q_\bullet(A)$ is a chain complex.
\end{proposition}

\begin{proof}
This is a standard consequence of the first cubical identity
\citep[\S1]{massey1980}. The computation is given in
\cref{app:cubical-chain}.
\end{proof}

The same definition applies to a semi-cubical object, since only the face
maps are used.

\subsection{From \texorpdfstring{$\N$}{N}-sets to semi-cubical sets}
\label{sec:N-to-cubical}

The $p$-cubes of \cref{def:p-cube} assemble into a semi-cubical set.

\begin{definition}\label{def:QX}
Let $(X,\N^n)$ be an $\N$-set. Define $\Qc(X)\colon\Bxinj^\op\to\Set$ by
$\Qc(X)_p=\{\cube{x}{I}:\ x\in X,\ I\subseteq[n],\ |I|=p\}$, with face maps,
for $I=\{i_1<\dots<i_p\}$,
\[
  d^1_j\cube{x}{I}=\cube{e_{i_j}\cdot x}{I\setminus\{i_j\}},\qquad
  d^0_j\cube{x}{I}=\cube{x}{I\setminus\{i_j\}}.
\]
\end{definition}

\begin{proposition}\label{prop:QX-semicubical}
$\Qc(X)$ is a semi-cubical set.
\end{proposition}

\begin{proof}
We check $d_i^\alpha d_j^\beta=d_{j-1}^\beta d_i^\alpha$ for $i<j$. On
$\cube{x}{I}$, the left side removes $i_j$ and then the $i$-th element of
$I\setminus\{i_j\}$, which is $i_i$ since $i<j$. The right side removes
$i_i$ and then the $(j-1)$-st element of $I\setminus\{i_i\}$, which is
$i_j$. Both sides therefore remove $\{i_i,i_j\}$, and the base vertices
agree because the actions of $e_{i_i}$ and $e_{i_j}$ commute.
\end{proof}

We identified $p$-forms with functions on $p$-cubes, that is, with
$\Hom_\Set(\Qc(X)_p,\Z)=\Hom_\Z(\Z\Qc(X)_p,\Z)$. So it is reasonable to
expect the cubical cochain complex of $\Qc(X)$ to be related to the de Rham
complex. Let $\Cc[\bullet](X)$ denote the cubical chain complex of $\Qc(X)$
with integer coefficients, the unnormalized chain complex of the
semi-cubical abelian group $\Z\Qc(X)$, and let $\Cs(X)$ denote its dual.

\begin{theorem}\label{thm:deRham-cubical}
Let $(X,\N^n)$ be an $\N$-set. The assignment $\omega\mapsto\tilde\omega$,
$\tilde\omega\cube{x}{I}=\omega_I(x)$, is an isomorphism of cochain
complexes
\[
  \Omega^\bullet(X,\N^n)\isoto\Cs(X).
\]
In particular the de Rham cohomology of an $\N$-set is the cubical
cohomology of $\Qc(X)$ with integer coefficients.
\end{theorem}

\begin{proof}
By \eqref{eq:cubical-differential} the differential of $\Cc[\bullet](X)$ is
\[
  \partial\cube{x}{\{i_1<\dots<i_p\}}=\sum_{j=1}^p(-1)^{j-1}\Bigl(
  \cube{e_{i_j}\cdot x}{I\setminus\{i_j\}}-\cube{x}{I\setminus\{i_j\}}\Bigr).
\]
Dualizing, a function $\omega$ on $p$-cubes is sent to
$\omega\circ\partial$, which on the $(p+1)$-cube
$\cube{x}{\{i_0<\dots<i_p\}}$ equals
$\sum_{j=0}^p(-1)^j\bigl(\omega\cube{e_{i_j}\cdot x}{I\setminus\{i_j\}}
-\omega\cube{x}{I\setminus\{i_j\}}\bigr)$ after reindexing. This is the
formula of \cref{prop:d-on-cubes} for $d\omega$. The map
$\omega\mapsto\tilde\omega$ is bijective in each degree by construction.
\end{proof}

Every cubical set has a geometric realization $|\Qc|$ in $\Top$, a CW
complex with one $p$-cell for each non-degenerate $p$-cube, and the cubical
cohomology of the cubical set is the cellular cohomology of its realization
\citep{kan1955,massey1980,bhs2011}, which agrees with singular cohomology
\citep{hatcher2002}. Consequently $\HdR{\bullet}(X,\N^n)\cong H^\bullet(|\Qc(X)|;\Z)$,
and the de Rham cohomology of an $\N$-set is the cohomology of a space. The
same holds for the normalized cohomology of \cref{sec:normalized}, with the
degenerate cubes of \cref{def:degenerate} omitted from the realization,
and this is used in \cref{sec:triangulated-torus} and \cref{ch:examples}.
%Claude note: the draft asked for a citation for "cubical cohomology is the
%cellular cohomology of the realization". I have cited Kan, Massey and
%Brown-Higgins-Sivera; please confirm the precise statements you want to
%point to. Note also that for a semi-cubical set, which has no degenerate
%cubes, every cube gives a cell.

\subsection{The Serre diagonal}
\label{sec:serre}

The isomorphism $\Omega^\bullet(X,\N^n)\cong\Cs(X)$ can induce the twisted
product on the cubical cochains. Cubical chains carry a diagonal introduced
by Serre~\citep{serre1951}, whose dual is a cup product on cochains, and we show
that the two products agree. We follow the conventions of
Rivera and Zeinalian~\citep{rivera-zeinalian2020}. Throughout, for
$I\sqcup J=K$ with $I$ and $J$ listed in increasing order,
$\varepsilon(I,J)\in\{\pm1\}$ denotes the shuffle sign defined by
$dX_I\wedge dX_J=\varepsilon(I,J)\,dX_K$. 

\begin{definition}[Serre diagonal]\label{def:serre}
  Let $\Cc[\bullet](Q)$ be the chains of a semi-cubical set. Then the Serre diagonal 
  $\Delta : \Cc[\bullet](Q) \rightarrow \Cc[\bullet](Q) \otimes \Cc[\bullet](Q) $ is induced by
  \[
  \Delta(\sigma) = \sum \epsilon(I,J) d^{1}_{i_{1}} ... d^{1}_{i_{p}}(\sigma) \otimes d^{0}_{j_{1}} ... d^{0}_{j_{q}} (\sigma)
  \]
where the sum runs over all $I = {i_{1},...,i_{p}}$ and $J = (j_{1},...,j_{q})$ such that $I \sqcup J = [n]$.

Let $(X,\N^n)$ be an $\N$-set and $\Cc[\bullet](X)$ its cubical chain
complex. Recall the shift $S_J$ and the
twisted product $(f\,dX_I)\cdot(g\,dX_J)=S_J(f)\,g\,dX_I\wedge dX_J$. The Serre diagonal, 
where $Q(X)$ is the semi-cubical set of $(X,\mathbb{N}^{n})$, becomes
$\Delta\colon\Cc[\bullet](X)\to\Cc[\bullet](X)\otimes\Cc[\bullet](X)$,
\[
  \Delta\cube{x}{K}=\sum_{I\sqcup J=K}\varepsilon(I,J)\,
  \cube{S_J\cdot x}{I}\otimes\cube{x}{J},
\]
the sum over all ordered pairs of disjoint subsets with union $K$.

\end{definition}
%Claude note: conventions for the Serre diagonal differ between sources (which
%face is "front" and which is "back", and where the shift sits). The formula
%here is the one from your draft and it matches the twisted product as
%stated; please check it against Rivera-Zeinalian, whom you cite for it.

The map $\Delta$ is a chain map, making $(\Cc[\bullet](X),\partial,\Delta)$
a differential graded coalgebra. Let $\smile$ be the product on
$\Cs(X)=\Hom(\Cc[\bullet](X),\Z)$ dual to $\Delta$, so that
$(\varphi\smile\psi)(c)=\sum\varphi(c')\psi(c'')$ when
$\Delta(c)=\sum c'\otimes c''$. This makes $(\Cs(X),d,\smile)$ a DGA.

\begin{theorem}\label{thm:serre-dga}
The isomorphism $\omega\mapsto\tilde\omega$ of \cref{thm:deRham-cubical} is
an isomorphism of differential graded algebras
$(\Omega^\bullet(X,\N^n),d,\cdot)\isoto(\Cs(X),d,\smile)$.
\end{theorem}

\begin{proof}
By \cref{thm:deRham-cubical} the map is a cochain map, so it remains to
show that it respects products.
Let $\omega\in\Omega^p$, $\eta\in\Omega^q$, and let $K\subseteq[n]$ with
$|K|=p+q$. On the cubical side,
\[
  (\tilde\omega\smile\tilde\eta)\cube{x}{K}=\sum_{I\sqcup J=K}\varepsilon(I,J)\,
  \tilde\omega\cube{S_J\cdot x}{I}\,\tilde\eta\cube{x}{J}
  =\sum_{\substack{I\sqcup J=K\\ |I|=p,\ |J|=q}}\varepsilon(I,J)\,
  \omega_I(S_J\cdot x)\,\eta_J(x),
\]
the terms with $|I|\neq p$ vanishing because $\omega$ is homogeneous. On
the form side, expanding $\omega=\sum_I\omega_I\,dX_I$ and
$\eta=\sum_J\eta_J\,dX_J$,
\[
  \omega\cdot\eta=\sum_{I,J}S_J(\omega_I)\,\eta_J\,dX_I\wedge dX_J
  =\sum_K\Bigl(\sum_{I\sqcup J=K}\varepsilon(I,J)\,S_J(\omega_I)\,\eta_J\Bigr)dX_K,
\]
since $dX_I\wedge dX_J$ vanishes unless $I$ and $J$ are disjoint and equals
$\varepsilon(I,J)\,dX_{I\sqcup J}$ when they are. Evaluating the coefficient
of $dX_K$ at $x$ gives $\widetilde{\omega\cdot\eta}\cube{x}{K}
=(\tilde\omega\smile\tilde\eta)\cube{x}{K}$.
\end{proof}

\begin{intuition}
The twisted product of \cref{sec:algebra-structure} was introduced to
obtain a Leibniz rule for the difference operators. \Cref{thm:serre-dga}
identifies it with the cup product of cubical cochains. The shift $S_J$ in
the definition of the product corresponds to the fact that the cup product
evaluates one factor on a front face of the cube and the other on the
opposite back face, whose base vertex is $S_J\cdot x$.
\end{intuition}

\section{Degeneracies and the normalized complex}
\label{sec:normalized}

\subsection{The left Kan extension}
\label{sec:kan-extension}

The inclusion $\Bxinj\to\Bx$ induces a restriction functor
$\PSh(\Bx)\to\PSh(\Bxinj)$, and its left adjoint, the left Kan extension
$J_!\colon\PSh(\Bxinj)\to\PSh(\Bx)$, freely adds degeneracies
\citep{husainov2013}. Applied to the semi-cubical set $\Qc(X)$ of an
$\N$-set it produces a cubical set which we also denote $\Qc(X)$ when no
confusion can arise. Its $p$-cubes are the formal pairs $\cube{x}{I}$ where
now $I=(i_1,\dots,i_p)$ is a $p$-tuple with entries in $[n]\cup\{*\}$ whose
entries in $[n]$ are distinct and increasing. We call $*$ the
\emph{degenerate coordinate}. The degeneracy $s_k\colon\Qc(X)_{p-1}\to\Qc(X)_p$
inserts $*$ in position $k$,
\[
  s_k\cube{x}{(i_1,\dots,i_{p-1})}=\cube{x}{(i_1,\dots,i_{k-1},*,i_k,\dots,i_{p-1})}.
\]
It is convenient to let the degenerate coordinate act trivially,
$e_*\cdot x=x$. Then the face maps are given by the same formula as before,
$d^0_k\cube{x}{I}=\cube{x}{I\setminus i_k}$ and
$d^1_k\cube{x}{I}=\cube{e_{i_k}\cdot x}{I\setminus i_k}$, where
$I\setminus i_k$ deletes the $k$-th entry.

\subsection{Functoriality}
\label{sec:functor-cubical}

We ask whether $(X,\N^n)\mapsto\Qc(X)$ is a functor from $\Nset$ to
$\PSh(\Bx)$. Given a morphism $(f,\phi)$, the natural candidate is
$\cube{x}{(i_1,\dots,i_p)}\mapsto\cube{f(x)}{(\phi(e_{i_1}),\dots,\phi(e_{i_p}))}$,
which only makes sense if each $\phi(e_i)$ is a generator or zero. This is
the same restriction that functoriality of forms imposed in
\cref{sec:functoriality}, and it leads again to the coordinate morphisms
$(f,\sigma)$ with $\sigma\colon[n]\to[m]\cup\{*\}$. If moreover $\sigma$ is
increasing on $\sigma^{-1}([m])$, then
$\cube{x}{(i_1,\dots,i_p)}\mapsto\cube{f(x)}{(\sigma(i_1),\dots,\sigma(i_p))}$
is a morphism of cubical sets.

Two features of coordinate morphisms are not captured by this. A coordinate
morphism may send two coordinates to the same one, and it need not be
increasing. Both can be accommodated by enlarging the box category, the
first by connections and the second by passing to symmetric cubical sets,
presheaves on the box category with permutations adjoined
\citep{isaacson2009}. Neither enlargement is needed here. Coordinate
morphisms act on the cochain DGAs by \cref{thm:functor} whether or not they
act on the cubical sets, and in the gluing constructions of \cref{ch:gluing}
the transition maps are isomorphisms of semi-cubical sets, for which the
coordinates can be relabelled so that $\sigma$ is increasing.

\subsection{The normalized complex of an \texorpdfstring{$\N$}{N}-set}
\label{sec:normalized-complex}

Let $Q$ be a cubical set. Its \emph{degenerate $p$-cubes} are those of the
form $s_iq'$ for some $i$ and some $q'\in Q_{p-1}$. They do not form a
sub-cubical set, since $d_i^\alpha s_i=\id$ shows that a face of a
degenerate cube need not be degenerate, but they do span a subcomplex of
the chain complex.

\begin{proposition}\label{prop:degenerate-subcomplex}
Let $A$ be a cubical object in an abelian category and let $D_n\subseteq A_n$
be the image of $\sum_is_i\colon\bigoplus_iA_{n-1}\to A_n$. Then
$\partial(D_n)\subseteq D_{n-1}$.
\end{proposition}

\begin{proof}
For $a\in A_{n-1}$ split $\partial s_ia=\sum_j(-1)^{j-1}(d^1_j-d^0_j)s_ia$
according to the position of $j$ relative to $i$. The term $j=i$ vanishes
since $d^1_is_i=d^0_is_i=\id$. For $j<i$ the identities
\eqref{eq:cubical-identities} give $d^\alpha_js_i=s_{i-1}d^\alpha_j$, and
for $j>i$ they give $d^\alpha_js_i=s_id^\alpha_{j-1}$. Hence
\[
  \partial s_ia=s_{i-1}\Bigl(\sum_{j<i}(-1)^{j-1}(d^1_j-d^0_j)a\Bigr)
  +s_i\Bigl(\sum_{j>i}(-1)^{j-1}(d^1_{j-1}-d^0_{j-1})a\Bigr)\in D_{n-1}.\qedhere
\]
\end{proof}

\begin{definition}[Normalized complex]\label{def:normalized}
The \emph{normalized chain complex} of a cubical object $A$ is the quotient
$N_\bullet(A)=Q_\bullet(A)/D_\bullet(A)$.
\end{definition}

There is a short exact sequence of complexes
$0\to D_\bullet\to Q_\bullet\to N_\bullet\to 0$, and the degenerate
subcomplex is acyclic, so $Q_\bullet$ and $N_\bullet$ have the same
homology \citep{massey1980,bhs2011}. For the degeneracies freely added by
$J_!$ the normalized complex therefore computes the same cohomology as the
unnormalized one. The normalization becomes relevant when further
degeneracies are imposed.

\begin{intuition}[loops as degenerate cubes]
When a graph is modelled by an $\N$-set the coordinate fields may fix a
vertex, and each fixed point contributes a loop that the graph does not
have. The computations of \cref{sec:polar-grid} and \cref{sec:path-square}
show that these loops carry cohomology classes. In a cubical set a $1$-cube
whose two faces coincide is the image of a $0$-cube under a degeneracy, so
the cubes containing a fixed direction are declared degenerate and the
normalized complex is formed. The remainder of this section carries this
out.
\end{intuition}

\begin{definition}\label{def:collapsed-cubical}
Let $\Qc(X)$ be the cubical set of \cref{sec:kan-extension}. Let $\approx$
be the equivalence relation on each $\Qc(X)_p$ generated by
\[
  \cube{x}{(i_1,\dots,i_k,\dots,i_p)}\approx\cube{x}{(i_1,\dots,*,\dots,i_p)}
  \qquad\text{whenever }e_{i_k}\cdot x=x,
\]
and put $\overline{\Qc}(X)_p=\Qc(X)_p/\!\approx$.
\end{definition}

\begin{proposition}\label{prop:collapsed-cubical}
The face and degeneracy maps of $\Qc(X)$ descend to $\overline{\Qc}(X)$,
which is therefore a cubical set.
\end{proposition}

\begin{proof}
The degeneracies insert a dummy entry and are compatible with $\approx$. For
the faces, let $I'$ be $I$ with its $k$-th entry $i_k$ replaced by $*$,
where $e_{i_k}\cdot x=x$. Then $d^0_k\cube{x}{I}=\cube{x}{I\setminus i_k}
=d^0_k\cube{x}{I'}$, and
$d^1_k\cube{x}{I}=\cube{e_{i_k}\cdot x}{I\setminus i_k}=\cube{x}{I\setminus i_k}
=d^1_k\cube{x}{I'}$ because $e_{i_k}$ fixes $x$. For $j\neq k$ the faces
$d^\alpha_j\cube{x}{I}$ and $d^\alpha_j\cube{x}{I'}$ are again related by
$\approx$, since the base vertex is $x$ or $e_{i_j}\cdot x$ and
$e_{i_k}$ fixes $e_{i_j}\cdot x$ as well, the actions commuting.
\end{proof}

The $\N$-sets below serve as charts on simple graphs, which have no loops,
and a cube with an index fixing its base vertex corresponds to a loop of the
induced graph. All such cubes are declared degenerate.

\begin{definition}[Degenerate cubes of an $\N$-set]\label{def:degenerate}
Let $(X,\N^n)$ be an $\N$-set. A cube $\cube{x}{I}$ is \emph{degenerate} if
some entry of $I$ is the dummy coordinate $*$, or if some entry $i\in I$
satisfies $e_i\cdot x=x$. Otherwise it is \emph{non-degenerate}.
\end{definition}

Only the second clause affects the cohomology, since the cubes with a dummy
entry are removed by the normalization of \cref{def:normalized} without
changing it. From now on we therefore work with the semi-cubical set
$\Qc(X)$ of \cref{def:QX}. Let $\Dc[p](X)\subseteq\Cc[p](X)$ be the span of the
degenerate $p$-cubes $\cube{x}{I}$, $I\subseteq[n]$, in the sense of the
second clause.

\begin{proposition}\label{prop:D-subcomplex}
$\Dc[\bullet](X)$ is a subcomplex of $\Cc[\bullet](X)$.
\end{proposition}

\begin{proof}
Let $\cube{x}{I}$ be degenerate with $e_i\cdot x=x$ for some $i\in I$. In
$\partial\cube{x}{I}$, the two faces obtained by removing $i$ are
$\cube{e_i\cdot x}{I\setminus\{i\}}=\cube{x}{I\setminus\{i\}}$ and they
cancel because they are subtracted from each other in the sum. Every other face still contains $i$ in its index set, and its base
vertex is $x$ or $e_j\cdot x$, both of which are fixed by $e_i$. So every
face surviving in the same is degenerate.
\end{proof}

\begin{definition}[Normalized de Rham complex]\label{def:normalized-deRham}
Let $(X,\N^n)$ be an $\N$-set. Its \emph{normalized chain complex} is
$\Nc[\bullet](X)=\Cc[\bullet](X)/\Dc[\bullet](X)$, the free abelian group on
the non-degenerate cubes with the induced differential. Its
\emph{normalized cochain complex}, or normalized de Rham complex, is
\[
  \Ns(X)=\Hom\bigl(\Nc[\bullet](X),\Z\bigr)\cong
  \bigl\{\omega\in\Omega^\bullet(X,\N^n):\ \omega_I(x)=0\text{ if }
  e_i\cdot x=x\text{ for some }i\in I\bigr\},
\]
the forms vanishing on degenerate cubes, with the de Rham differential. Its
cohomology is the \emph{normalized de Rham cohomology} $\HN{p}(X,\N^n)$.
\end{definition}

Dualizing the short exact sequence
$0\to\Dc[\bullet](X)\to\Cc[\bullet](X)\to\Nc[\bullet](X)\to0$ of free
abelian groups gives a short exact sequence of cochain complexes
\[
  0\to\Ns(X)\to\Cs(X)\xrightarrow{\ \rho\ }\Ds(X)\to 0,
\]
where $\Ds(X)=\Hom(\Dc[\bullet](X),\Z)$ is the complex of functions on
degenerate cubes and $\rho$ is restriction, and hence a long exact
sequence
\begin{equation}\label{eq:les}
  \cdots\to H^{p-1}(\Ds(X))\xrightarrow{\ \partial\ }\HN{p}(X)\to\HdR{p}(X)
  \xrightarrow{\ \rho\ }H^p(\Ds(X))\xrightarrow{\ \partial\ }\HN{p+1}(X)\to\cdots
\end{equation}
relating the normalized and unnormalized cohomologies.

\begin{proposition}\label{prop:H0-normalized}
$\HN{0}(X,\N^n)\cong\HdR{0}(X,\N^n)$.
\end{proposition}

\begin{proof}
There are no degenerate $0$-cubes, so $\Ds[0](X)=0$ and the sequence
\eqref{eq:les} begins $0\to\HN{0}\to\HdR{0}\to 0$.
\end{proof}

\section{Examples}
\label{sec:normalized-examples}

\subsection{The path graph and the square graph}
\label{sec:path-square-normalized}

\emph{The path graph.} Let $X=\{1,2,3\}$ with the single coordinate field
of \cref{sec:path-square}, so $n=1$ and $e_1\cdot1=2$, $e_1\cdot2=3$,
$e_1\cdot3=3$. The only fixed point is $3$, so the only degenerate cube is
$\cube{3}{\{1\}}$.

\begin{proposition}\label{prop:path-normalized}
$\HN{\bullet}(X,\N)=(\Z,0)$.
\end{proposition}

\begin{proof}
The normalized complex is $0\to N^0\to N^1\to 0$ with $N^0=\F_X\cong\Z^3$
and $N^1$ the functions on the non-degenerate $1$-cubes $\cube{1}{\{1\}}$
and $\cube{2}{\{1\}}$, so $N^1\cong\Z^2$. The differential is
$(a,b,c)\mapsto(b-a,c-b)$, which is surjective with kernel the constants.
\end{proof}

The normalized cohomology of the path graph is that of a point. The loop
at $3$, which carried the class in $\HdR{1}$ in \cref{sec:path-square}, is
degenerate.

\emph{The square graph.} Let $X=\{1,2\}\times\{1,2\}$ with the two
coordinate fields of \cref{sec:path-square}, $e_1$ moving right and $e_2$
moving up, each fixing the vertices from which the move would leave the
square,
\[
  e_1\cdot(1,j)=(2,j),\quad e_1\cdot(2,j)=(2,j),\quad
  e_2\cdot(i,1)=(i,2),\quad e_2\cdot(i,2)=(i,2).
\]
The fixed points are $\Fix(e_1)=\{(2,1),(2,2)\}$ and
$\Fix(e_2)=\{(1,2),(2,2)\}$. The degenerate $1$-cubes are the four loops
$\cube{(2,1)}{\{1\}}$, $\cube{(2,2)}{\{1\}}$, $\cube{(1,2)}{\{2\}}$,
$\cube{(2,2)}{\{2\}}$, and every $2$-cube is degenerate except
$\cube{(1,1)}{\{1,2\}}$. So the normalized complex is
\[
  0\to\Z^4\xrightarrow{\ d^0\ }\Z^4\xrightarrow{\ d^1\ }\Z\to 0,
\]
with $\Z^4$ vertices, $8-4=4$ non-degenerate $1$-cubes, and one
non-degenerate $2$-cube. The kernel of $d^0$ is the constants, since the
$\N$-set is transitive, so $\HN{0}=\Z$. The map $d^1$ sends
$\omega$ to $\omega\cube{(2,1)}{2}-\omega\cube{(1,1)}{2}
-\omega\cube{(1,2)}{1}+\omega\cube{(1,1)}{1}$, which is surjective, so
$\HN{2}=0$. The quotient $N^1/\im d^0$ is the first simplicial cohomology of
the square graph, which is $\Z$ by \cref{prop:simplicial-H0}, and $d^1$
induces a surjection from it onto $\Z$, hence an isomorphism. So
$\HN{1}=\ker d^1/\im d^0=0$ and $\HN{\bullet}(X,\N^2)=(\Z,0,0)$.

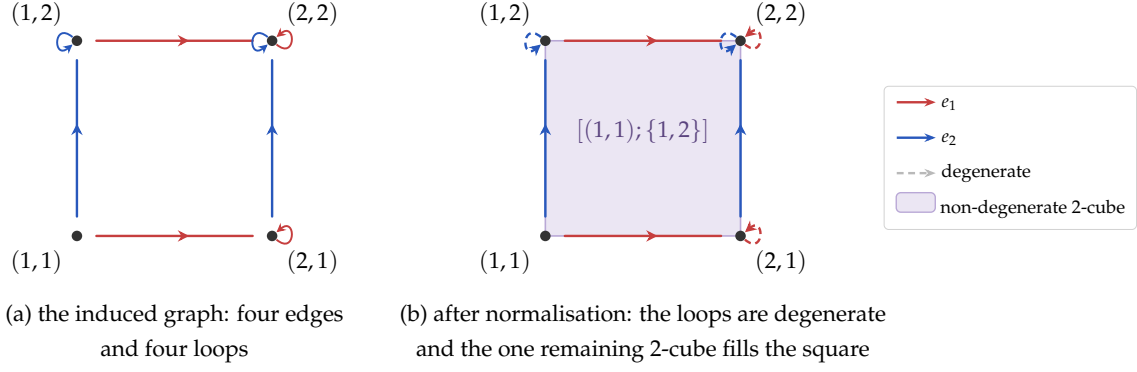
\begin{figure}[htb]
\centering
\resizebox{\textwidth}{!}{\begingroup
\newcommand{\board}[1]{%
  \node[vtx] (a) at (0,0) {};  \node[vtx] (b) at (3,0) {};
  \node[vtx] (c) at (0,3) {};  \node[vtx] (d) at (3,3) {};
  \node[font=\small,below left =1pt of a] {$(1,1)$};
  \node[font=\small,below right=1pt of b] {$(2,1)$};
  \node[font=\small,above left =1pt of c] {$(1,2)$};
  \node[font=\small,above right=1pt of d] {$(2,2)$};
  \draw[ar,e1col] (a) -- (b);
  \draw[ar,e1col] (c) -- (d);
  \draw[ar,e2col] (a) -- (c);
  \draw[ar,e2col] (b) -- (d);
  \draw[#1,e1col] (b) to[out=-40,in=40,looseness=9] (b);
  \draw[#1,e1col] (d) to[out=-40,in=40,looseness=9] (d);
  \draw[#1,e2col] (c) to[out=140,in=220,looseness=9] (c);
  \draw[#1,e2col] (d) to[out=140,in=220,looseness=9] (d);
}
\begin{tikzpicture}[x=1cm,y=1cm]
\begin{scope}
  \board{lp}
  \node[font=\small,align=center] at (1.5,-1.5)
    {(a) the induced graph: four edges\\ and four loops};
\end{scope}
\begin{scope}[xshift=7.2cm]
  \fill[fillcol!16] (0,0) -- (3,0) -- (3,3) -- (0,3) -- cycle;
  \draw[fillcol!55,line width=0.7pt] (0,0) -- (3,0) -- (3,3) -- (0,3) -- cycle;
  \node[fillcol!70!black,font=\small] at (1.5,1.62) {$[(1,1);\{1,2\}]$};
  \board{dlp}
  \node[font=\small,align=center] at (1.5,-1.5)
    {(b) after normalisation: the loops are degenerate\\
     and the one remaining $2$-cube fills the square};
\end{scope}
\begin{scope}[xshift=12.4cm,yshift=1.2cm]
  \node[anchor=west,font=\scriptsize,align=left,inner sep=5pt,
        draw=gray!35,rounded corners=2pt] at (0,0) {%
    \tikz{\draw[e1col,line width=1.05pt,-{Stealth[length=5pt]}](0,0)--(0.6,0);}~~$e_1$\\[3pt]
    \tikz{\draw[e2col,line width=1.05pt,-{Stealth[length=5pt]}](0,0)--(0.6,0);}~~$e_2$\\[3pt]
    \tikz{\draw[gray!55,densely dashed,line width=1.05pt,-{Stealth[length=5pt]}](0,0)--(0.6,0);}%
      ~~degenerate\\[3pt]
    \tikz{\fill[fillcol!16](0,-0.12) rectangle (0.6,0.18);
          \draw[fillcol!55,line width=0.6pt](0,-0.12) rectangle (0.6,0.18);}%
      ~~non-degenerate $2$-cube};
\end{scope}
\end{tikzpicture}
\endgroup}
\caption{The square graph with its rank-two coordinate system, before and
after normalization. The four loops become degenerate and the single
non-degenerate $2$-cube fills the square.}
\label{fig:square-filled}
\end{figure}

\begin{intuition}
\Cref{fig:square-filled} shows the effect of normalization on the square
graph. Before normalization the induced graph has four loops and the
cohomology of a torus. After normalization the loops are degenerate and
the one non-degenerate $2$-cube $\cube{(1,1)}{\{1,2\}}$ is a $2$-cell whose
boundary is the square, so the normalized complex computes the cohomology
of a filled square. In general the rank of the coordinate system is the
dimension of the cells it supplies.
\end{intuition}

\subsection{The polar grid}
\label{sec:polar-normalized}

Let $\PG(m)$ be the polar grid of \cref{sec:polar-grid} with the coordinate
system $\{R,\Theta\}$. Both fix the origin and nothing else, so
$\Fix(R)=\Fix(\Theta)=\{\mathbf 0\}$. The degenerate $1$-cubes are
$\cube{\mathbf 0}{\{R\}}$ and $\cube{\mathbf 0}{\{\Theta\}}$, the only
degenerate $2$-cube is $\cube{\mathbf 0}{\{R,\Theta\}}$, and
\cref{def:normalized-deRham} gives
\[
  N^0=\F,\qquad N^1=\{g\,dR+h\,d\Theta:\ g(\mathbf 0)=h(\mathbf 0)=0\},\qquad
  N^2=\{f\,dR\wedge d\Theta:\ f(\mathbf 0)=0\}.
\]

\begin{proposition}\label{prop:polar-normalized}
$\HN{\bullet}(\PG(m),\{R,\Theta\})=(\Z,0,0)$.
\end{proposition}

\begin{proof}
We use the long exact sequence \eqref{eq:les}. The complex $\Ds(\PG(m))$ of
functions on degenerate cubes is $0\to0\to\Z^2\xrightarrow{d^1}\Z\to0$, and
$d^1$ is zero: for $f$ a function on the two degenerate $1$-cubes,
\[
  (d^1f)\cube{\mathbf 0}{\{R,\Theta\}}
  =\bigl(f\cube{R\cdot\mathbf 0}{\Theta}-f\cube{\mathbf 0}{\Theta}\bigr)
  -\bigl(f\cube{\Theta\cdot\mathbf 0}{R}-f\cube{\mathbf 0}{R}\bigr)=0,
\]
since $R\cdot\mathbf 0=\Theta\cdot\mathbf 0=\mathbf 0$. So
$H^\bullet(\Ds)=(0,\Z^2,\Z)$. By \cref{prop:polar-torus},
$\HdR{\bullet}(\PG(m))=(\Z,\Z^2,\Z)$, and \eqref{eq:les} reads
\[
  0\to\HN{0}\to\Z\to0\to\HN{1}\to\Z^2\xrightarrow{\rho^1}\Z^2\to\HN{2}\to
  \Z\xrightarrow{\rho^2}\Z\to0.
\]
The maps $\rho^p$ restrict a closed form to the degenerate cubes. In degree
$1$ the classes $\delta_{\mathbf 0}\,dR$ and $d\Theta$ of
\cref{prop:polar-torus} restrict to the two basis vectors of $\Ds[1]$, so
$\rho^1$ is an isomorphism. In degree $2$ the generator
$\delta_{\mathbf 0}\,dR\wedge d\Theta$ restricts to $1$ on
$\cube{\mathbf 0}{\{R,\Theta\}}$, so $\rho^2$ is an isomorphism. Exactness
then forces $\HN{1}=0$ and $\HN{2}=0$.
\end{proof}

The normalized cohomology of the polar grid is that of a disc. The
coordinate system supplies $2$-cells for the squares between consecutive
rings and for the triangles around the origin, and the loops at the origin,
which carried the classes of \cref{prop:polar-torus}, are degenerate. In
\cref{ex:globe-mv} two polar grids glued along their boundary circles give
a sphere.

\subsection{The \texorpdfstring{$2$}{2}-simplex}
\label{sec:two-simplex}

A coordinate system of rank $2$ can also fill a triangle. The triangle
appears at the centre of the polar grid, and it is the simplest case in
which a cube $\cube{x}{I}$ and its image in $X$ differ.

Let $X=\{x,a,b\}$ with two coordinate fields,
\[
  e_1\cdot x=a,\quad e_1\cdot a=a,\quad e_1\cdot b=b,\qquad
  e_2\cdot x=b,\quad e_2\cdot a=b,\quad e_2\cdot b=b.
\]
The two fields commute, since $e_1e_2\cdot x=b=e_2e_1\cdot x$,
$e_1e_2\cdot a=b=e_2e_1\cdot a$ and $b$ is fixed by both. The induced graph
is the triangle on $x,a,b$ with loops at $a$ and $b$.

\begin{proposition}\label{prop:simplex}
$\HN{\bullet}(X,\N^2)=(\Z,0,0)$.
\end{proposition}

\begin{proof}
The non-degenerate cubes are the three vertices, the $1$-cubes
$\cube{x}{\{1\}}$, $\cube{x}{\{2\}}$ and $\cube{a}{\{2\}}$, the other
$1$-cubes being loops, and the single $2$-cube $\cube{x}{\{1,2\}}$, since
$a$ and $b$ are fixed by $e_1$. Writing a $1$-cochain as
$(\alpha,\beta,\gamma)$ in this order, the complex is
\[
  0\to\Z^3\xrightarrow{d^0}\Z^3\xrightarrow{d^1}\Z\to0,
\]
with $d^0(f_x,f_a,f_b)=(f_a-f_x,\ f_b-f_x,\ f_b-f_a)$ and
$d^1(\alpha,\beta,\gamma)=\alpha-\beta+\gamma$.
Here $d^1$ is computed from \cref{prop:d-on-cubes}, the face
$\cube{b}{\{1\}}$ being degenerate and contributing $0$. So $d^1$ is
surjective, $\ker d^1=\{\beta=\alpha+\gamma\}\cong\Z^2$, and
$\im d^0=\ker d^1$ since $f_b-f_x=(f_a-f_x)+(f_b-f_a)$. Hence
$\HN{\bullet}=(\Z,0,0)$.
\end{proof}

The triangle is filled by the single $2$-cube $\cube{x}{\{1,2\}}$, as the
square was filled in \cref{sec:path-square-normalized}, one edge of the
square, the face $\cube{b}{\{1\}}$, being collapsed to a point. Collapsing
an edge of a square to obtain a triangle is the standard comparison between
cubical and simplicial sets, and it underlies the equivalence of their
homotopy theories \citep{cisinski2006}. The triangle is used as a chart in
\cref{sec:triangulated-torus} and \cref{ch:examples}, where triangulated
surfaces are covered by one such chart per triangle.
%Claude note: the draft cited "nLab, cubulation" for the Quillen equivalence
%between simplicial and cubical sets; I have cited Cisinski's Asterisque
%volume instead. Please confirm this is the reference you want.

\chapter{Gluing \texorpdfstring{$\N$}{N}-sets and sheaves}
\label{ch:gluing}

\section{The sheaf of differential forms}
\label{sec:sheaf-of-forms}

\subsection{The need for gluing: two charts on a sphere}
\label{sec:globe}

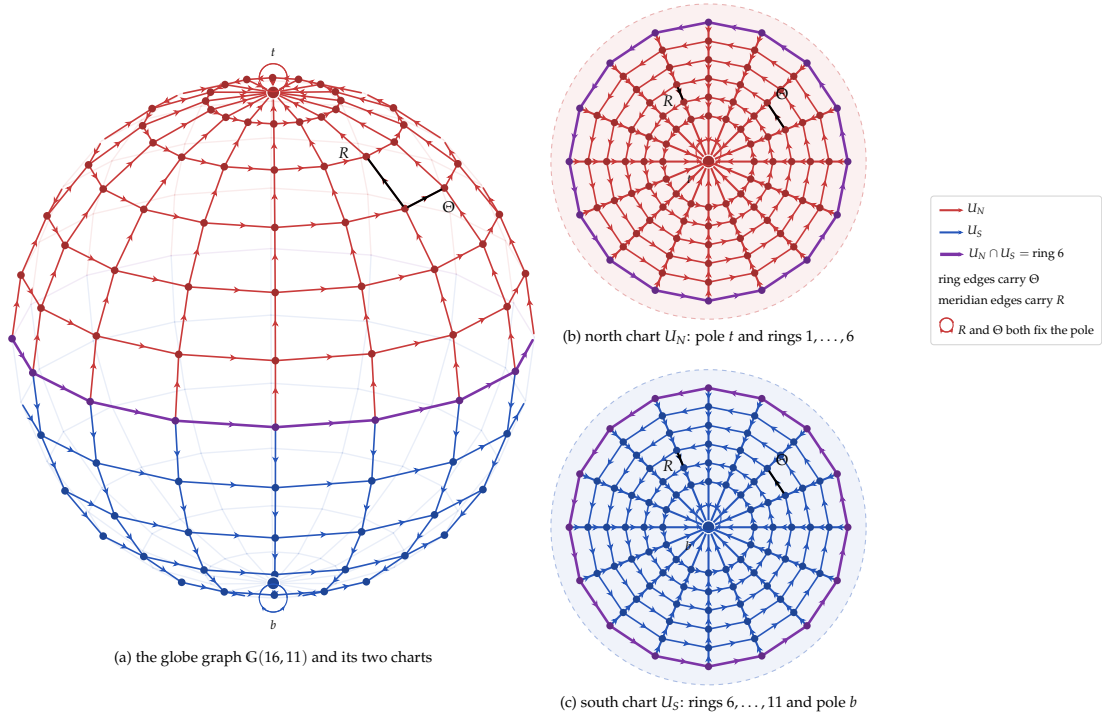
\begin{figure}[htb]
\centering
\resizebox{\textwidth}{!}{\begingroup
% ---------------- parameters ----------------------------------------------
\def\Nmer{16}      % meridians  n
\def\Mlat{11}      % latitude rings  m
\def\Eeq{6}        % the shared ring
\def\az{22}        % view azimuth
\def\el{20}        % view elevation
\def\Rad{6}        % radius of panel (a), in cm
\def\rmin{1.05}    % inner radius of the flat charts
\def\rmax{3.20}    % outer radius of the flat charts

% ---------------- projection ----------------------------------------------
% \sph{i}{u} sets \px,\py (screen, already scaled) and \pd (depth).
% u = 0 is the north pole t, u = \Mlat+1 is the south pole b.
\newcommand{\sph}[2]{%
  \pgfmathsetmacro{\phg}{#2*180/(\Mlat+1)}%
  \pgfmathsetmacro{\thg}{360*#1/\Nmer}%
  \pgfmathsetmacro{\XX}{sin(\phg)*cos(\thg)}%
  \pgfmathsetmacro{\YY}{sin(\phg)*sin(\thg)}%
  \pgfmathsetmacro{\ZZ}{cos(\phg)}%
  \pgfmathsetmacro{\px}{\Rad*(-\XX*sin(\az)+\YY*cos(\az))}%
  \pgfmathsetmacro{\py}{\Rad*(-(\XX*cos(\az)+\YY*sin(\az))*sin(\el)+\ZZ*cos(\el))}%
  \pgfmathsetmacro{\pd}{(\XX*cos(\az)+\YY*sin(\az))*cos(\el)+\ZZ*sin(\el)}%
}

\newif\iffrontpass

% \gedge{i1}{u1}{i2}{u2}{colour}{width} -- the arrow runs from 1 to 2.
\newcommand{\gedge}[6]{%
  \sph{#1}{#2}\edef\ax{\px}\edef\ay{\py}\edef\adep{\pd}%
  \sph{#3}{#4}\edef\bx{\px}\edef\by{\py}\edef\bdep{\pd}%
  \pgfmathsetmacro{\mdep}{0.5*(\adep+\bdep)}%
  \ifdim\mdep pt>0pt
    \iffrontpass \draw[#5,line width=#6,cedge] (\ax,\ay) -- (\bx,\by);\fi
  \else
    \iffrontpass\else
      \draw[#5!45,line width=0.9pt,line cap=round,opacity=0.22] (\ax,\ay) -- (\bx,\by);
    \fi
  \fi
}

% every edge of the globe, once
\newcommand{\globeedges}{%
  % ring edges: Theta, i -> i+1
  \foreach \u in {1,...,\Mlat}{%
    \pgfmathtruncatemacro{\isN}{\u<\Eeq}\pgfmathtruncatemacro{\isE}{\u==\Eeq}%
    \edef\ecol{\ifnum\isE=1 polE\else\ifnum\isN=1 polN\else polS\fi\fi}%
    \edef\ewid{\ifnum\isE=1 1.7pt\else 1.05pt\fi}%
    \foreach \i in {0,...,\the\numexpr\Nmer-1\relax}{%
      \pgfmathtruncatemacro{\j}{mod(\i+1,\Nmer)}%
      \expandafter\gedge\expandafter{\i}{\u}{\j}{\u}{\ecol}{\ewid}}}%
  % meridian edges: R, towards the chart's own pole
  \foreach \u in {0,...,\Mlat}{%
    \pgfmathtruncatemacro{\up}{\u+1}%
    \pgfmathtruncatemacro{\isN}{\up<\Eeq+1}%
    \edef\ecol{\ifnum\isN=1 polN\else polS\fi}%
    \foreach \i in {0,...,\the\numexpr\Nmer-1\relax}{%
      \ifnum\isN=1 \gedge{\i}{\up}{\i}{\u}{\ecol}{1.05pt}%
      \else        \gedge{\i}{\u}{\i}{\up}{\ecol}{1.05pt}\fi}}%
}

% ---------------- a flattened polar chart ----------------------------------
% \cappanel{N or S}{cx}{cy}.  Both caps are drawn looking DOWN THE POLAR AXIS
% FROM THE NORTH, so Theta has the same sense in the two pictures and the
% shared ring looks identical in both.
\newcommand{\cappanel}[3]{%
  \def\which{#1}%
  \ifx\which\Nletter
    \pgfmathtruncatemacro{\kk}{\Eeq-1}\def\ccol{polN}\def\cpole{t}%
  \else
    \pgfmathtruncatemacro{\kk}{\Mlat-\Eeq}\def\ccol{polS}\def\cpole{b}%
  \fi
  \begin{scope}[shift={(#2,#3)}]
    \pgfmathsetmacro{\rout}{\rmax+0.42}
    \fill[\ccol!7] (0,0) circle (\rout);
    \draw[\ccol!40,dashed,line width=0.7pt] (0,0) circle (\rout);
    % rings: Theta
    \foreach \a in {0,...,\kk}{%
      \pgfmathsetmacro{\rr}{\rmin+(\rmax-\rmin)*\a/\kk}%
      \pgfmathtruncatemacro{\isout}{\a==\kk}%
      \edef\ecol{\ifnum\isout=1 polE\else\ccol\fi}%
      \edef\ewid{\ifnum\isout=1 1.7pt\else 1.05pt\fi}%
      \foreach \i in {0,...,\the\numexpr\Nmer-1\relax}{%
        \pgfmathsetmacro{\ta}{360*\i/\Nmer}\pgfmathsetmacro{\tb}{360*(\i+1)/\Nmer}%
        \draw[\ecol,line width=\ewid,cedge]
          (\ta:\rr) -- (\tb:\rr);}}
    % spokes: R, inward to the chart's own pole
    \foreach \a in {1,...,\kk}{%
      \pgfmathsetmacro{\ro}{\rmin+(\rmax-\rmin)*\a/\kk}%
      \pgfmathsetmacro{\ri}{\rmin+(\rmax-\rmin)*(\a-1)/\kk}%
      \foreach \i in {0,...,\the\numexpr\Nmer-1\relax}{%
        \pgfmathsetmacro{\ta}{360*\i/\Nmer}%
        \draw[\ccol,line width=1.35pt,cedge] (\ta:\ro) -- (\ta:\ri);}}
    \foreach \i in {0,...,\the\numexpr\Nmer-1\relax}{%
      \pgfmathsetmacro{\ta}{360*\i/\Nmer}%
      \draw[\ccol,line width=1.35pt,cedge,shorten >=3.4pt] (\ta:\rmin) -- (0,0);}
    % vertices
    \foreach \a in {0,...,\kk}{%
      \pgfmathsetmacro{\rr}{\rmin+(\rmax-\rmin)*\a/\kk}%
      \pgfmathtruncatemacro{\isout}{\a==\kk}%
      \edef\ecol{\ifnum\isout=1 polE\else\ccol\fi}%
      \foreach \i in {0,...,\the\numexpr\Nmer-1\relax}{%
        \pgfmathsetmacro{\ta}{360*\i/\Nmer}%
        \fill[\ecol!80!black] (\ta:\rr) circle (2.4pt);}}
    \fill[white] (0,0) circle (4.6pt);
    \fill[\ccol!80!black] (0,0) circle (3.8pt);
    \draw[poleloop,\ccol] (-0.30,0.30) arc (162:-150:0.30);
    \draw[poleloop,\ccol] (0.30,0.30) arc (18:330:0.30);
    \node[font=\scriptsize,anchor=north east] at (-0.26,-0.18) {$\cpole$};
    % labelled Theta and R
    \pgfmathsetmacro{\rmid}{\rmin+(\rmax-\rmin)*int(\kk/2)/\kk}
    \pgfmathsetmacro{\rmidb}{\rmin+(\rmax-\rmin)*(int(\kk/2)-1)/\kk}
    \pgfmathsetmacro{\ta}{360*1/\Nmer}\pgfmathsetmacro{\tb}{360*2/\Nmer}
    \draw[callout] (\ta:\rmid) -- (\tb:\rmid);
    \node[font=\small] at ($(\tb:\rmid)+(0.34,0.20)$) {$\Theta$};
    \pgfmathsetmacro{\tc}{360*5/\Nmer}
    \draw[callout] (\tc:\rmid) -- (\tc:\rmidb);
    \node[font=\small] at ($(\tc:\rmidb)+(-0.36,0.04)$) {$R$};
  \end{scope}
}
\def\Nletter{N}

% ===========================================================================
\begin{tikzpicture}[x=1cm,y=1cm,line join=round]

% ---- panel (a): the globe -------------------------------------------------
\frontpassfalse \globeedges          % far side, faint
\frontpasstrue  \globeedges          % near side, with arrowheads

\foreach \u in {1,...,\Mlat}{%
  \pgfmathtruncatemacro{\isN}{\u<\Eeq}\pgfmathtruncatemacro{\isE}{\u==\Eeq}%
  \edef\vcol{\ifnum\isE=1 polE\else\ifnum\isN=1 polN\else polS\fi\fi}%
  \foreach \i in {0,...,\the\numexpr\Nmer-1\relax}{%
    \sph{\i}{\u}%
    \ifdim\pd pt>0pt \fill[\vcol!80!black] (\px,\py) circle (2.5pt);\fi}}

\sph{0}{0}\fill[white] (\px,\py) circle (4.6pt);
          \fill[polN!80!black] (\px,\py) circle (3.8pt);
          \pgfmathsetmacro{\yy}{\py+0.34}
          \draw[poleloop,polN] (\px-0.32,\yy) arc (180:-145:0.32);
          \draw[poleloop,polN] (\px+0.32,\yy) arc (0:325:0.32);
          \node[font=\scriptsize,anchor=south] at (\px,\yy+0.36) {$t$};
\pgfmathtruncatemacro{\bpole}{\Mlat+1}
\sph{0}{\bpole}\fill[white] (\px,\py) circle (4.6pt);
          \fill[polS!80!black] (\px,\py) circle (3.8pt);
          \pgfmathsetmacro{\yy}{\py-0.34}
          \draw[poleloop,polS] (\px-0.32,\yy) arc (180:-145:0.32);
          \draw[poleloop,polS] (\px+0.32,\yy) arc (0:325:0.32);
          \node[font=\scriptsize,anchor=north] at (\px,\yy-0.36) {$b$};

% labelled Theta and R on the sphere itself
\sph{3}{3}\edef\ax{\px}\edef\ay{\py}
\sph{4}{3}\edef\bx{\px}\edef\by{\py}
\draw[callout] (\ax,\ay) -- (\bx,\by);
\node[font=\small] at ($(\bx,\by)+(0.10,-0.42)$) {$\Theta$};
\sph{3}{3}\edef\ax{\px}\edef\ay{\py}
\sph{3}{2}\edef\bx{\px}\edef\by{\py}
\draw[callout] (\ax,\ay) -- (\bx,\by);
\node[font=\small] at ($(\bx,\by)+(-0.52,0.10)$) {$R$};

\node[font=\small] at (0,-7.35)
  {(a) the globe graph $\mathbb{G}(\Nmer,\Mlat)$ and its two charts};

% ---- panels (b), (c): the flattened charts --------------------------------
\cappanel{N}{10.0}{4.05}
\node[font=\small] at (10.0,0.0)
  {(b) north chart $U_N$: pole $t$ and rings $1,\dots,\Eeq$};

\cappanel{S}{10.0}{-4.35}
\node[font=\small] at (10.0,-8.40)
  {(c) south chart $U_S$: rings $\Eeq,\dots,\Mlat$ and pole $b$};

% ---- legend ---------------------------------------------------------------
\node[font=\scriptsize,anchor=west,align=left,inner sep=5pt,
      draw=gray!35,rounded corners=2pt,fill=white] at (15.1,1.55) {%
  \tikz{\draw[polN,line width=1.35pt,-{Stealth[length=3.4pt]}] (0,0)--(0.55,0);}~~$U_N$\\[2pt]
  \tikz{\draw[polS,line width=1.35pt,-{Stealth[length=3.4pt]}] (0,0)--(0.55,0);}~~$U_S$\\[2pt]
  \tikz{\draw[polE,line width=2.2pt,-{Stealth[length=3.8pt]}] (0,0)--(0.55,0);}%
    ~~$U_N\cap U_S=$ ring $\Eeq$\\[5pt]
  ring edges carry $\Theta$\\[1pt]
  meridian edges carry $R$\\[5pt]
  \tikz{\draw[polN,line width=0.75pt,-{Stealth[length=3.2pt]}]
        (-0.16,0.16) arc (162:-150:0.16);\draw[polN,line width=0.75pt,-{Stealth[length=3.2pt]}] (0.16,0.16) arc (18:330:0.16);}~~~$R$ and $\Theta$ both fix the pole};

\end{tikzpicture}
\endgroup}
\caption{The globe graph $\mathbb G(16,11)$ and its cover by two polar
charts, $U_N$ in red and $U_S$ in blue, overlapping in the ring $6$ drawn
in purple. Ring edges carry $\Theta$ and meridian edges carry $R$, which
runs towards the pole of its own chart. Both coordinates fix the pole of
their chart, so each pole carries two loops.}
\label{fig:globe}
\end{figure}

Consider the \emph{globe graph} $\mathbb G(n,m)$ of \cref{fig:globe}(a),
with $n$ meridians and $m$ latitude rings together with a north pole $t$
and a south pole $b$. It is a discrete sphere. The north pole together with
the rings $1,\dots,E$ is a truncated polar grid, and so is the south pole
together with the rings $E,\dots,m$. Each carries the coordinate system
$\{R,\Theta\}$ of \cref{sec:polar-grid}, with $R$ running inwards towards
its own pole and $\Theta$ rotating the rings, and the two charts overlap in
the ring $E$, on which they agree as far as $\Theta$ is concerned.
\Cref{fig:globe}(b,c) shows the two charts flattened out.

A single polar chart cannot cover the globe graph. In a polar chart the
radial field runs towards the pole of the chart, and a chart containing
both poles would have to run towards the north on the upper half and
towards the south on the lower half, whereas a coordinate field has a
single value on each meridian edge. At least two charts are therefore
needed, together with a rule for gluing them. With the definitions of this
chapter, \cref{ex:globe-mv} shows that the normalized cohomology of the
globe graph with this atlas is $(\Z,0,\Z)$.

\subsection{What to glue over}
\label{sec:what-to-glue}

When gluing local pieces $(U_i)_{i\in I}$ to construct a space one has
smaller pieces $U_{ij}\subseteq U_i$ and isomorphisms $U_{ij}\cong U_{ji}$
along which to identify. In the example above the charts are glued along
the ring $E$. One could define a sub-$\N$-set of $(X,\N^n)$ to be a subset
$V\subseteq X$ and glue over subsets. Gluing over a subset, however,
identifies every edge spanned by that subset, and the following example
shows that this is too restrictive.

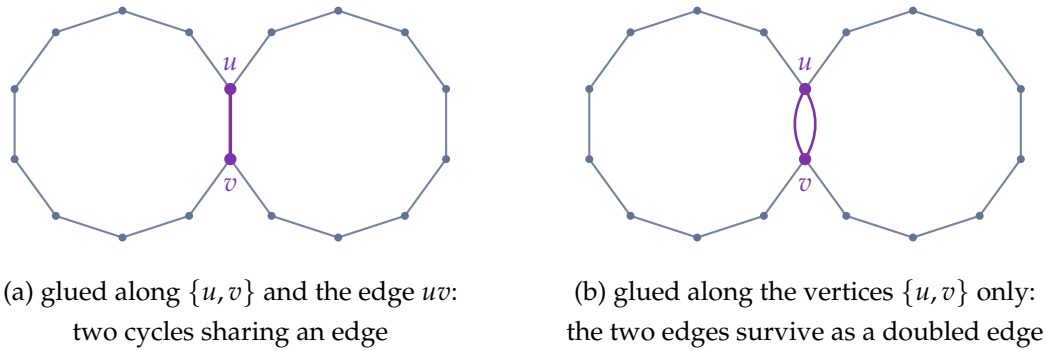
\begin{figure}[htb]
\centering
% Two copies of C_10 glued along two adjacent vertices u,v.
% (a) gluing over the sub-object spanned by u, v and the edge uv: the two
%     edges are identified and the result is two cycles sharing an edge.
% (b) gluing over the 0-dimensional sub-object {u,v} only: the two edges
%     survive as a doubled edge.
\begin{tikzpicture}[x=1cm,y=1cm]
\newcommand{\twocycles}[1]{%
  % circles of radius 1.5 centred at (-+1.4266,0); they meet at (0,+-0.4635),
  % which are vertices of both cycles (angles +-18 degrees).
  \foreach \i in {0,...,9}{
    \pgfmathsetmacro{\ang}{18+36*\i}
    \coordinate (L\i) at ($(-1.4266,0)+(\ang:1.5)$);
    \pgfmathsetmacro{\angb}{198+36*\i}
    \coordinate (R\i) at ($(1.4266,0)+(\angb:1.5)$);}
  % L0 = R9 = u (top), L9 = R0 = v (bottom)
  \foreach \i in {0,...,8}{\pgfmathtruncatemacro{\j}{\i+1}
    \draw[gline,line width=0.8pt] (L\i) -- (L\j);
    \draw[gline,line width=0.8pt] (R\i) -- (R\j);}
  #1
  \foreach \i in {1,...,8}{\node[gvtx] at (L\i) {}; \node[gvtx] at (R\i) {};}
  \node[vtx,fill=cE] at (L0) {}; \node[vtx,fill=cE] at (L9) {};
  \node[font=\small,cE] at ($(L0)+(0,0.34)$) {$u$};
  \node[font=\small,cE] at ($(L9)+(0,-0.34)$) {$v$};
}
\begin{scope}
  \twocycles{\draw[cE,line width=1.3pt] (L0) -- (L9);}
  \node[font=\small,align=center] at (0,-2.5) {(a) glued along $\{u,v\}$ and the edge $uv$:\\ two cycles sharing an edge};
\end{scope}
\begin{scope}[xshift=7.6cm]
  \twocycles{\draw[cE,line width=1.0pt] (L0) to[bend right=30] (L9); \draw[cE,line width=1.0pt] (L0) to[bend left=30] (L9);}
  \node[font=\small,align=center] at (0,-2.5) {(b) glued along the vertices $\{u,v\}$ only:\\ the two edges survive as a doubled edge};
\end{scope}
\end{tikzpicture}
\caption{Two copies of $C_{10}$ glued along two adjacent vertices. The
result depends on whether the edge between them is glued as well.}
\label{fig:glue-c10}
\end{figure}

\begin{example}\label{ex:two-cycles}
Take two copies of the cycle graph $C_{10}$ and choose in each a pair of
adjacent vertices $u,v$. If we glue the copies along $\{u,v\}$ together
with the edge $uv$, the two edges are identified and the result is two
cycles sharing an edge, as in \cref{fig:glue-c10}(a). If we glue along the
vertices $\{u,v\}$ only, the two edges remain distinct and the result has a
doubled edge, as in \cref{fig:glue-c10}(b). The two results have different
first cohomology, of ranks $2$ and $3$.
\end{example}

The overlaps are therefore taken to be sub-objects of the cubical set, which
record the identified cubes of \textit{every} dimension and not only the vertices.

\begin{definition}[Cubical subset]\label{def:cubical-subset}
A \emph{cubical subset} of a semi-cubical set $Q$ is a sub-functor
$S\subseteq Q$, that is, a family of subsets $S_p\subseteq Q_p$ closed
under the face maps.
\end{definition}

Cubical subsets can also be described as the open sets of a topology. Sheaves
on discrete spaces such as simplicial complexes and posets are developed
extensively by Curry~\citep{curry2014}, and we use the notions developed there.

\begin{definition}[Face poset]\label{def:face-poset}
Let $Q$ be a semi-cubical set. Its \emph{face poset} $\Fp(Q)$ is the set of
all cubes of $Q$, of all dimensions, with $c\preceq c'$ if $c$ can be
obtained from $c'$ by a finite sequence of face maps. We write
$U^{(p)}=U\cap Q_p$ for the set of $p$-cubes of a subset $U\subseteq\Fp(Q)$.
\end{definition}

For an $\N$-set, $\cube{x}{\{i\}}\preceq\cube{x}{\{i,j\}}$ while
$\cube{x}{\{i,j\}}\not\preceq\cube{x}{\{i\}}$. Had we used the cubical set
with degeneracies of \cref{sec:kan-extension}, and allowed degeneracy maps
as well as faces in the definition, the relation would only be a preorder,
since $\cube{x}{\{i,j\}}$ and $\cube{x}{(i,j,*)}$ would each be obtained
from the other. We work with the semi-cubical set, where $\preceq$ is a
partial order.

\begin{definition}[Alexandrov topology]\label{def:alexandrov}
The \emph{Alexandrov topology} on a preordered set $(A,\preceq)$ is the
topology whose open sets are the downward closed subsets $U\subseteq A$,
those with $a\in U$ whenever $a\preceq b$ and $b\in U$
\citep{curry2014}.
\end{definition}

Arbitrary unions and intersections of downward closed sets are downward
closed, so this is a topology. A subset of $\Fp(Q)$ is open exactly when it
is closed under faces, so the open subsets of $\Fp(Q)$ are exactly the cubical
subsets of $Q$, and the intersection of two cubical subsets is a cubical
subset. Every cube $q$ has a smallest open neighbourhood $U_q$, the set of
all its iterated faces. Differential forms are functions on cubes, that is,
on the points of $\Fp(Q)$, and can therefore be described as a sheaf.

\subsection{The sheaves of forms}
\label{sec:sheaves}

Fix an $\N$-set $(X,\N^n)$, let $\Qc(X)$ be its semi-cubical set and
$\Fp(X)=\Fp(\Qc(X))$ its face poset with the Alexandrov topology.

\begin{definition}[The three sheaves]\label{def:three-sheaves}
For $U\subseteq\Fp(X)$ open put
\begin{itemize}
\item $\Cc^{p}(U)=\Fun(U^{(p)},\Z)$, the functions on the $p$-cubes of $U$;
\item $\Dc^{p}(U)\subseteq\Cc^{p}(U)$, the functions supported on the
  degenerate $p$-cubes of $U$, those $\cube{x}{I}\in U$ with $e_i\cdot x=x$
  for some $i\in I$;
\item $\Nc^{p}(U)=\Cc^{p}(U)/\Dc^{p}(U)$, equivalently the functions on the
  non-degenerate $p$-cubes of $U$.
\end{itemize}
The restriction maps are restriction of functions.
\end{definition}

Each of these is a sheaf of abelian groups on $\Fp(X)$. The sheaf axiom for
$\Cc^{p}$ is the statement that a function on $\bigcup_\alpha U_\alpha^{(p)}$
is the same as a compatible family of functions on the $U_\alpha^{(p)}$,
and $\Dc^{p}$ and $\Nc^{p}$ are the sub- and quotient sheaves of functions
supported on, respectively defined on, the degenerate and non-degenerate
cubes. Since an open set is closed under faces, the de Rham differential
$d\omega=\omega\circ\partial$ of \cref{thm:deRham-cubical} makes sense on
$\Cs(U)$, it preserves $\Ds(U)$ by
\cref{prop:D-subcomplex}, and it commutes with restriction. So
$\Cs$, $\Ds$ and $\Ns$ are sheaves of cochain complexes on $\Fp(X)$, and
there is a short exact sequence of sheaves of complexes
\[
  0\to\Ns\to\Cs\to\Ds\to0.
\]
The global sections $\Cs(\Fp(X))$ and $\Ns(\Fp(X))$ are the de Rham complex
and the normalized de Rham complex of \cref{ch:cubical}. The cup product of
\cref{sec:serre} also makes sense on $\Cs(U)$, since it only involves faces
of the cube it is evaluated on, so $\Cs$ is a sheaf of DGAs.
%Claude note: the draft states that D and N are sheaves of DGAs as well and
%that the sequence is one of DGAs. I do not think this is right: D is not an
%ideal for the cup product, because a face of a non-degenerate cube can be
%degenerate. In the 2-simplex of Section 4.3.3 the cube [x;{1,2}] is
%non-degenerate while its face [b;{1}] is a loop, so the cup product of a
%form supported on degenerate cubes with an arbitrary form can be non-zero
%on [x;{1,2}]. I have therefore only claimed complexes of sheaves, which is
%all that is used later. Please check.

\begin{remark}[Why cubes and not points]\label{rem:cubes-not-points}
Here the discrete construction departs from the smooth one. An open subset
of a manifold is a set of points, and it determines which infinitesimal
cubes lie in it, since an infinitesimal translate of a point of $U$ is
still in $U$. This is not rigorous, but the intuition is that since $U$ is open, it consists of 
all directions and all "infinitesimal translations". In the smooth setting an open set and the cubical sub-object
it spans therefore carry the same information. In the discrete setting they
do not, as \cref{ex:two-cycles} shows, two cubical subsets having the same
$0$-cubes. For this reason the points of the space $\Fp(X)$ are the cubes
rather than the vertices.
\end{remark}

\section{Atlases of \texorpdfstring{$\N$}{N}-sets}
\label{sec:atlas}

\subsection{Gluing \texorpdfstring{$\N$}{N}-sets}
\label{sec:gluing-N-sets}

We continue the dictionary with manifolds. A manifold structure on a set
$M$ is given by charts, bijections $\varphi\colon U\to\varphi(U)\subseteq\R^n$
onto open sets, covering $M$, with smooth transition maps
$\varphi_\alpha\circ\varphi_\beta^{-1}$ on the overlaps
$\varphi_\beta(U_\alpha\cap U_\beta)$. Equivalently, a manifold is given by
a family of open sets $U_\alpha\subseteq\R^n$, open subsets
$U_{\alpha\beta}\subseteq U_\alpha$ with $U_{\alpha\alpha}=U_\alpha$, and
diffeomorphisms $\varphi_{\alpha\beta}\colon U_{\alpha\beta}\to U_{\beta\alpha}$
with $\varphi_{\alpha\alpha}=\id$ and the cocycle condition
$\varphi_{\beta\gamma}\circ\varphi_{\alpha\beta}=\varphi_{\alpha\gamma}$ on
$U_{\alpha\beta}\cap U_{\alpha\gamma}$. The manifold is then
$M=\bigsqcup_\alpha U_\alpha/\!\sim$, with $\sim$ generated by
$x\sim\varphi_{\alpha\beta}(x)$, and the cocycle condition ensures that
each $U_\alpha$ maps injectively into $M$. We do the same with the
semi-cubical sets of $\N$-sets, gluing along cubical subsets.

\begin{definition}[Discrete coordinatized space]\label{def:dcs}
Let $(X_\alpha,\N^{n_\alpha})_{\alpha\in A}$ be a family of $\N$-sets and
write $U_\alpha=\Fp(\Qc(X_\alpha))$. \emph{Gluing data} for the family
consists of open subsets $U_{\alpha\beta}\subseteq U_\alpha$ and
isomorphisms of semi-cubical sets
$\tau_{\alpha\beta}\colon U_{\alpha\beta}\to U_{\beta\alpha}$ such that
\begin{enumerate}
\item $U_{\alpha\alpha}=U_\alpha$ and $\tau_{\alpha\alpha}=\id$;
\item $\tau_{\alpha\beta}(U_{\alpha\beta}\cap U_{\alpha\gamma})
       =U_{\beta\alpha}\cap U_{\beta\gamma}$ and
       $\tau_{\alpha\gamma}=\tau_{\beta\gamma}\circ\tau_{\alpha\beta}$ on
       $U_{\alpha\beta}\cap U_{\alpha\gamma}$, for all $\alpha,\beta,\gamma$.
\end{enumerate}
%Claude note: the first half of (ii) is implicit in the draft (it is needed
%for the cocycle condition to make sense, and is standard); I have made it
%explicit.
The resulting \emph{discrete coordinatized space} $X$ is the semi-cubical
set $\Qc(X)$ with $\Qc(X)_p=\bigsqcup_\alpha U_\alpha^{(p)}/\!\sim$, where
$\sim$ is the equivalence relation generated by $c\sim\tau_{\alpha\beta}(c)$
for $c\in U_{\alpha\beta}$, and with face maps induced from those of the
$U_\alpha$. We say that $X$ \emph{has the atlas} $(X_\alpha,\N^{n_\alpha})$,
and we write $\Fp(X)$ for its face poset with the Alexandrov topology.
\end{definition}

The face maps are well defined because the $\tau_{\alpha\beta}$ are
morphisms of semi-cubical sets: if $c\sim\tau_{\alpha\beta}(c)$ then
$d_i^\epsilon(\tau_{\alpha\beta}(c))=\tau_{\alpha\beta}(d_i^\epsilon c)\sim
d_i^\epsilon c$. The cocycle condition ensures, exactly as for manifolds,
that each $U_\alpha$ maps injectively into $\Fp(X)$, and its image is open
since it is closed under faces. So $\Fp(X)=\bigcup_\alpha U_\alpha$ is
covered by the open sets $U_\alpha$, each of which is the face poset of an
$\N$-set.

\begin{lemma}\label{lem:degenerate-intrinsic}
A cube $\cube{x}{I}$ of an $\N$-set is degenerate if and only if two of its
opposite faces coincide, $d^0_k\cube{x}{I}=d^1_k\cube{x}{I}$ for some $k$.
Consequently isomorphisms of semi-cubical sets between face posets of
$\N$-sets preserve degenerate cubes, and the degenerate cubes of a discrete
coordinatized space are well defined.
\end{lemma}

\begin{proof}
$d^1_k\cube{x}{I}=\cube{e_{i_k}\cdot x}{I\setminus\{i_k\}}$ equals
$d^0_k\cube{x}{I}=\cube{x}{I\setminus\{i_k\}}$ if and only if
$e_{i_k}\cdot x=x$.
\end{proof}

\begin{intuition}
An atlas of $\N$-sets is the discrete counterpart of an atlas of charts.
Each chart is a subgraph on which $\N^{n_\alpha}$ acts, with coordinate
lines the orbits of the generators. The overlaps are cubical subsets, which
specify the shared edges and squares as well as the shared vertices, and
the transition maps are isomorphisms of semi-cubical sets, which may
relabel the coordinates of one chart in terms of another. The ranks
$n_\alpha$ of the charts need not agree. The globe of \cref{fig:globe} has
two charts of rank $2$ glued along a cubical subset of rank $1$, the ring
with its $\Theta$-edges.
\end{intuition}

\subsection{Gluing sheaves}
\label{sec:gluing-sheaves}

Given gluing data $(U_\alpha,U_{\alpha\beta},\varphi_{\alpha\beta})$ for a
manifold $M$, the sheaf of forms on $M$ is obtained by gluing the sheaves
$\Omega^\bullet_\alpha$ on the $U_\alpha$. The transition maps induce
isomorphisms $\varphi_{\alpha\beta}^*$ of the restricted sheaves over the
overlaps, satisfying the cocycle condition, and for $U\subseteq M$ open
\[
  \Omega^\bullet(U)=\Bigl\{(\omega_\alpha)_\alpha:\ \omega_\alpha\in
  \Omega^\bullet_\alpha(U\cap U_\alpha),\
  \varphi_{\alpha\beta}^*\omega_\beta=\omega_\alpha\text{ on }U\cap U_\alpha\cap U_\beta\Bigr\}
\]
consists of the families of local forms that agree on overlaps.

Gluing in the discrete case is simpler, because differential forms are
functions on the points of the space, and gluing sheaves of functions
produces the sheaf of functions on the glued space.

\begin{proposition}\label{prop:glued-sheaves}
Let $X$ be a discrete coordinatized space with atlas
$(X_\alpha,\N^{n_\alpha})$. Gluing the sheaves $\Cs$, $\Ds$, $\Ns$ of
\cref{def:three-sheaves} along the $\tau_{\alpha\beta}$ produces sheaves of
cochain complexes $\Cs$, $\Ds$, $\Ns$ on $\Fp(X)$, and these are the
sheaves of the glued semi-cubical set: $\Cc^{p}(U)$ is the set of functions
on the $p$-cubes of $U$ and $\Nc^{p}(U)$ the set of functions on the
non-degenerate ones.
\end{proposition}

\begin{proof}
By definition of gluing,
$\Cc^{p}(U)=\{(f_\alpha)\in\prod_\alpha\Fun(U^{(p)}\cap U_\alpha^{(p)},\Z):\
f_\beta\circ\tau_{\alpha\beta}=f_\alpha\text{ on }U\cap U_{\alpha\beta}\}$.
Since $\Fp(X)=\bigsqcup U_\alpha/\!\sim$ with $c\sim\tau_{\alpha\beta}(c)$,
such a family is exactly a function on $U^{(p)}$, namely $c\mapsto f_\alpha(c)$
for $c\in U_\alpha$, which is well defined by the compatibility condition.
The same argument applies to $\Ds$ and $\Ns$, using
\cref{lem:degenerate-intrinsic} to see that the notion of degenerate cube
is compatible with the gluing.
\end{proof}

We define the \emph{normalized de Rham cohomology} of a discrete
coordinatized space $X$ to be $\HN{p}(X)=H^p(\Ns(\Fp(X)))$, the cohomology
of the global sections of the sheaf of normalized forms. For a single chart
this is the normalized cohomology of \cref{def:normalized-deRham}.

\section{Sheaf cohomology}
\label{sec:sheaf-cohomology}

\subsection{The constant sheaf and its resolution}
\label{sec:resolution}

For a manifold the de Rham cohomology is the sheaf cohomology of the
constant sheaf. Recall the definition. If $\mathcal G$ is a sheaf of abelian
groups on a topological space $Y$ then
$H^i(Y,\mathcal G)=R^i\Gamma(Y,-)(\mathcal G)$ is the $i$-th right derived
functor of global sections, computed from an injective resolution
$0\to\mathcal G\to\mathcal I^\bullet$ as $H^i(\Gamma(Y,\mathcal I^\bullet))$.
Any resolution by acyclic sheaves, sheaves with vanishing higher cohomology,
computes the same groups, and flasque sheaves are acyclic
\citep[III.2]{hartshorne1977}.

\begin{theorem}[de Rham]\label{thm:smooth-deRham}
For a smooth manifold $M$, $\HdR{p}(M;\R)\cong H^p(M,\underline{\R})$.
\end{theorem}

\begin{proof}[Sketch]
The sequence of sheaves $0\to\underline\R\to\Omega^0\to\Omega^1\to\cdots$
is exact, because every point has a neighbourhood $U\cong\R^n$ on which the
sequence of sections is exact by the Poincar\'e lemma, and exactness of
sheaves is checked on stalks. Each $\Omega^p$ is acyclic because smooth
functions admit partitions of unity. So the complex of global sections
computes $H^\bullet(M,\underline\R)$. Full accounts are in
Warner~\citep{warner1983} and Bott and Tu~\citep{bott-tu1982}.
\end{proof}

We consider the analogous complex in the discrete case. Let $X$ be a
discrete coordinatized space and consider the complex of sheaves on
$\Fp(X)$
\[
  0\to\mathcal K\to\Nc^{0}\xrightarrow{d}\Nc^{1}\xrightarrow{d}\Nc^{2}\to\cdots,
\]
where $\mathcal K=\ker(\Nc^{0}\to\Nc^{1})$.

\begin{proposition}\label{prop:kernel-constant}
$\mathcal K$ is the constant sheaf $\underline\Z$ on $\Fp(X)$.
\end{proposition}

\begin{proof}
For $U$ open, $\Nc^{0}(U)$ is the functions on the $0$-cubes of $U$, and
$(df)\cube{x}{i}=f(e_i\cdot x)-f(x)$, so $\mathcal K(U)$ is the
functions on $U^{(0)}$ with $f(e_i\cdot x)=f(x)$ for every $1$-cube
$\cube{x}{i}\in U$. The constant sheaf is
$\underline\Z(U)=\{f\colon U\to\Z\text{ locally constant}\}$. A locally
constant function is constant on the smallest open neighbourhood $U_q$ of
each cube $q\in U$, so it is determined by its values on $0$-cubes. This is because if $q \in U$, 
$U_{q}$ consists of all the faces of $q$, including the $0-$cubes of $q$. In particular, 
if $\cube{x}{i} \in U$, then the $0-$cubes $x, e_{i} \cdot x \in U$, giving $f(x) = f(e_{i} \cdot x)$. 
Thus, $f \vert_{U^{(0)}} \in \mathcal K$ and hence $\underline\Z(U) \subseteq \mathcal K$. Conversely, let $f \in \mathcal K(U)$. Then, constant along the $1$-cubes of $U$ extends uniquely to a locally constant
function on $U$, since all vertices of a cube $q\in U$ are joined by
$1$-cubes in $U_q$. So $\mathcal K(U)=\underline\Z(U)$.
\end{proof}

Exactness of a sequence of sheaves is checked on stalks, and on an
Alexandrov space the stalk of a sheaf at $q$ is its group of sections over
the smallest open neighbourhood $U_q$. If $q=\cube{x}{I}$ is a $p$-cube of
the chart $X_\alpha$, then $U_q$ is the set of all faces of $q$,
\begin{equation}\label{eq:Uq}
  U_q=\bigl\{\cube{S_J\cdot x}{I\setminus K}:\ J\subseteq K\subseteq I\bigr\},
\end{equation}
and the sequence of stalks at $q$ is
\begin{equation}\label{eq:stalk-sequence}
  0\to\Z\to\Nc^{0}(U_q)\to\Nc^{1}(U_q)\to\Nc^{2}(U_q)\to\cdots,
\end{equation}
where $\underline\Z(U_q)=\Z$ because $U_q$ is connected. The sequence is
exact at $\Nc^{0}$ by \cref{prop:kernel-constant}, so it is exact if and
only if the normalized cochain complex of the cubical subset $U_q$ is
acyclic in positive degrees.

\subsection{The basic shapes of \texorpdfstring{$\N$}{N}-sets}
\label{sec:basic-shapes}

In general \eqref{eq:stalk-sequence} is not exact, because $U_q$ is the
image of the standard $p$-cube in $X$ and not necessarily a copy of it.

\begin{figure}[htb]
\centering
\resizebox{\textwidth}{!}{% Two cubes whose face-closure U_q is not contractible.
% (a) q=[0;{1,2}] in Cay(Z/3,{+1,+2}): the square has corners 0,1,2,0 and its
%     1-skeleton is a wedge of two circles.
% (b) q=[0;{1,2,3}] in Cay(Z/7,{+1,+2,+3}): the corner 1+2 and the corner 3
%     of the standard cube are the same vertex of Z/7.
\begin{tikzpicture}[x=1cm,y=1cm]
% ---------------- (a) -------------------------------------------------------
\begin{scope}
  % the abstract square
  \node[vtx,label={[font=\small]below left:$0$}] (a) at (0,0) {};
  \node[vtx,label={[font=\small]below right:$1$}] (b) at (2.2,0) {};
  \node[vtx,label={[font=\small]above left:$2$}] (c) at (0,2.2) {};
  \node[vtx,label={[font=\small]above right:$0$}] (d) at (2.2,2.2) {};
  \fill[fillcol!14] (0,0) rectangle (2.2,2.2);
  \draw[ar,cA] (a) -- (b); \draw[ar,cA] (c) -- (d);
  \draw[ar,cB] (a) -- (c); \draw[ar,cB] (b) -- (d);
  \node[font=\scriptsize,fillcol!70!black] at (1.1,1.1) {$\cube{0}{\{1,2\}}$};
  \draw[->,gray,shorten >=3pt,shorten <=3pt] (3.0,1.1) -- (4.2,1.1);
  % its image in Z/3
  \begin{scope}[xshift=5.9cm,yshift=1.1cm]
    \node[vtx,label={[font=\small]below:$0$}] (o) at (0,0) {};
    \node[vtx,label={[font=\small]right:$1$}] (p) at (1.6,0.9) {};
    \node[vtx,label={[font=\small]left:$2$}] (q) at (-1.6,0.9) {};
    \draw[ar,cA] (o) to[bend left=22] (p);   % 0 -> 1  (+1)
    \draw[ar,cB] (p) to[bend left=22] (o);   % 1 -> 0  (+2)
    \draw[ar,cB] (o) to[bend right=22] (q);  % 0 -> 2  (+2)
    \draw[ar,cA] (q) to[bend right=22] (o);  % 2 -> 0  (+1)
  \end{scope}
  \node[font=\small,align=center] at (3.6,-1.3) {(a) $q=\cube{0}{\{1,2\}}$ in $\Cay(\Z/3,\{+1,+2\})$: the two corners\\ marked $0$ coincide, and $U_q$ is a wedge of two circles};
\end{scope}
% ---------------- (b) -------------------------------------------------------
\begin{scope}[xshift=10.6cm,yshift=-0.3cm]
  \def\ox{1.1}\def\oy{0.75}
  \coordinate (v0) at (0,0);        \coordinate (v1) at (2.4,0);
  \coordinate (v2) at (0,2.4);      \coordinate (v12) at (2.4,2.4);
  \coordinate (v3) at (\ox,\oy);    \coordinate (v13) at (2.4+\ox,\oy);
  \coordinate (v23) at (\ox,2.4+\oy); \coordinate (v123) at (2.4+\ox,2.4+\oy);
  \fill[fillcol!10] (v0) -- (v1) -- (v12) -- (v2) -- cycle;
  \fill[fillcol!10] (v3) -- (v13) -- (v123) -- (v23) -- cycle;
  \draw[ar,cA] (v0) -- (v1); \draw[ar,cA] (v2) -- (v12); \draw[ar,cA] (v3) -- (v13); \draw[ar,cA] (v23) -- (v123);
  \draw[ar,cB] (v0) -- (v2); \draw[ar,cB] (v1) -- (v12); \draw[ar,cB] (v3) -- (v23); \draw[ar,cB] (v13) -- (v123);
  \draw[ar,cC] (v0) -- (v3); \draw[ar,cC] (v1) -- (v13); \draw[ar,cC] (v2) -- (v23); \draw[ar,cC] (v12) -- (v123);
  \node[vtx,label={[font=\small]below left:$0$}] at (v0) {};
  \node[vtx,label={[font=\small]below right:$1$}] at (v1) {};
  \node[vtx,label={[font=\small]left:$2$}] at (v2) {};
  \node[vtx,fill=idc,label={[font=\small,idc]above left:$3$}] at (v12) {};
  \node[vtx,fill=idc,label={[font=\small,idc]below right:$3$}] at (v3) {};
  \node[vtx,label={[font=\small]right:$4$}] at (v13) {};
  \node[vtx,label={[font=\small]above:$5$}] at (v23) {};
  \node[vtx,label={[font=\small]above right:$6$}] at (v123) {};
  \node[font=\small,align=center] at (1.75,-1.3) {(b) $q=\cube{0}{\{1,2,3\}}$ in $\Cay(\Z/7,\{+1,+2,+3\})$: the corners\\ $1+2$ and $3$ coincide, so $U_q$ has seven vertices, not eight};
\end{scope}
\end{tikzpicture}}
\caption{Two cubes $q$ whose smallest open neighbourhood $U_q$ is not acyclic. In each case
distinct corners of the standard cube map to the same vertex.}
\label{fig:wedge-cubes}
\end{figure}
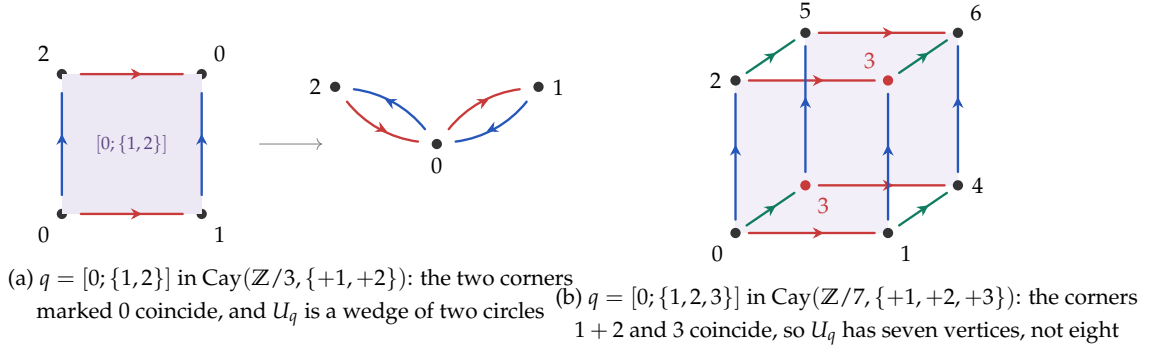

\begin{example}\label{ex:non-regular}
Consider $\Cay(\Z/3,\{+1,+2\})$ as the $\N$-set $(\Z/3,\N^2)$ with $e_1$
adding $1$ and $e_2$ adding $2$, and take the $2$-cube $q=\cube{0}{\{1,2\}}$.
Its four corners are $0$, $1$, $2$ and $1+2=0$, so two corners coincide,
and its $1$-skeleton is a wedge of two circles, \cref{fig:wedge-cubes}(a).
The normalized complex of $U_q$ has $3$ vertices, $4$ edges and $1$ square,
none of them degenerate, and a direct computation gives
$H^1(\Ns(U_q))=\Z$. Similarly, in $\Cay(\Z/7,\{+1,+2,+3\})=K_7$ the
$3$-cube $q=\cube{0}{\{1,2,3\}}$ has the corners $1+2$ and $3$ in common,
\cref{fig:wedge-cubes}(b), and $U_q$ has $7$ vertices, $12$ edges, $6$
squares and $1$ cube with $H^1(\Ns(U_q))=\Z$. Neither is acyclic. Both
groups are explained by \cref{prop:violation}, see \cref{ex:violation}.
%Claude note: the draft says these "can be shown to be non-contractible"
%without details. I verified both cohomology groups by computing the
%complexes directly (H^1 = Z in both cases, H^2 = 0 for Z/3, H^2 = H^3 = 0
%for Z/7); the wedge-of-circles picture gives H^1 of rank 2 for the
%1-skeleton, reduced to rank 1 by the single square.
\end{example}

In all the examples we have studied so far, the polar grid, the square and
the triangle among them, the smallest open neighbourhood $U_q$ of every
cube $q=\cube{x}{I}$ has been acyclic, even in the triangle where two
corners coincide. We give this
property a name, and in \cref{sec:regularity-criterion} we determine
exactly which $\N$-sets have it.

\begin{definition}[Regular $\N$-set]\label{def:regular}
Write $\HN{p}(U)=H^p(\Ns(U))$ for an open subset $U\subseteq\Fp(X)$. We
call $U$ \emph{contractible} if $\HN{0}(U)=\Z$ and $\HN{p}(U)=0$ for
$p>0$. An $\N$-set $(X,\N^n)$ is \emph{regular} if the smallest open neighbourhood $U_q$
of every cube $q=\cube{x}{I}$, with $x\in X$ and $I\subseteq[n]$, is
contractible. Since $U_q$ is connected, this is the same as asking that
the normalized complex of $U_q$ is acyclic in positive degrees,
$H^p(\Ns(U_q))=0$ for $p>0$. An atlas is regular if all its charts are.
\end{definition}
%Claude note: "contractible" is used here in the cohomological sense of
%your draft (H^0 = Z and no higher cohomology). It says nothing about the
%homotopy type of a geometric realization; if you prefer, "acyclic" would
%avoid that reading.

\begin{theorem}\label{thm:sheaf-coh}
Let $X$ be a discrete coordinatized space with a regular atlas
$(X_\alpha,\N^{n_\alpha})$. Then the normalized de Rham cohomology of $X$
is the sheaf cohomology of the constant sheaf on its face poset,
\[
  \HN{p}(X)\cong H^p\bigl(\Fp(X),\underline\Z\bigr)\qquad\text{for all }p.
\]
\end{theorem}

\begin{proof}
By \cref{prop:kernel-constant} and regularity, the sequence of stalks
\eqref{eq:stalk-sequence} is exact at every $q$, so
$0\to\underline\Z\to\Nc^{0}\to\Nc^{1}\to\cdots$ is an exact sequence of
sheaves on $\Fp(X)$, a resolution of $\underline\Z$. Each $\Nc^{p}$ is
flasque: for $V\subseteq U$ open, the restriction
$\Nc^{p}(U)\to\Nc^{p}(V)$ is restriction of functions on non-degenerate
$p$-cubes, which is surjective since a function on $V^{(p)}$ extends by
zero. Flasque sheaves are acyclic, so the resolution computes sheaf
cohomology, and the complex of global sections is $\Ns(\Fp(X))$, whose
cohomology is $\HN{\bullet}(X)$ by definition.
\end{proof}

\begin{intuition}
For regular atlases the normalized de Rham cohomology is an invariant of
the topological space $\Fp(X)$, which is the discrete form of de Rham's
theorem. In the proof, regularity plays the role of the Poincar\'e lemma
and flasqueness plays the role of partitions of unity. Flasqueness holds
because a function on a subset of a discrete set extends by zero.
\end{intuition}

\subsection{A criterion for regularity}
\label{sec:regularity-criterion}

%Claude note: this subsection is new. It is your Regularity.pdf (Condition
%(*), the up-set lemma and its corollary) typeset in your notation, followed
%by the converse (Proposition "violation") and the regularity criterion.
%Departures from the draft are flagged in further notes below.

\Cref{thm:sheaf-coh} applies only to regular atlases, so we now determine
which $\N$-sets are regular. Throughout this subsection $(X,\N^n)$ is an
$\N$-set and $q=\cube{x}{I}$ is a cube with $I=\{i_1,\dots,i_p\}$. By
\eqref{eq:Uq} the smallest open neighbourhood $U_q\subseteq\Fp(X)$ of $q$ is the poset with
elements $\cube{S_J\cdot x}{I\setminus K}$ for $J\subseteq K\subseteq I$.
Writing $L=I\setminus K$, its $r$-cubes are
\[
  U_q^{(r)}=\bigl\{\cube{S_J\cdot x}{L}:\ |L|=r,\ J\subseteq I,\ J\cap L=\varnothing\bigr\},
\]
and in particular its $0$-cubes are the vertices $S_J\cdot x$ for
$J\subseteq I$. This matches the observation that a vertex of the standard
cube $\Bx^p$ is determined by a subset of $I$. Let $2^I$ be the set of
subsets of $I$, the vertices of $\Bx^p$, and let
\[
  \phi\colon 2^I\to U_q^{(0)},\qquad J\mapsto S_J\cdot x.
\]
Since the coordinate fields can have fixed points and need not be
injective, $\phi$ need not be injective either. It can identify different
corners of the standard cube, and this is why the sets $U_q$ in
\cref{ex:non-regular} are not contractible. A cube of $U_q$ is determined
by its direction set and its base vertex, so $U_q$ depends only on which
corners $\phi$ identifies and on which corners are fixed by which
coordinate fields.

\begin{definition}[Condition $(\ast)$]\label{def:condition-star}
A cube $q=\cube{x}{I}$ satisfies \emph{condition $(\ast)$} if, whenever
$S_K\cdot x=S_{K'}\cdot x=z$ for subsets $K,K'\subseteq I$, which may be
empty, we have $e_i\cdot z=z$ for every $i\in K\triangle K'$.
\end{definition}

In other words, if two corners $K$ and $K'$ land on the same vertex,
$\phi(K)=\phi(K')$, then the cube is already collapsed at that vertex in
every direction in which the two corners differ. If $\phi$ is injective
the condition holds trivially. If $q$ satisfies $(\ast)$ then so does every
face of $q$. Indeed the corners of a face $\cube{S_J\cdot x}{L}$ are the
vertices $S_{J\cup A}\cdot x$ for $A\subseteq L$, which are corners of $q$,
and $(J\cup A)\triangle(J\cup A')=A\triangle A'$.

\begin{example}\label{ex:condition-star}
In the $2$-simplex of \cref{sec:two-simplex} the cube $\cube{x}{\{1,2\}}$
has corners $\phi(\varnothing)=x$, $\phi(\{1\})=a$ and
$\phi(\{2\})=\phi(\{1,2\})=b$. The only two corners that coincide differ in
the direction $1$, and $e_1\cdot b=b$, so $(\ast)$ holds. In
\cref{ex:non-regular} the cube $\cube{0}{\{1,2\}}$ has
$\phi(\varnothing)=\phi(\{1,2\})=0$, but $e_1\cdot0=1\neq0$, so $(\ast)$
fails. Here one of the two corners is the empty set, which is why empty
sets are allowed in \cref{def:condition-star}.
\end{example}

\begin{intuition}
Condition $(\ast)$ says that corners of a cube may only meet by collapsing.
In the $2$-simplex the corners $\{2\}$ and $\{1,2\}$ meet because the edge
$\cube{b}{\{1\}}$ joining them is a loop, which the normalization removes,
and the square becomes a triangle, which is still a disc. In
$\Cay(\Z/3,\{+1,+2\})$ the opposite corners $\varnothing$ and $\{1,2\}$ are
identified without any edge being collapsed, and identifying two opposite
corners of a square produces a circle.
\end{intuition}

We show that $(\ast)$ implies contractibility by induction on the
dimension $|I|$ of the cube. The induction needs pieces of $U_q$ that are
slightly more general than the sets $U_q$. An \emph{up-set} of $2^I$ is a
subset $\cU\subseteq2^I$ such that $J\in\cU$ and $J\subseteq J'\subseteq I$
imply $J'\in\cU$. For an up-set $\cU$ we define the subposet
\[
  P_{[x;I]}(\cU)=\bigl\{\cube{S_J\cdot x}{L}:\ J\in\cU,\ L\subseteq I\setminus J\bigr\}
  \subseteq U_q .
\]
It is an open subset of $\Fp(X)$, that is, it is closed under faces,
because $\cU$ is an up-set. The two faces of $\cube{S_J\cdot x}{L}$ in a
direction $l\in L$ are $\cube{S_J\cdot x}{L\setminus\{l\}}$ and
$\cube{S_{J\cup\{l\}}\cdot x}{L\setminus\{l\}}$, and $J\cup\{l\}\in\cU$.
For the up-set of all subsets we recover $P_{[x;I]}(2^I)=U_q$.

\begin{construction}\label{con:vertical}
Let $\cU\subseteq2^I$ be a non-empty up-set and fix a direction $i\in I$,
for instance $i=i_1$. Put $I'=I\setminus\{i\}$ and
\[
  \cU_0=\{J\in\cU:\ i\notin J\},\qquad
  \cU_1=\{J\setminus\{i\}:\ i\in J\in\cU\},
\]
which are up-sets of $2^{I'}$. If $J\in\cU_0$ then $J\cup\{i\}\in\cU$, so
$\cU_0\subseteq\cU_1$, and $\cU_1\neq\varnothing$ because $I\in\cU$. We
write $P=P_{[x;I]}(\cU)$ and
\begin{align*}
  B&=P_{[x;I']}(\cU_0) &&\text{(bottom)},\\
  T&=P_{[e_i\cdot x;I']}(\cU_1) &&\text{(top)},\\
  V&=\bigl\{\cube{S_J\cdot x}{L}\in P:\ i\in L\bigr\} &&\text{(vertical)},\\
  F&=\bigl\{\cube{y}{L}\in B:\ e_i\cdot y=y\bigr\},
\end{align*}
and we let $\cU_{\mathrm{fix}}=\{J\subseteq I':\ e_i\cdot S_J\cdot x=S_J\cdot x\}$.
\end{construction}

%Claude note: the draft fixes the direction i_1. The proof of
%Proposition "violation" below needs to choose the direction, so the
%construction allows any i in I; for i = i_1 nothing changes except the
%sign in the map psi of Lemma "vertical".
%
%Claude note: the draft describes B and T as P(U) intersected with
%U_[x;I\{i_1}] and U_[e_{i_1}x;I\{i_1}]. These intersections can be
%larger than B and T. Take X={x,a} with e_1x=e_2x=a and a fixed by both,
%q=[x;{1,2}], i=1 and U={{1},{1,2}}. Then U_0 is empty, so B is empty, but
%P(U) and U_[x;{2}] share the vertex a. The proof only
%uses B and T as defined here, so I have defined them directly.

\begin{lemma}\label{lem:vertical}
The following hold in the situation of \cref{con:vertical}, whether or not
$q$ satisfies $(\ast)$.
\begin{enumerate}
\item $P=B\cup T\cup V$, the set $V$ is disjoint from $B\cup T$, and $B$,
  $T$ and $B\cup T$ are open.
\item $\cU_{\mathrm{fix}}$ is an up-set of $2^{I'}$,
  $F=P_{[x;I']}(\cU_0\cap\cU_{\mathrm{fix}})$, and $F\subseteq B\cap T$. In
  particular $F$ is open.
\item Let $\Ns(V)$ denote the kernel of the restriction
  $\Ns(P)\to\Ns(B\cup T)$. Then
  \[
    0\to\Ns(V)\to\Ns(P)\to\Ns(B\cup T)\to0
  \]
  is a short exact sequence of cochain complexes, and $\Nc^{r}(V)$ is the
  group of functions on the non-degenerate $r$-cubes of $V$.
\item Let $\Ns(B,F)$ be the kernel of the restriction $\Ns(B)\to\Ns(F)$.
  Then the maps
  \[
    \psi\colon\Nc^{r-1}(B,F)\to\Nc^{r}(V),\qquad
    (\psi\omega)\cube{y}{M\cup\{i\}}=(-1)^{m(M)}\,\omega\cube{y}{M},
  \]
  where $m(M)$ is the number of elements of $M$ greater than $i$, form an
  isomorphism of cochain complexes $\Ns(B,F)[-1]\cong\Ns(V)$. Hence
  $\HN{r}(V)\cong\HN{r-1}(B,F)$ for all $r$.
\end{enumerate}
\end{lemma}

Here and below $\HN{r}(U)=H^r(\Ns(U))$ for an open set $U$, and
$\HN{r}(V)$ and $\HN{r}(B,F)$ denote the cohomology of the complexes
$\Ns(V)$ and $\Ns(B,F)$. The set $V$ is not open. It is the complement of
the open set $B\cup T$ in $P$, and $\Ns(V)$ is the complex of cochains on
$P$ supported on $V$.

\begin{proof}
(i) A cube of $P$ is $\cube{S_J\cdot x}{L}$ with $J\in\cU$ and
$J\cap L=\varnothing$. If $i\in L$ it lies in $V$. If $i\notin L$ and
$i\notin J$, then $J\in\cU_0$ and $L\subseteq I'\setminus J$, so the cube
lies in $B$. If $i\notin L$ and $i\in J$, then
$S_J\cdot x=S_{J\setminus\{i\}}\cdot(e_i\cdot x)$ with
$J\setminus\{i\}\in\cU_1$ and $L\subseteq I'\setminus(J\setminus\{i\})$, so
the cube lies in $T$. Conversely $B\subseteq P$ because
$\cU_0\subseteq\cU$, and a cube $\cube{S_J\cdot(e_i\cdot x)}{L}$ of $T$
equals $\cube{S_{J\cup\{i\}}\cdot x}{L}$ with $J\cup\{i\}\in\cU$ and
$L\subseteq I\setminus(J\cup\{i\})$, so $T\subseteq P$. The cubes of $V$
have $i$ in their direction set and those of $B\cup T$ do not, so $V$ is
disjoint from $B\cup T$. Finally $B$ and $T$ are open, being up-set pieces
of the cubes $\cube{x}{I'}$ and $\cube{e_i\cdot x}{I'}$, and a union of
open sets is open.

(ii) If $J\subseteq J'\subseteq I'$ and $e_i\cdot S_J\cdot x=S_J\cdot x$,
then $e_i\cdot S_{J'}\cdot x=S_{J'\setminus J}\cdot e_i\cdot S_J\cdot x
=S_{J'}\cdot x$ because the coordinate fields commute, so
$\cU_{\mathrm{fix}}$ is an up-set. A cube $\cube{y}{L}\in B$ has
$y=S_J\cdot x$ for some $J\in\cU_0$ with $L\subseteq I'\setminus J$, and
$e_i\cdot y=y$ holds exactly when $J\in\cU_{\mathrm{fix}}$, a condition
that depends only on the vertex $y$. Hence
$F=P_{[x;I']}(\cU_0\cap\cU_{\mathrm{fix}})$, which is open because
$\cU_0\cap\cU_{\mathrm{fix}}$ is an up-set. If $\cube{y}{L}\in F$ with
$y=S_J\cdot x$ as above, then $J\cup\{i\}\in\cU$, so $J\in\cU_1$, and
$S_J\cdot(e_i\cdot x)=e_i\cdot y=y$. Hence
$\cube{y}{L}=\cube{S_J\cdot(e_i\cdot x)}{L}\in T$, and $F\subseteq B\cap T$.

(iii) Since $B\cup T$ is open, the value of $d\omega$ on a cube of
$B\cup T$ only involves the values of $\omega$ on faces of that cube, and
these lie in $B\cup T$. So restriction commutes with $d$. It is surjective
because a function on the non-degenerate cubes of $B\cup T$ extends by zero
to $P$. Its kernel is a subcomplex, and by (i) it consists of the functions
on the non-degenerate cubes of $P$ that vanish outside $V$.

(iv) A cube $\cube{y}{L}$ lies in $V$ exactly when $L=M\cup\{i\}$ with
$i\notin M$ and $y=S_J\cdot x$ for some $J\in\cU_0$ with
$J\cap M=\varnothing$, that is, exactly when $\cube{y}{M}\in B$. So
$\cube{y}{M\cup\{i\}}\mapsto\cube{y}{M}$ is a bijection from $V$ to $B$.
The cube $\cube{y}{M\cup\{i\}}$ is non-degenerate if and only if
$e_l\cdot y\neq y$ for all $l\in M$ and $e_i\cdot y\neq y$, that is, if and
only if $\cube{y}{M}$ is a non-degenerate cube of $B$ that does not lie in
$F$. Since $\Nc^{r-1}(B,F)$ is the group of functions on the non-degenerate
$(r-1)$-cubes of $B$ that vanish on $F$, each $\psi$ is an isomorphism of
groups.

It remains to check that $\psi$ commutes with the differentials. Let
$\omega\in\Nc^{r-1}(B,F)$ and let $c=\cube{y}{M\cup\{i\}}$ be a
non-degenerate cube of $V$ with $|M|=r$. We compute $(d\psi\omega)(c)$ with
\cref{prop:d-on-cubes}. The two faces of $c$ in the direction $i$ lie in
$B\cup T$, where $\psi\omega$ vanishes. Let $l$ be the $k$-th element of
$M$, counting from $0$. Its position in $M\cup\{i\}$ is $k+[l>i]$, where
$[l>i]$ is $1$ if $l>i$ and $0$ otherwise, and the two faces of $c$ in the
direction $l$ are $\cube{e_l\cdot y}{(M\setminus\{l\})\cup\{i\}}$ and
$\cube{y}{(M\setminus\{l\})\cup\{i\}}$. Since
$m(M\setminus\{l\})=m(M)-[l>i]$, the term of $l$ in $(d\psi\omega)(c)$ is
\[
  (-1)^{k+[l>i]}(-1)^{m(M)-[l>i]}
  \Bigl(\omega\cube{e_l\cdot y}{M\setminus\{l\}}-\omega\cube{y}{M\setminus\{l\}}\Bigr),
\]
which is $(-1)^{m(M)}$ times the term of $l$ in $(d\omega)\cube{y}{M}$.
Faces that are degenerate or that correspond to cubes of $F$ contribute $0$
on both sides. Summing over $l\in M$ gives
$(d\psi\omega)(c)=(-1)^{m(M)}(d\omega)\cube{y}{M}=(\psi\,d\omega)(c)$.
\end{proof}

%Claude note: the draft defines psi without the sign (-1)^{m(M)}. For
%i = i_1 every element of M is greater than i, so m(M)=|M|, and without the
%sign one gets d(psi omega) = -psi(d omega). The unsigned map is still an
%isomorphism of complexes after composing with (-1)^{degree}, so the
%cohomology statement in the draft is correct; the sign just makes psi
%itself a chain map.

\begin{proposition}\label{prop:upset-contractible}
If $q=\cube{x}{I}$ satisfies $(\ast)$, then for every non-empty up-set
$\cU\subseteq2^I$ the open set $P_{[x;I]}(\cU)$ is contractible, that is,
$\HN{0}(P_{[x;I]}(\cU))=\Z$ and $\HN{r}(P_{[x;I]}(\cU))=0$ for $r\geq1$.
\end{proposition}

%Claude note: the draft says "for any up-set". The empty up-set gives the
%empty set, whose H^0 is 0, so the up-set has to be non-empty.

\begin{proof}
We induct on $p=|I|$. If $p=0$, then $\cU=\{\varnothing\}$ and
$P_{[x;I]}(\cU)$ is the single vertex $x$, so the claim holds. Let
$p\geq1$, and assume the claim for every cube of dimension at most $p-1$
that satisfies $(\ast)$ and every non-empty up-set of it. Let $\cU$ be a
non-empty up-set, fix $i\in I$ and use the notation of
\cref{con:vertical}. The faces $\cube{x}{I'}$ and $\cube{e_i\cdot x}{I'}$
satisfy $(\ast)$ and have dimension $p-1$, and $\cU_0$, $\cU_1$ and
$\cU_0\cap\cU_{\mathrm{fix}}$ are up-sets. By the induction hypothesis and
\cref{lem:vertical}(ii), each of $B$, $T$ and $F$ is contractible when it
is non-empty.

We first show that $F=B\cap T$. This is the only place where $(\ast)$ is
used. By \cref{lem:vertical}(ii), $F\subseteq B\cap T$. Conversely let
$\cube{y}{L}\in B\cap T$. Then $y=S_J\cdot x$ with $J\in\cU_0$, and
$y=S_K\cdot(e_i\cdot x)=S_{K\cup\{i\}}\cdot x$ with $K\in\cU_1$. Since
$i\notin J$ we have $i\in J\triangle(K\cup\{i\})$, so $(\ast)$ gives
$e_i\cdot y=y$, and $\cube{y}{L}\in F$.

If $\cU_0=\varnothing$, then $B=\varnothing$, so $V=\varnothing$ by
\cref{lem:vertical}(iv), and $P_{[x;I]}(\cU)=T$ is contractible by the
induction hypothesis. So assume $\cU_0\neq\varnothing$. Then $B$ and $T$
are non-empty and contractible. Since $F=B\cap T$ there is a
Mayer--Vietoris short exact sequence
\[
  0\to\Ns(B\cup T)\to\Ns(B)\oplus\Ns(T)\to\Ns(F)\to0,
\]
with maps $\omega\mapsto(\omega|_B,\omega|_T)$ and
$(\beta,\tau)\mapsto\beta|_F-\tau|_F$. It is exact because all three are
complexes of functions on non-degenerate cubes and $B\cap T=F$.

\emph{Case 1: $F=\varnothing$.} Then $B$ and $T$ are disjoint, so
$\HN{0}(B\cup T)=\HN{0}(B)\oplus\HN{0}(T)=\Z^2$ and $\HN{r}(B\cup T)=0$ for
$r\geq1$. By \cref{lem:vertical}(iv),
$\HN{r}(V)\cong\HN{r-1}(B,\varnothing)=\HN{r-1}(B)$, which is $\Z$ for
$r=1$ and $0$ otherwise. The short exact sequence of
\cref{lem:vertical}(iii) gives the long exact sequence
\[
  0\to\HN{0}(P_{[x;I]}(\cU))\to\Z^2\xrightarrow{\ \delta\ }\Z
  \to\HN{1}(P_{[x;I]}(\cU))\to0,
\]
and $\HN{r}(P_{[x;I]}(\cU))=0$ for $r\geq2$, since it lies between
$\HN{r}(V)=0$ and $\HN{r}(B\cup T)=0$. To compute $\delta$, let $f$ be the
function equal to $a$ on the vertices of $B$ and to $b$ on the vertices of
$T$. Every vertex of $P$ lies in $B\cup T$, so $f$ is a $0$-cochain on $P$,
and $df$ vanishes on $B\cup T$. On a non-degenerate vertical edge
$\cube{y}{\{i\}}$ we have $y\in B$ and $e_i\cdot y\in T$, so
$(df)\cube{y}{\{i\}}=f(e_i\cdot y)-f(y)=b-a$. Under $\psi$ this is the
constant function $b-a$ on $B$, so $\delta(a,b)=b-a$. This map is
surjective with kernel $\Z$, hence $\HN{0}(P_{[x;I]}(\cU))=\Z$ and
$\HN{1}(P_{[x;I]}(\cU))=0$.

\emph{Case 2: $F\neq\varnothing$.} Then $F$ is contractible, and the long
exact sequence of the Mayer--Vietoris sequence begins
\[
  0\to\HN{0}(B\cup T)\to\Z^2\to\Z\to\HN{1}(B\cup T)\to0,
\]
with $\HN{r}(B\cup T)=0$ for $r\geq2$. The map $\Z^2\to\Z$ is
$(a,b)\mapsto a-b$, which is surjective, so $\HN{0}(B\cup T)=\Z$ and
$\HN{1}(B\cup T)=0$. The restriction $\HN{0}(B)\to\HN{0}(F)$ is the
identity of $\Z$, since both groups consist of the constant functions, so
the long exact sequence of $0\to\Ns(B,F)\to\Ns(B)\to\Ns(F)\to0$ gives
$\HN{\bullet}(B,F)=0$. By \cref{lem:vertical}(iv) $\HN{\bullet}(V)=0$, and
the long exact sequence of \cref{lem:vertical}(iii) gives
$\HN{\bullet}(P_{[x;I]}(\cU))\cong\HN{\bullet}(B\cup T)$, which is $\Z$ in
degree $0$ and $0$ in higher degrees.
\end{proof}

%Claude note: the draft also treats the case U_1 = empty. It cannot occur,
%since a non-empty up-set contains I, and I contains i. In the case
%U_0 = empty the draft writes P(U) = U_{[e_{i_1}x; I\{i_1}]}; in general it is
%the up-set piece P_{[e_i x; I']}(U_1), which is contractible by the same
%induction hypothesis.

\begin{corollary}\label{cor:star-contractible}
If $q=\cube{x}{I}$ satisfies $(\ast)$, then $U_q$ is contractible, that is,
$\HN{0}(U_q)=\Z$ and $\HN{r}(U_q)=0$ for $r\geq1$.
\end{corollary}

\begin{proof}
$U_q=P_{[x;I]}(2^I)$, the up-set being all of $2^I$.
\end{proof}

\begin{intuition}[bottom, top and vertical cubes]
Choosing a direction $i$ cuts the cube into a bottom face, a top face and
the cubes that contain the direction $i$. The vertical cubes form a copy of
the bottom $B$ times an interval, whose lower end is the bottom and whose
upper end is glued to the top by $e_i$. Where $e_i$ fixes a vertex the
vertical cubes above it are degenerate and the cylinder is pinched, and this
happens exactly over $F$. Condition $(\ast)$ says that the bottom and the
top only meet where the cylinder is pinched. The pieces are then glued along
contractible pieces, and the result is contractible.
\end{intuition}

The converse of \cref{cor:star-contractible} fails for a single cube, but
it holds once all faces are taken into account.

\begin{example}\label{ex:cone}
Let $X=\Z/3\sqcup\{\ast\}$ with three coordinate fields. Let $e_1$ add $1$
and $e_2$ add $2$ on $\Z/3$, both fixing $\ast$, and let $e_3$ send every
point to $\ast$. The fields commute, since $e_3e_j$ and $e_je_3$ both send
every point to $\ast$. The cube $q=\cube{0}{\{1,2,3\}}$ fails $(\ast)$,
because its face $\cube{0}{\{1,2\}}$ is the cube of \cref{ex:non-regular}.
Nevertheless a direct computation gives $\HN{\bullet}(U_q)=(\Z,0,0,0)$.
Geometrically, every corner containing the direction $3$ is sent to
$\ast$, and $U_q$ is a cone on the circle $U_{\cube{0}{\{1,2\}}}$ with apex
$\ast$.
\end{example}

\begin{proposition}\label{prop:violation}
Suppose that $q=\cube{x}{I}$ does not satisfy $(\ast)$. Then some face $q'$
of $q$ has $\HN{1}(U_{q'})\neq0$, so $U_{q'}$ is not contractible.

More precisely, call a pair $(K_1,K_2)$ of subsets of $I$ a
\emph{violating pair} of $q$ if $S_{K_1}\cdot x=S_{K_2}\cdot x=z$ and
$e_i\cdot z\neq z$ for some $i\in K_1\triangle K_2$. Choose a violating
pair with $|K_1\triangle K_2|$ as small as possible, and put
\[
  y=S_{K_1\cap K_2}\cdot x,\qquad J=K_1\triangle K_2,\qquad
  q'=\cube{y}{J},\qquad z=S_{K_1}\cdot x .
\]
Choose $i\in J$ with $e_i\cdot z\neq z$ and let
\[
  W=\bigl\{S_K\cdot y:\ K\subseteq J\setminus\{i\},\
  S_K\cdot y=S_{J\setminus K}\cdot y,\ e_i\cdot S_K\cdot y\neq S_K\cdot y\bigr\}.
\]
Then $W\neq\varnothing$, $\HN{0}(U_{q'})=\Z$, $\HN{1}(U_{q'})\cong\Z^{|W|}$
and $\HN{r}(U_{q'})=0$ for $r\geq2$.
\end{proposition}

\begin{proof}
Put $C=K_1\cap K_2$, $A_1=K_1\setminus C$ and $A_2=K_2\setminus C$. Then
$A_1$ and $A_2$ are disjoint, $J=A_1\cup A_2\subseteq I\setminus C$, so
$q'=\cube{S_C\cdot x}{J}$ is a face of $q$, and
$S_{A_1}\cdot y=S_{A_2}\cdot y=z$. Exchanging $K_1$ and $K_2$ if necessary,
we may assume that $i\in A_1$.

\emph{Step 1. Every violating pair of $q'$ has symmetric difference $J$, and
every proper face of $q'$ satisfies $(\ast)$.} A violating pair $(A,A')$ of
$q'$, with $A,A'\subseteq J$, gives the violating pair $(C\cup A,C\cup A')$
of $q$, which has the same symmetric difference $A\triangle A'\subseteq J$.
By minimality $|A\triangle A'|\geq|J|$, so $A\triangle A'=J$. A proper face
of $q'$ is a cube $\cube{S_D\cdot y}{L}$ with $L\subsetneq J$, and a
violating pair of it would give a violating pair $(D\cup A,D\cup A')$ of
$q'$ whose symmetric difference is contained in $L$, which is impossible.

\emph{Step 2. The decomposition.} We apply \cref{con:vertical} to $q'$,
with the up-set $\cU=2^J$ and the direction $i$, and put
$J'=J\setminus\{i\}$. Then $\cU_0=\cU_1=2^{J'}$, so
$B=U_{\cube{y}{J'}}$ and $T=U_{\cube{e_i\cdot y}{J'}}$. Both are of the
form $U_{q''}$ for a proper face $q''$ of $q'$, so they are contractible by Step 1 and
\cref{cor:star-contractible}. Likewise $F=P_{[y;J']}(\cU_{\mathrm{fix}})$
is contractible if it is non-empty, by \cref{prop:upset-contractible}
applied to the face $\cube{y}{J'}$. \Cref{lem:vertical} does not use
$(\ast)$, so it applies.

\emph{Step 3. $B\cap T=F\sqcup W$, where the elements of $W$ are regarded as
$0$-cubes.} By \cref{lem:vertical}(ii), $F\subseteq B\cap T$. If
$w=S_K\cdot y\in W$, then $w$ is a vertex of $B$ because $K\subseteq J'$,
and it is a vertex of $T$ because
$w=S_{J\setminus K}\cdot y=S_{J'\setminus K}\cdot(e_i\cdot y)$, as
$i\in J\setminus K$. It does not lie in $F$ because $e_i\cdot w\neq w$.
Conversely let $\cube{w}{L}\in(B\cap T)\setminus F$. Then
$w=S_K\cdot y=S_{K'}\cdot(e_i\cdot y)=S_{K'\cup\{i\}}\cdot y$ with
$K,K'\subseteq J'\setminus L$, and $e_i\cdot w\neq w$. So
$(K,K'\cup\{i\})$ is a violating pair of $q'$, and by Step 1 its symmetric
difference is $J$. Both sets are disjoint from $L$, which forces
$L=\varnothing$ and $K'\cup\{i\}=J\setminus K$, so $w\in W$. Finally
$z\in W$, taking $K=A_2$. Indeed $A_2\subseteq J'$ because $i\in A_1$, and
$S_{A_2}\cdot y=z=S_{A_1}\cdot y=S_{J\setminus A_2}\cdot y$.

\emph{Step 4. Cohomology.} The Mayer--Vietoris sequence
$0\to\Ns(B\cup T)\to\Ns(B)\oplus\Ns(T)\to\Ns(B\cap T)\to0$ is exact as in
the proof of \cref{prop:upset-contractible}. Each vertex is an open set,
so $B\cap T$ is the disjoint union of the open sets $F$ and $W$, and
$\Ns(B\cap T)=\Ns(F)\oplus\Z^W$ with $\Z^W$ in degree $0$.

If $F\neq\varnothing$, then $\HN{\bullet}(B,F)=0$ as in Case 2 of the proof
of \cref{prop:upset-contractible}, so $\HN{\bullet}(V)=0$ and
$\HN{\bullet}(U_{q'})\cong\HN{\bullet}(B\cup T)$. The Mayer--Vietoris
sequence gives
\[
  0\to\HN{0}(B\cup T)\to\Z^2\xrightarrow{\ \mu\ }\Z^{1+|W|}\to\HN{1}(B\cup T)\to0
\]
and $\HN{r}(B\cup T)=0$ for $r\geq2$, where $\mu(a,b)=(a-b)(1,\dots,1)$.
The kernel of $\mu$ is $\Z$, and its cokernel is
$\Z^{1+|W|}/\Z(1,\dots,1)\cong\Z^{|W|}$ because $(1,\dots,1)$ extends to a
basis of $\Z^{1+|W|}$.

If $F=\varnothing$, then $B\cap T=W$, and the same argument gives
$\HN{0}(B\cup T)=\Z$, $\HN{1}(B\cup T)\cong\Z^{|W|-1}$ and
$\HN{r}(B\cup T)=0$ for $r\geq2$. As in Case 1 of the proof of
\cref{prop:upset-contractible}, $\HN{r}(V)\cong\HN{r-1}(B)$ is $\Z$ for
$r=1$ and $0$ otherwise. This time the connecting map
$\delta\colon\HN{0}(B\cup T)\to\HN{1}(V)$ is zero. The group
$\HN{0}(B\cup T)=\Z$ is generated by the constant function $1$, and $d1$
vanishes on every vertical edge, since $1-1=0$. The long exact sequence of
\cref{lem:vertical}(iii) therefore gives $\HN{0}(U_{q'})=\Z$,
$\HN{r}(U_{q'})=0$ for $r\geq2$, and a short exact sequence
\[
  0\to\Z\to\HN{1}(U_{q'})\to\Z^{|W|-1}\to0,
\]
which splits because $\Z^{|W|-1}$ is free. So $\HN{1}(U_{q'})\cong\Z^{|W|}$
in both cases.
\end{proof}

\begin{example}\label{ex:violation}
In $\Cay(\Z/3,\{+1,+2\})$ the cube $q=\cube{0}{\{1,2\}}$ has the violating
pair $(\varnothing,\{1,2\})$, so $q'=q$, $y=z=0$ and $J=\{1,2\}$. With
$i=2$ we find $W=\{0\}$, because $S_\varnothing\cdot0=0=S_{\{1,2\}}\cdot0$
while $S_{\{1\}}\cdot0=1\neq2=S_{\{2\}}\cdot0$. In $K_7$ the cube
$q=\cube{0}{\{1,2,3\}}$ has the violating pair $(\{1,2\},\{3\})$, which is
the only pair of corners that coincide, so again $q'=q$. With $i=3$ we find
$W=\{3\}$. In both cases \cref{prop:violation} gives $\HN{1}(U_q)=\Z$, in
agreement with \cref{ex:non-regular}.
\end{example}

We can now characterize regular $\N$-sets.

\begin{theorem}[Regularity criterion]\label{thm:regularity-criterion}
For an $\N$-set $(X,\N^n)$ the following are equivalent.
\begin{enumerate}
\item $(X,\N^n)$ is regular.
\item Every cube $\cube{x}{I}$ of $X$ satisfies $(\ast)$.
\item For every $y\in X$ and all disjoint subsets $A,A'\subseteq[n]$, which
  may be empty, with $S_A\cdot y=S_{A'}\cdot y$, the vertex $S_A\cdot y$ is
  fixed by $e_i$ for every $i\in A\cup A'$.
\end{enumerate}
\end{theorem}

\begin{proof}
(ii)$\Rightarrow$(i) is \cref{cor:star-contractible}. For
(i)$\Rightarrow$(ii), suppose some cube $q$ does not satisfy $(\ast)$. By
\cref{prop:violation} some face $q'$ of $q$ has $\HN{1}(U_{q'})\neq0$, and
$q'$ is itself a cube of $X$, so $X$ is not regular.

For (ii)$\Rightarrow$(iii), apply $(\ast)$ to the cube
$\cube{y}{A\cup A'}$ with $K=A$ and $K'=A'$, noting that
$A\triangle A'=A\cup A'$ because $A$ and $A'$ are disjoint. For
(iii)$\Rightarrow$(ii), let $S_K\cdot x=S_{K'}\cdot x=z$ with
$K,K'\subseteq I$. Put $y=S_{K\cap K'}\cdot x$, $A=K\setminus K'$ and
$A'=K'\setminus K$. These are disjoint with $A\cup A'=K\triangle K'$, and
since the coordinate fields commute,
$S_A\cdot y=S_K\cdot x=z=S_{K'}\cdot x=S_{A'}\cdot y$. So (iii) gives
$e_i\cdot z=z$ for every $i\in K\triangle K'$.
\end{proof}

Condition (iii) is the easiest to check in practice, since it involves one
vertex at a time and finitely many pairs $(A,A')$. For $\N$-sets without
fixed points it says that distinct corners of every cube are distinct
vertices.

\begin{corollary}\label{cor:regular-translation}
Let $G$ be an abelian group, let $g_1,\dots,g_n\in G$ be non-zero, and let
$(G,\N^n)$ be the $\N$-set in which $e_i$ adds $g_i$. Then $(G,\N^n)$ is
regular if and only if
\[
  \sum_{i=1}^n\varepsilon_ig_i\neq0\qquad
  \text{for every non-zero }\varepsilon\in\{-1,0,1\}^n .
\]
In particular a homogeneous coordinate system on a Cayley graph
$\Cay(A,\{a_1,\dots,a_n\})$, as in \cref{thm:cayley}, is regular if and
only if its generators satisfy this condition.
\end{corollary}

\begin{proof}
No $e_i$ has a fixed point, since $g_i\neq0$. So condition (iii) of
\cref{thm:regularity-criterion} says that $S_A\cdot y\neq S_{A'}\cdot y$
for all disjoint $A,A'\subseteq[n]$ that are not both empty. Since
$S_A\cdot y=y+\sum_{i\in A}g_i$, this means
$\sum_{i\in A}g_i\neq\sum_{i\in A'}g_i$, that is,
$\sum_i\varepsilon_ig_i\neq0$ where $\varepsilon_i$ is $1$ on $A$, $-1$ on
$A'$ and $0$ elsewhere. Every non-zero $\varepsilon\in\{-1,0,1\}^n$ arises
in this way, with $A=\{i:\varepsilon_i=1\}$ and
$A'=\{i:\varepsilon_i=-1\}$.
\end{proof}

In additive combinatorics a family with this property is called
\emph{dissociated}. The two Cayley graphs of \cref{ex:non-regular} are not
regular, because $1+2=0$ in $\Z/3$ and $1+2-3=0$ in $\Z/7$. Neither is the
single chart of rank $3$ on the triangulated torus of
\cref{sec:triangulated-torus}, because $(1,0)+(0,1)-(1,1)=0$. This is a
further reason to replace it by the atlas of $2$-simplex charts used there.
On the other hand the cycle $\Cay(\Z/m,\{1\})$ and the square torus
$\Cay(\Z/n\times\Z/m,\{(1,0),(0,1)\})$ with $n,m\geq2$ are regular.
%Claude note: the term "dissociated" is standard in additive combinatorics
%(for example Tao and Vu, Additive Combinatorics, CUP 2006). I have not added
%a citation because refs.bib has no entry for it; add one if you want the
%reference.

\subsection{Mayer--Vietoris and \v{C}ech--de Rham}
\label{sec:mayer-vietoris}

Let $X$ be a discrete coordinatized space, let $\cU=(U_\alpha)_{\alpha\in A}$
be a family of open subsets of $\Fp(X)$ and $U=\bigcup_\alpha U_\alpha$.
Write $U_{\alpha_0\cdots\alpha_q}=U_{\alpha_0}\cap\dots\cap U_{\alpha_q}$.
The \emph{generalized Mayer--Vietoris sequence} of $\Nc^{p}$ is
\begin{equation}\label{eq:gmv}\small
  0\to\Nc^{p}(U)\xrightarrow{r}\prod_\alpha\Nc^{p}(U_\alpha)
  \xrightarrow{\delta}\prod_{\alpha_0,\alpha_1}\Nc^{p}(U_{\alpha_0\alpha_1})
  \xrightarrow{\delta}\cdots\xrightarrow{\delta}
  \prod_{\alpha_0,\dots,\alpha_q}\!\Nc^{p}(U_{\alpha_0\cdots\alpha_q})
  \xrightarrow{\delta}\cdots,
\end{equation}
where the products are over all tuples of indices, $r$ is restriction and
$(\delta\varphi)_{\alpha_0\cdots\alpha_{q+1}}
=\sum_{k=0}^{q+1}(-1)^k\varphi_{\alpha_0\cdots\widehat{\alpha_k}\cdots\alpha_{q+1}}
\big|_{U_{\alpha_0\cdots\alpha_{q+1}}}$.

\begin{proposition}\label{prop:gmv}
The sequence \eqref{eq:gmv} is exact.
\end{proposition}

\begin{proof}
Exactness at $\Nc^{p}(U)$ and at $\prod_\alpha\Nc^{p}(U_\alpha)$ is the sheaf
axiom. For the rest we use a discrete partition of unity. For every
non-degenerate $p$-cube $c\in U$ choose an index $\alpha(c)$ with
$c\in U_{\alpha(c)}$, and define
\[
  K\colon\prod\Nc^{p}(U_{\alpha_0\cdots\alpha_q})\to
  \prod\Nc^{p}(U_{\alpha_0\cdots\alpha_{q-1}}),\qquad
  (K\varphi)_{\alpha_0\cdots\alpha_{q-1}}(c)=\varphi_{\alpha(c)\alpha_0\cdots\alpha_{q-1}}(c),
\]
which makes sense because $c\in U_{\alpha_0\cdots\alpha_{q-1}}$ and
$c\in U_{\alpha(c)}$ imply $c\in U_{\alpha(c)\alpha_0\cdots\alpha_{q-1}}$. A
direct computation gives
\begin{align*}
  (\delta K\varphi)_{\alpha_0\cdots\alpha_q}(c)
  &=\sum_{k=0}^q(-1)^k\varphi_{\alpha(c)\alpha_0\cdots\widehat{\alpha_k}\cdots\alpha_q}(c),\\
  (K\delta\varphi)_{\alpha_0\cdots\alpha_q}(c)
  &=\varphi_{\alpha_0\cdots\alpha_q}(c)
  -\sum_{k=0}^q(-1)^k\varphi_{\alpha(c)\alpha_0\cdots\widehat{\alpha_k}\cdots\alpha_q}(c),
\end{align*}
so $\delta K+K\delta=\id$ in degrees $q\geq1$, and every cocycle
$\varphi=\delta(K\varphi)$ is a coboundary.
\end{proof}

\begin{intuition}[partitions of unity]
The smooth proof of the Mayer--Vietoris principle \citep[Prop.~8.5]{bott-tu1982}
uses a partition of unity $\{\rho_\alpha\}$ subordinate to the cover and
the operator $(K\varphi)_{\alpha_0\cdots\alpha_{q-1}}
=\sum_\alpha\rho_\alpha\,\varphi_{\alpha\alpha_0\cdots\alpha_{q-1}}$. In the
discrete setting the partition of unity is a family of indicator functions.
The choice $c\mapsto\alpha(c)$ partitions the cubes of $U$ among the
$U_\alpha$, and $\rho_\alpha$ is the indicator function of the cubes
assigned to $\alpha$. The functions $\rho_\alpha$ sum to $1$ and each is
supported in $U_\alpha$, and the smooth proof applies without change. The
same property of functions on discrete sets, extension by zero, is what
made the sheaves $\Nc^{p}$ flasque in \cref{thm:sheaf-coh}.
\end{intuition}

For two open sets $U$ and $V$ with $\Fp(X)=U\cup V$, the short exact
sequence $0\to\Ns(U\cup V)\to\Ns(U)\oplus\Ns(V)\to\Ns(U\cap V)\to0$ of
cochain complexes gives the usual long exact sequence.

\begin{corollary}[Mayer--Vietoris]\label{cor:mv}
If $\Fp(X)=U\cup V$ with $U$ and $V$ open, there is a long exact sequence
\[
  \cdots\to\HN{p}(X)\to\HN{p}(U)\oplus\HN{p}(V)\to\HN{p}(U\cap V)
  \to\HN{p+1}(X)\to\cdots.
\]
\end{corollary}

\begin{definition}[\v{C}ech--de Rham complex]\label{def:cech-deRham}
For an open cover $\cU$ of $\Fp(X)$ indexed by a totally ordered set, let
\[
  K^{p,q}=\prod_{\alpha_0<\dots<\alpha_q}\Nc^{p}(U_{\alpha_0\cdots\alpha_q}),
\]
a double complex with the de Rham differential $d$ in the $p$ direction
and the \v{C}ech differential $\delta$ in the $q$ direction. The
\emph{\v{C}ech--de Rham complex} of $\cU$ is the total complex, with
$D=d+(-1)^p\delta$ in bidegree $(p,q)$.
\end{definition}

\begin{theorem}\label{thm:cech-deRham}
The restriction map $\Ns(\Fp(X))\to K^{\bullet,0}\subseteq\Tot K$ induces an
isomorphism $\HN{\bullet}(X)\cong H^\bullet(\Tot K,D)$.
\end{theorem}

\begin{proof}
By \cref{prop:gmv}, applied to the alternating cochains, each row
$0\to\Nc^{p}(\Fp(X))\to K^{p,0}\to K^{p,1}\to\cdots$ is exact. The acyclic
assembly lemma \citep[2.7.3]{weibel1994}, or the generalized
Mayer--Vietoris principle of Bott and Tu~\citep[Prop.~8.8]{bott-tu1982}, then identifies
the cohomology of the total complex with that of the augmentation column
$\Ns(\Fp(X))$.
\end{proof}

Filtering the double complex by $q$ gives a spectral sequence with
$E_1^{p,q}=\prod\HN{p}(U_{\alpha_0\cdots\alpha_q})$ converging to
$\HN{\bullet}(X)$. In this sense the \v{C}ech--de Rham complex computes the
cohomology of a graph from the cohomology of its charts and their overlaps,
and \cref{cor:mv} is the case of two charts.

\begin{example}[The globe]\label{ex:globe-mv}
Consider the globe graph of \cref{fig:globe}, with $\Fp(X)=U_N\cup U_S$
the two polar charts. The overlap $U_N\cap U_S$ is the ring $E$ with its
$\Theta$-edges, a cycle graph with its rotation, so
$\HN{\bullet}(U_N\cap U_S)=(\Z,\Z)$ by \cref{ex:cycle}, there being no
degenerate cubes. Each chart is a truncated polar grid, and the argument of
\cref{prop:polar-normalized} applies verbatim to give
$\HN{\bullet}(U_N)=\HN{\bullet}(U_S)=(\Z,0,0)$. The Mayer--Vietoris
sequence reads
\[
  0\to\HN{0}(X)\to\Z\oplus\Z\to\Z\to\HN{1}(X)\to0\oplus0\to\Z\to\HN{2}(X)\to0.
\]
The map $\Z\oplus\Z\to\Z$ sends $(a,b)$ to the difference $a-b$ of the two
restrictions to the overlap. It is surjective with kernel $\Z$, so
$\HN{0}(X)=\Z$, exactness at the fourth term gives $\HN{1}(X)=0$, and the
last three terms give $\HN{2}(X)=\Z$. The globe graph with this atlas has
the cohomology of the $2$-sphere.
\end{example}

\subsection{The triangulated torus}
\label{sec:triangulated-torus}

%Claude note: this section is new, added at your request. It carries out the
%computation announced in Section 2.3.1.

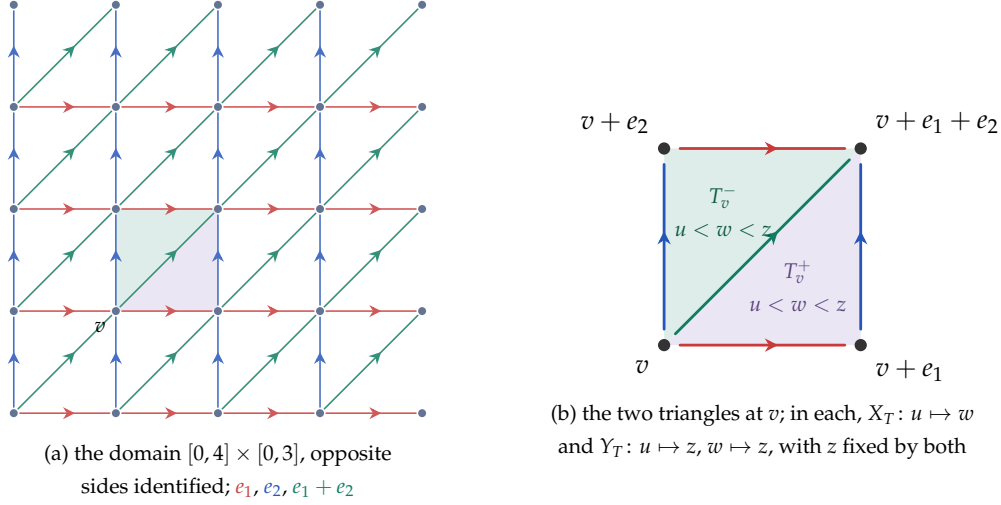
\begin{figure}[htb]
\centering
% Fundamental domain of Cay(Z/4 x Z/3, {(1,0),(0,1),(1,1)}) with the arrows
% of the Cayley graph (left), and the two triangles T_v^+ and T_v^- at a
% vertex v with their 2-simplex charts (right).
\begin{tikzpicture}[x=1cm,y=1cm]
\begin{scope}
  \def\s{1.35}
  % the two highlighted triangles at v=(1,1)
  \fill[fillcol!16] (1*\s,1*\s) -- (2*\s,1*\s) -- (2*\s,2*\s) -- cycle;
  \fill[cC!14] (1*\s,1*\s) -- (1*\s,2*\s) -- (2*\s,2*\s) -- cycle;
  % arrows: horizontal e1 (cA), vertical e2 (cB), diagonal e1+e2 (cC)
  \foreach \i in {0,...,3}{\foreach \j in {0,...,3}{
    \pgfmathtruncatemacro{\ii}{\i+1}\pgfmathtruncatemacro{\jj}{\j+1}
    \draw[cA!85,line width=0.6pt,cedge] (\i*\s,\j*\s) -- (\ii*\s,\j*\s);
    \draw[cB!85,line width=0.6pt,cedge] (\i*\s,\j*\s) -- (\i*\s,\jj*\s);
    \draw[cC!85,line width=0.6pt,cedge] (\i*\s,\j*\s) -- (\ii*\s,\jj*\s);}}
  \foreach \i in {0,...,4}{\foreach \j in {0,...,4}{\node[gvtx] at (\i*\s,\j*\s) {};}}
  \node[font=\scriptsize,below left] at (1*\s,1*\s) {$v$};
  \node[font=\scriptsize,align=center] at (2*\s,-0.75) {(a) the domain $[0,4]\times[0,3]$, opposite\\ sides identified; $\textcolor{cA}{e_1}$, $\textcolor{cB}{e_2}$, $\textcolor{cC}{e_1+e_2}$};
\end{scope}
\begin{scope}[xshift=8.6cm,yshift=0.9cm]
  \def\s{2.6}
  \fill[fillcol!16] (0,0) -- (\s,0) -- (\s,\s) -- cycle;
  \fill[cC!14] (0,0) -- (0,\s) -- (\s,\s) -- cycle;
  \draw[ar,cA] (0,0) -- (\s,0);
  \draw[ar,cB] (\s,0) -- (\s,\s);
  \draw[ar,cC] (0,0) -- (\s,\s);
  \draw[ar,cB] (0,0) -- (0,\s);
  \draw[ar,cA] (0,\s) -- (\s,\s);
  \node[vtx,label={[font=\small]below left:$v$}] at (0,0) {};
  \node[vtx,label={[font=\small]below right:$v+e_1$}] at (\s,0) {};
  \node[vtx,label={[font=\small]above left:$v+e_2$}] at (0,\s) {};
  \node[vtx,label={[font=\small]above right:$v+e_1+e_2$}] at (\s,\s) {};
  \node[font=\scriptsize,fillcol!70!black,align=center] at (0.68*\s,0.3*\s) {$T_v^+$\\ $u<w<z$};
  \node[font=\scriptsize,cC!70!black,align=center] at (0.3*\s,0.68*\s) {$T_v^-$\\ $u<w<z$};
  \node[font=\scriptsize,align=center] at (0.5*\s,-1.15) {(b) the two triangles at $v$; in each, $X_T\colon u\mapsto w$\\ and $Y_T\colon u\mapsto z$, $w\mapsto z$, with $z$ fixed by both};
\end{scope}
\end{tikzpicture}
\caption{A fundamental domain of the triangulated torus
$\Cay(\Z/4\times\Z/3,\{(1,0),(0,1),(1,1)\})$, with the arrows of the Cayley
graph, and the two triangles $T_v^{\pm}$ at a vertex $v$ with their
$2$-simplex charts. In each triangle the arrows order the vertices as
$u<w<z$, the field $X_T$ moves $u$ to $w$ and the field $Y_T$ moves $u$ and
$w$ to $z$.}
\label{fig:torus-tri-atlas}
\end{figure}

We return to the triangulated torus of \cref{sec:clique}. Let $n,m\geq3$
and let $G=\Cay(\Z/n\times\Z/m,\{(1,0),(0,1),(1,1)\})$, regarded as a
directed graph with an arrow $v\to v+s$ for each generator $s$. Write
$e_1=(1,0)$ and $e_2=(0,1)$. The graph has $nm$ vertices and $3nm$ edges,
and its triangles are the $2nm$ sets
\[
  T_v^+=\{v,\ v+e_1,\ v+e_1+e_2\},\qquad T_v^-=\{v,\ v+e_2,\ v+e_1+e_2\},
  \qquad v\in\Z/n\times\Z/m,
\]
which triangulate the torus, \cref{fig:torus-tri-atlas}. The single chart
of rank $3$ given by the generators has the cohomology of the $3$-torus by
\cref{cor:toral}. We replace it by an atlas with one chart of rank $2$ on
each triangle.

On each triangle the arrows of the Cayley graph form a transitive
tournament. In $T_v^+$ the arrows are $v\to v+e_1$, $v+e_1\to v+e_1+e_2$
and $v\to v+e_1+e_2$, so $v<v+e_1<v+e_1+e_2$, and similarly
$v<v+e_2<v+e_1+e_2$ in $T_v^-$. Each triangle $T=\{u<w<z\}$ therefore
carries the $2$-simplex chart of \cref{sec:two-simplex},
\[
  X_T\colon u\mapsto w,\ w\mapsto w,\ z\mapsto z,\qquad
  Y_T\colon u\mapsto z,\ w\mapsto z,\ z\mapsto z.
\]
The non-degenerate $1$-cubes of $(T,\N^2)$ are $\cube{u}{X}$, $\cube{u}{Y}$
and $\cube{w}{Y}$, with faces $(u,w)$, $(u,z)$ and $(w,z)$ in the order
$(d^0,d^1)$. These are the three arrows of $T$, each oriented as in the
Cayley graph, and the only non-degenerate $2$-cube is $\cube{u}{\{X,Y\}}$.

Two triangles $T$ and $T'$ sharing an edge $s\to t$ are glued along the
cubical subset $U_{TT'}=\{s,\ t,\ c\}\subseteq\Fp(\Qc(T))$, where $c$ is
the non-degenerate $1$-cube of $T$ with faces $d^0c=s$ and $d^1c=t$, and
$\tau_{TT'}$ sends $s,t,c$ to the corresponding cubes of $T'$. This is an
isomorphism of semi-cubical sets because the edge is oriented from $s$ to
$t$ in both triangles, so that $d^0$ and $d^1$ are preserved. Triangles
sharing only a vertex are glued along that vertex, and disjoint triangles
are not glued. The orientation of the Cayley graph is what makes the
transition maps compatible with the face maps. Every edge lies in exactly
two triangles and every vertex in six, and the identifications are those of
the simplicial complex, so the cocycle condition holds. Let $X$ denote the
resulting discrete coordinatized space.

\begin{proposition}\label{prop:triangulated-torus}
$\HN{\bullet}(X)\cong(\Z,\Z^2,\Z)$, the cohomology of the $2$-torus.
\end{proposition}

\begin{proof}
By \cref{prop:glued-sheaves}, $\Ns(\Fp(X))$ is the complex of functions on
the non-degenerate cubes of $X$. These are the $nm$ vertices, the $3nm$
arrows and the $2nm$ cubes $\cube{u_T}{\{X_T,Y_T\}}$, one for each triangle
$T$. On a $1$-cochain $\omega$ the differential is, by
\cref{prop:d-on-cubes} with the degenerate face $\cube{z}{X_T}$
contributing $0$,
\[
  (d\omega)\cube{u}{\{X_T,Y_T\}}
  =\omega\cube{w}{Y_T}-\omega\cube{u}{Y_T}+\omega\cube{u}{X_T}
  =\omega(w\to z)-\omega(u\to z)+\omega(u\to w),
\]
and on a $0$-cochain it is $(df)(s\to t)=f(t)-f(s)$. These are the
coboundaries of the simplicial cochain complex of the triangulated torus,
with the triangles oriented by the order of their vertices. Hence
$\Ns(\Fp(X))$ is the simplicial cochain complex of a triangulation of
$T^2$, and its cohomology is $H^\bullet(T^2;\Z)=(\Z,\Z^2,\Z)$.
\end{proof}

The same conclusion follows from the geometric realization of
\cref{sec:N-to-cubical}. The realization of the non-degenerate part of
$\Qc(X)$ has one $0$-cell for each vertex, one $1$-cell for each arrow and
one $2$-cell for each triangle, the $2$-cell $\cube{u}{\{X_T,Y_T\}}$ being a
square with the edge $\cube{z}{X_T}$ collapsed, and this cell complex is the
triangulated torus. Each chart is regular. If $q$ is the
non-degenerate $2$-cube of $(T,\N^2)$, then $U_q$ contains all non-degenerate cubes of
the chart, and its normalized complex is acyclic in positive degrees by
\cref{prop:simplex}; for every other cube $q'$ the set $U_{q'}$ has, after
normalization, at most one non-degenerate $1$-cube and no $2$-cube, and is
acyclic. So \cref{thm:sheaf-coh} also identifies $\HN{\bullet}(X)$ with
$H^\bullet(\Fp(X),\underline\Z)$.

The contrast with the single chart of rank $3$ illustrates the dependence
of the cohomology on the atlas. The graph is the same, but the atlas of
$2nm$ charts of rank $2$ supplies one $2$-cell for each triangle, while the
single chart of rank $3$ supplies $2$-cells and $3$-cells indexed by the
generators and computes the cohomology of a quotient of $\Z^3$. The same
construction applies to any triangulated surface whose $1$-skeleton is
given an acyclic orientation on every triangle, and \cref{sec:ex-triangles}
applies it to the complete graphs $K_4$, $K_5$ and $K_6$.

\chapter{Examples}
\label{ch:examples}

%Claude note: this chapter has been rewritten from your rough version. It
%presents the figures in the folder ch6-figures with short explanations.
%The cohomology groups are the ones recorded in your figure files; I have
%checked each against the Euler characteristic of the normalized complex
%and the evident topology of the filled complex, and corrected two typos in
%the torus captions (Z3 and Z7 should read Z^13 and Z^17). The
%identification of H_N with the cohomology of the "filled" space is the
%sheaf-cohomology theorem of Chapter 5 together with the cellular
%description of the cubical set in Chapter 4, and it needs the atlases to
%be regular; I have not checked regularity of every chart below.

This chapter collects examples of graphs with atlases of $\N$-sets and
records the normalized de Rham cohomology $\HN{\bullet}$ of the resulting
discrete coordinatized spaces. Since the cohomology depends on the atlas,
several atlases on the same graph are compared. We first say what it means
for an atlas to be an atlas on a given graph.

\begin{definition}[Induced graph, coordinatized graph]\label{def:coordinatized-graph}
Let $X$ be a discrete coordinatized space with atlas
$(X_\alpha,\N^{n_\alpha})$. Its \emph{induced directed graph} has vertex set
$\Fp(X)^{(0)}$, the $0$-cubes, and a directed edge $x\to e_i\cdot x$ for
each non-degenerate $1$-cube $\cube{x}{i}$. A simple undirected graph $G$
is \emph{coordinatized} by $X$ if the induced directed graph, with its
directions forgotten and self-loops removed, is isomorphic to $G$.
\end{definition}

In the figures, each chart is a subgraph on which several commuting
coordinate fields are drawn in colour, one colour per field. A small loop
at a vertex indicates that the field of that colour fixes the vertex, so
that the corresponding $1$-cube is degenerate. Charts are glued along the
cubical subsets they have in common.

The Euler characteristic of
$\HN{\bullet}(X)$ is that of the normalized complex,
\[
  \chi=\sum_p(-1)^p\,\#\{\text{non-degenerate }p\text{-cubes of }X\},
\]
and we record the numbers of non-degenerate cubes as a cell count
$(c_0,c_1,c_2,\dots)$.

\section{One chart: Cayley graphs}
\label{sec:ex-cayley}

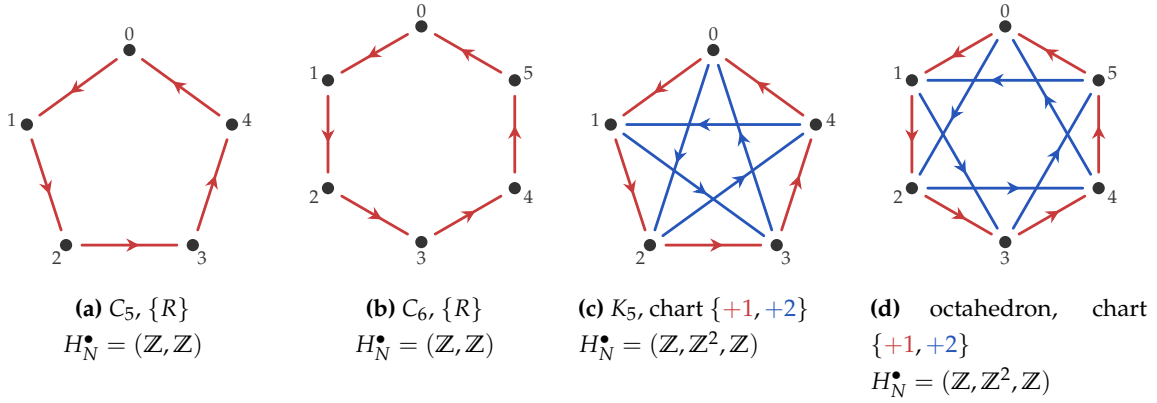
\begin{figure}[htb]
\centering
\renewcommand{\figscale}{0.62}
\begin{subfigure}[t]{0.24\textwidth}\centering
\begin{tikzpicture}[x=1cm,y=1cm,scale=\figscale]
  \draw[cA,ar] (0.0000,2.3000) -- (-2.1874,0.7107);
  \draw[cA,ar] (-2.1874,0.7107) -- (-1.3519,-1.8607);
  \draw[cA,ar] (-1.3519,-1.8607) -- (1.3519,-1.8607);
  \draw[cA,ar] (1.3519,-1.8607) -- (2.1874,0.7107);
  \draw[cA,ar] (2.1874,0.7107) -- (0.0000,2.3000);
  \node[vtx] at (0.0000,2.3000) {};
  \node[vtx] at (-2.1874,0.7107) {};
  \node[vtx] at (-1.3519,-1.8607) {};
  \node[vtx] at (1.3519,-1.8607) {};
  \node[vtx] at (2.1874,0.7107) {};
  \node[font=\tiny,gray!45!black] at (0.0000,2.6400) {0};
  \node[font=\tiny,gray!45!black] at (-2.5108,0.8158) {1};
  \node[font=\tiny,gray!45!black] at (-1.5518,-2.1358) {2};
  \node[font=\tiny,gray!45!black] at (1.5518,-2.1358) {3};
  \node[font=\tiny,gray!45!black] at (2.5108,0.8158) {4};
\end{tikzpicture}
\caption{$C_5$, $\{R\}$\\ $\HN{\bullet}=(\Z,\Z)$}
\end{subfigure}\hfill
\begin{subfigure}[t]{0.24\textwidth}\centering
\begin{tikzpicture}[x=1cm,y=1cm,scale=\figscale]
  \draw[cA,ar] (0.0000,2.3000) -- (-1.9919,1.1500);
  \draw[cA,ar] (-1.9919,1.1500) -- (-1.9919,-1.1500);
  \draw[cA,ar] (-1.9919,-1.1500) -- (-0.0000,-2.3000);
  \draw[cA,ar] (-0.0000,-2.3000) -- (1.9919,-1.1500);
  \draw[cA,ar] (1.9919,-1.1500) -- (1.9919,1.1500);
  \draw[cA,ar] (1.9919,1.1500) -- (0.0000,2.3000);
  \node[vtx] at (0.0000,2.3000) {};
  \node[vtx] at (-1.9919,1.1500) {};
  \node[vtx] at (-1.9919,-1.1500) {};
  \node[vtx] at (-0.0000,-2.3000) {};
  \node[vtx] at (1.9919,-1.1500) {};
  \node[vtx] at (1.9919,1.1500) {};
  \node[font=\tiny,gray!45!black] at (0.0000,2.6400) {0};
  \node[font=\tiny,gray!45!black] at (-2.2863,1.3200) {1};
  \node[font=\tiny,gray!45!black] at (-2.2863,-1.3200) {2};
  \node[font=\tiny,gray!45!black] at (-0.0000,-2.6400) {3};
  \node[font=\tiny,gray!45!black] at (2.2863,-1.3200) {4};
  \node[font=\tiny,gray!45!black] at (2.2863,1.3200) {5};
\end{tikzpicture}
\caption{$C_6$, $\{R\}$\\ $\HN{\bullet}=(\Z,\Z)$}
\end{subfigure}\hfill
\begin{subfigure}[t]{0.24\textwidth}\centering
\begin{tikzpicture}[x=1cm,y=1cm,scale=\figscale]
  \draw[cA,ar] (0.0000,2.3000) -- (-2.1874,0.7107);
  \draw[cA,ar] (-2.1874,0.7107) -- (-1.3519,-1.8607);
  \draw[cA,ar] (-1.3519,-1.8607) -- (1.3519,-1.8607);
  \draw[cA,ar] (1.3519,-1.8607) -- (2.1874,0.7107);
  \draw[cA,ar] (2.1874,0.7107) -- (0.0000,2.3000);
  \draw[cB,ar] (0.0000,2.3000) -- (-1.3519,-1.8607);
  \draw[cB,ar] (-2.1874,0.7107) -- (1.3519,-1.8607);
  \draw[cB,ar] (-1.3519,-1.8607) -- (2.1874,0.7107);
  \draw[cB,ar] (1.3519,-1.8607) -- (0.0000,2.3000);
  \draw[cB,ar] (2.1874,0.7107) -- (-2.1874,0.7107);
  \node[vtx] at (0.0000,2.3000) {};
  \node[vtx] at (-2.1874,0.7107) {};
  \node[vtx] at (-1.3519,-1.8607) {};
  \node[vtx] at (1.3519,-1.8607) {};
  \node[vtx] at (2.1874,0.7107) {};
  \node[font=\tiny,gray!45!black] at (0.0000,2.6400) {0};
  \node[font=\tiny,gray!45!black] at (-2.5108,0.8158) {1};
  \node[font=\tiny,gray!45!black] at (-1.5518,-2.1358) {2};
  \node[font=\tiny,gray!45!black] at (1.5518,-2.1358) {3};
  \node[font=\tiny,gray!45!black] at (2.5108,0.8158) {4};
\end{tikzpicture}
\caption{$K_5$, chart $\{\textcolor{cA}{+1},\textcolor{cB}{+2}\}$\\ $\HN{\bullet}=(\Z,\Z^2,\Z)$}
\end{subfigure}\hfill
\begin{subfigure}[t]{0.24\textwidth}\centering
\begin{tikzpicture}[x=1cm,y=1cm,scale=\figscale]
  \draw[cA,ar] (0.0000,2.3000) -- (-1.9919,1.1500);
  \draw[cA,ar] (-1.9919,1.1500) -- (-1.9919,-1.1500);
  \draw[cA,ar] (-1.9919,-1.1500) -- (-0.0000,-2.3000);
  \draw[cA,ar] (-0.0000,-2.3000) -- (1.9919,-1.1500);
  \draw[cA,ar] (1.9919,-1.1500) -- (1.9919,1.1500);
  \draw[cA,ar] (1.9919,1.1500) -- (0.0000,2.3000);
  \draw[cB,ar] (0.0000,2.3000) -- (-1.9919,-1.1500);
  \draw[cB,ar] (-1.9919,1.1500) -- (-0.0000,-2.3000);
  \draw[cB,ar] (-1.9919,-1.1500) -- (1.9919,-1.1500);
  \draw[cB,ar] (-0.0000,-2.3000) -- (1.9919,1.1500);
  \draw[cB,ar] (1.9919,-1.1500) -- (0.0000,2.3000);
  \draw[cB,ar] (1.9919,1.1500) -- (-1.9919,1.1500);
  \node[vtx] at (0.0000,2.3000) {};
  \node[vtx] at (-1.9919,1.1500) {};
  \node[vtx] at (-1.9919,-1.1500) {};
  \node[vtx] at (-0.0000,-2.3000) {};
  \node[vtx] at (1.9919,-1.1500) {};
  \node[vtx] at (1.9919,1.1500) {};
  \node[font=\tiny,gray!45!black] at (0.0000,2.6400) {0};
  \node[font=\tiny,gray!45!black] at (-2.2863,1.3200) {1};
  \node[font=\tiny,gray!45!black] at (-2.2863,-1.3200) {2};
  \node[font=\tiny,gray!45!black] at (-0.0000,-2.6400) {3};
  \node[font=\tiny,gray!45!black] at (2.2863,-1.3200) {4};
  \node[font=\tiny,gray!45!black] at (2.2863,1.3200) {5};
\end{tikzpicture}
\caption{octahedron, chart $\{\textcolor{cA}{+1},\textcolor{cB}{+2}\}$\\ $\HN{\bullet}=(\Z,\Z^2,\Z)$}
\end{subfigure}
\caption{Cayley graphs with their homogeneous chart. By \cref{cor:toral} a
single chart gives the cohomology of a torus.}
\label{fig:ex-cayley}
\end{figure}

By \cref{thm:cayley} the finite connected graphs carrying a single
homogeneous chart are the Cayley graphs of finite abelian groups, and by
Shapiro's lemma, their cohomology is that of an $n$-torus,
$n$ being the number of generators. Since the coordinates have no fixed
points there are no degenerate cubes and $\HN{\bullet}=\HdR{\bullet}$.
\Cref{fig:ex-cayley} shows four instances. The cycle graphs $C_5$ and $C_6$
carry only the rotation chart and have the cohomology of the circle. The
complete graph $K_5=\Cay(\Z/5,\{+1,+2\})$ and the octahedron graph, which is (surprisingly) 
isomorphic to $\Cay(\Z/6,\{+1,+2\})$, carry charts of rank $2$ and have the cohomology of
the $2$-torus.

\section{One chart: cones}
\label{sec:ex-cones}

A chart with fixed points can have the normalized cohomology of a point.
The simplest case is a cone. Let $G$ be a graph with an $\N$-set structure
$(V,\N^{n})$ and add a new vertex $0$ joined to every vertex of $G$. Let
$e_1$ send every vertex of $G$ to $0$ and fix $0$, and let the coordinates
of $G$ fix $0$. Then every $p$-cube of $G$ becomes the base of a
$(p+1)$-cube whose $e_1$-faces collapse to $0$, and the result is a filled
cone on $G$, with $\HN{\bullet}=(\Z,0,\dots,0)$.

\begin{figure}[htb]
\centering
\renewcommand{\figscale}{0.7}
\begin{subfigure}[t]{0.3\textwidth}\centering
\begin{tikzpicture}[x=1cm,y=1cm,scale=\figscale]
  \draw[cA,ar] (0.0000,2.1000) -- (0.0000,0.0000);
  \draw[cA,ar] (-1.8200,-1.0500) -- (0.0000,0.0000);
  \draw[cA,ar] (1.8200,-1.0500) -- (0.0000,0.0000);
  \draw[cB,ar] (0.0000,2.1000) -- (-1.8200,-1.0500);
  \draw[cB,ar] (-1.8200,-1.0500) -- (1.8200,-1.0500);
  \draw[cB,ar] (1.8200,-1.0500) -- (0.0000,2.1000);
  \draw[cA,lp] (-0.2088,0.1200) arc (150:-160:0.240);
  \draw[cB,lp] (-0.2088,0.1200) arc (150:-160:0.240);
  \node[vtx] at (0.0000,0.0000) {};
  \node[vtx] at (0.0000,2.1000) {};
  \node[vtx] at (-1.8200,-1.0500) {};
  \node[vtx] at (1.8200,-1.0500) {};
  \node[font=\tiny,gray!45!black] at (0.0000,0.0000) {0};
  \node[font=\tiny,gray!45!black] at (0.0000,2.4400) {1};
  \node[font=\tiny,gray!45!black] at (-2.1145,-1.2199) {2};
  \node[font=\tiny,gray!45!black] at (2.1145,-1.2199) {3};
\end{tikzpicture}
\caption{$K_4$: cone on a triangle\\ $\textcolor{cA}{e_1}$ to the apex, $\textcolor{cB}{e_2}$ rotates\\ $\HN{\bullet}=(\Z)$}
\end{subfigure}\hfill
\begin{subfigure}[t]{0.33\textwidth}\centering
\begin{tikzpicture}[x=1cm,y=1cm,scale=\figscale]
  \draw[cA,ar] (0.0000,2.2000) -- (0.0000,0.0000);
  \draw[cA,ar] (-2.2000,0.0000) -- (0.0000,0.0000);
  \draw[cA,ar] (-0.0000,-2.2000) -- (0.0000,0.0000);
  \draw[cA,ar] (2.2000,-0.0000) -- (0.0000,0.0000);
  \draw[cB,ar] (-2.2000,0.0000) -- (0.0000,2.2000);
  \draw[cB,ar] (-0.0000,-2.2000) -- (0.0000,2.2000);
  \draw[cB,ar] (2.2000,-0.0000) -- (0.0000,2.2000);
  \draw[cC,ar] (-2.2000,0.0000) -- (-0.0000,-2.2000);
  \draw[cC,ar] (-0.0000,-2.2000) -- (2.2000,-0.0000);
  \draw[cC,ar] (2.2000,-0.0000) -- (-2.2000,0.0000);
  \draw[cA,lp] (-0.2088,0.1200) arc (150:-160:0.240);
  \draw[cB,lp] (-0.2088,0.1200) arc (150:-160:0.240);
  \draw[cB,lp] (-0.2088,2.3200) arc (150:-160:0.240);
  \draw[cC,lp] (-0.2088,0.1200) arc (150:-160:0.240);
  \draw[cC,lp] (-0.2088,2.3200) arc (150:-160:0.240);
  \node[vtx] at (0.0000,0.0000) {};
  \node[vtx] at (0.0000,2.2000) {};
  \node[vtx] at (-2.2000,0.0000) {};
  \node[vtx] at (-0.0000,-2.2000) {};
  \node[vtx] at (2.2000,-0.0000) {};
  \node[font=\tiny,gray!45!black] at (0.3325,0.0709) {0};
  \node[font=\tiny,gray!45!black] at (0.0000,2.5400) {1};
  \node[font=\tiny,gray!45!black] at (-2.5400,0.0000) {2};
  \node[font=\tiny,gray!45!black] at (-0.0000,-2.5400) {3};
  \node[font=\tiny,gray!45!black] at (2.5400,-0.0000) {4};
\end{tikzpicture}
\caption{$K_5$: one chart with three fields\\ $\textcolor{cA}{e_1}$, $\textcolor{cB}{e_2}$, $\textcolor{cC}{e_3}$\\ $\HN{\bullet}=(\Z)$}
\end{subfigure}\hfill
\begin{subfigure}[t]{0.33\textwidth}\centering
\begin{tikzpicture}[x=1cm,y=1cm,scale=\figscale]
  \draw[cA,ar] (0.0000,2.3000) -- (0.0000,0.0000);
  \draw[cA,ar] (-2.1874,0.7107) -- (0.0000,0.0000);
  \draw[cA,ar] (-1.3519,-1.8607) -- (0.0000,0.0000);
  \draw[cA,ar] (1.3519,-1.8607) -- (0.0000,0.0000);
  \draw[cA,ar] (2.1874,0.7107) -- (0.0000,0.0000);
  \draw[cB,ar] (0.0000,2.3000) -- (-2.1874,0.7107);
  \draw[cB,ar] (-2.1874,0.7107) -- (-1.3519,-1.8607);
  \draw[cB,ar] (-1.3519,-1.8607) -- (1.3519,-1.8607);
  \draw[cB,ar] (1.3519,-1.8607) -- (2.1874,0.7107);
  \draw[cB,ar] (2.1874,0.7107) -- (0.0000,2.3000);
  \draw[cC,ar] (0.0000,2.3000) -- (-1.3519,-1.8607);
  \draw[cC,ar] (-2.1874,0.7107) -- (1.3519,-1.8607);
  \draw[cC,ar] (-1.3519,-1.8607) -- (2.1874,0.7107);
  \draw[cC,ar] (1.3519,-1.8607) -- (0.0000,2.3000);
  \draw[cC,ar] (2.1874,0.7107) -- (-2.1874,0.7107);
  \draw[cA,lp] (-0.2088,0.1200) arc (150:-160:0.240);
  \draw[cB,lp] (-0.2088,0.1200) arc (150:-160:0.240);
  \draw[cC,lp] (-0.2088,0.1200) arc (150:-160:0.240);
  \node[vtx] at (0.0000,0.0000) {};
  \node[vtx] at (0.0000,2.3000) {};
  \node[vtx] at (-2.1874,0.7107) {};
  \node[vtx] at (-1.3519,-1.8607) {};
  \node[vtx] at (1.3519,-1.8607) {};
  \node[vtx] at (2.1874,0.7107) {};
  \node[font=\tiny,gray!45!black] at (0.0000,0.3400) {0};
  \node[font=\tiny,gray!45!black] at (0.0000,2.6400) {1};
  \node[font=\tiny,gray!45!black] at (-2.5108,0.8158) {2};
  \node[font=\tiny,gray!45!black] at (-1.5518,-2.1358) {3};
  \node[font=\tiny,gray!45!black] at (1.5518,-2.1358) {4};
  \node[font=\tiny,gray!45!black] at (2.5108,0.8158) {5};
\end{tikzpicture}
\caption{$K_6$: cone on $\Cay(\Z/5,\{+1,+2\})$\\ $\textcolor{cA}{e_1}$ to the apex, $\textcolor{cB}{e_2}$, $\textcolor{cC}{e_3}$\\ $\HN{\bullet}=(\Z)$}
\end{subfigure}
\caption{Cones. A coordinate field pointing at an apex makes the chart
contractible. The loops mark the fixed points, whose cubes are degenerate.}
\label{fig:ex-cones}
\end{figure}
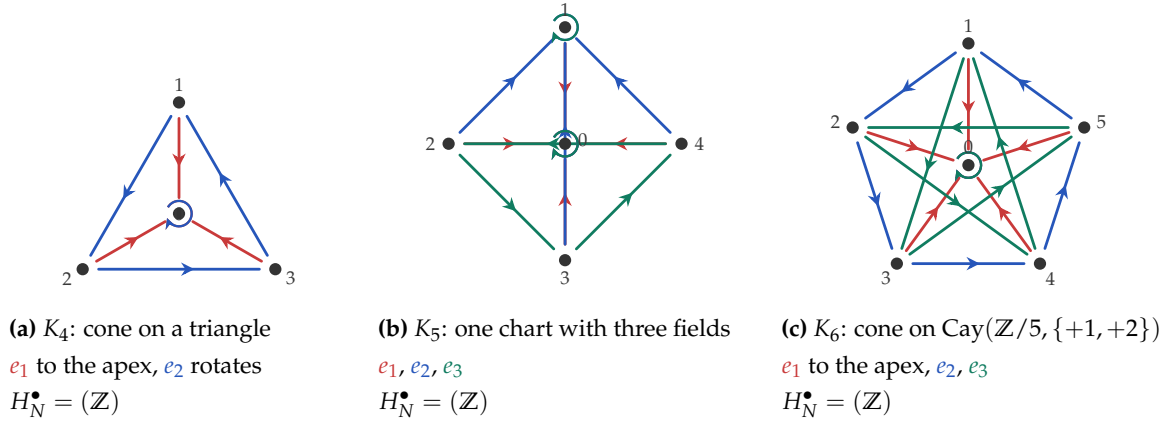

\Cref{fig:ex-cones} shows three cones. On $K_4$, the field $e_1$ sends the
three outer vertices to the apex and $e_2$ rotates them, giving a cone on a
triangle with cells $(4,6,3)$ and $\chi=1$. On $K_5$ the apex is joined to a
$K_4$ which is itself a cone, with three fields and cells $(5,10,9,3)$. On
$K_6$ the apex is joined to $K_5=\Cay(\Z/5,\{+1,+2\})$ with its torus
chart, and the cone on a torus is contractible.

\section{Atlases of triangles on complete graphs}
\label{sec:ex-triangles}

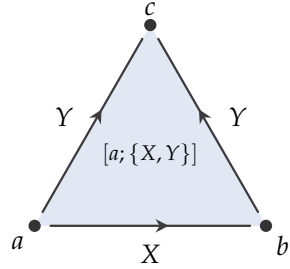
\begin{figure}[htb]
\centering
\renewcommand{\figscale}{0.9}
\begin{tikzpicture}[x=1cm,y=1cm,scale=\figscale]
  \fill[fillc] (0,0) -- (3.4,0) -- (1.7,2.95) -- cycle;
  \draw[grar] (0,0) -- (3.4,0);
  \draw[grar] (3.4,0) -- (1.7,2.95);
  \draw[grar] (0,0) -- (1.7,2.95);
  \node[vtx] at (0,0) {};  \node[font=\small,below left] at (0,0) {$a$};
  \node[vtx] at (3.4,0) {}; \node[font=\small,below right] at (3.4,0) {$b$};
  \node[vtx] at (1.7,2.95) {}; \node[font=\small,above] at (1.7,2.95) {$c$};
  \node[font=\small,below] at (1.7,-0.12) {$X$};
  \node[font=\small,left] at (0.72,1.60) {$Y$};
  \node[font=\small,right] at (2.70,1.60) {$Y$};
  \node[font=\scriptsize,align=center] at (1.7,1.05) {$[a;\{X,Y\}]$};
\end{tikzpicture}
\caption{The triangle chart on $T=\{a<b<c\}$. Here $X\colon a\mapsto b$,
$b\mapsto b$, $c\mapsto c$ and $Y\colon a\mapsto c$, $b\mapsto c$,
$c\mapsto c$. The fields commute, $\cube{c}{X}$ is the only degenerate face
of the square $\cube{a}{\{X,Y\}}$, and the $2$-cube has the three edges of
$T$ as its boundary.}
\label{fig:rp2-chart}
\end{figure}

The $2$-simplex of \cref{sec:two-simplex} is a chart of rank $2$ on a
triangle $T=\{a<b<c\}$, with $X\colon a\mapsto b$ fixing $b$ and $c$, and
$Y\colon a\mapsto c$, $b\mapsto c$ fixing $c$. Its only non-degenerate
$2$-cube $\cube{a}{\{X,Y\}}$ has the three edges of $T$ as its boundary, so
the chart fills the triangle, \cref{fig:rp2-chart}. As in
\cref{sec:triangulated-torus}, placing one such chart on each triangle of
a triangulated surface whose edges are oriented acyclically on every
triangle, and gluing adjacent charts along their common edges, produces an
atlas whose normalized complex is the simplicial cochain complex of the
surface.

In the complete graph $K_n$ every triple of vertices is a triangle. Orient
the edges of $K_n$ by the order $0<1<\dots<n-1$ of the vertices and make
every triangle a chart.
\Cref{fig:ex-complete} shows the results. On $K_4$ the four triangles form
the boundary of a tetrahedron, a sphere, with cells $(4,6,4)$ and
$\HN{\bullet}=(\Z,0,\Z)$. On $K_5$ the ten triangles form the $2$-skeleton
of the $4$-simplex, with cells $(5,10,10)$ and $\HN{\bullet}=(\Z,0,\Z^4)$,
a wedge of four $2$-spheres. On $K_6$ the twenty triangles give cells
$(6,15,20)$ and $\HN{\bullet}=(\Z,0,\Z^{10})$.

\begin{figure}[htb]
\centering
\renewcommand{\figscale}{0.7}
\begin{subfigure}[t]{0.3\textwidth}\centering
\begin{tikzpicture}[x=1cm,y=1cm,scale=\figscale]
  \fill[cA,opacity=0.13] (-1.9000,-1.1000) -- (1.9000,-1.1000) -- (0.0000,2.1000) -- cycle;
  \fill[cB,opacity=0.13] (-1.9000,-1.1000) -- (1.9000,-1.1000) -- (0.0000,0.0500) -- cycle;
  \fill[cC,opacity=0.13] (-1.9000,-1.1000) -- (0.0000,2.1000) -- (0.0000,0.0500) -- cycle;
  \fill[cD,opacity=0.13] (1.9000,-1.1000) -- (0.0000,2.1000) -- (0.0000,0.0500) -- cycle;
  \draw[cA,ar] (-1.9000,-1.1000) -- (1.9000,-1.1000);
  \draw[cA,ar] (-1.9000,-1.1000) -- (0.0000,2.1000);
  \draw[cA,ar] (-1.9000,-1.1000) -- (0.0000,0.0500);
  \draw[cA,ar] (1.9000,-1.1000) -- (0.0000,2.1000);
  \draw[cA,ar] (1.9000,-1.1000) -- (0.0000,0.0500);
  \draw[cA,ar] (0.0000,2.1000) -- (0.0000,0.0500);
  \node[vtx] at (-1.9000,-1.1000) {};
  \node[vtx] at (1.9000,-1.1000) {};
  \node[vtx] at (0.0000,2.1000) {};
  \node[vtx] at (0.0000,0.0500) {};
  \node[font=\tiny,gray!45!black] at (-2.1951,-1.2689) {0};
  \node[font=\tiny,gray!45!black] at (2.1951,-1.2689) {1};
  \node[font=\tiny,gray!45!black] at (0.0000,2.4400) {2};
  \node[font=\tiny,gray!45!black] at (0.0000,0.3900) {3};
\end{tikzpicture}
\caption{$K_4$, four triangular charts\\ cells $(4,6,4)$, $\HN{\bullet}=(\Z,0,\Z)$}
\end{subfigure}\hfill
\begin{subfigure}[t]{0.33\textwidth}\centering
\begin{tikzpicture}[x=1cm,y=1cm,scale=\figscale]
  \draw[cA,ar] (0.0000,2.3000) -- (-2.1874,0.7107);
  \draw[cA,ar] (0.0000,2.3000) -- (-1.3519,-1.8607);
  \draw[cA,ar] (0.0000,2.3000) -- (1.3519,-1.8607);
  \draw[cA,ar] (0.0000,2.3000) -- (2.1874,0.7107);
  \draw[cA,ar] (-2.1874,0.7107) -- (-1.3519,-1.8607);
  \draw[cA,ar] (-2.1874,0.7107) -- (1.3519,-1.8607);
  \draw[cA,ar] (-2.1874,0.7107) -- (2.1874,0.7107);
  \draw[cA,ar] (-1.3519,-1.8607) -- (1.3519,-1.8607);
  \draw[cA,ar] (-1.3519,-1.8607) -- (2.1874,0.7107);
  \draw[cA,ar] (1.3519,-1.8607) -- (2.1874,0.7107);
  \node[vtx] at (0.0000,2.3000) {};
  \node[vtx] at (-2.1874,0.7107) {};
  \node[vtx] at (-1.3519,-1.8607) {};
  \node[vtx] at (1.3519,-1.8607) {};
  \node[vtx] at (2.1874,0.7107) {};
  \node[font=\tiny,gray!45!black] at (0.0000,2.6400) {0};
  \node[font=\tiny,gray!45!black] at (-2.5108,0.8158) {1};
  \node[font=\tiny,gray!45!black] at (-1.5518,-2.1358) {2};
  \node[font=\tiny,gray!45!black] at (1.5518,-2.1358) {3};
  \node[font=\tiny,gray!45!black] at (2.5108,0.8158) {4};
\end{tikzpicture}
\caption{$K_5$, one chart per triangle\\ cells $(5,10,10)$, $\HN{\bullet}=(\Z,0,\Z^4)$}
\end{subfigure}\hfill
\begin{subfigure}[t]{0.33\textwidth}\centering
\begin{tikzpicture}[x=1cm,y=1cm,scale=\figscale]
  \draw[cA,ar] (0.0000,2.3000) -- (-1.9919,1.1500);
  \draw[cA,ar] (0.0000,2.3000) -- (-1.9919,-1.1500);
  \draw[cA,ar] (0.0000,2.3000) -- (-0.0000,-2.3000);
  \draw[cA,ar] (0.0000,2.3000) -- (1.9919,-1.1500);
  \draw[cA,ar] (0.0000,2.3000) -- (1.9919,1.1500);
  \draw[cA,ar] (-1.9919,1.1500) -- (-1.9919,-1.1500);
  \draw[cA,ar] (-1.9919,1.1500) -- (-0.0000,-2.3000);
  \draw[cA,ar] (-1.9919,1.1500) -- (1.9919,-1.1500);
  \draw[cA,ar] (-1.9919,1.1500) -- (1.9919,1.1500);
  \draw[cA,ar] (-1.9919,-1.1500) -- (-0.0000,-2.3000);
  \draw[cA,ar] (-1.9919,-1.1500) -- (1.9919,-1.1500);
  \draw[cA,ar] (-1.9919,-1.1500) -- (1.9919,1.1500);
  \draw[cA,ar] (-0.0000,-2.3000) -- (1.9919,-1.1500);
  \draw[cA,ar] (-0.0000,-2.3000) -- (1.9919,1.1500);
  \draw[cA,ar] (1.9919,-1.1500) -- (1.9919,1.1500);
  \node[vtx] at (0.0000,2.3000) {};
  \node[vtx] at (-1.9919,1.1500) {};
  \node[vtx] at (-1.9919,-1.1500) {};
  \node[vtx] at (-0.0000,-2.3000) {};
  \node[vtx] at (1.9919,-1.1500) {};
  \node[vtx] at (1.9919,1.1500) {};
  \node[font=\tiny,gray!45!black] at (0.0000,2.6400) {0};
  \node[font=\tiny,gray!45!black] at (-2.2863,1.3200) {1};
  \node[font=\tiny,gray!45!black] at (-2.2863,-1.3200) {2};
  \node[font=\tiny,gray!45!black] at (-0.0000,-2.6400) {3};
  \node[font=\tiny,gray!45!black] at (2.2863,-1.3200) {4};
  \node[font=\tiny,gray!45!black] at (2.2863,1.3200) {5};
\end{tikzpicture}
\caption{$K_6$, one chart per triangle\\ cells $(6,15,20)$, $\HN{\bullet}=(\Z,0,\Z^{10})$}
\end{subfigure}
\caption{Complete graphs with all triangles filled. Edges are oriented by
the order of the vertices and each shaded triangle is one chart.}
\label{fig:ex-complete}
\end{figure}
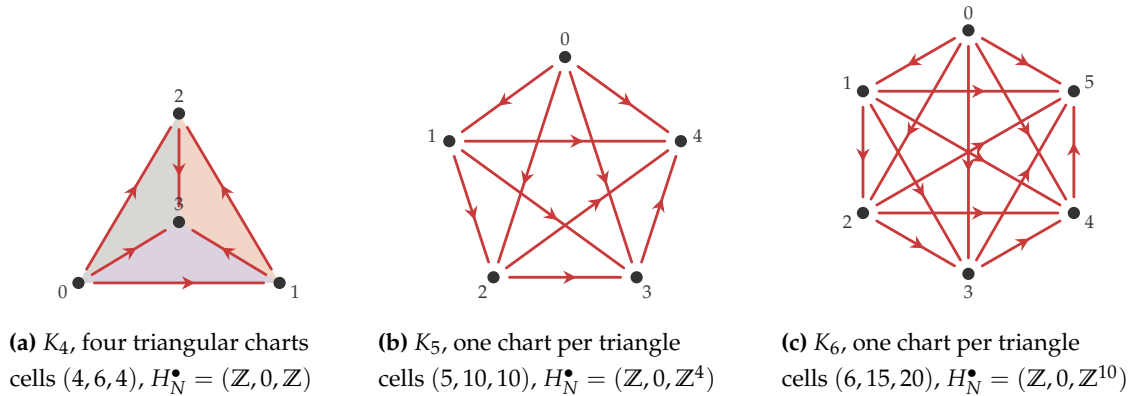

\begin{figure}[htb]
\centering
\renewcommand{\figscale}{0.78}
\begin{tikzpicture}[x=1cm,y=1cm,scale=\figscale]
  \fill[fillc] (0.0000,4.2000) -- (3.6373,2.1000) -- (0.7750,1.3423) -- cycle;
  \fill[fillc] (3.6373,2.1000) -- (0.7750,1.3423) -- (0.7750,-1.3423) -- cycle;
  \fill[fillc] (3.6373,2.1000) -- (3.6373,-2.1000) -- (0.7750,-1.3423) -- cycle;
  \fill[fillc] (3.6373,-2.1000) -- (0.0000,-4.2000) -- (0.7750,-1.3423) -- cycle;
  \fill[fillc] (0.0000,-4.2000) -- (0.7750,-1.3423) -- (-1.5500,0.0000) -- cycle;
  \fill[fillc] (0.0000,-4.2000) -- (-3.6373,-2.1000) -- (-1.5500,0.0000) -- cycle;
  \fill[fillc] (-3.6373,-2.1000) -- (-3.6373,2.1000) -- (-1.5500,0.0000) -- cycle;
  \fill[fillc] (-3.6373,2.1000) -- (-1.5500,0.0000) -- (0.7750,1.3423) -- cycle;
  \fill[fillc] (-3.6373,2.1000) -- (0.0000,4.2000) -- (0.7750,1.3423) -- cycle;
  \fill[fillc] (0.7750,1.3423) -- (0.7750,-1.3423) -- (-1.5500,0.0000) -- cycle;
  \draw[idar] (0.0000,4.2000) -- (3.6373,2.1000);
  \draw[idc,line width=1.1pt] (1.9187,3.3232) -- (1.7187,2.9768);
  \draw[grar] (3.6373,2.1000) -- (0.7750,1.3423);
  \draw[grar] (0.0000,4.2000) -- (0.7750,1.3423);
  \draw[grar] (0.7750,1.3423) -- (0.7750,-1.3423);
  \draw[grar] (3.6373,2.1000) -- (0.7750,-1.3423);
  \draw[idar] (3.6373,2.1000) -- (3.6373,-2.1000);
  \draw[idc,line width=1.1pt] (3.8373,0.1000) -- (3.4373,0.1000);
  \draw[idc,line width=1.1pt] (3.8373,-0.1000) -- (3.4373,-0.1000);
  \draw[grar] (3.6373,-2.1000) -- (0.7750,-1.3423);
  \draw[idar] (0.0000,-4.2000) -- (3.6373,-2.1000);
  \draw[idc,line width=1.1pt] (1.5454,-3.0768) -- (1.7454,-3.4232);
  \draw[idc,line width=1.1pt] (1.7187,-2.9768) -- (1.9187,-3.3232);
  \draw[idc,line width=1.1pt] (1.8919,-2.8768) -- (2.0919,-3.2232);
  \draw[grar] (0.0000,-4.2000) -- (0.7750,-1.3423);
  \draw[grar] (0.7750,-1.3423) -- (-1.5500,0.0000);
  \draw[grar] (0.0000,-4.2000) -- (-1.5500,0.0000);
  \draw[idar] (0.0000,-4.2000) -- (-3.6373,-2.1000);
  \draw[idc,line width=1.1pt] (-1.9187,-3.3232) -- (-1.7187,-2.9768);
  \draw[grar] (-3.6373,-2.1000) -- (-1.5500,0.0000);
  \draw[idar] (-3.6373,-2.1000) -- (-3.6373,2.1000);
  \draw[idc,line width=1.1pt] (-3.8373,-0.1000) -- (-3.4373,-0.1000);
  \draw[idc,line width=1.1pt] (-3.8373,0.1000) -- (-3.4373,0.1000);
  \draw[grar] (-3.6373,2.1000) -- (-1.5500,0.0000);
  \draw[grar] (0.7750,1.3423) -- (-1.5500,0.0000);
  \draw[grar] (-3.6373,2.1000) -- (0.7750,1.3423);
  \draw[idar] (0.0000,4.2000) -- (-3.6373,2.1000);
  \draw[idc,line width=1.1pt] (-1.5454,3.0768) -- (-1.7454,3.4232);
  \draw[idc,line width=1.1pt] (-1.7187,2.9768) -- (-1.9187,3.3232);
  \draw[idc,line width=1.1pt] (-1.8919,2.8768) -- (-2.0919,3.2232);
  \node[font=\scriptsize,edgec] at (1.4708,2.5474) {$124$};
  \node[font=\scriptsize,edgec] at (1.7291,0.7000) {$245$};
  \node[font=\scriptsize,edgec] at (2.6832,-0.4474) {$235$};
  \node[font=\scriptsize,edgec] at (1.4708,-2.5474) {$135$};
  \node[font=\scriptsize,edgec] at (-0.2583,-1.8474) {$156$};
  \node[font=\scriptsize,edgec] at (-1.7291,-2.1000) {$126$};
  \node[font=\scriptsize,edgec] at (-2.9415,0.0000) {$236$};
  \node[font=\scriptsize,edgec] at (-1.4708,1.1474) {$346$};
  \node[font=\scriptsize,edgec] at (-0.9541,2.5474) {$134$};
  \node[font=\scriptsize,edgec] at (0.0000,0.0000) {$456$};
  \node[vtx] at (0.0000,4.2000) {};
  \node[font=\small,idc] at (0.0000,4.6200) {$1$};
  \node[vtx] at (3.6373,2.1000) {};
  \node[font=\small,idc] at (4.0010,2.3100) {$2$};
  \node[vtx] at (3.6373,-2.1000) {};
  \node[font=\small,idc] at (4.0010,-2.3100) {$3$};
  \node[vtx] at (0.0000,-4.2000) {};
  \node[font=\small,idc] at (0.0000,-4.6200) {$1$};
  \node[vtx] at (-3.6373,-2.1000) {};
  \node[font=\small,idc] at (-4.0010,-2.3100) {$2$};
  \node[vtx] at (-3.6373,2.1000) {};
  \node[font=\small,idc] at (-4.0010,2.3100) {$3$};
  \node[vtx] at (0.7750,1.3423) {};
  \node[font=\small,black] at (0.9850,1.7061) {$4$};
  \node[vtx] at (0.7750,-1.3423) {};
  \node[font=\small,black] at (0.9850,-1.7061) {$5$};
  \node[vtx] at (-1.5500,0.0000) {};
  \node[font=\small,black] at (-1.9700,0.0000) {$6$};
\end{tikzpicture}
\caption{The projective plane as ten triangular charts on $K_6$. The
hexagon is a fundamental domain, with boundary edges carrying the same
number of ticks identified, so its six corners are the three vertices
$1,2,3$ and every edge drawn is an edge of $K_6$. Cells $(6,15,10)$ and
$\HN{\bullet}=(\Z,0,\Z/2)$.}
\label{fig:rp2-k6}
\end{figure}
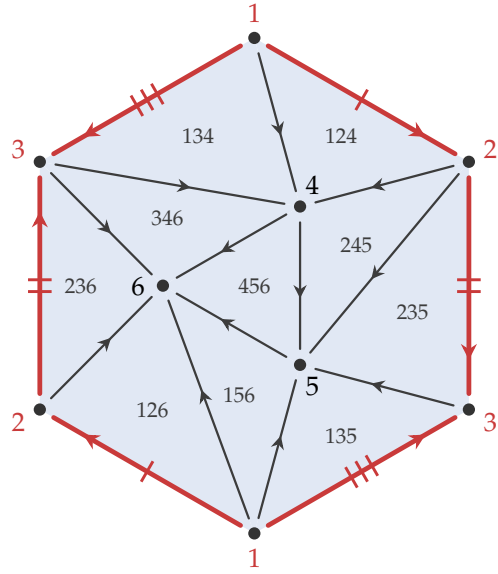

A subset of the triangles may be used instead. The six-vertex
triangulation of the projective plane uses all fifteen edges of $K_6$ and
ten of its twenty triangles, \cref{fig:rp2-k6}. Making these ten triangles
the charts gives cells $(6,15,10)$ and $\HN{\bullet}=(\Z,0,\Z/2)$, the
integral cohomology of $\R P^2$. The graph $K_6$ therefore carries an atlas
with the cohomology of a point, one with the cohomology of a wedge of ten
spheres, and one with the cohomology of the projective plane. The torsion
in the last case cannot occur for a single chart, whose cohomology groups are free
by \cref{cor:toral}.

\section{The same graph with different atlases}
\label{sec:ex-spectrum}

The dependence of the cohomology on the atlas distinguishes this theory
from the intrinsic theories of \cref{sec:other-theories}. In each of the
following families the graph is fixed and the atlas varies.

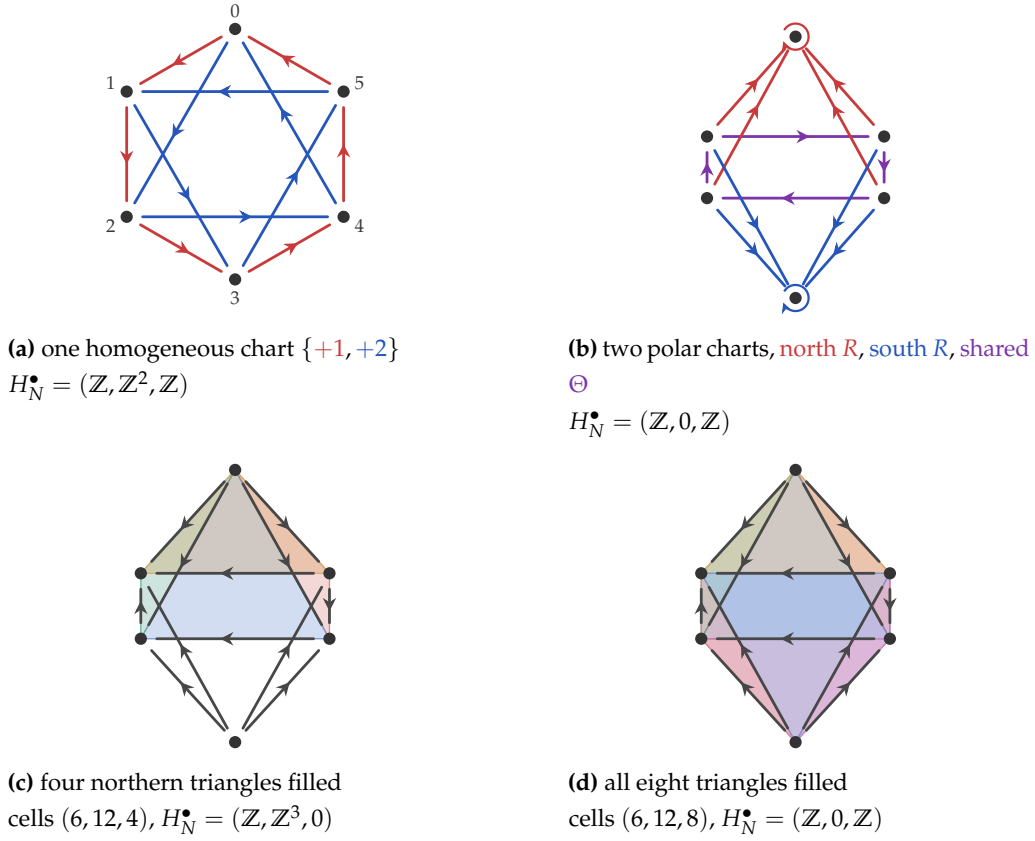
\begin{figure}[htb]
\centering
\renewcommand{\figscale}{0.72}
\begin{subfigure}[t]{0.4\textwidth}\centering
\begin{tikzpicture}[x=1cm,y=1cm,scale=\figscale]
  \draw[cA,ar] (0.0000,2.3000) -- (-1.9919,1.1500);
  \draw[cA,ar] (-1.9919,1.1500) -- (-1.9919,-1.1500);
  \draw[cA,ar] (-1.9919,-1.1500) -- (-0.0000,-2.3000);
  \draw[cA,ar] (-0.0000,-2.3000) -- (1.9919,-1.1500);
  \draw[cA,ar] (1.9919,-1.1500) -- (1.9919,1.1500);
  \draw[cA,ar] (1.9919,1.1500) -- (0.0000,2.3000);
  \draw[cB,ar] (0.0000,2.3000) -- (-1.9919,-1.1500);
  \draw[cB,ar] (-1.9919,1.1500) -- (-0.0000,-2.3000);
  \draw[cB,ar] (-1.9919,-1.1500) -- (1.9919,-1.1500);
  \draw[cB,ar] (-0.0000,-2.3000) -- (1.9919,1.1500);
  \draw[cB,ar] (1.9919,-1.1500) -- (0.0000,2.3000);
  \draw[cB,ar] (1.9919,1.1500) -- (-1.9919,1.1500);
  \node[vtx] at (0.0000,2.3000) {};
  \node[vtx] at (-1.9919,1.1500) {};
  \node[vtx] at (-1.9919,-1.1500) {};
  \node[vtx] at (-0.0000,-2.3000) {};
  \node[vtx] at (1.9919,-1.1500) {};
  \node[vtx] at (1.9919,1.1500) {};
  \node[font=\tiny,gray!45!black] at (0.0000,2.6400) {0};
  \node[font=\tiny,gray!45!black] at (-2.2863,1.3200) {1};
  \node[font=\tiny,gray!45!black] at (-2.2863,-1.3200) {2};
  \node[font=\tiny,gray!45!black] at (-0.0000,-2.6400) {3};
  \node[font=\tiny,gray!45!black] at (2.2863,-1.3200) {4};
  \node[font=\tiny,gray!45!black] at (2.2863,1.3200) {5};
\end{tikzpicture}
\caption{one homogeneous chart $\{\textcolor{cA}{+1},\textcolor{cB}{+2}\}$\\ $\HN{\bullet}=(\Z,\Z^2,\Z)$}
\end{subfigure}\hspace{0.08\textwidth}
\begin{subfigure}[t]{0.4\textwidth}\centering
\begin{tikzpicture}[x=1cm,y=1cm,scale=\figscale]
  \draw[cE,ar] (1.6263,0.5657) -- (1.6263,-0.5657);
  \draw[cE,ar] (1.6263,-0.5657) -- (-1.6263,-0.5657);
  \draw[cE,ar] (-1.6263,-0.5657) -- (-1.6263,0.5657);
  \draw[cE,ar] (-1.6263,0.5657) -- (1.6263,0.5657);
  \draw[cA,ar] (1.6263,0.5657) -- (0.0000,2.4000);
  \draw[cA,ar] (1.6263,-0.5657) -- (0.0000,2.4000);
  \draw[cA,ar] (-1.6263,-0.5657) -- (0.0000,2.4000);
  \draw[cA,ar] (-1.6263,0.5657) -- (0.0000,2.4000);
  \draw[cB,ar] (1.6263,0.5657) -- (0.0000,-2.4000);
  \draw[cB,ar] (1.6263,-0.5657) -- (0.0000,-2.4000);
  \draw[cB,ar] (-1.6263,-0.5657) -- (0.0000,-2.4000);
  \draw[cB,ar] (-1.6263,0.5657) -- (0.0000,-2.4000);
  \draw[cA,lp] (-0.2088,2.5200) arc (150:-160:0.240);
  \draw[cB,lp] (-0.2088,-2.2800) arc (150:-160:0.240);
  \node[vtx] at (0.0000,2.4000) {};
  \node[vtx] at (0.0000,-2.4000) {};
  \node[vtx] at (1.6263,0.5657) {};
  \node[vtx] at (1.6263,-0.5657) {};
  \node[vtx] at (-1.6263,-0.5657) {};
  \node[vtx] at (-1.6263,0.5657) {};
\end{tikzpicture}
\caption{two polar charts, \textcolor{cA}{north $R$}, \textcolor{cB}{south $R$}, \textcolor{cE}{shared $\Theta$}\\ $\HN{\bullet}=(\Z,0,\Z)$}
\end{subfigure}\\[8pt]
\begin{subfigure}[t]{0.4\textwidth}\centering
\begin{tikzpicture}[x=1cm,y=1cm,scale=\figscale]
  \fill[cA,opacity=0.20] (1.7324,0.6010) -- (1.7324,-0.6010) -- (0.0000,2.5000) -- cycle;
  \draw[cA!60,line width=0.5pt] (1.7324,0.6010) -- (1.7324,-0.6010) -- (0.0000,2.5000) -- cycle;
  \fill[cB,opacity=0.20] (1.7324,-0.6010) -- (-1.7324,-0.6010) -- (0.0000,2.5000) -- cycle;
  \draw[cB!60,line width=0.5pt] (1.7324,-0.6010) -- (-1.7324,-0.6010) -- (0.0000,2.5000) -- cycle;
  \fill[cC,opacity=0.20] (-1.7324,-0.6010) -- (-1.7324,0.6010) -- (0.0000,2.5000) -- cycle;
  \draw[cC!60,line width=0.5pt] (-1.7324,-0.6010) -- (-1.7324,0.6010) -- (0.0000,2.5000) -- cycle;
  \fill[cD,opacity=0.20] (-1.7324,0.6010) -- (1.7324,0.6010) -- (0.0000,2.5000) -- cycle;
  \draw[cD!60,line width=0.5pt] (-1.7324,0.6010) -- (1.7324,0.6010) -- (0.0000,2.5000) -- cycle;
  \draw[black!72,ar] (0.0000,-2.5000) -- (1.7324,0.6010);
  \draw[black!72,ar] (0.0000,2.5000) -- (1.7324,0.6010);
  \draw[black!72,ar] (1.7324,0.6010) -- (1.7324,-0.6010);
  \draw[black!72,ar] (1.7324,0.6010) -- (-1.7324,0.6010);
  \draw[black!72,ar] (0.0000,-2.5000) -- (1.7324,-0.6010);
  \draw[black!72,ar] (0.0000,2.5000) -- (1.7324,-0.6010);
  \draw[black!72,ar] (1.7324,-0.6010) -- (-1.7324,-0.6010);
  \draw[black!72,ar] (0.0000,-2.5000) -- (-1.7324,-0.6010);
  \draw[black!72,ar] (0.0000,2.5000) -- (-1.7324,-0.6010);
  \draw[black!72,ar] (-1.7324,-0.6010) -- (-1.7324,0.6010);
  \draw[black!72,ar] (0.0000,-2.5000) -- (-1.7324,0.6010);
  \draw[black!72,ar] (0.0000,2.5000) -- (-1.7324,0.6010);
  \node[vtx] at (0.0000,2.5000) {};
  \node[vtx] at (0.0000,-2.5000) {};
  \node[vtx] at (1.7324,0.6010) {};
  \node[vtx] at (1.7324,-0.6010) {};
  \node[vtx] at (-1.7324,-0.6010) {};
  \node[vtx] at (-1.7324,0.6010) {};
\end{tikzpicture}
\caption{four northern triangles filled\\ cells $(6,12,4)$, $\HN{\bullet}=(\Z,\Z^3,0)$}
\end{subfigure}\hspace{0.08\textwidth}
\begin{subfigure}[t]{0.4\textwidth}\centering
\begin{tikzpicture}[x=1cm,y=1cm,scale=\figscale]
  \fill[cA,opacity=0.20] (1.7324,0.6010) -- (1.7324,-0.6010) -- (0.0000,2.5000) -- cycle;
  \draw[cA!60,line width=0.5pt] (1.7324,0.6010) -- (1.7324,-0.6010) -- (0.0000,2.5000) -- cycle;
  \fill[cB,opacity=0.20] (1.7324,-0.6010) -- (-1.7324,-0.6010) -- (0.0000,2.5000) -- cycle;
  \draw[cB!60,line width=0.5pt] (1.7324,-0.6010) -- (-1.7324,-0.6010) -- (0.0000,2.5000) -- cycle;
  \fill[cC,opacity=0.20] (-1.7324,-0.6010) -- (-1.7324,0.6010) -- (0.0000,2.5000) -- cycle;
  \draw[cC!60,line width=0.5pt] (-1.7324,-0.6010) -- (-1.7324,0.6010) -- (0.0000,2.5000) -- cycle;
  \fill[cD,opacity=0.20] (-1.7324,0.6010) -- (1.7324,0.6010) -- (0.0000,2.5000) -- cycle;
  \draw[cD!60,line width=0.5pt] (-1.7324,0.6010) -- (1.7324,0.6010) -- (0.0000,2.5000) -- cycle;
  \fill[cE,opacity=0.20] (1.7324,0.6010) -- (1.7324,-0.6010) -- (0.0000,-2.5000) -- cycle;
  \draw[cE!60,line width=0.5pt] (1.7324,0.6010) -- (1.7324,-0.6010) -- (0.0000,-2.5000) -- cycle;
  \fill[cF,opacity=0.20] (1.7324,-0.6010) -- (-1.7324,-0.6010) -- (0.0000,-2.5000) -- cycle;
  \draw[cF!60,line width=0.5pt] (1.7324,-0.6010) -- (-1.7324,-0.6010) -- (0.0000,-2.5000) -- cycle;
  \fill[cA,opacity=0.20] (-1.7324,-0.6010) -- (-1.7324,0.6010) -- (0.0000,-2.5000) -- cycle;
  \draw[cA!60,line width=0.5pt] (-1.7324,-0.6010) -- (-1.7324,0.6010) -- (0.0000,-2.5000) -- cycle;
  \fill[cB,opacity=0.20] (-1.7324,0.6010) -- (1.7324,0.6010) -- (0.0000,-2.5000) -- cycle;
  \draw[cB!60,line width=0.5pt] (-1.7324,0.6010) -- (1.7324,0.6010) -- (0.0000,-2.5000) -- cycle;
  \draw[black!72,ar] (0.0000,-2.5000) -- (1.7324,0.6010);
  \draw[black!72,ar] (0.0000,2.5000) -- (1.7324,0.6010);
  \draw[black!72,ar] (1.7324,0.6010) -- (1.7324,-0.6010);
  \draw[black!72,ar] (1.7324,0.6010) -- (-1.7324,0.6010);
  \draw[black!72,ar] (0.0000,-2.5000) -- (1.7324,-0.6010);
  \draw[black!72,ar] (0.0000,2.5000) -- (1.7324,-0.6010);
  \draw[black!72,ar] (1.7324,-0.6010) -- (-1.7324,-0.6010);
  \draw[black!72,ar] (0.0000,-2.5000) -- (-1.7324,-0.6010);
  \draw[black!72,ar] (0.0000,2.5000) -- (-1.7324,-0.6010);
  \draw[black!72,ar] (-1.7324,-0.6010) -- (-1.7324,0.6010);
  \draw[black!72,ar] (0.0000,-2.5000) -- (-1.7324,0.6010);
  \draw[black!72,ar] (0.0000,2.5000) -- (-1.7324,0.6010);
  \node[vtx] at (0.0000,2.5000) {};
  \node[vtx] at (0.0000,-2.5000) {};
  \node[vtx] at (1.7324,0.6010) {};
  \node[vtx] at (1.7324,-0.6010) {};
  \node[vtx] at (-1.7324,-0.6010) {};
  \node[vtx] at (-1.7324,0.6010) {};
\end{tikzpicture}
\caption{all eight triangles filled\\ cells $(6,12,8)$, $\HN{\bullet}=(\Z,0,\Z)$}
\end{subfigure}
\caption{Four atlases on the octahedron graph.}
\label{fig:ex-oct}
\end{figure}

\emph{The octahedron.} \Cref{fig:ex-oct} shows four atlases. As the Cayley
graph $\Cay(\Z/6,\{+1,+2\})$ it carries a single homogeneous chart with the
cohomology of a torus. Viewed as a globe with two poles and one equatorial
square, it carries two polar charts of rank $2$, each a cone on the
equator, glued along the equatorial cycle. Mayer--Vietoris, exactly as in
\cref{ex:globe-mv}, gives $\HN{\bullet}=(\Z,0,\Z)$. Filling only the four
northern triangles with triangle charts gives a disc with a hollow southern
cap, cells $(6,12,4)$ and $\HN{\bullet}=(\Z,\Z^3,0)$, and filling all eight
triangles gives the sphere again, cells $(6,12,8)$.

\begin{figure}[htb]
\centering
\renewcommand{\figscale}{0.72}
\begin{subfigure}[t]{0.32\textwidth}\centering
\begin{tikzpicture}[x=1cm,y=1cm,scale=\figscale]
  \draw[black!72,ar] (0.0000,0.0000) -- (1.2500,0.9500);
  \draw[black!72,ar] (0.0000,0.0000) -- (0.0000,2.6000);
  \draw[black!72,ar] (1.2500,0.9500) -- (1.2500,3.5500);
  \draw[black!72,ar] (0.0000,2.6000) -- (1.2500,3.5500);
  \draw[black!72,ar] (0.0000,0.0000) -- (2.6000,0.0000);
  \draw[black!72,ar] (2.6000,0.0000) -- (3.8500,0.9500);
  \draw[black!72,ar] (2.6000,0.0000) -- (2.6000,2.6000);
  \draw[black!72,ar] (1.2500,0.9500) -- (3.8500,0.9500);
  \draw[black!72,ar] (3.8500,0.9500) -- (3.8500,3.5500);
  \draw[black!72,ar] (0.0000,2.6000) -- (2.6000,2.6000);
  \draw[black!72,ar] (2.6000,2.6000) -- (3.8500,3.5500);
  \draw[black!72,ar] (1.2500,3.5500) -- (3.8500,3.5500);
  \node[vtx] at (0.0000,0.0000) {};
  \node[vtx] at (1.2500,0.9500) {};
  \node[vtx] at (0.0000,2.6000) {};
  \node[vtx] at (1.2500,3.5500) {};
  \node[vtx] at (2.6000,0.0000) {};
  \node[vtx] at (3.8500,0.9500) {};
  \node[vtx] at (2.6000,2.6000) {};
  \node[vtx] at (3.8500,3.5500) {};
\end{tikzpicture}
\caption{$Q_3$, rank-one charts only\\ cells $(8,12,0)$, $\HN{\bullet}=(\Z,\Z^5,0)$}
\end{subfigure}\hfill
\begin{subfigure}[t]{0.32\textwidth}\centering
\begin{tikzpicture}[x=1cm,y=1cm,scale=\figscale]
  \fill[cA,opacity=0.20] (0.0000,0.0000) -- (0.0000,2.6000) -- (1.2500,3.5500) -- (1.2500,0.9500) -- cycle;
  \draw[cA!60,line width=0.5pt] (0.0000,0.0000) -- (0.0000,2.6000) -- (1.2500,3.5500) -- (1.2500,0.9500) -- cycle;
  \fill[cB,opacity=0.20] (2.6000,0.0000) -- (2.6000,2.6000) -- (3.8500,3.5500) -- (3.8500,0.9500) -- cycle;
  \draw[cB!60,line width=0.5pt] (2.6000,0.0000) -- (2.6000,2.6000) -- (3.8500,3.5500) -- (3.8500,0.9500) -- cycle;
  \fill[cC,opacity=0.20] (0.0000,0.0000) -- (2.6000,0.0000) -- (3.8500,0.9500) -- (1.2500,0.9500) -- cycle;
  \draw[cC!60,line width=0.5pt] (0.0000,0.0000) -- (2.6000,0.0000) -- (3.8500,0.9500) -- (1.2500,0.9500) -- cycle;
  \fill[cD,opacity=0.20] (0.0000,2.6000) -- (2.6000,2.6000) -- (3.8500,3.5500) -- (1.2500,3.5500) -- cycle;
  \draw[cD!60,line width=0.5pt] (0.0000,2.6000) -- (2.6000,2.6000) -- (3.8500,3.5500) -- (1.2500,3.5500) -- cycle;
  \draw[black!72,ar] (0.0000,0.0000) -- (1.2500,0.9500);
  \draw[black!72,ar] (0.0000,0.0000) -- (0.0000,2.6000);
  \draw[black!72,ar] (1.2500,0.9500) -- (1.2500,3.5500);
  \draw[black!72,ar] (0.0000,2.6000) -- (1.2500,3.5500);
  \draw[black!72,ar] (0.0000,0.0000) -- (2.6000,0.0000);
  \draw[black!72,ar] (2.6000,0.0000) -- (3.8500,0.9500);
  \draw[black!72,ar] (2.6000,0.0000) -- (2.6000,2.6000);
  \draw[black!72,ar] (1.2500,0.9500) -- (3.8500,0.9500);
  \draw[black!72,ar] (3.8500,0.9500) -- (3.8500,3.5500);
  \draw[black!72,ar] (0.0000,2.6000) -- (2.6000,2.6000);
  \draw[black!72,ar] (2.6000,2.6000) -- (3.8500,3.5500);
  \draw[black!72,ar] (1.2500,3.5500) -- (3.8500,3.5500);
  \node[vtx] at (0.0000,0.0000) {};
  \node[vtx] at (1.2500,0.9500) {};
  \node[vtx] at (0.0000,2.6000) {};
  \node[vtx] at (1.2500,3.5500) {};
  \node[vtx] at (2.6000,0.0000) {};
  \node[vtx] at (3.8500,0.9500) {};
  \node[vtx] at (2.6000,2.6000) {};
  \node[vtx] at (3.8500,3.5500) {};
\end{tikzpicture}
\caption{$Q_3$, four side faces\\ cells $(8,12,4)$, $\HN{\bullet}=(\Z,\Z,0)$}
\end{subfigure}\hfill
\begin{subfigure}[t]{0.32\textwidth}\centering
\begin{tikzpicture}[x=1cm,y=1cm,scale=\figscale]
  \fill[cA,opacity=0.20] (0.0000,0.0000) -- (0.0000,2.6000) -- (1.2500,3.5500) -- (1.2500,0.9500) -- cycle;
  \draw[cA!60,line width=0.5pt] (0.0000,0.0000) -- (0.0000,2.6000) -- (1.2500,3.5500) -- (1.2500,0.9500) -- cycle;
  \fill[cB,opacity=0.20] (2.6000,0.0000) -- (2.6000,2.6000) -- (3.8500,3.5500) -- (3.8500,0.9500) -- cycle;
  \draw[cB!60,line width=0.5pt] (2.6000,0.0000) -- (2.6000,2.6000) -- (3.8500,3.5500) -- (3.8500,0.9500) -- cycle;
  \fill[cC,opacity=0.20] (0.0000,0.0000) -- (2.6000,0.0000) -- (3.8500,0.9500) -- (1.2500,0.9500) -- cycle;
  \draw[cC!60,line width=0.5pt] (0.0000,0.0000) -- (2.6000,0.0000) -- (3.8500,0.9500) -- (1.2500,0.9500) -- cycle;
  \fill[cD,opacity=0.20] (0.0000,2.6000) -- (2.6000,2.6000) -- (3.8500,3.5500) -- (1.2500,3.5500) -- cycle;
  \draw[cD!60,line width=0.5pt] (0.0000,2.6000) -- (2.6000,2.6000) -- (3.8500,3.5500) -- (1.2500,3.5500) -- cycle;
  \fill[cE,opacity=0.20] (0.0000,0.0000) -- (2.6000,0.0000) -- (2.6000,2.6000) -- (0.0000,2.6000) -- cycle;
  \draw[cE!60,line width=0.5pt] (0.0000,0.0000) -- (2.6000,0.0000) -- (2.6000,2.6000) -- (0.0000,2.6000) -- cycle;
  \fill[cF,opacity=0.20] (1.2500,0.9500) -- (3.8500,0.9500) -- (3.8500,3.5500) -- (1.2500,3.5500) -- cycle;
  \draw[cF!60,line width=0.5pt] (1.2500,0.9500) -- (3.8500,0.9500) -- (3.8500,3.5500) -- (1.2500,3.5500) -- cycle;
  \draw[black!72,ar] (0.0000,0.0000) -- (1.2500,0.9500);
  \draw[black!72,ar] (0.0000,0.0000) -- (0.0000,2.6000);
  \draw[black!72,ar] (1.2500,0.9500) -- (1.2500,3.5500);
  \draw[black!72,ar] (0.0000,2.6000) -- (1.2500,3.5500);
  \draw[black!72,ar] (0.0000,0.0000) -- (2.6000,0.0000);
  \draw[black!72,ar] (2.6000,0.0000) -- (3.8500,0.9500);
  \draw[black!72,ar] (2.6000,0.0000) -- (2.6000,2.6000);
  \draw[black!72,ar] (1.2500,0.9500) -- (3.8500,0.9500);
  \draw[black!72,ar] (3.8500,0.9500) -- (3.8500,3.5500);
  \draw[black!72,ar] (0.0000,2.6000) -- (2.6000,2.6000);
  \draw[black!72,ar] (2.6000,2.6000) -- (3.8500,3.5500);
  \draw[black!72,ar] (1.2500,3.5500) -- (3.8500,3.5500);
  \node[vtx] at (0.0000,0.0000) {};
  \node[vtx] at (1.2500,0.9500) {};
  \node[vtx] at (0.0000,2.6000) {};
  \node[vtx] at (1.2500,3.5500) {};
  \node[vtx] at (2.6000,0.0000) {};
  \node[vtx] at (3.8500,0.9500) {};
  \node[vtx] at (2.6000,2.6000) {};
  \node[vtx] at (3.8500,3.5500) {};
\end{tikzpicture}
\caption{$Q_3$, all six faces\\ cells $(8,12,6)$, $\HN{\bullet}=(\Z,0,\Z)$}
\end{subfigure}\\[8pt]
\begin{subfigure}[t]{0.4\textwidth}\centering
\begin{tikzpicture}[x=1cm,y=1cm,scale=\figscale]
  \draw[black!72,ar] (0.0000,2.6000) -- (0.0000,1.3500);
  \draw[black!72,ar] (0.0000,1.3500) -- (-1.2839,0.4172);
  \draw[black!72,ar] (0.0000,1.3500) -- (1.2839,0.4172);
  \draw[black!72,ar] (0.0000,2.6000) -- (-2.4727,0.8034);
  \draw[black!72,ar] (-2.4727,0.8034) -- (-1.2839,0.4172);
  \draw[black!72,ar] (-2.4727,0.8034) -- (-1.5282,-2.1034);
  \draw[black!72,ar] (-1.2839,0.4172) -- (-0.7935,-1.0922);
  \draw[black!72,ar] (-1.5282,-2.1034) -- (-0.7935,-1.0922);
  \draw[black!72,ar] (-1.5282,-2.1034) -- (1.5282,-2.1034);
  \draw[black!72,ar] (-0.7935,-1.0922) -- (0.7935,-1.0922);
  \draw[black!72,ar] (1.5282,-2.1034) -- (0.7935,-1.0922);
  \draw[black!72,ar] (0.7935,-1.0922) -- (1.2839,0.4172);
  \draw[black!72,ar] (0.0000,2.6000) -- (2.4727,0.8034);
  \draw[black!72,ar] (1.5282,-2.1034) -- (2.4727,0.8034);
  \draw[black!72,ar] (2.4727,0.8034) -- (1.2839,0.4172);
  \node[vtx] at (0.0000,2.6000) {};
  \node[vtx] at (0.0000,1.3500) {};
  \node[vtx] at (-2.4727,0.8034) {};
  \node[vtx] at (-1.2839,0.4172) {};
  \node[vtx] at (-1.5282,-2.1034) {};
  \node[vtx] at (-0.7935,-1.0922) {};
  \node[vtx] at (1.5282,-2.1034) {};
  \node[vtx] at (0.7935,-1.0922) {};
  \node[vtx] at (2.4727,0.8034) {};
  \node[vtx] at (1.2839,0.4172) {};
\end{tikzpicture}
\caption{$C_5\times K_2$, rank-one charts only\\ cells $(10,15,0)$, $\HN{\bullet}=(\Z,\Z^6,0)$}
\end{subfigure}\hspace{0.08\textwidth}
\begin{subfigure}[t]{0.4\textwidth}\centering
\begin{tikzpicture}[x=1cm,y=1cm,scale=\figscale]
  \fill[cA,opacity=0.20] (0.0000,2.6000) -- (-2.4727,0.8034) -- (-1.2839,0.4172) -- (0.0000,1.3500) -- cycle;
  \draw[cA!60,line width=0.5pt] (0.0000,2.6000) -- (-2.4727,0.8034) -- (-1.2839,0.4172) -- (0.0000,1.3500) -- cycle;
  \fill[cB,opacity=0.20] (-2.4727,0.8034) -- (-1.5282,-2.1034) -- (-0.7935,-1.0922) -- (-1.2839,0.4172) -- cycle;
  \draw[cB!60,line width=0.5pt] (-2.4727,0.8034) -- (-1.5282,-2.1034) -- (-0.7935,-1.0922) -- (-1.2839,0.4172) -- cycle;
  \fill[cC,opacity=0.20] (-1.5282,-2.1034) -- (1.5282,-2.1034) -- (0.7935,-1.0922) -- (-0.7935,-1.0922) -- cycle;
  \draw[cC!60,line width=0.5pt] (-1.5282,-2.1034) -- (1.5282,-2.1034) -- (0.7935,-1.0922) -- (-0.7935,-1.0922) -- cycle;
  \fill[cD,opacity=0.20] (1.5282,-2.1034) -- (2.4727,0.8034) -- (1.2839,0.4172) -- (0.7935,-1.0922) -- cycle;
  \draw[cD!60,line width=0.5pt] (1.5282,-2.1034) -- (2.4727,0.8034) -- (1.2839,0.4172) -- (0.7935,-1.0922) -- cycle;
  \fill[cE,opacity=0.20] (2.4727,0.8034) -- (0.0000,2.6000) -- (0.0000,1.3500) -- (1.2839,0.4172) -- cycle;
  \draw[cE!60,line width=0.5pt] (2.4727,0.8034) -- (0.0000,2.6000) -- (0.0000,1.3500) -- (1.2839,0.4172) -- cycle;
  \draw[black!72,ar] (0.0000,2.6000) -- (0.0000,1.3500);
  \draw[black!72,ar] (0.0000,1.3500) -- (-1.2839,0.4172);
  \draw[black!72,ar] (0.0000,1.3500) -- (1.2839,0.4172);
  \draw[black!72,ar] (0.0000,2.6000) -- (-2.4727,0.8034);
  \draw[black!72,ar] (-2.4727,0.8034) -- (-1.2839,0.4172);
  \draw[black!72,ar] (-2.4727,0.8034) -- (-1.5282,-2.1034);
  \draw[black!72,ar] (-1.2839,0.4172) -- (-0.7935,-1.0922);
  \draw[black!72,ar] (-1.5282,-2.1034) -- (-0.7935,-1.0922);
  \draw[black!72,ar] (-1.5282,-2.1034) -- (1.5282,-2.1034);
  \draw[black!72,ar] (-0.7935,-1.0922) -- (0.7935,-1.0922);
  \draw[black!72,ar] (1.5282,-2.1034) -- (0.7935,-1.0922);
  \draw[black!72,ar] (0.7935,-1.0922) -- (1.2839,0.4172);
  \draw[black!72,ar] (0.0000,2.6000) -- (2.4727,0.8034);
  \draw[black!72,ar] (1.5282,-2.1034) -- (2.4727,0.8034);
  \draw[black!72,ar] (2.4727,0.8034) -- (1.2839,0.4172);
  \node[vtx] at (0.0000,2.6000) {};
  \node[vtx] at (0.0000,1.3500) {};
  \node[vtx] at (-2.4727,0.8034) {};
  \node[vtx] at (-1.2839,0.4172) {};
  \node[vtx] at (-1.5282,-2.1034) {};
  \node[vtx] at (-0.7935,-1.0922) {};
  \node[vtx] at (1.5282,-2.1034) {};
  \node[vtx] at (0.7935,-1.0922) {};
  \node[vtx] at (2.4727,0.8034) {};
  \node[vtx] at (1.2839,0.4172) {};
\end{tikzpicture}
\caption{$C_5\times K_2$, five side squares\\ cells $(10,15,5)$, $\HN{\bullet}=(\Z,\Z,0)$}
\end{subfigure}
\caption{The cube graph and the pentagonal prism with progressively more
squares filled.}
\label{fig:ex-cube}
\end{figure}
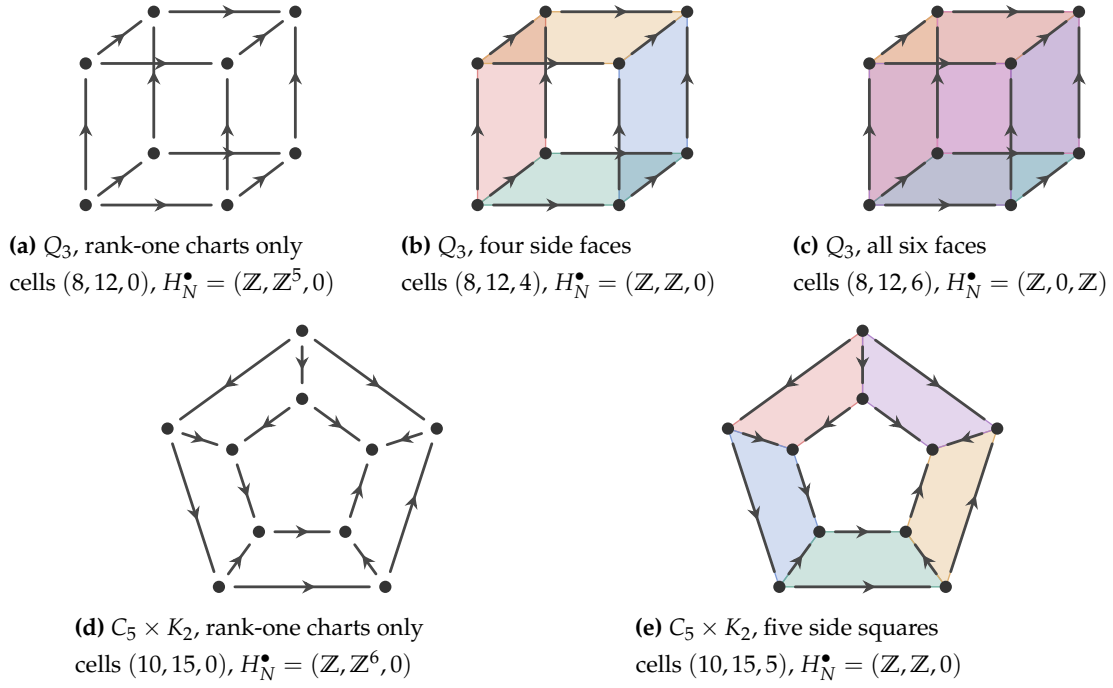

\emph{The cube and the prism.} The graph $Q_3$ of the cube has $8$ vertices
and $12$ edges. With charts of rank $1$ only, no square is filled and
$\HN{\bullet}=(\Z,\Z^5,0)$ is the cohomology of the graph itself. Filling the
four side faces with charts of rank $2$ gives a cylinder, cells $(8,12,4)$
and $\HN{\bullet}=(\Z,\Z,0)$, and filling all six faces gives a sphere,
cells $(8,12,6)$. The prism $C_5\times K_2$ behaves the same way, with
$\HN{\bullet}=(\Z,\Z^6,0)$ for the bare graph and $(\Z,\Z,0)$ for the
cylinder of five side squares. See \cref{fig:ex-cube}.

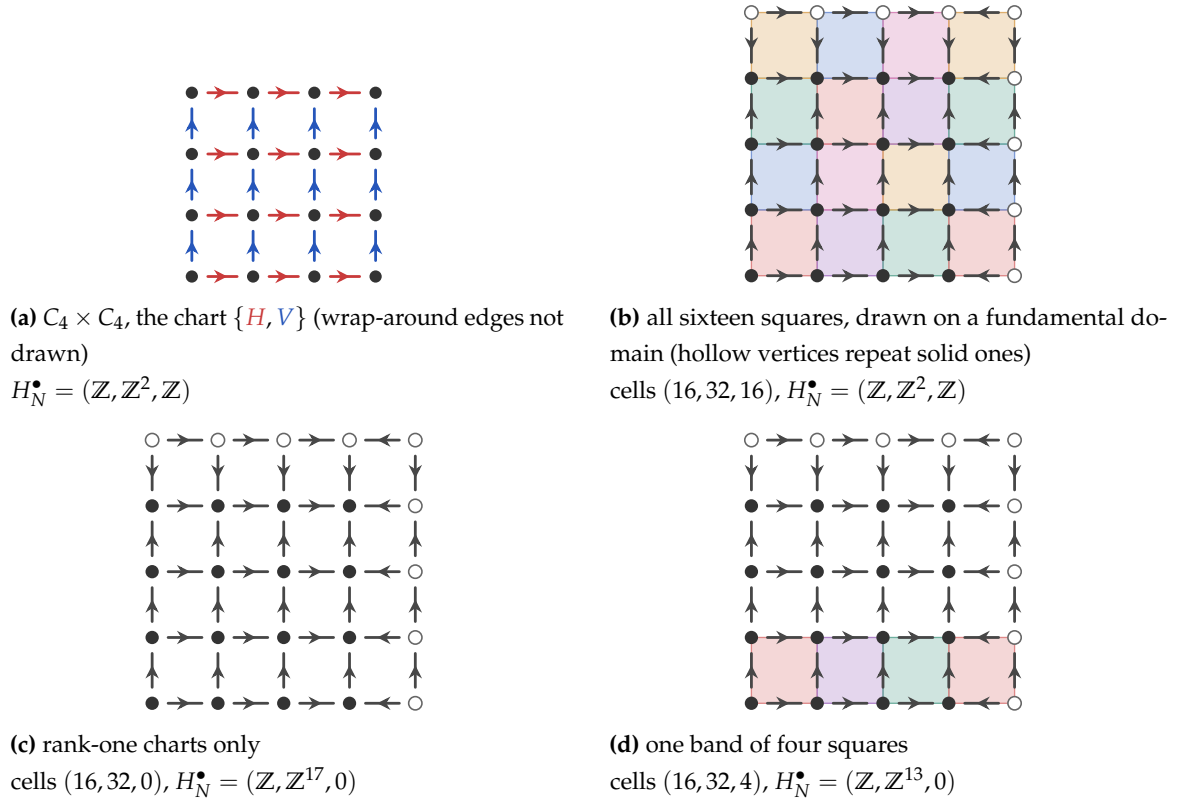
\begin{figure}[htb]
\centering
\renewcommand{\figscale}{0.6}
\begin{subfigure}[t]{0.48\textwidth}\centering
\begin{tikzpicture}[x=1cm,y=1cm,scale=\figscale]
  \draw[cA,ar] (0.0000,0.0000) -- (1.3500,0.0000);
  \draw[cA,ar] (0.0000,1.3500) -- (1.3500,1.3500);
  \draw[cA,ar] (0.0000,2.7000) -- (1.3500,2.7000);
  \draw[cA,ar] (0.0000,4.0500) -- (1.3500,4.0500);
  \draw[cA,ar] (1.3500,0.0000) -- (2.7000,0.0000);
  \draw[cA,ar] (1.3500,1.3500) -- (2.7000,1.3500);
  \draw[cA,ar] (1.3500,2.7000) -- (2.7000,2.7000);
  \draw[cA,ar] (1.3500,4.0500) -- (2.7000,4.0500);
  \draw[cA,ar] (2.7000,0.0000) -- (4.0500,0.0000);
  \draw[cA,ar] (2.7000,1.3500) -- (4.0500,1.3500);
  \draw[cA,ar] (2.7000,2.7000) -- (4.0500,2.7000);
  \draw[cA,ar] (2.7000,4.0500) -- (4.0500,4.0500);
  \draw[cB,ar] (0.0000,0.0000) -- (0.0000,1.3500);
  \draw[cB,ar] (0.0000,1.3500) -- (0.0000,2.7000);
  \draw[cB,ar] (0.0000,2.7000) -- (0.0000,4.0500);
  \draw[cB,ar] (1.3500,0.0000) -- (1.3500,1.3500);
  \draw[cB,ar] (1.3500,1.3500) -- (1.3500,2.7000);
  \draw[cB,ar] (1.3500,2.7000) -- (1.3500,4.0500);
  \draw[cB,ar] (2.7000,0.0000) -- (2.7000,1.3500);
  \draw[cB,ar] (2.7000,1.3500) -- (2.7000,2.7000);
  \draw[cB,ar] (2.7000,2.7000) -- (2.7000,4.0500);
  \draw[cB,ar] (4.0500,0.0000) -- (4.0500,1.3500);
  \draw[cB,ar] (4.0500,1.3500) -- (4.0500,2.7000);
  \draw[cB,ar] (4.0500,2.7000) -- (4.0500,4.0500);
  \node[vtx] at (0.0000,0.0000) {};
  \node[vtx] at (0.0000,1.3500) {};
  \node[vtx] at (0.0000,2.7000) {};
  \node[vtx] at (0.0000,4.0500) {};
  \node[vtx] at (1.3500,0.0000) {};
  \node[vtx] at (1.3500,1.3500) {};
  \node[vtx] at (1.3500,2.7000) {};
  \node[vtx] at (1.3500,4.0500) {};
  \node[vtx] at (2.7000,0.0000) {};
  \node[vtx] at (2.7000,1.3500) {};
  \node[vtx] at (2.7000,2.7000) {};
  \node[vtx] at (2.7000,4.0500) {};
  \node[vtx] at (4.0500,0.0000) {};
  \node[vtx] at (4.0500,1.3500) {};
  \node[vtx] at (4.0500,2.7000) {};
  \node[vtx] at (4.0500,4.0500) {};
\end{tikzpicture}
\caption{$C_4\times C_4$, the chart $\{\textcolor{cA}{H},\textcolor{cB}{V}\}$ (wrap-around edges not drawn)\\ $\HN{\bullet}=(\Z,\Z^2,\Z)$}
\end{subfigure}\hfill
\begin{subfigure}[t]{0.48\textwidth}\centering
\begin{tikzpicture}[x=1cm,y=1cm,scale=\figscale]
  \fill[cA,opacity=0.20] (0.000,0.000) -- (1.450,0.000) -- (1.450,1.450) -- (0.000,1.450) -- cycle;
  \draw[cA!60,line width=0.5pt] (0.000,0.000) -- (1.450,0.000) -- (1.450,1.450) -- (0.000,1.450) -- cycle;
  \fill[cB,opacity=0.20] (0.000,1.450) -- (1.450,1.450) -- (1.450,2.900) -- (0.000,2.900) -- cycle;
  \draw[cB!60,line width=0.5pt] (0.000,1.450) -- (1.450,1.450) -- (1.450,2.900) -- (0.000,2.900) -- cycle;
  \fill[cC,opacity=0.20] (0.000,2.900) -- (1.450,2.900) -- (1.450,4.350) -- (0.000,4.350) -- cycle;
  \draw[cC!60,line width=0.5pt] (0.000,2.900) -- (1.450,2.900) -- (1.450,4.350) -- (0.000,4.350) -- cycle;
  \fill[cD,opacity=0.20] (0.000,4.350) -- (1.450,4.350) -- (1.450,5.800) -- (0.000,5.800) -- cycle;
  \draw[cD!60,line width=0.5pt] (0.000,4.350) -- (1.450,4.350) -- (1.450,5.800) -- (0.000,5.800) -- cycle;
  \fill[cE,opacity=0.20] (1.450,0.000) -- (2.900,0.000) -- (2.900,1.450) -- (1.450,1.450) -- cycle;
  \draw[cE!60,line width=0.5pt] (1.450,0.000) -- (2.900,0.000) -- (2.900,1.450) -- (1.450,1.450) -- cycle;
  \fill[cF,opacity=0.20] (1.450,1.450) -- (2.900,1.450) -- (2.900,2.900) -- (1.450,2.900) -- cycle;
  \draw[cF!60,line width=0.5pt] (1.450,1.450) -- (2.900,1.450) -- (2.900,2.900) -- (1.450,2.900) -- cycle;
  \fill[cA,opacity=0.20] (1.450,2.900) -- (2.900,2.900) -- (2.900,4.350) -- (1.450,4.350) -- cycle;
  \draw[cA!60,line width=0.5pt] (1.450,2.900) -- (2.900,2.900) -- (2.900,4.350) -- (1.450,4.350) -- cycle;
  \fill[cB,opacity=0.20] (1.450,4.350) -- (2.900,4.350) -- (2.900,5.800) -- (1.450,5.800) -- cycle;
  \draw[cB!60,line width=0.5pt] (1.450,4.350) -- (2.900,4.350) -- (2.900,5.800) -- (1.450,5.800) -- cycle;
  \fill[cC,opacity=0.20] (2.900,0.000) -- (4.350,0.000) -- (4.350,1.450) -- (2.900,1.450) -- cycle;
  \draw[cC!60,line width=0.5pt] (2.900,0.000) -- (4.350,0.000) -- (4.350,1.450) -- (2.900,1.450) -- cycle;
  \fill[cD,opacity=0.20] (2.900,1.450) -- (4.350,1.450) -- (4.350,2.900) -- (2.900,2.900) -- cycle;
  \draw[cD!60,line width=0.5pt] (2.900,1.450) -- (4.350,1.450) -- (4.350,2.900) -- (2.900,2.900) -- cycle;
  \fill[cE,opacity=0.20] (2.900,2.900) -- (4.350,2.900) -- (4.350,4.350) -- (2.900,4.350) -- cycle;
  \draw[cE!60,line width=0.5pt] (2.900,2.900) -- (4.350,2.900) -- (4.350,4.350) -- (2.900,4.350) -- cycle;
  \fill[cF,opacity=0.20] (2.900,4.350) -- (4.350,4.350) -- (4.350,5.800) -- (2.900,5.800) -- cycle;
  \draw[cF!60,line width=0.5pt] (2.900,4.350) -- (4.350,4.350) -- (4.350,5.800) -- (2.900,5.800) -- cycle;
  \fill[cA,opacity=0.20] (4.350,0.000) -- (5.800,0.000) -- (5.800,1.450) -- (4.350,1.450) -- cycle;
  \draw[cA!60,line width=0.5pt] (4.350,0.000) -- (5.800,0.000) -- (5.800,1.450) -- (4.350,1.450) -- cycle;
  \fill[cB,opacity=0.20] (4.350,1.450) -- (5.800,1.450) -- (5.800,2.900) -- (4.350,2.900) -- cycle;
  \draw[cB!60,line width=0.5pt] (4.350,1.450) -- (5.800,1.450) -- (5.800,2.900) -- (4.350,2.900) -- cycle;
  \fill[cC,opacity=0.20] (4.350,2.900) -- (5.800,2.900) -- (5.800,4.350) -- (4.350,4.350) -- cycle;
  \draw[cC!60,line width=0.5pt] (4.350,2.900) -- (5.800,2.900) -- (5.800,4.350) -- (4.350,4.350) -- cycle;
  \fill[cD,opacity=0.20] (4.350,4.350) -- (5.800,4.350) -- (5.800,5.800) -- (4.350,5.800) -- cycle;
  \draw[cD!60,line width=0.5pt] (4.350,4.350) -- (5.800,4.350) -- (5.800,5.800) -- (4.350,5.800) -- cycle;
  \draw[black!72,ar] (0.000,0.000) -- (1.450,0.000);
  \draw[black!72,ar] (0.000,0.000) -- (0.000,1.450);
  \draw[black!72,ar] (0.000,1.450) -- (1.450,1.450);
  \draw[black!72,ar] (0.000,1.450) -- (0.000,2.900);
  \draw[black!72,ar] (0.000,2.900) -- (1.450,2.900);
  \draw[black!72,ar] (0.000,2.900) -- (0.000,4.350);
  \draw[black!72,ar] (0.000,4.350) -- (1.450,4.350);
  \draw[black!72,ar] (0.000,5.800) -- (0.000,4.350);
  \draw[black!72,ar] (0.000,5.800) -- (1.450,5.800);
  \draw[black!72,ar] (1.450,0.000) -- (2.900,0.000);
  \draw[black!72,ar] (1.450,0.000) -- (1.450,1.450);
  \draw[black!72,ar] (1.450,1.450) -- (2.900,1.450);
  \draw[black!72,ar] (1.450,1.450) -- (1.450,2.900);
  \draw[black!72,ar] (1.450,2.900) -- (2.900,2.900);
  \draw[black!72,ar] (1.450,2.900) -- (1.450,4.350);
  \draw[black!72,ar] (1.450,4.350) -- (2.900,4.350);
  \draw[black!72,ar] (1.450,5.800) -- (1.450,4.350);
  \draw[black!72,ar] (1.450,5.800) -- (2.900,5.800);
  \draw[black!72,ar] (2.900,0.000) -- (4.350,0.000);
  \draw[black!72,ar] (2.900,0.000) -- (2.900,1.450);
  \draw[black!72,ar] (2.900,1.450) -- (4.350,1.450);
  \draw[black!72,ar] (2.900,1.450) -- (2.900,2.900);
  \draw[black!72,ar] (2.900,2.900) -- (4.350,2.900);
  \draw[black!72,ar] (2.900,2.900) -- (2.900,4.350);
  \draw[black!72,ar] (2.900,4.350) -- (4.350,4.350);
  \draw[black!72,ar] (2.900,5.800) -- (2.900,4.350);
  \draw[black!72,ar] (2.900,5.800) -- (4.350,5.800);
  \draw[black!72,ar] (5.800,0.000) -- (4.350,0.000);
  \draw[black!72,ar] (4.350,0.000) -- (4.350,1.450);
  \draw[black!72,ar] (5.800,1.450) -- (4.350,1.450);
  \draw[black!72,ar] (4.350,1.450) -- (4.350,2.900);
  \draw[black!72,ar] (5.800,2.900) -- (4.350,2.900);
  \draw[black!72,ar] (4.350,2.900) -- (4.350,4.350);
  \draw[black!72,ar] (5.800,4.350) -- (4.350,4.350);
  \draw[black!72,ar] (4.350,5.800) -- (4.350,4.350);
  \draw[black!72,ar] (5.800,5.800) -- (4.350,5.800);
  \draw[black!72,ar] (5.800,0.000) -- (5.800,1.450);
  \draw[black!72,ar] (5.800,1.450) -- (5.800,2.900);
  \draw[black!72,ar] (5.800,2.900) -- (5.800,4.350);
  \draw[black!72,ar] (5.800,5.800) -- (5.800,4.350);
  \node[circle,fill=black!80,inner sep=0pt,minimum size=5pt] at (0.000,0.000) {};
  \node[circle,fill=black!80,inner sep=0pt,minimum size=5pt] at (0.000,1.450) {};
  \node[circle,fill=black!80,inner sep=0pt,minimum size=5pt] at (0.000,2.900) {};
  \node[circle,fill=black!80,inner sep=0pt,minimum size=5pt] at (0.000,4.350) {};
  \node[circle,draw=black!60,fill=white,line width=.6pt,inner sep=0pt,minimum size=5pt] at (0.000,5.800) {};
  \node[circle,fill=black!80,inner sep=0pt,minimum size=5pt] at (1.450,0.000) {};
  \node[circle,fill=black!80,inner sep=0pt,minimum size=5pt] at (1.450,1.450) {};
  \node[circle,fill=black!80,inner sep=0pt,minimum size=5pt] at (1.450,2.900) {};
  \node[circle,fill=black!80,inner sep=0pt,minimum size=5pt] at (1.450,4.350) {};
  \node[circle,draw=black!60,fill=white,line width=.6pt,inner sep=0pt,minimum size=5pt] at (1.450,5.800) {};
  \node[circle,fill=black!80,inner sep=0pt,minimum size=5pt] at (2.900,0.000) {};
  \node[circle,fill=black!80,inner sep=0pt,minimum size=5pt] at (2.900,1.450) {};
  \node[circle,fill=black!80,inner sep=0pt,minimum size=5pt] at (2.900,2.900) {};
  \node[circle,fill=black!80,inner sep=0pt,minimum size=5pt] at (2.900,4.350) {};
  \node[circle,draw=black!60,fill=white,line width=.6pt,inner sep=0pt,minimum size=5pt] at (2.900,5.800) {};
  \node[circle,fill=black!80,inner sep=0pt,minimum size=5pt] at (4.350,0.000) {};
  \node[circle,fill=black!80,inner sep=0pt,minimum size=5pt] at (4.350,1.450) {};
  \node[circle,fill=black!80,inner sep=0pt,minimum size=5pt] at (4.350,2.900) {};
  \node[circle,fill=black!80,inner sep=0pt,minimum size=5pt] at (4.350,4.350) {};
  \node[circle,draw=black!60,fill=white,line width=.6pt,inner sep=0pt,minimum size=5pt] at (4.350,5.800) {};
  \node[circle,draw=black!60,fill=white,line width=.6pt,inner sep=0pt,minimum size=5pt] at (5.800,0.000) {};
  \node[circle,draw=black!60,fill=white,line width=.6pt,inner sep=0pt,minimum size=5pt] at (5.800,1.450) {};
  \node[circle,draw=black!60,fill=white,line width=.6pt,inner sep=0pt,minimum size=5pt] at (5.800,2.900) {};
  \node[circle,draw=black!60,fill=white,line width=.6pt,inner sep=0pt,minimum size=5pt] at (5.800,4.350) {};
  \node[circle,draw=black!60,fill=white,line width=.6pt,inner sep=0pt,minimum size=5pt] at (5.800,5.800) {};
\end{tikzpicture}
\caption{all sixteen squares, drawn on a fundamental domain (hollow vertices repeat solid ones)\\ cells $(16,32,16)$, $\HN{\bullet}=(\Z,\Z^2,\Z)$}
\end{subfigure}\\[8pt]
\begin{subfigure}[t]{0.48\textwidth}\centering
\begin{tikzpicture}[x=1cm,y=1cm,scale=\figscale]
  \draw[black!72,ar] (0.000,0.000) -- (1.450,0.000);
  \draw[black!72,ar] (0.000,0.000) -- (0.000,1.450);
  \draw[black!72,ar] (0.000,1.450) -- (1.450,1.450);
  \draw[black!72,ar] (0.000,1.450) -- (0.000,2.900);
  \draw[black!72,ar] (0.000,2.900) -- (1.450,2.900);
  \draw[black!72,ar] (0.000,2.900) -- (0.000,4.350);
  \draw[black!72,ar] (0.000,4.350) -- (1.450,4.350);
  \draw[black!72,ar] (0.000,5.800) -- (0.000,4.350);
  \draw[black!72,ar] (0.000,5.800) -- (1.450,5.800);
  \draw[black!72,ar] (1.450,0.000) -- (2.900,0.000);
  \draw[black!72,ar] (1.450,0.000) -- (1.450,1.450);
  \draw[black!72,ar] (1.450,1.450) -- (2.900,1.450);
  \draw[black!72,ar] (1.450,1.450) -- (1.450,2.900);
  \draw[black!72,ar] (1.450,2.900) -- (2.900,2.900);
  \draw[black!72,ar] (1.450,2.900) -- (1.450,4.350);
  \draw[black!72,ar] (1.450,4.350) -- (2.900,4.350);
  \draw[black!72,ar] (1.450,5.800) -- (1.450,4.350);
  \draw[black!72,ar] (1.450,5.800) -- (2.900,5.800);
  \draw[black!72,ar] (2.900,0.000) -- (4.350,0.000);
  \draw[black!72,ar] (2.900,0.000) -- (2.900,1.450);
  \draw[black!72,ar] (2.900,1.450) -- (4.350,1.450);
  \draw[black!72,ar] (2.900,1.450) -- (2.900,2.900);
  \draw[black!72,ar] (2.900,2.900) -- (4.350,2.900);
  \draw[black!72,ar] (2.900,2.900) -- (2.900,4.350);
  \draw[black!72,ar] (2.900,4.350) -- (4.350,4.350);
  \draw[black!72,ar] (2.900,5.800) -- (2.900,4.350);
  \draw[black!72,ar] (2.900,5.800) -- (4.350,5.800);
  \draw[black!72,ar] (5.800,0.000) -- (4.350,0.000);
  \draw[black!72,ar] (4.350,0.000) -- (4.350,1.450);
  \draw[black!72,ar] (5.800,1.450) -- (4.350,1.450);
  \draw[black!72,ar] (4.350,1.450) -- (4.350,2.900);
  \draw[black!72,ar] (5.800,2.900) -- (4.350,2.900);
  \draw[black!72,ar] (4.350,2.900) -- (4.350,4.350);
  \draw[black!72,ar] (5.800,4.350) -- (4.350,4.350);
  \draw[black!72,ar] (4.350,5.800) -- (4.350,4.350);
  \draw[black!72,ar] (5.800,5.800) -- (4.350,5.800);
  \draw[black!72,ar] (5.800,0.000) -- (5.800,1.450);
  \draw[black!72,ar] (5.800,1.450) -- (5.800,2.900);
  \draw[black!72,ar] (5.800,2.900) -- (5.800,4.350);
  \draw[black!72,ar] (5.800,5.800) -- (5.800,4.350);
  \node[circle,fill=black!80,inner sep=0pt,minimum size=5pt] at (0.000,0.000) {};
  \node[circle,fill=black!80,inner sep=0pt,minimum size=5pt] at (0.000,1.450) {};
  \node[circle,fill=black!80,inner sep=0pt,minimum size=5pt] at (0.000,2.900) {};
  \node[circle,fill=black!80,inner sep=0pt,minimum size=5pt] at (0.000,4.350) {};
  \node[circle,draw=black!60,fill=white,line width=.6pt,inner sep=0pt,minimum size=5pt] at (0.000,5.800) {};
  \node[circle,fill=black!80,inner sep=0pt,minimum size=5pt] at (1.450,0.000) {};
  \node[circle,fill=black!80,inner sep=0pt,minimum size=5pt] at (1.450,1.450) {};
  \node[circle,fill=black!80,inner sep=0pt,minimum size=5pt] at (1.450,2.900) {};
  \node[circle,fill=black!80,inner sep=0pt,minimum size=5pt] at (1.450,4.350) {};
  \node[circle,draw=black!60,fill=white,line width=.6pt,inner sep=0pt,minimum size=5pt] at (1.450,5.800) {};
  \node[circle,fill=black!80,inner sep=0pt,minimum size=5pt] at (2.900,0.000) {};
  \node[circle,fill=black!80,inner sep=0pt,minimum size=5pt] at (2.900,1.450) {};
  \node[circle,fill=black!80,inner sep=0pt,minimum size=5pt] at (2.900,2.900) {};
  \node[circle,fill=black!80,inner sep=0pt,minimum size=5pt] at (2.900,4.350) {};
  \node[circle,draw=black!60,fill=white,line width=.6pt,inner sep=0pt,minimum size=5pt] at (2.900,5.800) {};
  \node[circle,fill=black!80,inner sep=0pt,minimum size=5pt] at (4.350,0.000) {};
  \node[circle,fill=black!80,inner sep=0pt,minimum size=5pt] at (4.350,1.450) {};
  \node[circle,fill=black!80,inner sep=0pt,minimum size=5pt] at (4.350,2.900) {};
  \node[circle,fill=black!80,inner sep=0pt,minimum size=5pt] at (4.350,4.350) {};
  \node[circle,draw=black!60,fill=white,line width=.6pt,inner sep=0pt,minimum size=5pt] at (4.350,5.800) {};
  \node[circle,draw=black!60,fill=white,line width=.6pt,inner sep=0pt,minimum size=5pt] at (5.800,0.000) {};
  \node[circle,draw=black!60,fill=white,line width=.6pt,inner sep=0pt,minimum size=5pt] at (5.800,1.450) {};
  \node[circle,draw=black!60,fill=white,line width=.6pt,inner sep=0pt,minimum size=5pt] at (5.800,2.900) {};
  \node[circle,draw=black!60,fill=white,line width=.6pt,inner sep=0pt,minimum size=5pt] at (5.800,4.350) {};
  \node[circle,draw=black!60,fill=white,line width=.6pt,inner sep=0pt,minimum size=5pt] at (5.800,5.800) {};
\end{tikzpicture}
\caption{rank-one charts only\\ cells $(16,32,0)$, $\HN{\bullet}=(\Z,\Z^{17},0)$}
\end{subfigure}\hfill
\begin{subfigure}[t]{0.48\textwidth}\centering
\begin{tikzpicture}[x=1cm,y=1cm,scale=\figscale]
  \fill[cA,opacity=0.20] (0.000,0.000) -- (1.450,0.000) -- (1.450,1.450) -- (0.000,1.450) -- cycle;
  \draw[cA!60,line width=0.5pt] (0.000,0.000) -- (1.450,0.000) -- (1.450,1.450) -- (0.000,1.450) -- cycle;
  \fill[cE,opacity=0.20] (1.450,0.000) -- (2.900,0.000) -- (2.900,1.450) -- (1.450,1.450) -- cycle;
  \draw[cE!60,line width=0.5pt] (1.450,0.000) -- (2.900,0.000) -- (2.900,1.450) -- (1.450,1.450) -- cycle;
  \fill[cC,opacity=0.20] (2.900,0.000) -- (4.350,0.000) -- (4.350,1.450) -- (2.900,1.450) -- cycle;
  \draw[cC!60,line width=0.5pt] (2.900,0.000) -- (4.350,0.000) -- (4.350,1.450) -- (2.900,1.450) -- cycle;
  \fill[cA,opacity=0.20] (4.350,0.000) -- (5.800,0.000) -- (5.800,1.450) -- (4.350,1.450) -- cycle;
  \draw[cA!60,line width=0.5pt] (4.350,0.000) -- (5.800,0.000) -- (5.800,1.450) -- (4.350,1.450) -- cycle;
  \draw[black!72,ar] (0.000,0.000) -- (1.450,0.000);
  \draw[black!72,ar] (0.000,0.000) -- (0.000,1.450);
  \draw[black!72,ar] (0.000,1.450) -- (1.450,1.450);
  \draw[black!72,ar] (0.000,1.450) -- (0.000,2.900);
  \draw[black!72,ar] (0.000,2.900) -- (1.450,2.900);
  \draw[black!72,ar] (0.000,2.900) -- (0.000,4.350);
  \draw[black!72,ar] (0.000,4.350) -- (1.450,4.350);
  \draw[black!72,ar] (0.000,5.800) -- (0.000,4.350);
  \draw[black!72,ar] (0.000,5.800) -- (1.450,5.800);
  \draw[black!72,ar] (1.450,0.000) -- (2.900,0.000);
  \draw[black!72,ar] (1.450,0.000) -- (1.450,1.450);
  \draw[black!72,ar] (1.450,1.450) -- (2.900,1.450);
  \draw[black!72,ar] (1.450,1.450) -- (1.450,2.900);
  \draw[black!72,ar] (1.450,2.900) -- (2.900,2.900);
  \draw[black!72,ar] (1.450,2.900) -- (1.450,4.350);
  \draw[black!72,ar] (1.450,4.350) -- (2.900,4.350);
  \draw[black!72,ar] (1.450,5.800) -- (1.450,4.350);
  \draw[black!72,ar] (1.450,5.800) -- (2.900,5.800);
  \draw[black!72,ar] (2.900,0.000) -- (4.350,0.000);
  \draw[black!72,ar] (2.900,0.000) -- (2.900,1.450);
  \draw[black!72,ar] (2.900,1.450) -- (4.350,1.450);
  \draw[black!72,ar] (2.900,1.450) -- (2.900,2.900);
  \draw[black!72,ar] (2.900,2.900) -- (4.350,2.900);
  \draw[black!72,ar] (2.900,2.900) -- (2.900,4.350);
  \draw[black!72,ar] (2.900,4.350) -- (4.350,4.350);
  \draw[black!72,ar] (2.900,5.800) -- (2.900,4.350);
  \draw[black!72,ar] (2.900,5.800) -- (4.350,5.800);
  \draw[black!72,ar] (5.800,0.000) -- (4.350,0.000);
  \draw[black!72,ar] (4.350,0.000) -- (4.350,1.450);
  \draw[black!72,ar] (5.800,1.450) -- (4.350,1.450);
  \draw[black!72,ar] (4.350,1.450) -- (4.350,2.900);
  \draw[black!72,ar] (5.800,2.900) -- (4.350,2.900);
  \draw[black!72,ar] (4.350,2.900) -- (4.350,4.350);
  \draw[black!72,ar] (5.800,4.350) -- (4.350,4.350);
  \draw[black!72,ar] (4.350,5.800) -- (4.350,4.350);
  \draw[black!72,ar] (5.800,5.800) -- (4.350,5.800);
  \draw[black!72,ar] (5.800,0.000) -- (5.800,1.450);
  \draw[black!72,ar] (5.800,1.450) -- (5.800,2.900);
  \draw[black!72,ar] (5.800,2.900) -- (5.800,4.350);
  \draw[black!72,ar] (5.800,5.800) -- (5.800,4.350);
  \node[circle,fill=black!80,inner sep=0pt,minimum size=5pt] at (0.000,0.000) {};
  \node[circle,fill=black!80,inner sep=0pt,minimum size=5pt] at (0.000,1.450) {};
  \node[circle,fill=black!80,inner sep=0pt,minimum size=5pt] at (0.000,2.900) {};
  \node[circle,fill=black!80,inner sep=0pt,minimum size=5pt] at (0.000,4.350) {};
  \node[circle,draw=black!60,fill=white,line width=.6pt,inner sep=0pt,minimum size=5pt] at (0.000,5.800) {};
  \node[circle,fill=black!80,inner sep=0pt,minimum size=5pt] at (1.450,0.000) {};
  \node[circle,fill=black!80,inner sep=0pt,minimum size=5pt] at (1.450,1.450) {};
  \node[circle,fill=black!80,inner sep=0pt,minimum size=5pt] at (1.450,2.900) {};
  \node[circle,fill=black!80,inner sep=0pt,minimum size=5pt] at (1.450,4.350) {};
  \node[circle,draw=black!60,fill=white,line width=.6pt,inner sep=0pt,minimum size=5pt] at (1.450,5.800) {};
  \node[circle,fill=black!80,inner sep=0pt,minimum size=5pt] at (2.900,0.000) {};
  \node[circle,fill=black!80,inner sep=0pt,minimum size=5pt] at (2.900,1.450) {};
  \node[circle,fill=black!80,inner sep=0pt,minimum size=5pt] at (2.900,2.900) {};
  \node[circle,fill=black!80,inner sep=0pt,minimum size=5pt] at (2.900,4.350) {};
  \node[circle,draw=black!60,fill=white,line width=.6pt,inner sep=0pt,minimum size=5pt] at (2.900,5.800) {};
  \node[circle,fill=black!80,inner sep=0pt,minimum size=5pt] at (4.350,0.000) {};
  \node[circle,fill=black!80,inner sep=0pt,minimum size=5pt] at (4.350,1.450) {};
  \node[circle,fill=black!80,inner sep=0pt,minimum size=5pt] at (4.350,2.900) {};
  \node[circle,fill=black!80,inner sep=0pt,minimum size=5pt] at (4.350,4.350) {};
  \node[circle,draw=black!60,fill=white,line width=.6pt,inner sep=0pt,minimum size=5pt] at (4.350,5.800) {};
  \node[circle,draw=black!60,fill=white,line width=.6pt,inner sep=0pt,minimum size=5pt] at (5.800,0.000) {};
  \node[circle,draw=black!60,fill=white,line width=.6pt,inner sep=0pt,minimum size=5pt] at (5.800,1.450) {};
  \node[circle,draw=black!60,fill=white,line width=.6pt,inner sep=0pt,minimum size=5pt] at (5.800,2.900) {};
  \node[circle,draw=black!60,fill=white,line width=.6pt,inner sep=0pt,minimum size=5pt] at (5.800,4.350) {};
  \node[circle,draw=black!60,fill=white,line width=.6pt,inner sep=0pt,minimum size=5pt] at (5.800,5.800) {};
\end{tikzpicture}
\caption{one band of four squares\\ cells $(16,32,4)$, $\HN{\bullet}=(\Z,\Z^{13},0)$}
\end{subfigure}
\caption{Four atlases on the discrete torus $C_4\times C_4$.}
\label{fig:ex-torus}
\end{figure}

\emph{The discrete torus.} The graph $C_4\times C_4$ has $16$ vertices and
$32$ edges. With its product shift chart $\{H,V\}$ of rank $2$ it is a
single homogeneous chart with $\HN{\bullet}=(\Z,\Z^2,\Z)$, which is also
what one gets by filling all sixteen squares with separate charts. With
rank-one charts only, no square is filled and $\HN{\bullet}=(\Z,\Z^{17},0)$.
Filling one band of four squares gives an annulus attached to the rest of
the graph, cells $(16,32,4)$ and $\HN{\bullet}=(\Z,\Z^{13},0)$. See
\cref{fig:ex-torus}.

\begin{figure}[htb]
\centering
\renewcommand{\figscale}{0.75}
\begin{subfigure}[t]{0.45\textwidth}\centering
\begin{tikzpicture}[x=1cm,y=1cm,scale=\figscale]
  \draw[cA,ar] (0.0000,1.2500) -- (0.0000,0.0000);
  \draw[cA,ar] (-1.0825,0.6250) -- (0.0000,0.0000);
  \draw[cA,ar] (-1.0825,-0.6250) -- (0.0000,0.0000);
  \draw[cA,ar] (-0.0000,-1.2500) -- (0.0000,0.0000);
  \draw[cA,ar] (1.0825,-0.6250) -- (0.0000,0.0000);
  \draw[cA,ar] (1.0825,0.6250) -- (0.0000,0.0000);
  \draw[cA,ar] (0.0000,2.5000) -- (0.0000,1.2500);
  \draw[cA,ar] (-2.1651,1.2500) -- (-1.0825,0.6250);
  \draw[cA,ar] (-2.1651,-1.2500) -- (-1.0825,-0.6250);
  \draw[cA,ar] (-0.0000,-2.5000) -- (-0.0000,-1.2500);
  \draw[cA,ar] (2.1651,-1.2500) -- (1.0825,-0.6250);
  \draw[cA,ar] (2.1651,1.2500) -- (1.0825,0.6250);
  \draw[cB,ar] (0.0000,1.2500) -- (-1.0825,0.6250);
  \draw[cB,ar] (-1.0825,0.6250) -- (-1.0825,-0.6250);
  \draw[cB,ar] (-1.0825,-0.6250) -- (-0.0000,-1.2500);
  \draw[cB,ar] (-0.0000,-1.2500) -- (1.0825,-0.6250);
  \draw[cB,ar] (1.0825,-0.6250) -- (1.0825,0.6250);
  \draw[cB,ar] (1.0825,0.6250) -- (0.0000,1.2500);
  \draw[cB,ar] (0.0000,2.5000) -- (-2.1651,1.2500);
  \draw[cB,ar] (-2.1651,1.2500) -- (-2.1651,-1.2500);
  \draw[cB,ar] (-2.1651,-1.2500) -- (-0.0000,-2.5000);
  \draw[cB,ar] (-0.0000,-2.5000) -- (2.1651,-1.2500);
  \draw[cB,ar] (2.1651,-1.2500) -- (2.1651,1.2500);
  \draw[cB,ar] (2.1651,1.2500) -- (0.0000,2.5000);
  \draw[cA,lp] (-0.2088,0.1200) arc (150:-160:0.240);
  \draw[cB,lp] (-0.2088,0.1200) arc (150:-160:0.240);
  \node[vtx] at (0.0000,0.0000) {};
  \node[vtx] at (0.0000,1.2500) {};
  \node[vtx] at (-1.0825,0.6250) {};
  \node[vtx] at (-1.0825,-0.6250) {};
  \node[vtx] at (-0.0000,-1.2500) {};
  \node[vtx] at (1.0825,-0.6250) {};
  \node[vtx] at (1.0825,0.6250) {};
  \node[vtx] at (0.0000,2.5000) {};
  \node[vtx] at (-2.1651,1.2500) {};
  \node[vtx] at (-2.1651,-1.2500) {};
  \node[vtx] at (-0.0000,-2.5000) {};
  \node[vtx] at (2.1651,-1.2500) {};
  \node[vtx] at (2.1651,1.2500) {};
\end{tikzpicture}
\caption{one chart $\{\textcolor{cA}{R},\textcolor{cB}{\Theta}\}$\\ $\HN{\bullet}=(\Z,0,0)$}
\end{subfigure}\hspace{0.05\textwidth}
\begin{subfigure}[t]{0.45\textwidth}\centering
\begin{tikzpicture}[x=1cm,y=1cm,scale=\figscale]
  \draw[cA,ar] (0.0000,2.5000) -- (0.0000,1.2500);
  \draw[cA,ar] (-2.1651,1.2500) -- (-1.0825,0.6250);
  \draw[cA,ar] (-2.1651,-1.2500) -- (-1.0825,-0.6250);
  \draw[cA,ar] (-0.0000,-2.5000) -- (-0.0000,-1.2500);
  \draw[cA,ar] (2.1651,-1.2500) -- (1.0825,-0.6250);
  \draw[cA,ar] (2.1651,1.2500) -- (1.0825,0.6250);
  \draw[cB,ar] (0.0000,1.2500) -- (-1.0825,0.6250);
  \draw[cB,ar] (-1.0825,0.6250) -- (-1.0825,-0.6250);
  \draw[cB,ar] (-1.0825,-0.6250) -- (-0.0000,-1.2500);
  \draw[cB,ar] (-0.0000,-1.2500) -- (1.0825,-0.6250);
  \draw[cB,ar] (1.0825,-0.6250) -- (1.0825,0.6250);
  \draw[cB,ar] (1.0825,0.6250) -- (0.0000,1.2500);
  \draw[cB,ar] (0.0000,2.5000) -- (-2.1651,1.2500);
  \draw[cB,ar] (-2.1651,1.2500) -- (-2.1651,-1.2500);
  \draw[cB,ar] (-2.1651,-1.2500) -- (-0.0000,-2.5000);
  \draw[cB,ar] (-0.0000,-2.5000) -- (2.1651,-1.2500);
  \draw[cB,ar] (2.1651,-1.2500) -- (2.1651,1.2500);
  \draw[cB,ar] (2.1651,1.2500) -- (0.0000,2.5000);
  \draw[cC,ar] (0.0000,1.2500) -- (0.0000,0.0000);
  \draw[cC,ar] (-1.0825,0.6250) -- (0.0000,0.0000);
  \draw[cC,ar] (-1.0825,-0.6250) -- (0.0000,0.0000);
  \draw[cC,ar] (-0.0000,-1.2500) -- (0.0000,0.0000);
  \draw[cC,ar] (1.0825,-0.6250) -- (0.0000,0.0000);
  \draw[cC,ar] (1.0825,0.6250) -- (0.0000,0.0000);
  \draw[cA,lp] (-0.2088,0.1200) arc (150:-160:0.240);
  \draw[cA,lp] (-0.2088,1.3700) arc (150:-160:0.240);
  \draw[cA,lp] (-1.2913,0.7450) arc (150:-160:0.240);
  \draw[cA,lp] (-1.2913,-0.5050) arc (150:-160:0.240);
  \draw[cA,lp] (-0.2088,-1.1300) arc (150:-160:0.240);
  \draw[cA,lp] (0.8737,-0.5050) arc (150:-160:0.240);
  \draw[cA,lp] (0.8737,0.7450) arc (150:-160:0.240);
  \draw[cB,lp] (-0.2088,0.1200) arc (150:-160:0.240);
  \draw[cC,lp] (-0.2088,0.1200) arc (150:-160:0.240);
  \draw[cC,lp] (-0.2088,2.6200) arc (150:-160:0.240);
  \draw[cC,lp] (-2.3739,1.3700) arc (150:-160:0.240);
  \draw[cC,lp] (-2.3739,-1.1300) arc (150:-160:0.240);
  \draw[cC,lp] (-0.2088,-2.3800) arc (150:-160:0.240);
  \draw[cC,lp] (1.9563,-1.1300) arc (150:-160:0.240);
  \draw[cC,lp] (1.9563,1.3700) arc (150:-160:0.240);
  \node[vtx] at (0.0000,0.0000) {};
  \node[vtx] at (0.0000,1.2500) {};
  \node[vtx] at (-1.0825,0.6250) {};
  \node[vtx] at (-1.0825,-0.6250) {};
  \node[vtx] at (-0.0000,-1.2500) {};
  \node[vtx] at (1.0825,-0.6250) {};
  \node[vtx] at (1.0825,0.6250) {};
  \node[vtx] at (0.0000,2.5000) {};
  \node[vtx] at (-2.1651,1.2500) {};
  \node[vtx] at (-2.1651,-1.2500) {};
  \node[vtx] at (-0.0000,-2.5000) {};
  \node[vtx] at (2.1651,-1.2500) {};
  \node[vtx] at (2.1651,1.2500) {};
\end{tikzpicture}
\caption{annulus chart $\{\textcolor{cA}{R},\textcolor{cB}{\Theta}\}$ plus \textcolor{cC}{spoke chart}\\ cells $(13,24,6)$, $\HN{\bullet}=(\Z,\Z^6)$}
\end{subfigure}
\caption{Two atlases on the polar disc with two rings.}
\label{fig:ex-polar}
\end{figure}
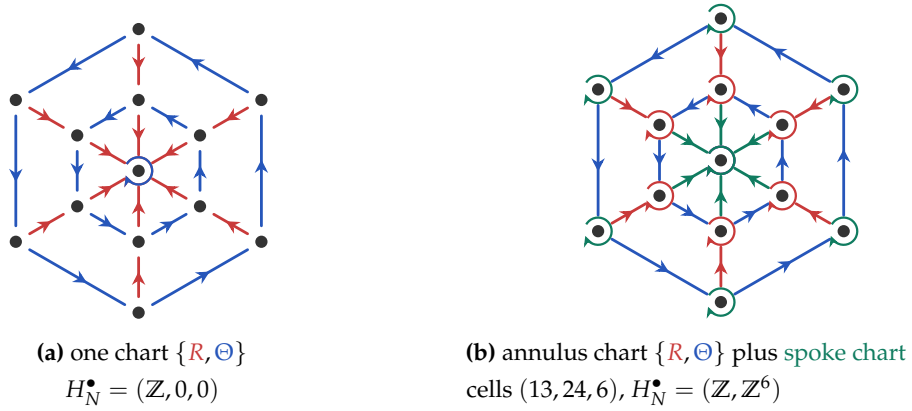

\emph{The polar disc.} Finally, \cref{fig:ex-polar} shows the polar grid
with two rings of six vertices. With the single chart $\{R,\Theta\}$ of
\cref{sec:polar-normalized} it is a disc, $\HN{\bullet}=(\Z,0,0)$. With an
annulus chart of rank $2$ on the two rings, in which $R$ stops on the inner
ring, glued to a spoke chart of rank $1$ sending the inner ring to the
centre, the six squares of the annulus are filled but the six triangles
around the centre are not, and $\HN{\bullet}=(\Z,\Z^6)$. The cells are
$(13,24,6)$, and the six surviving classes are carried by the six unfilled
triangles.

\begin{intuition}
A chart of rank $1$ supplies edges, a chart of rank $2$ supplies squares
or triangles, and a chart of rank $3$ supplies other shapes (we don't have a list of such; to produce one for every rank $r$ is a future goal of this project). The graph fixes the
$1$-skeleton and the atlas determines the higher cells, so the normalized
cohomology of a coordinatized graph is the cohomology of a cell complex
with the graph as $1$-skeleton, and different atlases on the same graph
correspond to different cell complexes. The cohomology is an invariant of
the atlas, not of the graph. Since a finite graph can only have finitely many coordinate systems, where each coordinate field produces a new edge,
we can define the set of its cohomologies as the \textit{de Rham cohomology spectrum} of the graph, which is an invariant. 
\end{intuition}

\chapter*{Conclusion}
\markboth{Conclusion}{Conclusion}
\label{ch:conclusion}

%Claude note: the draft has an empty "Conclusion" heading. This summary and
%the list of open questions are written from the body of the dissertation
%and from claims in your original introduction that are not proved in the
%body (the cohomology spectrum, the stack point of view). Please edit freely.

This dissertation suggests that \(\mathbb N\)-sets provide a useful 
local model for discrete differential geometry. 
They come naturally equipped with a de Rham complex, 
built from the difference operators associated with the action of \(\mathbb N^n\). 
The discrete nature of the construction appears in the twisted Leibniz rule: 
unlike the classical de Rham complex, whose wedge product is graded-commutative, 
the discrete complex carries a twisted, generally non-commutative product.

The relation with cubical sets gives a more geometric interpretation of this construction. 
A point of an \(\mathbb N\)-set together with a collection of coordinate directions determines a cube.
This is also what makes the passage from local to global possible. \(\mathbb N\)-sets 
themselves are not preserved under gluing in the way required for an atlas, 
whereas cubical sets are. The coordinate charts therefore serve as local descriptions, 
while the associated cubical sets provide the objects that can actually be glued. 
From this point of view, the role of \(\mathbb N\)-sets is to place familiar 
cubical sets in a language resembling differential geometry: 
coordinate directions give discrete derivatives, 
commuting directions give higher-dimensional cubes, 
and existing cubical constructions supply degeneracies, products and cohomology. 
One of the main conclusions of the dissertation is therefore that relatively 
little new formal machinery is needed once the correct local model has been chosen.

There are several directions in which this picture remains incomplete.

\begin{enumerate}
\item \emph{Regularity.} The comparison with sheaf cohomology in Theorem 5.3.5 
requires the cubes in each chart to be contractible. The examples considered
here show both possibilities: the polar grid, square and triangle are regular, 
while the Cayley graph examples of Section 5.3.3 are not. 
A combinatorial characterization of regular \(\mathbb N\)-sets would 
therefore clarify which local models can occur in regular atlases.

\item \emph{The set of cohomologies of a graph.} A finite graph admits only 
finitely many coordinate atlases, and hence only finitely many normalized de Rham 
cohomologies. The collection of all such cohomologies may therefore be regarded as
 an invariant of the underlying graph. This raises several natural questions. 
 How much information about the graph is retained by this collection? 
 Can two non-isomorphic graphs have the same collection of normalized cohomologies?
 More generally, does the collection carry additional structure coming from refinements, inclusions or transformations between atlases?

\item \emph{Descent.} The discrete coordinatized spaces of \cref{ch:gluing} are obtained by 
gluing \(\mathbb N\)-sets along cubical subsets by isomorphisms satisfying a cocycle condition. 
This is already close in form to a descent datum. A more categorical formulation would identify the site and 
fibered category in which such an atlas defines a global object. 

This would also connect the construction with the description of differentiable 
stacks by Lie groupoids \citep{behrend-xu2011}, of which the action
categories of \cref{sec:action-categories} are a discrete analogue: their objects are points and their morphisms record the action of coordinate directions. In a more general discrete coordinatized space, one might similarly regard vertices as objects and directed coordinate steps as morphisms, with the higher-dimensional cubical data recording relations among compositions. Such a formulation could make the gluing procedure more conceptual and clarify which parts of the present theory are instances of more 
general descent constructions.

\item \emph{Comparison with path homology.} Path cochains and the de Rham
complex are both differential graded algebras attached to a graph, and the
two appear to agree in low degrees for graphs with few coincidences among
walks. No precise comparison has been made.

\item \emph{Schemes.} A more speculative direction is to ask whether the same philosophy can be applied to discrete schemes. Graphs admit very little structure before a coordinate system is chosen, whereas schemes already come equipped with a sheaf of rings and with strong local algebraic structure. At the same time, many schemes of arithmetic interest are discrete in a geometric sense, so the usual differential intuition coming from smooth manifolds is not directly available.
The generalized calculus philosophy of this dissertation suggests looking for local symmetries or endomorphisms that can play the role of coordinate directions.
Difference operators associated with such maps would then 
interact directly with the structure sheaf, 
giving an entirely algebraic version of the construction. 
The relevant notion of a coordinate chart would no longer 
be imposed only on an underlying set or graph, 
but would have to respect the local ringed-space structure. 
Whether this leads to a useful discrete differential calculus on 
arithmetic or finite schemes remains open, but it provides a natural setting in which to test how 
far the ideas developed here extend beyond graphs.
\end{enumerate}

\appendix
\chapter{Proofs involving sign bookkeeping}
\label{app:proofs}

This appendix collects the proofs from \cref{ch:coordinates} and
\cref{ch:cubical} that consist mainly of expanding both sides of an
identity and keeping track of signs.

\section{The de Rham complex of a product}
\label{app:product}

We complete the proof of \cref{prop:product-dga}. Let $(X,\N^n)$ and
$(Y,\N^m)$ be finite $\N$-sets, $Z=X\times Y$ with the action of
$\N^{n+m}$, and
\[
  \Phi\colon\Omega^\bullet(X,\N^n)\otimes_\Z\Omega^\bullet(Y,\N^m)\to
  \Omega^\bullet(Z,\N^{n+m}),\qquad
  \Phi(f\,dX_I\otimes g\,dY_J)=(f\otimes g)\,dZ_{I\cup J[n]},
\]
where $(f\otimes g)(x,y)=f(x)g(y)$ and $J[n]=\{j+n:\ j\in J\}$. We showed in
the text that $\Phi$ is an isomorphism of graded $\Z$-modules. Two
observations about $Z$ are used below. The shift in a direction
$K=I\cup J[n]$ acts componentwise, $S^Z_K(f\otimes g)=S^X_I(f)\otimes S^Y_J(g)$,
and the difference operators satisfy
$\Delta^Z_i(f\otimes g)=\Delta^X_i(f)\otimes g$ for $i\leq n$ and
$\Delta^Z_{n+j}(f\otimes g)=f\otimes\Delta^Y_j(g)$. Also
$dZ_I\wedge dZ_{J[n]}=dZ_{I\cup J[n]}$ with no sign, since every index of $I$
is smaller than every index of $J[n]$.

%Claude note: the draft's appendix stopped in the middle of this argument
%("Claude: complete this proof"); the rest is supplied here.

\emph{$\Phi$ is multiplicative.} Take $f_1\,dX_{I_1}\otimes g_1\,dY_{J_1}$
and $f_2\,dX_{I_2}\otimes g_2\,dY_{J_2}$ of bidegrees $(a_1,b_1)$ and
$(a_2,b_2)$. By \cref{def:tensor-dga} and the twisted product,
\begin{multline*}
  (f_1\,dX_{I_1}\otimes g_1\,dY_{J_1})\cdot(f_2\,dX_{I_2}\otimes g_2\,dY_{J_2})\\
  =(-1)^{a_2b_1}\,S^X_{I_2}(f_1)f_2\,dX_{I_1}\wedge dX_{I_2}\otimes
  S^Y_{J_2}(g_1)g_2\,dY_{J_1}\wedge dY_{J_2},
\end{multline*}
whose image under $\Phi$ is
\[
  (-1)^{a_2b_1}\varepsilon(I_1,I_2)\varepsilon(J_1,J_2)\,
  \bigl(S^X_{I_2}(f_1)f_2\otimes S^Y_{J_2}(g_1)g_2\bigr)\,
  dZ_{(I_1\cup I_2)\cup(J_1\cup J_2)[n]},
\]
where $dX_{I_1}\wedge dX_{I_2}=\varepsilon(I_1,I_2)\,dX_{I_1\cup I_2}$ and
similarly for $J$, both sides being zero unless the sets are disjoint. On
the other side, with $K_r=I_r\cup J_r[n]$,
\[
  \Phi(f_1\,dX_{I_1}\otimes g_1\,dY_{J_1})\cdot\Phi(f_2\,dX_{I_2}\otimes g_2\,dY_{J_2})
  =S^Z_{K_2}(f_1\otimes g_1)\,(f_2\otimes g_2)\,dZ_{K_1}\wedge dZ_{K_2}.
\]
Here $S^Z_{K_2}(f_1\otimes g_1)(f_2\otimes g_2)=S^X_{I_2}(f_1)f_2\otimes S^Y_{J_2}(g_1)g_2$,
and
\begin{align*}
  dZ_{K_1}\wedge dZ_{K_2}&=dZ_{I_1}\wedge dZ_{J_1[n]}\wedge dZ_{I_2}\wedge dZ_{J_2[n]}
  =(-1)^{a_2b_1}\,dZ_{I_1}\wedge dZ_{I_2}\wedge dZ_{J_1[n]}\wedge dZ_{J_2[n]}\\
  &=(-1)^{a_2b_1}\varepsilon(I_1,I_2)\varepsilon(J_1,J_2)\,dZ_{(I_1\cup I_2)\cup(J_1\cup J_2)[n]},
\end{align*}
moving the $b_1$ factors of $dZ_{J_1[n]}$ past the $a_2$ factors of
$dZ_{I_2}$. The two expressions agree, and $\Phi(1\otimes1)=1$.

\emph{$\Phi$ is a cochain map.} With $\omega=f\,dX_I\otimes g\,dY_J$ of
bidegree $(a,b)$,
\[
  d\omega=\sum_{i=1}^n\Delta^X_i(f)\,dX_i\wedge dX_I\otimes g\,dY_J
  +(-1)^a\sum_{j=1}^mf\,dX_I\otimes\Delta^Y_j(g)\,dY_j\wedge dY_J,
\]
and applying $\Phi$ gives
\[
  \sum_{i=1}^n\bigl(\Delta^X_i(f)\otimes g\bigr)\,dZ_i\wedge dZ_I\wedge dZ_{J[n]}
  +(-1)^a\sum_{j=1}^m\bigl(f\otimes\Delta^Y_j(g)\bigr)\,dZ_I\wedge dZ_{n+j}\wedge dZ_{J[n]}.
\]
Moving $dZ_{n+j}$ to the front of the second sum costs $(-1)^a$, which
cancels the sign in front, and the two sums combine into
$\sum_{k=1}^{n+m}\Delta^Z_k(f\otimes g)\,dZ_k\wedge dZ_{I\cup J[n]}=d\,\Phi(\omega)$.
\qed

\section{Functoriality of the pullback}
\label{app:functor}

We prove \cref{thm:functor}. Let $F=(f,\sigma)\colon(X,\N^n)\to(Y,\N^m)$
be a coordinate morphism, with associated monoid map $\phi$ given by
$\phi(e_j)=e_{\sigma(j)}$ if $\sigma(j)\in[m]$ and $\phi(e_j)=0$ if
$\sigma(j)=*$. We show in turn that $F^*$ is well defined, that it is a
morphism of graded algebras, that it is a cochain map, and that
$F\mapsto F^*$ is functorial.

\emph{Step 1: shifts.} Set $S^Y_*=\id$ and $\Delta^Y_*=0$. Then for every
$j\in[n]$ and $g\in\F_Y$,
\begin{equation}\label{eq:app-shift}
  S^X_j(F^*g)=F^*\bigl(S^Y_{\sigma(j)}g\bigr),\qquad
  \Delta^X_j(F^*g)=F^*\bigl(\Delta^Y_{\sigma(j)}g\bigr).
\end{equation}
Indeed $(S^X_jF^*g)(x)=g(f(e_j\cdot x))=g(\phi(e_j)\cdot f(x))$. If
$\sigma(j)\in[m]$ this is $g(e_{\sigma(j)}\cdot f(x))=(F^*S^Y_{\sigma(j)}g)(x)$,
and if $\sigma(j)=*$ then $\phi(e_j)=0$ and it is $g(f(x))=(F^*S^Y_*g)(x)$.
Subtracting $F^*g$ gives the second identity. In particular
$\Delta^X_j\circ F^*=0$ when $\sigma(j)=*$. For $L\subseteq[n]$ put
$S^X_L=\prod_{l\in L}S^X_l$, and similarly $S^Y_J$ for $J\subseteq[m]$.
Iterating \eqref{eq:app-shift}, for any $L$ on which $\sigma$ is injective
with $\sigma(L)\subseteq[m]$,
\begin{equation}\label{eq:app-shift-L}
  S^X_L\circ F^*=F^*\circ S^Y_{\sigma(L)}.
\end{equation}

\emph{Step 2: the map on the exterior algebra.} Let
$\lambda\colon\Z^m\to\Z^n$ be the linear map with
$\lambda(dY_i)=\sum_{j:\ \sigma(j)=i}dX_j$, the empty sum being $0$, and let
$\Lam\lambda\colon\Lam^\bullet\Z^m\to\Lam^\bullet\Z^n$ be the induced
map of exterior algebras, which is multiplicative by construction. We
define
\[
  F^*(g\,dY_I)=F^*(g)\,\Lam^p\lambda(dY_I),\qquad g\in\F_Y,\ |I|=p,
\]
that is, $F^*=(-\circ f)\otimes\Lam^p\lambda$ on
$\Omega^p(Y,\N^m)=\F_Y\otimes\Lam^p\Z^m$. This is well defined and
$\Z$-linear. Writing out the wedge for $I=\{i_1<\dots<i_p\}$,
\[
  \Lam^p\lambda(dY_I)=\lambda(dY_{i_1})\wedge\dots\wedge\lambda(dY_{i_p})
  =\sum_{(j_1,\dots,j_p)}dX_{j_1}\wedge\dots\wedge dX_{j_p},
\]
the sum over tuples with $\sigma(j_a)=i_a$. The fibres $\sigma^{-1}(i_a)$
are pairwise disjoint, so the $j_a$ are distinct and the tuple determines
the set $K=\{j_1,\dots,j_p\}$, on which $\sigma$ restricts to a bijection
onto $I$. Conversely such a $K$ determines the tuple. Sorting the tuple
into increasing order contributes the sign $\varepsilon_\sigma(K)$ of the
permutation, and
\begin{equation}\label{eq:app-lambda}
  \Lam^p\lambda(dY_I)=\sum_{K}\varepsilon_\sigma(K)\,dX_K,
\end{equation}
the sum over the $K\subseteq[n]$ with $|K|=p$ for which
$\sigma|_K\colon K\to I$ is a bijection. In particular the coefficient of
$dX_K$ is zero unless $\sigma|_K$ is a bijection onto $I$, and
$\Lam^0\lambda(1)=1$ so $F^*(1)=1$.

\emph{Step 3: $F^*$ is a morphism of graded algebras.} Both sides of
$F^*(\omega\cdot\eta)=F^*(\omega)\cdot F^*(\eta)$ are $\Z$-bilinear, so it
suffices to take $\omega=g\,dY_I$ and $\eta=h\,dY_J$. By definition of the
twisted product, $\omega\cdot\eta=S^Y_J(g)h\,dY_I\wedge dY_J$, so using
that $F^*$ is a ring homomorphism on functions and $\Lam\lambda$ is
multiplicative,
\[
  F^*(\omega\cdot\eta)=F^*(S^Y_Jg)\,F^*(h)\,\Lam\lambda(dY_I)\wedge\Lam\lambda(dY_J).
\]
For the other side write $\Lam\lambda(dY_I)=\sum_Kc_K\,dX_K$ and
$\Lam\lambda(dY_J)=\sum_Ld_L\,dX_L$. Then
\[
  F^*(\omega)\cdot F^*(\eta)=\sum_{K,L}c_Kd_L\,S^X_L(F^*g)\,F^*(h)\,dX_K\wedge dX_L.
\]
The two agree provided $S^X_L(F^*g)=F^*(S^Y_Jg)$ whenever $d_L\neq0$. By
Step~2, $d_L\neq0$ forces $\sigma|_L\colon L\to J$ to be a bijection, so
$\sigma$ is injective on $L$ with $\sigma(L)=J\subseteq[m]$, and
\eqref{eq:app-shift-L} applies. Together with $F^*(1)=1$, $F^*$ is a
homomorphism of graded algebras.

\emph{Step 4: $F^*$ is a cochain map.} Let $\omega=g\,dY_I$ with $|I|=p$
and abbreviate $\Theta=\Lam^p\lambda(dY_I)\in\Lam^p\Z^n$, so that
$F^*\omega=F^*(g)\,\Theta$. Since $\Theta$ has constant coefficients,
\[
  d_X(F^*\omega)=\sum_{j=1}^n\Delta^X_j(F^*g)\,dX_j\wedge\Theta.
\]
On the other side $d_Y\omega=\sum_{i=1}^m\Delta^Y_i(g)\,dY_i\wedge dY_I$,
and applying $F^*$, which is multiplicative by Step~3,
\[
  F^*(d_Y\omega)=\sum_{i=1}^mF^*\bigl(\Delta^Y_ig\bigr)\,\lambda(dY_i)\wedge\Theta
  =\sum_{i=1}^m\ \sum_{j:\ \sigma(j)=i}F^*\bigl(\Delta^Y_ig\bigr)\,dX_j\wedge\Theta.
\]
Every $j$ with $\sigma(j)\neq*$ occurs in exactly one inner sum and no $j$
with $\sigma(j)=*$ occurs at all. Using
$F^*(\Delta^Y_{\sigma(j)}g)=\Delta^X_j(F^*g)$ from Step~1,
\[
  F^*(d_Y\omega)=\sum_{j:\ \sigma(j)\neq*}\Delta^X_j(F^*g)\,dX_j\wedge\Theta,
  \qquad
  d_X(F^*\omega)-F^*(d_Y\omega)=\sum_{j:\ \sigma(j)=*}\!\Delta^X_j(F^*g)\,dX_j\wedge\Theta,
\]
and the last sum vanishes because $\Delta^X_j\circ F^*=0$ when
$\sigma(j)=*$.
Since the $g\,dY_I$ span $\Omega^p(Y,\N^m)$ and both sides are additive,
$d_X\circ F^*=F^*\circ d_Y$.

\emph{Step 5: functoriality.} Let $F=(f,\sigma)\colon(X,\N^n)\to(Y,\N^m)$
and $G=(g,\tau)\colon(Y,\N^m)\to(W,\N^l)$ be coordinate morphisms. Then
$G\circ F=(g\circ f,\tau\circ\sigma)$ with $\tau(*)=*$ is again a
coordinate morphism. On functions,
$(G\circ F)^*(h)=h\circ g\circ f=F^*(G^*h)$. On $\Lam^1$,
\[
  \lambda_\sigma(\lambda_\tau(dW_a))=\sum_{i:\ \tau(i)=a}\ \sum_{j:\ \sigma(j)=i}dX_j
  =\sum_{j:\ \tau(\sigma(j))=a}dX_j=\lambda_{\tau\sigma}(dW_a),
\]
the middle equality because a $j$ with $\sigma(j)=*$ has $\tau(\sigma(j))=*$
and is counted on neither side. Hence $\lambda_{\tau\sigma}=\lambda_\sigma\circ\lambda_\tau$
and, taking exterior powers, $\Lam\lambda_{\tau\sigma}=\Lam\lambda_\sigma\circ\Lam\lambda_\tau$.
Since $(G\circ F)^*=(-\circ gf)\otimes\Lam\lambda_{\tau\sigma}$ and
$F^*\circ G^*=((-\circ f)\circ(-\circ g))\otimes(\Lam\lambda_\sigma\circ\Lam\lambda_\tau)$,
the two agree. For $F=\id$ we have $\sigma=\id$, $\lambda_\sigma=\id$ and
$F^*=\id$. Therefore $\Omega^\bullet\colon\CoordSet^\op\to\DGA$ is a
functor.\qed

\section{The differential as a function on cubes}
\label{app:d-on-cubes}

We prove \cref{prop:d-on-cubes}. Let $\omega=\sum_{|I|=p}\omega_I\,dX_I$.
Then
\[
  d\omega=\sum_{i=1}^n\sum_{|I|=p}\Delta_i(\omega_I)\,dX_i\wedge dX_I
  =\sum_{|K|=p+1}(d\omega)_K\,dX_K,
\]
and we must identify $(d\omega)_K$. A pair $(i,I)$ contributes to
$dX_K$ exactly when $i\notin I$ and $\{i\}\cup I=K$. If
$K=\{k_0<k_1<\dots<k_p\}$ the contributing pairs are
$(k_j,K\setminus\{k_j\})$ for $0\leq j\leq p$. Since $k_j$ is in position
$j+1$ of $K$, antisymmetry of the wedge gives
$dX_{k_j}\wedge dX_{K\setminus\{k_j\}}=(-1)^jdX_K$, hence
\[
  (d\omega)_K=\sum_{j=0}^p(-1)^j\Delta_{k_j}\bigl(\omega_{K\setminus\{k_j\}}\bigr).
\]
Evaluating at $x$ and using
$\Delta_{k_j}(\omega_{K\setminus\{k_j\}})(x)
=\omega_{K\setminus\{k_j\}}(e_{k_j}\cdot x)-\omega_{K\setminus\{k_j\}}(x)$
gives the formula of the proposition, with $I=K$ and $i_j=k_j$.\qed

\section{The cubical chain complex}
\label{app:cubical-chain}

We prove \cref{prop:cubical-chain}, that
$\partial_n=\sum_{i=1}^n(-1)^{i-1}(d^1_i-d^0_i)$ satisfies
$\partial_{n}\partial_{n+1}=0$. Write
$\partial^\epsilon_n=\sum_{i=1}^n(-1)^{i-1}d^\epsilon_i$ for
$\epsilon\in\{0,1\}$, so that $\partial_n=\partial^1_n-\partial^0_n$.

%Claude note: the draft's proof ends with "...easy to show. this is true";
%the computation is completed here.

For $\alpha,\beta\in\{0,1\}$,
\[
  \partial^\alpha_n\partial^\beta_{n+1}
  =\sum_{i=1}^n\sum_{j=1}^{n+1}(-1)^{i+j}d^\alpha_id^\beta_j
  =\sum_{1\leq i<j\leq n+1}(-1)^{i+j}d^\alpha_id^\beta_j
  +\sum_{1\leq j\leq i\leq n}(-1)^{i+j}d^\alpha_id^\beta_j.
\]
In the first sum the cubical identity $d^\alpha_id^\beta_j=d^\beta_{j-1}d^\alpha_i$
for $i<j$ applies. Substituting $j'=i$ and $i'=j-1$, so that
$1\leq j'\leq i'\leq n$ and $(-1)^{i+j}=-(-1)^{i'+j'}$, the first sum
becomes $-\sum_{1\leq j'\leq i'\leq n}(-1)^{i'+j'}d^\beta_{i'}d^\alpha_{j'}$.
Hence
\[
  \partial^\alpha_n\partial^\beta_{n+1}
  =\sum_{1\leq j\leq i\leq n}(-1)^{i+j}\bigl(d^\alpha_id^\beta_j-d^\beta_id^\alpha_j\bigr).
\]
For $\alpha=\beta$ the right-hand side vanishes, so
$\partial^1\partial^1=\partial^0\partial^0=0$. For $\alpha\neq\beta$ the
expression is antisymmetric in $(\alpha,\beta)$, so
$\partial^1_n\partial^0_{n+1}=-\partial^0_n\partial^1_{n+1}$. Therefore
\[
  \partial_n\partial_{n+1}
  =\partial^1\partial^1-\partial^1\partial^0-\partial^0\partial^1+\partial^0\partial^0
  =-\partial^1\partial^0-\partial^0\partial^1=0.
\]
\qed

\backmatter
\chapter*{Acknowledgements}
\markboth{Acknowledgements}{Acknowledgements}

I would like to thank my supervisor Prof.~Kobi Kremnizer for his continuous
guidance and support, and my external supervisor Dr.~Xerxes Arsiwalla for numerous
discussions about the nature of discrete geometry.

\bibliographystyle{alpha}
\bibliography{refs}

\end{document}